\documentclass{article}
\usepackage[english]{babel}
\usepackage[T1]{fontenc}
\usepackage[utf8]{inputenc}
\usepackage{amsthm}
\usepackage{amsmath}
\usepackage{amssymb}

\usepackage{palatino}
\usepackage{graphicx}
\usepackage{float}
\usepackage{pdflscape}
\usepackage{upgreek}
\usepackage[colorlinks=true, citecolor =blue]{hyperref}
\usepackage{dsfont}
\usepackage[all]{xy}
\usepackage{array}
\usepackage[top=3cm,bottom=3cm,left=1.5cm,right=1.5cm]{geometry}
\usepackage{enumitem}
\usepackage{mathtools}
\usepackage{subcaption}
\usepackage{authblk}
\usepackage{bbm}
\usepackage{bm}
\usepackage{mathrsfs}
\newcommand{\bl}{\Big\langle}
\newcommand{\br}{\Big\rangle}

\newcommand{\blue}{\textcolor{blue}}
\newcommand{\logles}{\lesssim_{\textup{(log)}}\ }

\newcommand{\p}{\partial}
\newcommand{\eps}{\varepsilon}
\newcommand{\doh}{\mathbf{O}_{\eps,\mathbf{h}}^N}

\newtheorem{theorem}{Theorem}[section]
\newtheorem{lemma}[theorem]{Lemma}
\newtheorem{definition}[theorem]{Definition}
\newtheorem{remark}[theorem]{Remark}
\newtheorem{proposition}[theorem]{Proposition}

\newtheorem{corollary}[theorem]{Corollary}

\newtheorem{conditionbase}{Condition}

\newcommand{\lin}{\textup{lin}}
\usepackage{tikz}
\usetikzlibrary{arrows.meta,positioning,shapes.geometric,fit,calc}
\newcommand{\op}{\textup{Op}}
\numberwithin{equation}{section}
\newcommand{\lomega}{\overleftarrow{\omega}}
\newcommand{\romega}{\overrightarrow{\omega}}

\renewcommand{\SS}{\textup{SS}}
\newcommand{\tr}{\dagger}

\usepackage{tikz-cd}

\newcommand{\avedisk}{\mathbf{D}_{\ell}^N}
\newcommand{\avediskplus}{\mathbf{D}_{2\ell-1}^N}
\newcommand{\local}{\bm{\mu}^N_}
\newcommand{\err}{\textbf{\textup{Err}}}
\newcommand{\ddiv}{\nabla_{\textup{disc}}\cdot}
\newcommand{\lden}{\mathbf{L}_{\textup{den}}^N}
\newcommand{\lcur}{\mathbf{L}_{\textup{cur}}^N}
\newcommand{\errtotal}{\err(\mathfrak{b},\delta,\eps,N)}
\newcommand{\tildeS}{\widetilde{\mathbf{S}}_{\delta,\eps}^N}

\author{Chenjiayue Qi\thanks{chenjiayue@ihes.fr}}
\affil{Institut des Hautes Études Scientifiques (IHES), Université Paris-Saclay, 91440, Bures-sur-Yvette, France}

\date{\today}

\begin{document}

\title{Singular Vortex in the Large Deviation of Empirical Current for Symmetric Simple Exclusion Process}

\maketitle

\begin{abstract}
We prove a large deviation estimate for singular current observables of the two-dimensional symmetric simple exclusion process. In particular, the observable has a non-standard large deviation scaling $N^2(\log N)^{-1}$ due to the presence of singular vortices, whereby trajectories rotating near the origin point yield an enhanced contribution to the current observable. The proof combines the heuristics of Bodineau–Derrida–Lebowitz \cite{Bodineau_Derrida_Lebowitz_2008}, the averaging method of Jara–Menezes \cite{Jara_Menezes_2018}, and novel cancellation and multi-scale renormalization arguments developed by this paper.
\end{abstract}

\paragraph{Ackownledgement:} The author is very grateful to Thierry Bodineau for all the enlightening discussions.

\tableofcontents

\section{Introduction}\label{sec:introduction}

Interacting particle systems give microscopic models for non-equilibrium statistical mechanics. For diffusive systems such as the symmetric simple exclusion process (SSEP), the hydrodynamic limit describes the typical macroscopic evolution through a conservation law \cite{kipnis1989hydrodynamics,kipnis1998scaling}
\begin{equation}\label{eq:hydro_limit_intro}
    \p_t\rho_t=-\nabla\cdot q_t,
\end{equation}
where $\rho$ is the macroscopic particle density and $q_t$ is the macroscopic instantaneous current. For SSEP, the current $q_t$ in the hydrodynamic limit \eqref{eq:hydro_limit_intro} is
\begin{equation*}
q_t=-\frac{1}{2}\nabla\rho_t,
\end{equation*}
so that the density solves the heat equation $\p_t\rho_t=\frac{1}{2}\Delta\rho_t$.

Large deviation principles complement this law by quantifying the cost of atypical density-current trajectories. We refer to \cite{bertini2007stochastic,derrida2007non} for general background on non-equilibrium large deviations. At the macroscopic level, the large deviation functional assigns a cost functional to density and current profiles $(\rho_t,q_t)_{t\geq 0}$, see \cite{bertini2002macroscopic,bertini2007large,bertini2015macroscopic}.

In this paper, we work with the two-dimensional SSEP on the discrete disk $\mathds{B}^N$ with diffusive scaling. The contribution is 
\begin{itemize}
    \item We identify a large-deviation regime for a two-dimensional winding-current observable with speed $N^2/\log N$, due to the presence of a singular vortex. Several large deviation estimates (with non-matching upper and lower bounds) and large deviation principles (with matching bounds) are proved under different conditions, by a multi-scale renormalization argument developped in this paper.
\end{itemize}
\begin{definition}
For $N>1$, define
\begin{equation}\label{eq:discrete_and_continuum_disk}
    \mathds{B}^N:=B_1\cap \frac{1}{N}\mathds{Z}^2,\qquad
    B_1:=\{x\in\mathds{R}^2:|x|\leq 1\}.
\end{equation}
The \underline{generator $L_{\SS}^N$ of symmetric simple exclusion process} is
\begin{equation}\label{eq:SSEP_generator}
L_{\SS}^Nf(\eta)=\frac{N^2}{2}\sum_{x\in\mathds{B}^N}\sum_{e\in\{e_1,e_2\}}\Big(f(\eta^{x,x+e})-f(\eta)\Big),
\end{equation}
where $e_1$ is the horizontal unit edge going right and $e_2$ is the vertical unit edge going up, and
\begin{equation*}
\eta_z^{x,y}=\left\{
\begin{array}{ll}
\eta_z    &z\not\in\{x,y\},  \\
\eta_y     &z=x,\\
\eta_x  & z=y.
\end{array}
\right.
\end{equation*}
Whenever a formula contains a term involving $x\pm e$, it is understood that the corresponding bond is included only when both endpoints belong to $\mathds{B}^N$. We will not repeat this boundary restriction in later formulas. Given a density profile $\rho$ as a function of space, the corresponding \underline{local equilibrium  probability measure} $d\local{\rho}$ is defined as
\begin{equation}\label{eq:def_of_local_equilibrium}
d\local{\rho}:=\bigotimes_{x\in\mathds{B}^N} \textup{Bern}\big(\rho_x\big),
\end{equation}
with $\textup{Bern}(\rho_x)$ being the Bernoulli probability measure with parameter $\rho_x$. Given a function of space $\gamma$ at the initial time, the measure of the SSEP is defined as the local equilibrium measure $d\local{\gamma}$ with density $\gamma_x$ at point $x\in\mathds{B}^N$. We denote the law of the SSEP trajectory starting from $d\local{\gamma}$ as $d\mathds{P}_{\gamma,\SS}^N$, and the projection of $d\mathds{P}_{\gamma,\SS}^N$ at time $t$ as $d\mathds{P}_{\gamma,\SS,t}^N$. In particular, we have
\begin{equation*}
d\mathds{P}_{\gamma,\SS,0}^N=d\local{\gamma}.
\end{equation*}
Given a microscopic trajectory $\{\eta_{t,x}\}_{t\geq 0,x\in\mathds{B}^N}$, the empirical density is
\begin{equation}\label{eq:empirical_density}
\rho_t^N:=\frac{1}{N^2}\sum_{x\in\mathds{B}^N}\eta_{t,x}\delta_{x}.
\end{equation}
For $t\geq 0$, $x\in\mathds{B}^N$, and $e\in\{e_1,e_2\}$, define the \underline{normalized net current} across the oriented bond $(x,x+e)$ by
\begin{equation}\label{eq:net_flow_number}
\begin{split}
\mathcal{Q}_{t,x,e}^{N}:=N^{-1}&\Big(\textup{number of jumps $x\mapsto x+e$ during $[0,t]$}-\textup{number of jumps $x+e\mapsto x$ during $[0,t]$}\Big).
\end{split}
\end{equation}
The associated time-integrated empirical current is the vector-valued measure
\begin{equation}\label{eq:empirical_current}
\mathcal{Q}_{t,e}^{N}:=\frac{1}{N^2}\sum_{x\in \mathds{B}^N}\mathcal{Q}_{t,x,e}^{N}\delta_{x}.
\end{equation}
\end{definition}
\subsection{Motivation: Integrated Current Crossing a Segment}

Large deviations for the density of SSEP was first derived by \cite{kipnis1989hydrodynamics} and for the current by \cite{bertini2007large}.  Intuitively, it can be understood as the following \cite{bertini2015macroscopic,derrida2007non,bodineau2007cumulants}, a profile  $(\rho_t,q_t)_{0\leq t\leq T}$ satisfying the continuity equation
\begin{equation}\label{eq;forced_equation_intro}
    \p_t\rho_t=-\nabla\cdot q_t
\end{equation}
has dynamical cost
\begin{equation}\label{eq:MFT_cost_intro}
    \frac{1}{2}\int_0^T\int_{B_1}
    \frac{\left|q_t+\frac{1}{2}\nabla\rho_t\right|^2}{\rho_t(1-\rho_t)}
    dxdt,
\end{equation}
up to the usual initial and boundary contributions. More precisely, the probability of the empirical density and current $(\rho^N,q^N)$ being close to an atypical profile $(\rho_t,q_t)_{0\leq t\leq T}$ is formally given by
\begin{equation}\label{eq:MFT_exponential_weight}
\mathds{P}_{\gamma,\SS}^N\big[(\rho^N,q^N)\approx(\rho_t,q_t)_{0\leq t\leq T}\big]\approx\exp\Bigg(-N^2 \frac{1}{2}\int_0^T\int_{B_1}\frac{\left|q_t+\frac{1}{2}\nabla\rho_t\right|^2}{\rho_t(1-\rho_t)}dxdt-N^2\textup{Ent}_{\gamma}(\rho_0)\Bigg).
\end{equation}
Here $\textup{Ent}_{\gamma}(\rho_0)$ is an entropy functional representing the initial time large deviation cost. The typical hydrodynamic current $q_t=-\frac{1}{2}\nabla\rho_t$ is exactly the zero-cost profile of \eqref{eq:MFT_cost_intro}. 

The motivation of this paper is to study local current observables in two dimensions.  An example, in the spirit of the formulation of \cite{Bodineau_Derrida_Lebowitz_2008}, is the integrated current crossing a fixed line segment.  To keep the picture concrete, take the horizontal segment with $w\in (0,1)$
\begin{equation*}
    \Gamma_0:=\{(u,0): -w\leq u\leq w\}\subset B_1,
\end{equation*}
with upward normal $b_2=(0,1)$.  At the continuum level, the integrated current crossing this segment is measured by
\begin{equation*}
    \int_0^T\int_{-w}^{w} q_t(u,0)\cdot b_2dudt .
\end{equation*}
At the microscopic level, the same observable is represented by the normalized net number of upward jumps across all the vertical bonds which intersect the segment. The notion of $\mathcal{Q}^N$ is introduced in \eqref{eq:net_flow_number}
\begin{equation}\label{eq:microscopic_crossing}
    \frac{1}{N}\sum_{|i|\leq wN}\mathcal{Q}_{T,(i/N,0),e_2}^{N},\quad e_2:=N^{-1}b_2.
\end{equation}
The large-deviation question is to understand the probability that this crossing current is atypically large, say of size $JT$.

\begin{figure}[h]
    \centering
    \begin{tikzpicture}[scale=1.05]
        \draw[thin, gray] (-3.4,0) -- (3.4,0);
        \draw[very thick, blue] (-1.7,0) -- (1.7,0);
        \fill (-1.7,0) circle (0.045);
        \fill (1.7,0) circle (0.045);
        \node[below] at (-1.7,-0.05) {$\partial_-\Gamma_0$};
        \node[below] at (1.7,-0.05) {$\partial_+\Gamma_0$};
        \draw[->, thick] (0,-0.55) -- (0,0.55);
        \node[right] at (0.08,0.36) {$b_2$};
        \draw[->, thick, red] (-1.7,0.65) arc (90:430:0.65);
        \draw[->, thick, red] (1.7,-0.65) arc (-90:-430:0.65);
        \node[red] at (-2.55,0.95) {local loop};
        \node[red] at (2.55,-0.95) {local loop};
    \end{tikzpicture}
    \caption{Endpoint-vortex for a horizontal segment $\Gamma_0$ with upward normal.  A current through $\Gamma_0$ can be generated by microscopic loops around the two endpoints.  Each loop contributes to the signed upward crossing number of $\Gamma_0$, while its spatial support can be much smaller than the macroscopic length of the segment.}
    \label{fig:bdl_endpoint_vortices}
\end{figure}
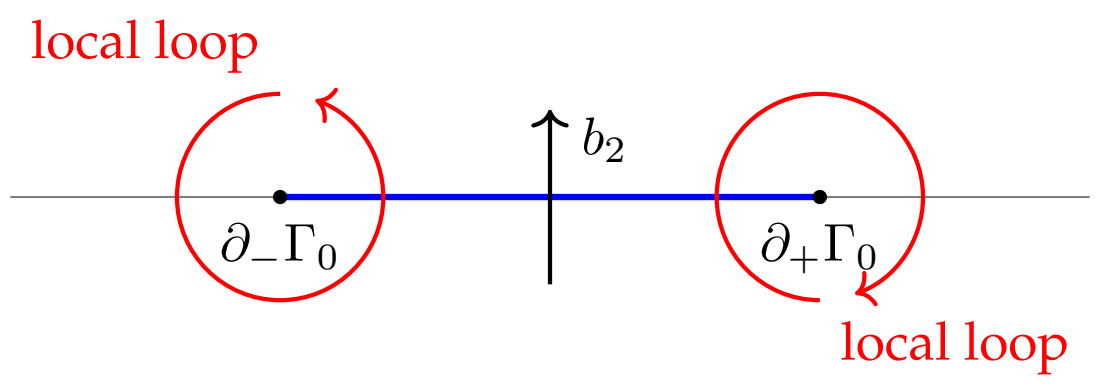

The key observation of \cite{Bodineau_Derrida_Lebowitz_2008} is that the efficient way to produce such an atypical partial current is not to push particles upward uniformly across the whole segment.  Instead, the current concentrates near the two endpoints of $\Gamma_0$, forming a pair of counter-rotating vortices, see Figure \ref{fig:bdl_endpoint_vortices}. At the particle level this means that, particles rotating faster around endpoints can contribute significantly to the observable \eqref{eq:microscopic_crossing}, while only traveling through a very small distance.

This mechanism is exactly where the large deviation functional \eqref{eq:MFT_cost_intro} becomes singular.  If the density is kept essentially constant and a vortex of circulation $J$ is supported on an annulus $\varepsilon\leq r\leq 1$ around an endpoint, its current profile has a rotational form
\begin{equation}\label{eq:singular_vortex_1.1}
q_\varepsilon(x)\simeq\frac{J}{\log(1/\varepsilon)}\nabla^{\perp}\log|x|=\frac{J}{\log(1/\varepsilon)}\frac{1}{|x|}\frac{(-x_2,x_1)}{|x|}.
\end{equation}
Assuming that the density $\rho_t$ is kept constant, the current part of \eqref{eq:MFT_cost_intro} then scales as
\begin{equation}\label{eq:bdl_vortex_cost_intro}
\int_{\varepsilon\leq r\leq 1}\frac{|q_\varepsilon|^2}{2\sigma(\rho)}dx\asymp\frac{J^2}{\sigma(\rho)\log(1/\varepsilon)}.
\end{equation}
Thus the formal continuum cost decreases as the vortex core is squeezed to smaller scales.  In the continuum variational problem, where no microscopic cutoff is present, this produces a degeneracy of the smooth-profile cost.  On a lattice of mesh size $N^{-1}$, however, the collapse stops at $\varepsilon\simeq N^{-1}$, and \eqref{eq:bdl_vortex_cost_intro} leaves the logarithmic correction $(\log N)^{-1}$. Plugging this back into the exponential weight \eqref{eq:MFT_exponential_weight} gives a probability of the form
\begin{equation}\label{eq:non_classical_rate_1.1}
    \exp\left\{-N^2 T\cdot O\left(\frac{J^2}{\sigma(\rho)\log N}\right)\right\},
\end{equation}
which is precisely the large-deviation scale $N^2(\log N)^{-1}$, in contrast with the standard $N^2$ scaling in \eqref{eq:MFT_exponential_weight}.

\subsection{Singular Vortex and Large Deviations of Winding Number}

This mechanism is the conceptual starting point of the present work. To simplify the mathematical analysis of the singular logarithmic large deviation scaling, we focus on a single vortex appearing in the study of integrated currents crossing a segment, see Figure \ref{fig:bdl_endpoint_vortices}. The resulting observable is the winding number around a point: instead of studying how much current crosses a segment, we study how much empirical current circulates around the origin.

\begin{definition}
The single-vortex reduction is encoded by testing the empirical current against the singular circulation field
\begin{equation*}
    \mathcal{V}(x):=\nabla^\perp\log|x|=\frac{(-x_2,x_1)}{|x|^2}.
\end{equation*}
The winding number is
\begin{equation}\label{eq:winding_number_intro}
\langle \mathcal{Q}_T^N, \mathcal{V} \rangle
:= \frac{1}{N^2} \sum_{x \in \mathds{B}^N} \mathcal{Q}_T^N(x) \cdot \mathcal{V}(x),
\end{equation}
with the harmless convention that the value of $\mathcal{V}$ at the origin is fixed arbitrarily, for instance as $0$. According to the gauge change in \cite{Bodineau_Derrida_Lebowitz_2008}, the winding number is the same as the current crossing a segment, in the long time limit.
\end{definition}
For the singular vortex profile just described in \eqref{eq:singular_vortex_1.1}, the inner product with the winding field is of order $O(1)$. Indeed, if for a certain parameter $J\in \mathds{R}$
\begin{equation*}
    q^N(x)=\frac{J}{2\pi\log N}\nabla^\perp\log|x|,
\end{equation*}
then, at the continuum level
\begin{equation*}
\int_{N^{-1}\leq |x|\leq 1} q^N(x)\cdot\mathcal{V}(x)dx=\frac{J}{2\pi\log N}\int_{N^{-1}}^1\int_0^{2\pi}\frac{1}{r^2}rd\theta dr=J.
\end{equation*}
Near the origin, a particle can accumulate a large angular rotation while traveling only a small physical distance. Hence an optimal current producing an order-one winding concentrates around the origin and has the scale
\begin{equation}\label{eq:singular_intro}
    |q^N(x)|\asymp \frac{1}{|x|\log N}.
\end{equation}
Thus the microscopic cutoff at distance $N^{-1}$ predicts the non-standard large-deviation speed $\frac{N^2}{\log N}$ as in \eqref{eq:non_classical_rate_1.1}. At the heuristic level, we can use a cost functional with rate $N^2(\log N)^{-1}$ to describe the large deviations of the winding number
\begin{equation*}
\mathds{P}^N\Big[\textup{winding number $\sim J$}\Big]
\approx
\exp\Big(-\frac{N^2}{\log N}\textup{cost}(J)\Big).
\end{equation*}
In Section \ref{sec:specific_main_result}, this heuristic statement is expanded into rigorous forms, where the cost functional is expressed via \eqref{eq:MFT_cost_intro}.

The vortex mechanism affects not only the current but also the density. Since the current part of the macroscopic cost in \eqref{eq:MFT_cost_intro} has the form $|q|^2/(2\sigma(\rho))$, once a large current is forced near the vortex core it is favorable to increase the local mobility by tilting the density toward a more mobile value. For SSEP the mobility $\sigma(\rho)=\rho(1-\rho)$ is maximal at $\rho=1/2$.

This density response cannot occur on a macroscopic scale. A macroscopic deformation of the density has a cost of order $N^2$, which is much larger than $N^2/\log N$. The density must therefore tilt on a scale smaller than the macroscopic one, but still larger than the microscopic lattice spacing. This leads to density profiles which equal the boundary density $\gamma$ at $|x|=1$, move toward a preferred value $\zeta$ near the vortex core, and interpolate logarithmically between these values on an annulus. The family of density profiles $\mathfrak{u}^N[\gamma,\zeta;\alpha]$ defined in \eqref{eq:u_profile} makes this idea precise: the parameter $\alpha$ fixes a core scale $N^{-\alpha}$, and the logarithmic interpolation balances the gain in current cost against the gradient cost of changing the density.

The large-deviation problem studied here is therefore a coupled singular optimization problem. The current wants to form the vortex profile \eqref{eq:singular_intro}; the density varies near the same singular region in order to increase mobility; and both effects occur across a continuum of scales between $N^{-1}$ and $1$. The main results of this paper identify the corresponding large-deviation bounds with rate $\frac{N^2}{\log N}$ for the winding number.

\textit{Organization of the paper.}
In Section \ref{sec:specific_main_result} we state the main large-deviation estimates, and in Section \ref{sec:proof_strategy} we summarize the proof strategy. Section \ref{sec:multi-scale_replacement} proves the multi-scale super-exponential replacement lemma. Section \ref{sec:proof_upper} proves the upper bound. Section \ref{sec:lower_bound} proves the lower bound, relying on the relative-entropy and concentration estimates developed in Sections \ref{sec:relative_entropy}--\ref{sec:concentration_RN}.

\section{Main Results}\label{sec:main_result}

Let us first introduce some notations before stating the main results.

\begin{definition}\label{def:spaces}
Define $B_1$ as the closed ball of radius $1$, and for arbitrary $N>1$ define $\mathds{B}^N$ as the discretization of $B_1$
\begin{equation}\label{eq:def_of_B1_BN}
B_1:=\Big\{x\in \mathds{R}^2\Big| |x|\leq 1\Big\},\quad \mathds{B}^N:=B_1\cap \frac{1}{N}\mathds{Z}^2.
\end{equation}
For any $N>1$, define the discrete set of particle configurations and current configurations
\begin{equation}\label{eq:def_of_M01_MR2}
\mathbf{M}_{\{0,1\}}^N:=\Big\{\eta\Big|\ \eta:\mathds{B}^N\rightarrow \{0,1\}\Big\},\quad \mathbf{M}_{\mathds{R}^2}^N:=\Big\{q\Big|\ q:\mathds{B}^N\rightarrow \mathds{R}^2\Big\}.
\end{equation}
For any $T>0$, define 
\begin{equation*}
D\big([0,T];\mathbf{M}_{\{0,1\}}^N\times \mathbf{M}_{\mathds{R}^2}^N\big),\quad D\big([0,T];\mathbf{M}_{\{0,1\}}^N\big),\quad D\big([0,T]; \mathbf{M}_{\mathds{R}^2}^N\big),
\end{equation*}
as the cad-lag spaces with value in the corresponding spaces. For any positive integer $k\geq 1$, define $\mathcal{C}^k\big(B_1;(0,1)\big)$ as the space of $k$th-differentiable functions on $B_1$ with value in $(0,1)$. The space $\mathcal{C}^k\big(B_1;\mathds{R}^2\big)$ can be defined similarly. Define $\mathcal{M}_{\mathds{R}^2}(B_1)$ as the space of $\mathds{R}^2$-valued finite Borel measure on $B_1$. Given any measure $q\in \mathcal{M}_{\mathds{R}^2}$ and any $\mathds{R}^2$-valued function $h$, define $\langle q,h \rangle$ as the integral of $h$ with respect to $q$
\begin{equation}\label{eq:bracket}
\langle q,h\rangle:=\int_{B_1}h(x)dq(x).
\end{equation}
\end{definition}
The bold font ($\mathbf{R}$ and $\mathbf{M}$ etc.) refers to discrete objects, and the calligraphy font ($\mathcal{C}^1$ and $\mathcal{M}$ etc.) usually refers to continuum objects.
 
To study the large deviations of the winding number $\langle \mathcal{Q}_T^N,\mathcal{V} \rangle$ defined in \eqref{eq:winding_number_intro}, we need to introduce a proper modified dynamics. Given a forcing field $\mathds{F}$ defined on $\mathds{B}^N$
\begin{equation*}
\mathds{F}_{x}:=(\mathds{F}_{x,e_1},\mathds{F}_{x,e_2}),\quad x\in\mathds{B}^N,
\end{equation*}
the generator of the modified dynamics is given as
\begin{equation}\label{eq:modified_generator}
\begin{split}
L_{\mathds{F}}^Nf(\eta)=\frac{N^2}{2}\sum_{x\in\mathds{B}^N}\sum_{e\in\{e_1,e_2\}}\Big[&\eta_x(1-\eta_{x+e})\big(f(\eta^{x,x+e})-f(\eta)\big)\exp\Big(\frac{\mathds{F}_{x,e}}{N\log N}\Big)\\
+&\eta_x
(1-\eta_{x-e})\big(f(\eta^{x,x-e})-f(\eta)\big)\exp\Big(-\frac{\mathds{F}_{x-e,e}}{N\log N}\Big)\Big].
\end{split}
\end{equation}
We denote the law of trajectories of this modified dynamics starting from $d\local{\gamma}$ as $d\mathds{P}_{\gamma,\mathds{F}}^N$, and the projection of $d\mathds{P}_{\gamma,\mathds{F}}^N$ at time $t$ as $d\mathds{P}_{\gamma,\mathds{F},t}^N$.

\begin{remark}
The generator \eqref{eq:modified_generator} differs from the forms commonly used for the weakly asymmetric simple exclusion process, in which no $\log N$ factor appears. Here the forcing field $\mathds{F}$ typically satisfies
\begin{equation*}
|\mathds{F}_x|\sim \frac{1}{|x|}, \quad x\in\mathds{B}^N.
\end{equation*}
In this case when $x$ is close to the origin, we have $|x|\sim \frac{1}{N}$ and the forcing in \eqref{eq:modified_generator} is much larger than for the standard large deviation theory, for example \cite{kipnis1989hydrodynamics}. The logarithmic factor in \eqref{eq:modified_generator} serves as a normalization, whose role will become clear later.
\end{remark}

\subsection{Statement of Main Results}\label{sec:specific_main_result}

In this subsection, we state the main results of our paper on the large deviation of winding number, for SSEP with constant initial density profile $\gamma\in (0,1)$. The large deviation upper and lower bounds involve the following functionals and profiles.

\begin{definition}\label{def:cost_functional_lower}
For arbitrary parameters $(\gamma,\zeta,\alpha,J,N)$ with
\begin{equation*}
\gamma\in (0,1),\quad \zeta\in (0,1),\quad \alpha\in (0,1),\quad J\in\mathds{R},\quad N>1,
\end{equation*}
define the \underline{density profile} $\mathfrak{u}^N[\gamma,\zeta;\alpha]$ and the \underline{current profile} $\mathfrak{q}^N[J]$ for $x\in B_1$

\begin{equation}\label{eq:u_profile}
\begin{aligned}
\mathfrak{u}^N[\gamma,\zeta;\alpha](x):=
\begin{cases}
\zeta+(\gamma-\zeta)\Big(1+\frac{\log|x|}{\alpha \log N}\Big)^3+\varphi^N(|x|) & N^{-\alpha}\leq |x|<1\\
\zeta & |x|<N^{-\alpha}
\end{cases}
\quad
\mathfrak{q}^N[J]:=\frac{J}{2\pi\log N}\nabla^{\perp}\log|x|\mathds{1}_{N^{-1}\leq |x|\leq 1},
\end{aligned}
\end{equation}
where $\{\varphi^N(r)\}$ is a sequence of smooth functions supported on $[\frac{1}{2},1]$ such that $\mathfrak{u}^N[\gamma,\zeta;\alpha](x)\equiv \gamma$ for $|x|\in [\frac{3}{4},1]$ and $\lVert \varphi^N\rVert_{\mathcal{C}^2}\leq C(\log N)^{-1}$. Consider a profile $(\rho,q)$. We define the \underline{cost functional} $\mathcal{C}^N$ as
\begin{equation}\label{eq:C_cost_functional}
\mathcal{C}^N(\rho,q):=\frac{1}{2}\int_{N^{-1}\leq |x|\leq 1}\frac{|q|^2}{\sigma(\rho)}\log Ndx+\frac{1}{8}\int_{N^{-1}\leq|x|\leq 1}\frac{|\nabla \rho|^2}{\sigma(\rho)}\log Ndx.
\end{equation}
We introduce the initial entropy functional as
\begin{equation*}
\textup{Ent}_{\gamma}(\rho)=\int_{|x|\leq 1}\rho\log\frac{\rho}{\gamma}dx+\int_{|x|\leq 1}(1-\rho)\log\frac{1-\rho}{1-\gamma}dx.
\end{equation*}
Define the $N$-independent function
\begin{equation}\label{eq:C_functional_parametrized}
\mathcal{C}(\gamma,\zeta,\alpha,J):=\lim_{N\rightarrow +\infty}\mathcal{C}^N\Big(\mathfrak{u}^N[\gamma,\zeta;\alpha],\mathfrak{q}^N[J]\Big),
\end{equation}
whose explicit expression is given in Proposition \ref{prop:C_cost_cubic}.
\end{definition}

\begin{figure}
    \centering
    \includegraphics[width=0.6\linewidth]{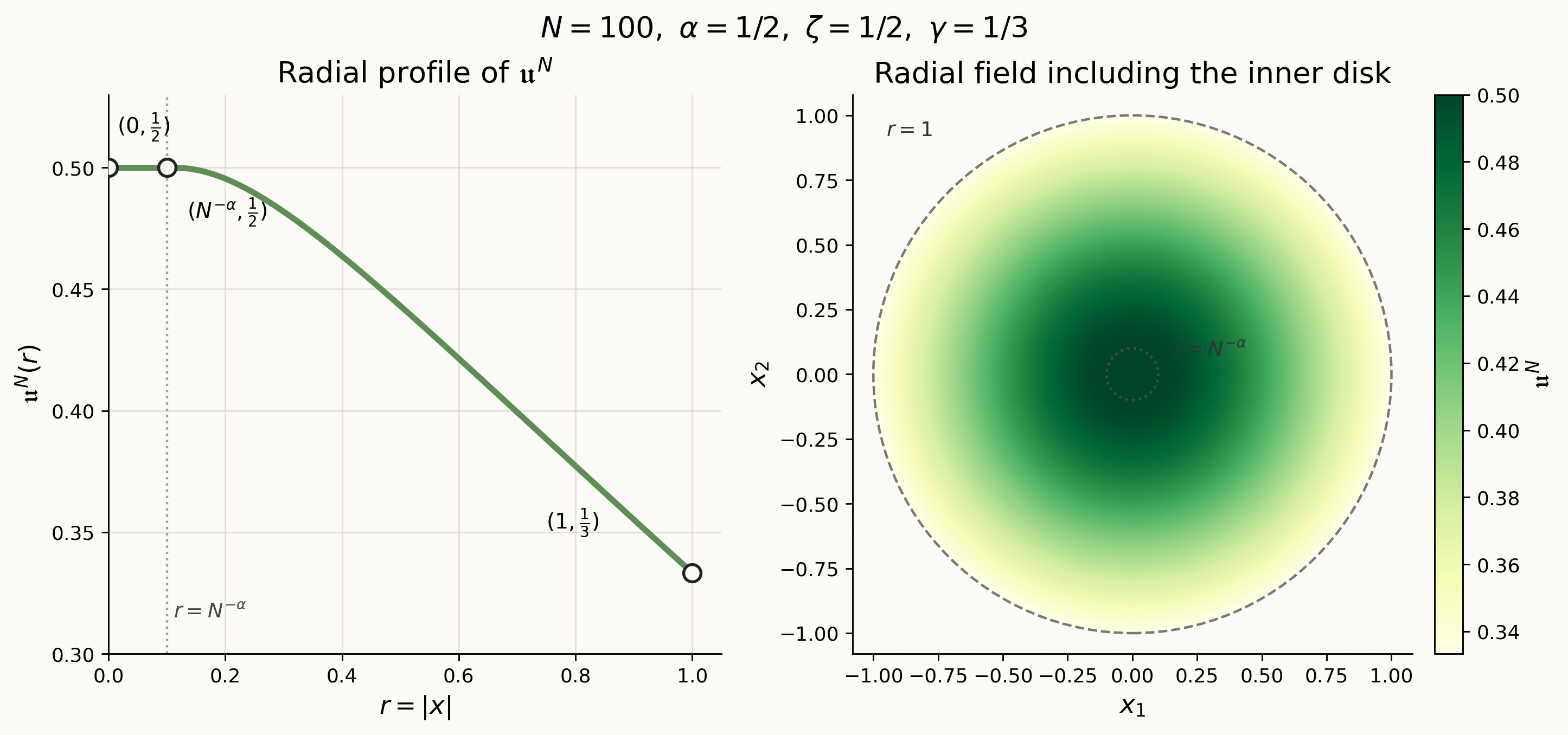}
    \caption{This is an illustration of the density profile $\mathfrak{u}^N[\gamma,\zeta;\alpha]$, with the boundary truncation $\varphi^N$ not pictured.}
    \label{fig:u_profile}
\end{figure}

In the definition \eqref{eq:u_profile} of the density profile $\mathfrak{u}^N[\gamma,\zeta;\alpha]$, $\gamma$ prescribes the density at the outer boundary $|x|=1$, while $\zeta$ is the constant density inside the inner disk $|x|\leq N^{-\alpha}$. The parameter $\alpha$ (together with $N$) determines the size of this disk. These density profiles are illustrated in Figure \ref{fig:u_profile}. The reason of introducing these profiles $\mathfrak{u}^N[\gamma,\zeta;\alpha]$ is that, in the study of large deviations, we need to optimize the cost functional $\mathcal{C}^N$ given in \eqref{eq:C_cost_functional}. When the current $q$ is fixed, pushing $\rho$ closer to $\frac{1}{2}$ increases the mobility $\sigma(\rho)$. This decreases the cost in \eqref{eq:C_cost_functional} related to the current $q$
\begin{equation*}
\log N\int_{N^{-1}\leq |x|\leq 1}\frac{|q|^2}{\sigma(\rho)}dx.
\end{equation*}
This tilting of the density profile can only be at a scale finer than the macroscopic scale, for example $N^{-\alpha}$. Otherwise the cost functional $\mathcal{C}^N$ will blow-up to infinity due to the factor $\log N$. Therefore, a log-interpolation between the outer-boundary value $\gamma$ and the tilted value $\zeta$ is needed, as in the definition \eqref{eq:u_profile} of $\mathfrak{u}^N[\gamma,\zeta;\alpha]$.

The main results of this paper have two folds:
\begin{itemize}
    \item First we establish the large deviation estimate of the event $\langle \mathcal{Q}_T^N,\mathcal{V}\rangle\geq JT$ as $N\rightarrow +\infty$, with parameter $J$ and time $T$ fixed. The convexity of the cost functional $\mathcal{C}^N$ in Definition \ref{def:cost_functional_lower} enables us to derive the large deviation upper bound as a time-independent variational formula. The large deviation lower bound is given by the cost functional of a specific family of profiles. Here for fixed $J$ and $T$, the current upper and lower bounds only match when $\gamma=\frac{1}{2}$, while for general $\gamma\in (0,1)$ they do not match. Then we show, as $a\rightarrow +\infty$ the lower and upper bounds of $\langle \mathcal{Q}_T^N,\mathcal{V}\rangle\geq aJT$ with fixed $J>0$ asymptotically match for general $\gamma\in (0,1)$.
\end{itemize}

\begin{theorem}\label{th:one-scale-deviation}
\textup{\textbf{[Singular Large Deviation Estimate]}} Take any parameters $(J,T,\gamma)$ with
\begin{equation*}
J>0,\quad T>0,\quad \gamma\in (0,1).
\end{equation*}
The large deviation upper bound for the winding number \eqref{eq:winding_number_intro} is
\begin{equation}\label{eq:LDP_upper_statement}
\limsup_{N\rightarrow +\infty}\frac{\log N}{N^2}\log\mathds{P}_{\gamma,\SS}^N\Big[\langle \mathcal{Q}_T^N,\mathcal{V}\rangle\geq JT\Big]\leq -F^{\uparrow}(J)T.
\end{equation}
Here the function $F^{\uparrow}(J)$ is defined as
\begin{equation}\label{eq:upper_rate}
F^{\uparrow}(J):=\inf_{( \rho,q)}\Bigg\{\frac{1}{2}\int_{e^{-1}\leq |x|\leq 1}\frac{|q|^2}{\sigma(\rho)}dx+\frac{1}{8}\int_{e^{-1}\leq|x|\leq 1}\frac{|\nabla \rho|^2}{\sigma(\rho)}dx\Bigg\},
\end{equation}
where the infimum over $(\rho,q)$ is taken within
\begin{equation*}
\rho\in \mathcal{C}^1\big(B_1;(0,1)\big),\quad \rho_x\equiv \gamma \textup{ on $x\in \p B_1$};\quad \textup{ $q$ is a function from $B_1$ to $\mathds{R}^2$ such that } \int \mathds{1}_{\{e^{-1}\leq |x|\leq 1\}}q\cdot\mathcal{V}dx=J.
\end{equation*}
We also have the large deviation lower bound
\begin{equation*}
\liminf_{N\rightarrow +\infty}\frac{\log N}{N^2}\log\mathds{P}_{\gamma,\SS}^N\Big[\langle \mathcal{Q}_T^N,\mathcal{V}\rangle\geq  JT\Big]\geq -F^{\downarrow}(J)T,
\end{equation*}
where the rate function $F^{\downarrow}(J)$ is explicitly given by the functional $\mathcal{C}$ in \eqref{eq:C_functional_parametrized}
\begin{equation*}
F^{\downarrow}(J)=\mathcal{C}\Big(\gamma,\frac{1}{2},\frac{1}{10}(1+|J|)^{-1},J\Big).
\end{equation*}
\end{theorem}

The two bounds in Theorem \ref{th:one-scale-deviation} do not match. The upper bound is already optimal, while the lower bound only uses the cost of an approximate optimizer in \eqref{eq:u_profile} due to the singularity of the profiles.

\underline{\textit{Explanation of the upper bound:}}

The vortex structure of the optimal profile is self-similar, see for example Lemma \ref{lem:C_self_similarity}. Therefore the microscopic cut-off $N^{-1}$ in \eqref{eq:C_cost_functional} can be replaced by any constant, here $e^{-1}$, in the limit $N\rightarrow +\infty$. In the study of large deviations of simple exclusion process in a fixed finite time window $[0,T]$, the large deviation upper and lower bounds are frequently expressed as variational problems over \textit{trajectories}, where the profiles may depend on time, for example see \cite[Theorem 4.1 and 5.3]{kipnis1989hydrodynamics}, \cite[Theorem 1]{benois1995large}, \cite[Theorem 2.4 and 2.6]{bodineau2012large}. In this paper, based on the convexity of the cost functional, we derive the variational formula over stationary profiles as the optimal upper bound. 

\underline{\textit{Explanation of the lower bound:}}

In Theorem \ref{th:one-scale-deviation}, the rate function $F^{\downarrow}(J)$ for the lower bound is given as the deviation cost of a test profile $(\rho,q)$ according to Definition \ref{def:cost_functional_lower}
\begin{equation}\label{eq:approx_optimizer}
\rho=\mathfrak{u}^N[\gamma,\frac{1}{2};\alpha],\quad q=\frac{J}{2\pi\log N}\nabla^{\perp}\log|x|\mathds{1}_{N^{-1}\leq |x|\leq 1},\quad \textup{with }\alpha=\frac{1}{10}\frac{1}{1+|J|},
\end{equation}
where the density profile is tilted to $\frac{1}{2}$ within the disk
\begin{equation*}
|x|\leq N^{-\alpha},\quad \alpha=\frac{1}{10}\frac{1}{1+|J|}.
\end{equation*}
The profile $(\rho,q)$ in \eqref{eq:approx_optimizer} can be seen as a kind of \textit{approximate optimizer}: when $|J|$ is large, the cost is improved by tilting the density profile from $\gamma$ to $\frac{1}{2}$.

For general initial data $\gamma$, the upper and lower bounds in Theorem \ref{th:one-scale-deviation} do not match. However by choosing the initial density profile $\gamma=\frac{1}{2}$, and therefore as the flat optimal density profile, Theorem \ref{th:one-scale-deviation} gives the following corollary with matching lower and upper bounds, therefore a large deviation principle.
\begin{corollary}
Fix $\gamma=\frac{1}{2}$. Take any time $T>0$, and any $J>0$. The large deviation upper and lower bounds match
\begin{equation*}
\begin{split}
&\lim_{N\rightarrow +\infty}\frac{\log N}{N^2}\log\mathds{P}_{\frac{1}{2},\SS}^N\Big[\langle \mathcal{Q}_T^N,\mathcal{V}\rangle\geq JT\Big]= -\pi^{-1} J^2T.
\end{split}
\end{equation*}
\end{corollary}

For general $\gamma\in (0,1)$, as $J\rightarrow +\infty$, the approximate optimizer in \eqref{eq:approx_optimizer} converges to the constant density profile $\frac{1}{2}$, with the highest possible mobility $\sigma(\rho)=\frac{1}{4}$. In this limit, we recover a full large deviation principle.

\begin{figure}[h]
    \centering
    \includegraphics[width=0.8\linewidth]{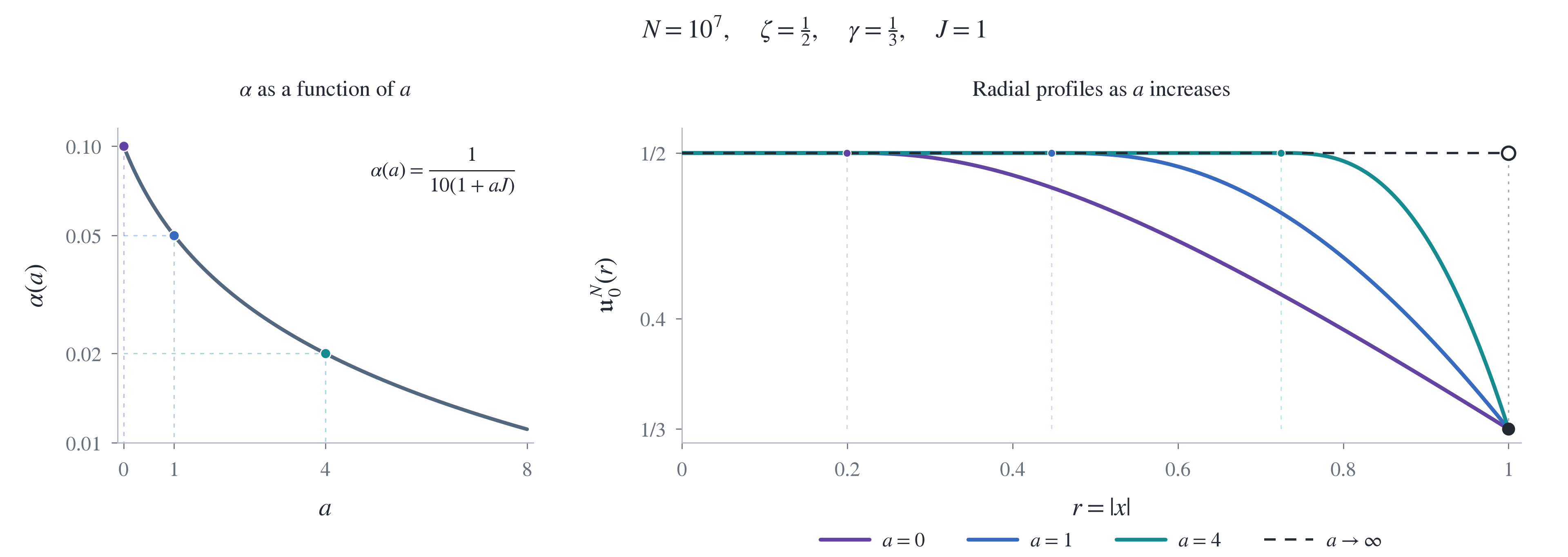}
    \caption{Behavior of the density profile $\mathfrak{u}^N$ in \eqref{eq:u_profile} for large current $J\rightarrow
    +\infty$.}
    \label{fig:large_J_behavior}
\end{figure}

\begin{corollary}\label{cor:two_scale}
\textup{\textbf{[Large Current Asymptotics]}} Consider any $\gamma\in (0,1)$. Take any $T>0$ and $J>0$. The large deviation upper and lower bounds match
\begin{equation*}
\limsup_{a\rightarrow +\infty}\limsup_{N\rightarrow +\infty}\frac{1}{a^2}\frac{\log N}{N^2}\log\mathds{P}_{\gamma,\SS}^N\Big[\langle \mathcal{Q}_T^N,\mathcal{V}\rangle\in [aJT,+\infty)\Big]= -\pi^{-1} TJ^2,
\end{equation*}

\end{corollary}

\subsection{Self-similarity of Large Deviation Costs}\label{sec:heuristic_explanation}

An important feature of the main results in Section \ref{sec:specific_main_result} is that, although as $N\rightarrow +\infty$ the density and current profiles become more and more singular around the origin, we can still derive $N$-independent large deviation upper and lower bounds. This is because the density-current profiles and the cost functionals have \textit{self-similarity} across all $N$. 

When we are evaluating the cost functional $\mathcal{C}^N$ in \eqref{eq:C_cost_functional}, we allow $N$ to be non-integer for convenience. Given $N\in (1,+\infty)$, define the logarithmic variable $s$ as
\begin{equation}\label{eq:log_variable_s}
s:=-\frac{\log |x|}{\log N},\quad |x|=N^{-s},
\end{equation}
which maps the interval $[N^{-1},1]$ of $|x|$ onto a fixed interval $[0,1]$ of $s$. 

\begin{remark}\label{rmk:vector_analysis}
Given any vector field $V_x$ on $\mathds{B}^N$ or $B_1$, it can be decomposed into the \underline{radial component} $V_{r}(x)\bm{r}_x$ and the \underline{rotation component} $V_{\theta}(x)\bm{\theta}_x$
\begin{equation}\label{eq:radial_rotation_decomposition}
V(x)=V_r(x)\bm{r}_x+V_{\theta}(x)\bm{\theta}_x,\quad \bm{r}_x:=\frac{x}{|x|},\quad \bm{\theta}_x:=\frac{(-x_2,x_1)}{|x|}.
\end{equation}
Sometimes we also use $V_x^{\textup{rot}}$ to denote the rotation component $V_{\theta}(x)\bm{\theta}_x$, and we say a vector field $V$ is \underline{'purely rotational'} if $V_{r}(x)\equiv 0$. The decomposition in \eqref{eq:radial_rotation_decomposition} can be applied to the current $q$, which gives the radial component $q_{r}$ and the rotation component $q_{\theta}$. We say a profile $(\rho,q)$ is \underline{radially symmetric}, if $(\rho,q_r,q_{\theta})$ depend only on $|x|$.
\end{remark}

Given a radially symmetric profile $(\rho,q)$, define $(\uprho,\mathtt{q}_r,\mathtt{q}_{\theta})$ as the functions of the variable $s$ in \eqref{eq:log_variable_s}
\begin{equation}\label{eq:self_similar_rho_q}
\mathtt{q}_{r}(s):=\Big(|x|\log N q_r(x)\Big)\Big|_{|x|=N^{-s}},\quad \mathtt{q}_{\theta}(s):=\Big(|x|\log N q_{\theta}(x)\Big)\Big|_{|x|=N^{-s}},\quad \uprho(s):=\rho(x)\big|_{|x|=N^{-s}}.
\end{equation}
In \eqref{eq:u_profile} we have introduced profiles $\mathfrak{u}^N[\gamma,\zeta;\alpha]$ and $\mathfrak{q}^N[J]$ as approximate optimal profiles for deviations of the winding number $\langle \mathcal{Q}_T^N,\mathcal{V}\rangle$. An observation is that, applying \eqref{eq:self_similar_rho_q} to the profiles $(\mathfrak{u}^N,\mathfrak{q}^N)$, they become \textit{$N$-independent} profiles (up to a vanishing error due to the boundary truncation $\varphi^N$ in \eqref{eq:u_profile}).

\begin{lemma}\label{lem:C_self_similarity}
For any $N\in (1,+\infty)$, any radially symmetric profiles $(\rho,q)$, the cost functional $\mathcal{C}^N(\rho,q)$ in \eqref{eq:C_cost_functional} can be expressed in an dimensional independent form
\begin{equation*}
\mathcal{C}^N(\rho,q)=\frac{1}{2}\int_0^12\pi\frac{|\mathtt{q}_r(s)|^2+|\mathtt{q}_{\theta}(s)|^2}{\sigma(\uprho(s))}ds+\frac{1}{8}\int_0^12\pi\Big(\frac{d}{ds}\uprho(s)\Big)^2\frac{1}{\sigma(\uprho(s))}ds=:\mathtt{C}(\uprho,\mathtt{q}),
\end{equation*}
with $(\mathtt{q}_r,\mathtt{q}_{\theta},\uprho)$ defined in \eqref{eq:self_similar_rho_q}.

\end{lemma}
This lemma proves that the cost function $\mathcal{C}(\gamma,\zeta,\alpha,J)$ in \eqref{eq:C_functional_parametrized} is a $N$-independent function, which gives a large deviation lower bound. In the analysis of large deviation upper in Section \ref{sec:proof_upper}, this self-similarity also plays an important role.

\begin{proof}
Re-write \eqref{eq:C_cost_functional} in polar coordinate using the radial symmetric of $(\rho,q)$, where $(\rho,q_r,q_\theta)$ can be seen as functions of $r=|x|$ (notice the subscript $r$ in $q_r$ refers to radial component, not $r=|x|$)
\begin{equation*}
\mathcal{C}^N(\rho,q)=\frac{1}{2}\int_{N^{-1}}^1\frac{|q_r(r)|^2+|q_{\theta}(r)|^2}{\sigma(\rho(r))}2\pi r\log Ndr+\frac{1}{8}\int_{N^{-1}}^1\Big(\frac{d}{dr}\rho(r)\Big)^2\frac{1}{\sigma(\rho(r))}2\pi r\log Ndr
\end{equation*}
The change of variable in \eqref{eq:self_similar_rho_q} implies 
\begin{equation*}
dr=-r\log Nds,\quad \frac{d}{dr}\rho(r)=\frac{ds}{dr}\frac{d}{ds}\rho(s)=-\frac{1}{r}\frac{1}{\log N}\frac{d}{ds}\rho(s).
\end{equation*}
This allows us to write
\begin{equation*}
\mathcal{C}^N(\rho,q)=\frac{1}{2}\int_0^12\pi\frac{|q_r(s)|^2+|q_{\theta}(s)|^2}{\sigma(\rho(s))}\big(r\log N\big)^2ds+\frac{1}{8}2\pi\int_0^1\Big(\frac{d}{ds}\rho(s)\Big)^2\frac{1}{\sigma(\rho(s))}ds,
\end{equation*}
which concludes the proof.
\end{proof}

\subsection{Proof Strategy: Multi-scale Methods and Cancellations}\label{sec:proof_strategy}

The notations used in this section have been summarized in Definition \ref{def:spaces}, where the bold font ($\mathbf{R}$ and $\mathbf{M}$ etc.) refers to discrete objects, and the calligraphy font ($\mathcal{C}^1$ and $\mathcal{M}$ etc.) refers to continuum objects.

\underline{\textit{Upper bound via multi-scale methods:}} 

The proof of the upper bound consists of
\begin{itemize}
    \item modifying the classical replacement lemma \cite[Lemma 5.10]{kipnis1998scaling} into a multi-scale decomposition, where the averaging radius depends on $|x|$, see Lemma \ref{lem:replacement_lemma}.
    \item applying this multi-scale replacement Lemma \ref{lem:replacement_lemma} to control the microscopic oscillations of $(\rho^N,\mathcal{Q}^N)$ at different scales, in Section \ref{sec:proof_upper}. This is heuristically explained in the beginning of Section \ref{sec:proof_upper}.
\end{itemize}

To derive the large deviation upper bound, we first prove the so-called replacement lemma. We begin by decomposing the domain $B_1$ into a family of layers $\{\mathcal{L}_i\}_{0\leq i\leq \lfloor \log_2N\rfloor}$. For $i\geq 0$, the layer $\mathcal{L}_i$ is defined by
\begin{equation}\label{eq:layer}
\mathcal{L}_i:=\Big\{x\in B_1\Big| 2^{-(i+1)}< |x|\leq  2^{-i}\Big\}.
\end{equation}
When $|x|\leq N^{-1}$ there is no new layers, and therefore the number of layers for each $N$ is of order $O(\log N)$. The microscopic averaging length $\mathfrak{l}_x$ is then defined as proportional to the scale of $\mathcal{L}_i$
\begin{equation}\label{eq:averaging_radius}
\mathfrak{l}_x:=N|x|_{\infty}.
\end{equation}
Accordingly, the averaged configuration is given by
\begin{equation}\label{eq:spatial_average}
\overline{\eta}_{x}^{\eps \mathfrak{l}}:=\frac{1}{(1+2\eps \mathfrak{l}_x)^2}\sum_{|y-x|_{\infty}\leq \eps \mathfrak{l}_xN^{-1}}\eta_{y}\mathds{1}_{y\in \mathds{B}^N},
\end{equation}
where by definition $x$ and $y$ are allowed to be in different layers. For $x\in \mathds{B}^N$ with $x\in \mathcal{L}_0$, the configuration $\eta$ may only be partially averaged. This is because close to the boundary $\p B_1$ we do not have enough room to average. However the contribution from $\mathcal{L}_0$ in the cost functionals is of vanishing order, because the singular vortex is localized around origin at smaller and smaller scales as $N\rightarrow +\infty$. Therefore the corresponding terms can be discarded.

Take $\psi$ as a local function, which will take the following forms with $e\in\{e_1,e_2\}$
\begin{equation}\label{eq:local_functions}
\begin{split}
&\psi(\eta)=\eta_0(1-\eta_{e}),\textup{ or }\psi(\eta)=\eta_{e}(1-\eta_0),\textup{ or }\psi(\eta)=\eta_0,\\ &\psi(\overline{\eta}^{\eps\mathfrak{l}})=\overline{\eta}_0^{\eps\mathfrak{l}}(1-\overline{\eta}_{0}^{\eps\mathfrak{l}}),\textup{ or }\psi(\overline{\eta}^{\eps\mathfrak{l}})=\overline{\eta}_0^{\eps\mathfrak{l}}.
\end{split}
\end{equation}

For arbitrary $x\in\mathds{B}^N$, the local function $\tau_x\psi$ is defined as
\begin{equation*}
\tau_x\psi(\eta):=\psi(\eta_{\cdot+x}),\quad \tau_x\psi\big(\overline{\eta}^{\eps\mathfrak{l}}\big):=\psi\big(\overline{\eta}_{\cdot+x}^{\eps\mathfrak{l}}\big).
\end{equation*}

For any parameters $(\delta>0,\eps>0,T>0,N>1)$ and function $H$, define
\begin{equation}\label{eq:replacement_error_set}
\mathbf{B}_{\delta,\eps,H}^N\big(\psi\big):=\Bigg\{\eta\in D\big([0,T];\mathbf{M}_{\{0,1\}}^N\big)\Bigg|\quad \int_0^T\Big|\frac{1}{N^2}\sum_{x\in\mathds{B}^N}H_x\Big((\tau_x\psi)\big(\overline{\eta}_{t}^{\eps \mathfrak{l}}\big)-(\tau_x\psi)(\eta_{t})\Big)\Big|dt\leq \delta\Bigg\}.
\end{equation}
We will consider test functions $H$ diverging close to the center of $B_1$. 

Given a function $h$ on $\mathds{B}^N$, define $\nabla_{\textup{disc}}h$ as
\begin{equation}\label{eq:discrete_gradient}
\nabla_{\textup{disc}}h_x=\Big(N(h_{x+e_1}-h_x),N(h_{x+e_2}-h_x)\Big).
\end{equation}
\begin{lemma}\label{lem:replacement_lemma}
Take arbitrary parameters $\big(\gamma\in(0,1),T>0,\mathfrak{b}>0,\delta>0,\eps>0\big)$. Consider families of functions $\{H^N\}_{N<+\infty}$ and $\{H_{\textup{lin}}^N\}_{N<+\infty}$ where each function is defined on $\mathds{B}^N$ with
\begin{equation}\label{eq:singular_condition_for_H}
\begin{split}
&|H_x^N|\leq \mathfrak{b}|x|^{-2}(\log N)^{-1},\qquad |\nabla_{\textup{disc}}H_x^N|\leq \mathfrak{b}|x|^{-3}(\log N)^{-1},\\
&|H_{\textup{lin},x}^N|\leq \mathfrak{b}|x|^{-2}\mathds{1}_{|x|\geq N^{-1+\eps}},\qquad |\nabla_{\textup{disc}}H_{\textup{lin},x}^N|\leq \mathfrak{b}|x|^{-3}\mathds{1}_{|x|\geq N^{-1+\eps}}.
\end{split}
\end{equation}
Then the configuration $\{\eta_x\}_{x\in\mathds{B}^N}$ is well approximated by the space-averaged configuration $\{\overline{\eta}_x^{\eps \mathfrak{l}}\}_{x\in\mathds{B}^N}$, in a sense that there exists a function $F_{\mathfrak{b},\delta,\eps}(N)$ determined by $(\mathfrak{b},\delta,\eps)$
\begin{equation}\label{eq:replacement_lemma}
\begin{split}
&\frac{\log N}{N^2}\log\mathds{P}_{\gamma,\SS}^N\Big[\mathbf{B}_{\delta,\eps,H^N}^{N,\textup{c}}\Big(\eta_{x}(1-\eta_{x+e})\Big)\cup\mathbf{B}_{\delta,\eps,H^N}^{N,\textup{c}}\Big(\eta_{x+e}(1-\eta_{x})\Big) \cup \mathbf{B}_{\delta,\eps,H_{\textup{lin}}^N}^{N,\textup{c}}\big(\eta_{x}\big)\Big]\leq F_{\mathfrak{b},\delta,\eps}(N),\\
&\hspace{10cm}\textup{with}\quad \lim_{\eps\rightarrow 0}\lim_{N\rightarrow +\infty}F_{\mathfrak{b},\delta,\eps}(N)=-\infty.
\end{split}
\end{equation}
\end{lemma}
The averaging radius $\mathfrak{l}_x$ in \eqref{eq:averaging_radius} is motivated by the singular behavior of $ H_x^N $ near the origin. When $ |x| $ is small, finer control over the averaging error is required. This ensures that the averaging is performed over a correspondingly small neighborhood, which allows for a more precise estimate of the averaging error and compensates for the singularity of $ H_x^N $ near the origin.

This replacement lemma allows us to describe the oscillation of $\eta_x$ at different scales via a multi-scale method developed in Section \ref{sec:proof_upper}. Together they conclude the proof of the large deviation upper bound, as shown in Section \ref{sec:true_proof_upper}.

\underline{\textit{Lower bound, relative entropy, cancellation, and concentration:}}

In the proof of the large deviation lower bound, to consider the probability of $\{\langle \mathcal{Q}_T^N,\mathcal{V} \rangle\geq JT\}$, we tilt the measure
\begin{equation*}
\mathds{E}_{\gamma,\SS}^N\Big[\mathds{1}_{\langle \mathcal{Q}_T^N,\mathcal{V} \rangle\geq JT}\Big]=\mathds{E}_{\rho,\mathds{F}}^N\Bigg[\frac{d\mathds{P}_{\gamma,\SS}^N}{d\mathds{P}_{\rho,\mathds{F}}^N}\mathds{1}_{\langle \mathcal{Q}_T^N,\mathcal{V} \rangle\geq JT}\Bigg].
\end{equation*}
Here $\rho$ is an initial density profile possibly different from $\gamma$, and $\mathds{F}$ is a forcing field. To derive the lower bound, we need to show that the Radon-Nikodym derivative concentrates to a certain limit as we take $N\rightarrow +\infty$. Thus a proper description of the measure $\mathds{P}_{\rho,\mathds{F}}^N$ is necessary, and this description will be given by the control of the relative entropy. In Proposition \ref{prop:relative_entropy_vague} we give a statement of this control.
\begin{proposition}\label{prop:relative_entropy_vague}
Take arbitrary parameters $(\gamma,\zeta,\alpha,J)$ with 
\begin{equation*}
\gamma\in (0,1),\quad \zeta\in (0,1),\quad \alpha\in (0,1),\quad J\in\mathds{R}.
\end{equation*}
Also take $\rho$ to represent $\mathfrak{u}^N[\gamma,\zeta;\alpha]$ in \eqref{eq:u_profile}, and $\mathds{F}$ as the forcing field $\mathds{F}^N[\gamma,\zeta,\alpha,J]$ in \eqref{eq:relative_entropy_inproof}. For the law $d\mathds{P}_{\rho,\mathds{F}}^N$ of the modified dynamics \eqref{eq:modified_generator} with generator $L_{\mathds{F}}^{N}$, its relative entropy satisfies for any time $t\geq 0$
\begin{equation*}
\mathds{H}^N\Big(d\mathds{P}_{\rho,\mathds{F},t}^N\Big|d\local{\rho}\Big)\leq C_{\gamma,\zeta,\alpha,J,t}O\Big(N^{4\alpha}\Big).
\end{equation*}
\end{proposition}

With the entropy bound given by Proposition \ref{prop:relative_entropy_vague}, we can prove the concentration estimate, stated in later in Proposition \ref{prop:concentration}. These allow us to establish the large deviation lower bound in Theorem \ref{th:one-scale-deviation}.

\section{Proof of the Multi-scale Replacement Lemma \ref{lem:replacement_lemma}}\label{sec:multi-scale_replacement}

In this section, we give the proof of the super-exponential replacement Lemma \ref{lem:replacement_lemma} to replace quadratic-in-$\eta$ terms (like $\eta_{x}(1-\eta_{x+e})$) and linear-in-$\eta$ terms by their averages. In this section, we step by step localize the estimate of nonlinear terms until we can apply the classical methods in \cite{kipnis1989hydrodynamics,kipnis1998scaling}. The estimate of the linear term is the same, and we give a remark on it at the end of this section.

To complete the proof of the replacement lemma \ref{lem:replacement_lemma}, we will give the details on the proof of
\begin{equation}\label{eq:replacement_goal}
\lim_{\eps\rightarrow 0}\lim_{N\rightarrow +\infty}\frac{\log N}{N^2}\log \mathds{P}_{\gamma,\SS}^N\Big[\int_0^T\mathbf{V}_{\eps}^{N}(\eta_{s})ds>\delta\Big]=-\infty,\qquad\forall \delta>0,
\end{equation}
where we choose any $e\in\{e_1,e_2\}$, and the term $\mathbf{V}_{\eps}^N$ is defined as
\begin{equation*}
\mathbf{V}_{\eps}^{N}(\eta):=\frac{1}{N^2}\Bigg|\sum_{x\in\mathds{B}^N}H_x^N\Big(\sigma(\bar{\eta}_x^{\eps \mathfrak{l}})-\eta_x(1-\eta_{x+e})\Big)\Bigg|.
\end{equation*}
The control of replacement errors involving $\eta_{x+e}(1-\eta_{x})$ or $\eta_x$ follows the same method.

The main difficulty of proving \eqref{eq:replacement_goal} relative to the standard replacement lemma arises from the singular behavior of the coefficient $H_x^N$ around the center, see \eqref{eq:singular_condition_for_H}. 

\textbf{Step 1: Decompose the error into different layers.}

To control this singularity, the key ingredient is the layer decomposition introduced in Equation \eqref{eq:layer}. The domain $\mathds{B}^N$ is written as a disjoint union
\begin{equation*}
\mathds{B}^N=\bigcup_{0\leq i\leq \lfloor \log_2N\rfloor}\Big(\mathcal{L}_i\cap \mathds{B}^N\Big),
\end{equation*}
where the size of the layers decays exponentially with respect to the index $i$, see \eqref{eq:layer}. The total replacement error $\mathbf{V}_{\eps}^{N}$ can be decomposed into the summation of replacement errors within each layer $\mathcal{L}_i$. Recall that in the boundary layer $\mathcal{L}_0$ we do not average, as stated in \eqref{eq:averaging_radius}
\begin{equation}\label{eq:layer_error}
\mathbf{V}_{\eps}^{N}(\eta)\leq \frac{1}{\log N}\sum_{1\leq i\leq \lfloor \log_2N\rfloor}\mathbf{V}_{i,\eps}^{N},\qquad \mathbf{V}_{i,\eps}^{N}(\eta):=\frac{\log N}{N^2}\Bigg|\sum_{x\in\mathcal{L}_i\cap \mathds{B}^N}H_x^N\Big(\sigma(\bar{\eta}_x^{\eps \mathfrak{l}})-\eta_x(1-\eta_{x+e})\Big)\Bigg|.
\end{equation}
In the definition of $\mathbf{V}_{i,\eps}^{N}$, we added a factor $\log N$ whose role will be clear later in \eqref{eq:rescale_H}. For each layer $i$, the normalizing factor $\frac{1}{N^2}$ controls the singularity of $H_x^N$. The singularity condition \eqref{eq:singular_condition_for_H} for $H_x^N$ and the definition of $\mathcal{L}_i$ in \eqref{eq:layer} imply
\begin{equation}\label{eq:rescale_H}
\big|H_x^N\big|\leq |x|^{-2}\frac{\mathfrak{b}}{\log N}\leq \mathfrak{b}C2^{2i}\frac{1}{\log N},\quad\textup{ and futher imply}\quad \log N\Big|H_x^N(2^{-i})^2\Big|\leq \mathfrak{b}C,\quad \textup{for any $x\in\mathcal{L}_{i}\cap \mathds{B}^N$}.
\end{equation}
Thus the error $\mathbf{V}_{i,\eps}^{N}$ in \eqref{eq:layer_error} is uniformly bounded and it also equals to
\begin{equation}\label{eq:layer_error_normalized}
\begin{split}
\mathbf{V}_{i,\eps}^N:=\frac{1}{(2^{-i}N)^2}\Bigg|\sum_{x\in\mathcal{L}_i\cap \mathds{B}^N}h_x^N \Big(\sigma(\bar{\eta}_{x}^{\eps\mathfrak{l}})-\eta_x(1-\eta_{x+e})\Big)\Bigg|,\quad h_x^N:=\log N (2^{-i})^2 H_x^N,
\end{split}
\end{equation}
where the function $h_x^N$ is bounded from above by $\mathfrak{b}C$ according to \eqref{eq:rescale_H}.

\textbf{Step 2: Apply the Feynman-Kac formula.}

The decomposition in \eqref{eq:layer_error} implies
\begin{equation*}
\begin{split}
\log \mathds{P}_{\gamma,\SS}^N\Bigg[\int_0^T\mathbf{V}_{\eps}^{N}(\eta_{s})ds>\delta\Bigg]&\leq \log \mathds{P}_{\gamma,\SS}^N\Bigg[\sum_{i=1}^{\lfloor \log_2N\rfloor}\int_0^T\mathbf{V}_{i,\eps}^{N}(\eta_{s})ds>\delta\log N\Bigg]\\
&=\log \mathds{P}_{\gamma,\SS}^N\Bigg[\exp\Bigg(\frac{N^2}{\log N}\sum_{i=1}^{\lfloor \log_2N\rfloor}\int_0^T\mathbf{V}_{i,\eps}^{N}(\eta_{s})ds\Bigg)>\exp\Big(N^2\delta\Big)\Bigg].
\end{split}
\end{equation*}
Thus we use the Markov inequality and the Feynman-Kac formula in Lemma \ref{lem:feynman_kac_variational_app}
\begin{equation*}
\begin{split}
\frac{\log N}{N^2}\log \mathds{P}_{\gamma,\SS}^N\Big[\int_0^T\mathbf{V}_{\eps}^{N}(\eta_{s})ds>\delta\Big]&\leq \frac{\log N}{N^2}\log \mathds{E}_{\gamma,\SS}^N\Bigg[\exp\Bigg(\frac{N^2}{\log N}\sum_{i=1}^{\lfloor \log_2N\rfloor}\int_0^T\mathbf{V}_{i,\eps}^{N}(\eta_{s})ds\Bigg)\Bigg]-\delta\log N\\
&\leq T\sup_{f}\Bigg\{\local{\gamma}\Big[f\sum_{i=1}^{\lfloor \log_2N\rfloor}\mathbf{V}_{i,\eps}^{N}\Big]-\frac{\log N}{2N^2}\local{\gamma}\big[\Gamma_{\SS}^N(\sqrt{f})\big]\Bigg\}-\delta \log N,
\end{split}
\end{equation*}
where $f$ is taken over all density functions with respect to $\local{\gamma}$, and $\Gamma_{\SS}^N(\sqrt{f})$ is the Dirichlet form given as 
\begin{equation}\label{eq:full_dirichley}
\Gamma_{\SS}^N(\sqrt{f}):=\frac{N^2}{2}\sum_{x\in\mathds{B}^N}\Big(\sqrt{f(\eta^{x,x+e_1})}-\sqrt{f(\eta)}\Big)^2+\frac{N^2}{2}\sum_{x\in\mathds{B}^N}\Big(\sqrt{f(\eta^{x,x+e_2})}-\sqrt{f(\eta)}\Big)^2.
\end{equation}

To prove \eqref{eq:replacement_goal}, it is sufficient to consider only $1\leq i\leq \lfloor \log_2N\rfloor-\log\log N$, because the layers with $\lfloor \log_2N\rfloor-\log\log N\leq i\leq \lfloor \log_2N\rfloor$ only contributes a negligible order $\log\log N$. In the rest of the proof, we will neglect certain issues involving $i\geq \lfloor \log_2N\rfloor-\log\log N$.

Decompose the Dirichlet form $\Gamma_{\SS}^N(\sqrt{f})$ into different layers
\begin{equation*}
\local{\gamma}\big[\Gamma_{\SS}^N(\sqrt{f})\big]\geq \frac{1}{3}\sum_{1\leq i\leq \lfloor \log_2N\rfloor}\local{\gamma}\big[\Gamma_i^N(\sqrt{f})\big],
\end{equation*}
with $\Gamma_i^N(\sqrt{f})$ defined as
\begin{equation}\label{eq:Gamma_i}
\Gamma_i^N(\sqrt{f}):=\frac{N^2}{2}\sum_{x\in \mathcal{L}_i^{\sharp}}\Big(\sqrt{f(\eta^{x,x+e_1})}-\sqrt{f(\eta)}\Big)^2+\frac{N^2}{2}\sum_{x\in \mathcal{L}_i^{\sharp}}\Big(\sqrt{f(\eta^{x,x+e_2})}-\sqrt{f(\eta)}\Big)^2,\quad \mathcal{L}_i^{\sharp}:=\mathcal{L}_{i-1}\cup \mathcal{L}_i\cup \mathcal{L}_{i+1}.
\end{equation}
Heuristically, we will use $\Gamma_i^N$ to control the replacement error for $x\in \mathcal{L}_i\cap \mathds{B}^N$. The union of layers in the definition \eqref{eq:Gamma_i} of $\mathcal{L}_i^{\sharp}$ is because that, for $x\in \mathcal{L}_i\cap \mathds{B}^N$ its averaging radius may exceed layer $\mathcal{L}_i$. Now we derive
\begin{equation}\label{eq:replacement_event_into_layers}
\frac{\log N}{N^2}\log\mathds{P}_{\gamma,\SS}^N\Bigg[\int_0^T \mathbf{V}_{\eps}^{N}(\eta_s)ds>\delta\Bigg]\leq T\lfloor \log_2N\rfloor\sup_{1\leq i\leq \lfloor \log_2N\rfloor}\sup_{f_i}\Bigg\{\local{\gamma,i}\Big[\mathbf{V}_{i,\eps}^{N}f_i\Big]-\frac{\log N}{6N^2}\local{\gamma}\big[\Gamma_i^N(\sqrt{f_i})\big]\Bigg\}-\delta\log N,
\end{equation}
where $\local{\gamma,i}$ is the projection of $\local{\gamma}$ onto $\mathcal{L}_{i-1}\cup \mathcal{L}_{i}\cup \mathcal{L}_{i+1}$, and $f_i$ is taken over density functions with respect to $\local{\gamma,i}$. 

In the following proof, we continue to work with $1\leq i\leq \lfloor \log_2N\rfloor-\log\log N$. 

\textbf{Step 3: One-block and two-block decomposition}

The function $\mathbf{V}_{i,\eps}^{N}$ can be decomposed into the one-block error $\mathbf{V}_{i,\eps}^N[1]$ and the two-block error $\mathbf{V}_{i,\eps}^N[2]$. Choose a positive integer $k$, and define $\overline{\eta}_x^{k}$ as
\begin{equation}\label{eq:k_average}
\overline{\eta}_x^k:=\frac{1}{(1+2k)^2}\sum_{|y|_{\infty}\leq kN^{-1}}\eta_{x+y}.
\end{equation}
Recalling the function $h_x^N$ in \eqref{eq:layer_error_normalized}, the error $\sigma(\overline{\eta}_x^{\eps\mathfrak{l}})-\eta_x(1-\eta_{x+e})$ can be decomposed into $\sigma(\overline{\eta}_x^{\eps\mathfrak{l}})-\sigma(\overline{\eta}_x^{k})$ and $\sigma(\overline{\eta}_x^k)-\eta_x(1-\eta_{x+e})$
\begin{equation*}
\begin{split}
&\mathbf{V}_{i,\eps}^N\leq \underbrace{\Bigg|\frac{1}{(2^{-i}N)^2}\sum_{x\in \mathcal{L}_i\cap \mathds{B}^N}h_x^N\Big(\sigma(\overline{\eta}_x^{\eps\mathfrak{l}})-\sigma(\overline{\eta}_x^k)\Big)\Bigg|}_{\textup{\blue{error between averages at large and small scales}}}+\underbrace{\Bigg|\frac{1}{(2^{-i}N)^2}\sum_{x\in \mathcal{L}_i\cap \mathds{B}^N}h_x^N\Big(\sigma(\overline{\eta}_x^{k})-\eta_x(1-\eta_{x+e})\Big)\Bigg|}_{\textup{\blue{error between average at small scale and original configuration}}}.
\end{split}
\end{equation*}
By classical methods, for example see \cite[Proof of Theorem 2.1]{kipnis1989hydrodynamics}, this leads to the one-block and two-block decomposition
\begin{equation}\label{eq:one_two_block_errors}
\begin{split}
& \mathbf{V}_{i,\eps}^N\leq \mathbf{V}_{i,\eps}^N[1]+\mathbf{V}_{i,\eps}^N[2]+\mathfrak{b}C_{\eps}O\Big(\frac{k 2^i}{N}\Big)+C_{\mathfrak{b}}O(\eps),\\
&\mathbf{V}_{i,\eps}^N[1]:=\frac{\mathfrak{b}}{(2^{-i}N)^2}\sum_{x\in \mathcal{L}_i\cap\mathds{B}^N}\Bigg|\frac{1}{(1+2k)^2}\sum_{|y|_{\infty}\leq kN^{-1}}\eta_{x+y}(1-\eta_{x+y+e})-\sigma\Big(\sum_{|y|_{\infty}\leq kN^{-1}}\frac{1}{(1+2k)^2}\eta_{x+y}\Big)\Bigg|,\\
&\mathbf{V}_{i,\eps}^N[2]:=\frac{\mathfrak{b}C_{\eps}}{(2^{-i}N)^4}\sum_{x,x'\in \mathcal{L}_i\cap \mathds{B}^N}\Bigg|\frac{1}{(1+2k)^2}\sum_{|y|_{\infty}\leq kN^{-1}}\eta_{x+y}-\frac{1}{(1+2k)^2}\sum_{|y'|_{\infty}\leq kN^{-1}}\eta_{x'+y'}\Bigg|.
\end{split}
\end{equation}
The error term $\mathfrak{b}C_{\eps}O((k2^i)/N)+C_{\mathfrak{b}}O(\eps)$ is negligible as we take $N\rightarrow +\infty$ and $\eps\rightarrow 0$, therefore we omit it in the following computations. For each layer $\mathcal{L}_i$, reorganize the Dirichlet form and the block errors
\begin{equation}\label{eq:layer_one_and_two}
\begin{split}
&\sup_{f_i}\Bigg\{\local{\gamma,i}\Big[\mathbf{V}_{i,\eps}^Nf_i\Big]-\frac{1}{6}\local{\gamma}\Big[\frac{\log N}{N^2}\Gamma_{i}^N(\sqrt{f_i})\Big]\Bigg\}\\
&\leq \sup_{f_i}\Bigg\{\local{\gamma,i}\Big[\mathbf{V}_{i,\eps}^N[1]f_i\Big]+\local{\gamma,i}\Big[\mathbf{V}_{i,\eps}^N[2]f_i\Big]-\frac{1}{6}\local{\gamma}\Big[\frac{\log N}{N^2}\Gamma_{i}^N(\sqrt{f_i})\Big]\Bigg\}\\
&\leq \underbrace{\sup_{f_i}\Bigg\{\local{\gamma,i}\Big[\mathbf{V}_{i,\eps}^N[1]f_i\Big]-\frac{1}{12}\local{\gamma}\Big[\frac{\log N}{N^2}\Gamma_{i}^N(\sqrt{f_i})\Big]\Bigg\}}_{\textup{\blue{one-block estimate}}}+\underbrace{\sup_{f_i}\Bigg\{\local{\gamma,i}\Big[\mathbf{V}_{i,\eps}^N[2]f_i\Big]-\frac{1}{12}\local{\gamma}\Big[\frac{\log N}{N^2}\Gamma_{i}^N(\sqrt{f_i})\Big]\Bigg\}}_{\textup{\blue{two-block estimate}}}.
\end{split}
\end{equation}

\textbf{Step 4: Classical methods in \cite{kipnis1998scaling,kipnis1989hydrodynamics} for each layer $\mathcal{L}_i$.}

Applying Equation \eqref{eq:layer_one_and_two} to bound \eqref{eq:replacement_event_into_layers}, we derive
\begin{equation}\label{eq:replacement_proof_final}
\begin{split}
\frac{\log N}{N^2}\log\mathds{P}_{\gamma,\SS}^N&\Bigg[\int_0^T \mathbf{V}_{\eps}^{N}(\eta_s)ds>\delta\Bigg]\\
&\leq T\lfloor \log_2N\rfloor\sup_{1\leq i\leq \lfloor \log_2N\rfloor-\log\log N}\sup_{f_i}\Bigg\{\local{\gamma,i}\Big[\mathbf{V}_{i,\eps}^N[1]f_i\Big]-\frac{1}{12}\local{\gamma}\Big[\frac{\log N}{N^2}\Gamma_{i}^N(\sqrt{f_i})\Big]\Bigg\}\\
&+ T\lfloor \log_2N\rfloor\sup_{1\leq i\leq \lfloor \log_2N\rfloor-\log\log N}\sup_{f_i}\Bigg\{\local{\gamma,i}\Big[\mathbf{V}_{i,\eps}^N[2]f_i\Big]-\frac{1}{12}\local{\gamma}\Big[\frac{\log N}{N^2}\Gamma_{i}^N(\sqrt{f_i})\Big]\Bigg\}-\delta \log N.
\end{split}
\end{equation}
Since we only consider layer $\mathcal{L}_i$ with $1\leq i\leq \lfloor \log_2N\rfloor-\log \log N$, there is always $|\mathcal{L}_i\cap \mathds{B}^N|\geq O\big((\log N)^{2\log 2}\big)$. This means the sizes of the layers we consider go to $+\infty$ uniformly as $N\rightarrow +\infty$.

In summary, we are considering layers with sizes growing to $+\infty$, and function $h_x^N$ is uniformly bounded by the rescaling in \eqref{eq:layer_error_normalized}, and the one-block and two-block errors within a layer have the same form as in classical contexts. In addition, we have extra logarithmic term here in $\frac{\log N}{N^2}\Gamma_i^N(\sqrt{f})$ compared to just $\frac{1}{N^2}\Gamma_i^N(\sqrt{f})$. At this stage, we can apply the methods in \cite{kipnis1998scaling,kipnis1989hydrodynamics} to control the one-block and two-block terms by decomposing Dirichlet forms. This derives
\begin{equation*}
\limsup_{k\rightarrow +\infty}\limsup_{N\rightarrow +\infty}\sup_{1\leq i\leq \lfloor \log_2N\rfloor-\log\log N}\sup_{f_i}\Bigg\{\local{\gamma,i}\Big[\mathbf{V}_{i,\eps}^N[1]f_i\Big]-\frac{1}{12}\local{\gamma}\Big[\frac{\log N}{N^2}\Gamma_{i}^N(\sqrt{f_i})\Big]\Bigg\}\leq 0,
\end{equation*}
and the same inequality for $\mathbf{V}_{i,\eps}^N[2]$. This concludes the estimate of the nonlinear-in-$\eta$ terms. For the linear-in-$\eta$ term, the steps 1-3 can be performed in the same manner. At step 4, the test function $H_{\textup{lin}}^N$ has an extra $\log N$ than $H^N$ in \eqref{eq:singular_condition_for_H}, but this factor can be canceled by the extra $\log N$ factor in \eqref{eq:replacement_proof_final}. Therefore the classical methods still apply. Lemma \ref{lem:replacement_lemma} is proved.

\section{Large Deviation Upper Bound Via a multi-scale Argument}\label{sec:proof_upper}

Before giving a complete proof of the large deviation upper bound in Theorem \ref{th:one-scale-deviation}, we first present a \textit{heuristic} overview of the classical approach in order to emphasize the difficulties occuring from the singularity of the density and current profiles. The equations and notations below are schematic and are used solely to illustrate the main idea.

\underline{\textit{Classical methods for the large deviation upper bound.}}

In classical settings (for example \cite{kipnis1998scaling} Appendix 2.3), to derive the large deviation upper bound for a possibly non-compact closed set $\mathcal{S}$, one typically proceeds as follows.
One first prove exponential tightness by constructing a compact set $\mathcal{K}$ with
\begin{equation}\label{eq:compact_reduction}
\mathds{P}^N\big[\mathcal{S}\setminus\mathcal{K}\big]\ \textup{is exponentially small in an appropriate sense}.
\end{equation}
As a consequence, it suffices to establish the large deviation upper bound within $\mathcal{K}$. To this end, one introduces a family of open sets $\{\mathcal{O}_i\}_i$ forming an open cover of the compact set $\mathcal{K}$ with $\mathcal{K}\subset \bigcup_i\mathcal{O}_i$. In each open set $\mathcal{O}_i$, one proves an estimate of the following form, with $\mathcal{R}_{T,\mathcal{F}}(\rho,\mathcal{Q})$ as a cost functional
\begin{equation}\label{eq:local_upper_bound}
\limsup_{N\rightarrow +\infty}\frac{1}{N^2}\log\mathds{P}^N\Big[(\rho^N,\mathcal{Q}^N)\in\mathcal{O}_i\Big]
\leq
-\sup_{\mathcal{F}}\inf_{(\rho,\mathcal{Q})\in\mathcal{O}_i}\mathcal{R}_{T,\mathcal{F}}\Big(\rho,\mathcal{Q}\Big)
+\textup{small error}.
\end{equation}

The open cover $\{\mathcal{O}_i\}_i$ is chosen in such a way that, the oscillation of the cost functional $\mathcal{R}_{T,\mathcal{F}}(\rho,\mathcal{Q})$ is small enough within each $\mathcal{O}_i$. In this case, the order of the supremum and infimum can be exchanged up to a controllable error, namely
\begin{equation}\label{eq:minmax_exchange}
-\sup_{\mathcal{F}}\inf_{(\rho,\mathcal{Q})\in\mathcal{O}_i}\mathcal{R}_{T,\mathcal{F}}\Big(\rho,\mathcal{Q}\Big)\leq-\inf_{(\rho,\mathcal{Q})\in\mathcal{O}_i}\sup_{\mathcal{F}}\mathcal{R}_{T,\mathcal{F}}\Big(\rho,\mathcal{Q}\Big)
+\textup{small error}.
\end{equation}

Combining \eqref{eq:compact_reduction}, \eqref{eq:local_upper_bound}, and \eqref{eq:minmax_exchange}, one arrives at
\begin{equation*}
\begin{split}
\limsup_{N\rightarrow +\infty}\frac{1}{N^2}\log\mathds{P}^N\Big[(\rho^N,\mathcal{Q}^N)\in \mathcal{S}\Big]
&\leq
\limsup_{N\rightarrow +\infty}\frac{1}{N^2}\log\mathds{P}^N\Big[(\rho^N,\mathcal{Q}^N)\in \mathcal{K}\cap \mathcal{S}\Big]
+\textup{small error}\\
&\leq \max_{i}\limsup_{N\rightarrow +\infty}\frac{1}{N^2}\log\mathds{P}^N\Big[(\rho^N,\mathcal{Q}^N)\in \mathcal{O}_i\Big] +\textup{small error} \\ &\leq -\inf_{(\rho,\mathcal{Q})\in \mathcal{K}\cap \mathcal{S}}\sup_{\mathcal{F}}\mathcal{R}_{T,\mathcal{F}}(\rho,\mathcal{Q}) +\textup{small error}. 
\end{split} 
\end{equation*} 
Letting $\mathcal{K}$ approach $\mathcal{S}$ and the open cover $\{\mathcal{O}_i\}_i$ become increasingly accurate, the error terms in the above bounds vanish, leading to the desired large deviation upper bound.

\underline{\textit{Multi-scale method in the singular large deviation.}}

In the large deviation of winding numbers, for particles around the origin, even a small rotational motion can produce a significant contribution to the winding number. As a result, particles exhibit a strong tendency to concentrate around the origin and to rotate there with high current, as shown for example in \eqref{eq:u_profile}. This phenomenon makes it impossible to perform the classical methods of proving large deviation upper bound, where the compact set $\mathcal{K}$ and the open cover $\{\mathcal{O}_i\}_i$ are chosen independent of $N$. That $N$-independent construction of the open cover is not able to capture the behavior of $(\rho^N,\mathcal{Q}^N)$ at finer scales. Instead, we need to modify the classical method in a \textit{multi-scale manner}, where we can capture the behavior of the profiles $(\rho^N,\mathcal{Q}^N)$ at every scale, for example $N^{-\alpha}$ with $\alpha\in (0,1)$. Specifically, this requires us to construct the open cover $\{\mathcal{O}_i\}_i$ in an $N$-dependent manner. 

To this end, the spatial average $\overline{\eta}_x^{\eps\mathfrak{l}}$ has been defined in \eqref{eq:spatial_average}, capturing the behavior of $\eta$ at different scales. The replacement Lemma \ref{lem:replacement_lemma} is an example of this idea.
\begin{lemma}\label{lem:regularity_of_average}
The space average $\overline{\eta}_x^{\eps\mathfrak{l}}$ in \eqref{eq:spatial_average} has the following regularity: for arbitrary shift $\boldsymbol{\delta}\in\mathds{R}^2$ with $|\boldsymbol{\delta}|_{\infty}\leq N^{-1}$
\begin{equation}\label{eq:discrete_as_continuous_gradient}
\Big|\overline{\eta}_{x}^{\eps\mathfrak{l}}-\overline{\eta}_{x+\boldsymbol{\delta}}^{\eps \mathfrak{l}}\Big|\leq C\eps^{-1}\frac{1}{N|x|}.
\end{equation}  
\end{lemma}

\begin{proof}
From $\overline{\eta}_{x}^{\eps\mathfrak{l}}$ to $\overline{\eta}_{x+\boldsymbol{\delta}}^{\eps\mathfrak{l}}$ according to the definition \eqref{eq:spatial_average} of the space average, the averaging center (from $x$ to $x+\boldsymbol{\delta}$) and the averaging radius (from $\mathfrak{l}_x$ to $\mathfrak{l}_{x+\boldsymbol{\delta}}$) both change. Therefore the difference $\overline{\eta}_x^{\eps\mathfrak{l}}-\overline{\eta}_{x+\boldsymbol{\delta}}^{\eps\mathfrak{l}}$ can be decomposed into the error of changing averaging center and the error of changing averaging radius
\begin{equation}\label{eq:spt_average_error_decompose}
\begin{split}
\Big|\overline{\eta}_{x}^{\eps\mathfrak{l}}-\overline{\eta}_{x+\boldsymbol{\delta}}^{\eps\mathfrak{l}}\Big|&\leq \Bigg|\underbrace{\sum_{|y-x|_{\infty}\leq \eps \mathfrak{l}_xN^{-1}}\frac{1}{(1+2\eps \mathfrak{l}_x)^2}\eta_{y}}_{\blue{\textup{average at $x$ with radius $\mathfrak{l}_x$}}}-\underbrace{\sum_{|y'-(x+\boldsymbol{\delta})|_{\infty}\leq \eps \mathfrak{l}_xN^{-1}}\frac{1}{(1+2\eps \mathfrak{l}_x)^2}\eta_{y'}}_{\textup{\blue{average at $x+\boldsymbol{\delta}$ with radius $\mathfrak{l}_x$}}}\Bigg|\qquad\qquad\quad\textup{\blue{(Term I)}}\\
&\quad+\Bigg|\underbrace{\sum_{|y-(x+\boldsymbol{\delta})|_{\infty}\leq \eps \mathfrak{l}_xN^{-1}}\frac{1}{(1+2\eps \mathfrak{l}_{x})^2}\eta_{y}}_{\textup{\blue{average at $x+\boldsymbol{\delta}$ with radius $\mathfrak{l}_x$}}}-\underbrace{\sum_{|y'-(x+\boldsymbol{\delta})|_{\infty}\leq \eps \mathfrak{l}_{x+\boldsymbol{\delta}}N^{-1}}\frac{1}{(1+2\eps \mathfrak{l}_{x+\boldsymbol{\delta}})^2}\eta_{y'}}_{\textup{\blue{average at $x+\boldsymbol{\delta}$ with radius $\mathfrak{l}_{x+\boldsymbol{\delta}}$}}}\Bigg|\quad\textup{\blue{(Term II)}}
\end{split}
\end{equation}
Substracting the common factor $(1+2\eps\mathfrak{l}_x)^{-2}$, the term (I) is less than
\begin{equation*}
\leq \frac{1}{(1+2\eps \mathfrak{l}_x)^2}\Bigg|\sum_{|y-x|_{\infty}\leq \eps \mathfrak{l}_x N^{-1}}\eta_{y}-\sum_{|y'-(x+\boldsymbol{\delta})|_{\infty}\leq \eps \mathfrak{l}_x N^{-1}}\eta_{y'}\Bigg|.
\end{equation*}
The set $\{y\in \mathds{B}^N|\ |y-x|_{\infty}\leq \eps \mathfrak{l}_xN^{-1}\}$ and $\{y'\in \mathds{B}^N|\ |y'-(x+\boldsymbol{\delta})|_{\infty}\leq \eps \mathfrak{l}_xN^{-1}\}$ only differ by $C\eps \mathfrak{l}_x$ points. The fact that $\mathfrak{l}_x=N|x|_{\infty}$ shows the term (I) in \eqref{eq:spt_average_error_decompose} is further less than
\begin{equation*}
\leq \frac{C}{(1+2\eps \mathfrak{l}_x)^2}\eps\mathfrak{l}_x\leq C\eps^{-1}\mathfrak{l}_x^{-1}\leq C\eps^{-1}(N|x|)^{-1}.
\end{equation*}
The same upper bound is true for the term (II) in \eqref{eq:spt_average_error_decompose}, since $|\eps\mathfrak{l}_{x+\boldsymbol{\delta}}-\eps\mathfrak{l}_x|\leq 1$ for any $x\in B_1$, according to the definition \eqref{eq:averaging_radius} of $\mathfrak{l}_x:=N|x|_{\infty}$ and the fact that $|\boldsymbol{\delta}|_{\infty}\leq N^{-1}$. In this case $|\eps\mathfrak{l}_{x+\boldsymbol{\delta}}-\eps\mathfrak{l}_x|\leq N|\delta|_{\infty}\leq 1$. This concludes the proof of Equation \eqref{eq:discrete_as_continuous_gradient}.
\end{proof}

Here we further introduce an averaged configuration, where the average is performed in a multi-scale manner same as for $\overline{\eta}_x^{\eps\mathfrak{l}}$.
\begin{definition}
Take arbitrary $T>0$, $N>1$, and $\eps>0$. For any $x\in B_1$, define $[\rho_x^N]_{\eps}$ as the space-time averaged density associated with $\eta\in D\big([0,T];\mathbf{M}_{\{0,1\}}^N\big)$
\begin{equation}\label{eq:space_time_average}
[\rho_x^N]_{\eps}:=\frac{1}{T}\int_0^T\overline{\eta}_{t,x}^{\eps\mathfrak{l}}dt.
\end{equation}
\end{definition}
Also, we introduce some cost functionals
\begin{definition}\label{def:R_functionals}
For any $x\in\mathds{B}^N$ and $e\in\{e_1,e_2\}$, define
\begin{equation}\label{eq:sigmae}
\sigma_{x,e}(\eta):=\frac{1}{2}\Big(\eta_{x}(1-\eta_{x+e})+\eta_{x+e}(1-\eta_{x})\Big).
\end{equation}
and $\nabla_{\textup{disc}}\cdot \mathbf{F}$ as the discrete divergence operator
\begin{equation}\label{eq:discrete_divergence}
\nabla_{\textup{disc}}\cdot \mathbf{F}:=N\Big(\mathbf{F}_{x,e_1}+\mathbf{F}_{x,e_2}-\mathbf{F}_{x-e_1,e_1}-\mathbf{F}_{x-e_2,e_2}\Big).
\end{equation}
Recall the inner product $\langle \cdot,\cdot\rangle$ is defined in Definition \ref{def:spaces}. Given a forcing $\mathbf{F}$ on $\mathds{B}^N$ and a microscopic configuration $\eta\in D\big([0,T];\mathbf{M}_{\{0,1\}}^N\big)$, define $\mathbf{R}_{T,\mathbf{F}}^N$ as the cost functional of discrete profiles
\begin{equation}\label{eq:discrete_R_functional}
\begin{split}
&\mathbf{R}_{T,\mathbf{F}}^N(\rho^N,\mathcal{Q}^N)\\
&:=\frac{1}{N^2}\sum_{x\in\mathds{B}^N}\mathcal{Q}_{T,x}^N\cdot\mathbf{F}_x-\int_0^T\sum_{x\in\mathds{B}^N,e\in\{e_1,e_2\}}\frac{1}{2N^2}\frac{(\mathbf{F}_{x,e})^2}{\log N}\sigma_{x,e}(\eta_{t})dt-\int_0^T \sum_{x\in\mathds{B}^N}\frac{1}{2N^2}\eta_{t,x}\big(\ddiv\mathbf{F}_x\big) dt.
\end{split}
\end{equation}
In addition, consider a time-independent continuum profile $(\rho,q)$ with $\rho\in \mathcal{C}^1(B_1;[0,1])$ and $q\in \mathcal{M}_{\mathds{R}^2}(B_1)$. Take a continuum forcing $\mathcal{F}$ with $\mathcal{F}^{\textup{rot}}$ as its rotation component in Remark \ref{rmk:vector_analysis}. Define the cost functional of continuum profiles
\begin{equation}\label{eq:continuum_R_cost}
\begin{split}
&\mathcal{R}_{T,\mathcal{F}}^N(\rho,q)=T\int_{N^{-1}\leq |x|\leq 1} q_x\cdot\mathcal{F}_x^{\textup{rot}}dx-\frac{T}{2}\int_{N^{-1}\leq |x|\leq 1} \frac{|\mathcal{F}|^2}{\log N}\sigma(\rho_x)dx-\frac{T}{2}\int_{N^{-1}\leq |x|\leq 1} \rho_x\nabla\cdot\mathcal{F}dx.
\end{split}
\end{equation}
\end{definition}

\begin{remark}\label{rmk:error_term}
Starting from the next subsection, we frequently encounter generic constants and small error terms. For example, $o_N(1;\delta)$ refers to an $o_N(1)$ whose form depends on parameter $\delta$. When the structure of errors terms becomes complicated, we may use the notation $\err(a_m,...,a_1)$, which means the error term converges to $0$ as we take parameters $(a_m,...,a_1)$ to its limit in the correct order (from $a_1$ to $a_m$). The limit of a parameter $a_i$ is either $0$ or $+\infty$, which will be clear from the context
\begin{equation}\label{eq:composed_error}
\err(a_m,...,a_1)\rightarrow 0,\quad \textup{as first $a_1\rightarrow +\infty$ or $0$,\quad then $a_2\rightarrow +\infty$ or $0$,\quad ...\quad and finally $a_m\rightarrow +\infty$ or $0$.}
\end{equation}
\end{remark}

\subsection{Approximating the event $\langle \mathcal{Q}_T^N,\mathcal{V} \rangle\geq JT$}

In this section, we prove that the large deviation upper bound of $\langle \mathcal{Q}_T^N,\mathcal{V}\rangle\geq JT$ can be reduced to the large deviation of simpler events, in the spirit of the \textit{construction of compact sets} in classical settings as \eqref{eq:compact_reduction}.

Define truncation functions $\mathcal{B}_+^{\eps}$, $\mathcal{B}_{-}^{\eps}$, $\Psi^{\eps}$ as the following, and later we will use these functions to truncate the empirical current $\mathcal{Q}_T^N$ to simplify some analysis in Sections \ref{sec:coarse_graining}-\ref{sec:true_proof_upper}
\begin{equation}\label{eq:radial_truncation}
\mathcal{B}_{+,x}^{\eps}:=\mathds{1}_{N^{-\eps}\leq |x|\leq 1},\quad \mathcal{B}_{-,x}^{\eps}:=\mathds{1}_{|x|\leq N^{-1+\eps}},\quad \Psi_x^{\eps}:=1-\mathcal{B}_{+,x}^{\eps}-\mathcal{B}_{-,x}^{\eps}.
\end{equation}

Recall $[\rho_x^N]_{\eps}$ is the space-time average of $\eta$ defined in \eqref{eq:space_time_average}. Given parameters $(J,T,\gamma,\eps,\delta)$, we construct an approximation $\mathbf{S}_{\delta,\eps}^N$ for the set $\mathbf{S}^N$ of profiles $(\rho^N,\mathcal{Q}^N)$ with $\langle \mathcal{Q}_T^N,\mathcal{V} \rangle \geq JT$
\begin{equation}\label{eq:S_sets}
\begin{split}
&\mathbf{S}^N:=\Big\{(\rho^N,\mathcal{Q}^N)\in D\big([0,T],\mathbf{M}_{\{0,1\}}^N\times\mathbf{M}_{\mathds{R}^2}^N \big)\Big|\quad  \langle\mathcal{Q}_T^N,\mathcal{V}\rangle\geq JT\Big\},\\
&\mathbf{S}_{\delta,\eps}^N:=\Bigg\{(\rho^N,\mathcal{Q}^N)\in \mathbf{S}^N\Bigg|\quad  \underbrace{|\mathcal{Q}_{T,x}^N|\leq N^{10}\textup{ for $x\in\mathds{B}^N$}}_{\textup{\blue{no extreme current}}}; \underbrace{\Big|[\rho_x^N]_{\eps}-\gamma\Big|\leq \eps,\ \textup{for $x\in\mathds{B}^N$ with $\eps\leq |x|\leq 1-2\eps$}}_{\textup{\blue{density is fixed away from $0$}}};\\
&\hspace{8cm} \underbrace{\big|\big\langle \mathcal{Q}_{T}^N(\mathcal{B}_+^{\eps}+\mathcal{B}_{-}^{\eps}),\mathcal{V}\big\rangle\big|\leq \delta}_{\textup{\blue{currents close to the boundary or the origin are negligible}}}\Bigg\}.
\end{split}
\end{equation}

\begin{lemma}\label{lem:SN_approximation}
The set $\mathbf{S}^N$ can be approximated by $\mathbf{S}_{\delta,\eps}^N$ with super-exponential accuracy in the following sense
\begin{equation*}
\limsup_{\eps\rightarrow 0}\limsup_{N\rightarrow +\infty}\frac{\log N}{N^2}\log\mathds{P}_{\gamma,\SS}^N\Big[\mathbf{S}^N\setminus\mathbf{S}_{\delta,\eps}^N\Big]=-\infty,\quad \textup{for any }J>0,\quad T>0,\quad \gamma\in (0,1),\quad \delta>0.
\end{equation*}
\end{lemma}

\begin{proof}
Define the bad events $\mathbf{A}^N$, $\mathbf{B}^N$, and $\mathbf{C}^N$ as
\begin{equation*}
\begin{split}
&\mathbf{A}^N:=\Big\{(\rho^N,\mathcal{Q}^N)\in \mathbf{S}^N\Big|\quad \exists x\in\mathds{B}^N:\ |\mathcal{Q}_{T,x}^N|>N^{10}\Big\},\\
&\mathbf{B}^N:=\Big\{(\rho^N,\mathcal{Q}^N)\in \mathbf{S}^N\Big|\quad \exists x\in\mathds{B}^N\ \textup{with }\eps\leq |x|\leq 1-2\eps:\ |[\rho_x^N]_{\eps}-\gamma|>\eps\Big\},\\
&\mathbf{C}^N:=\Big\{(\rho^N,\mathcal{Q}^N)\in \mathbf{S}^N\Big|\ \big|\big\langle \mathcal{Q}_{T}^N\mathcal{B}_+^{\eps},\mathcal{V}\big\rangle\big|> \frac{\delta}{2},\quad\textup{or}\quad  \big|\big\langle \mathcal{Q}_{T}^N\mathcal{B}_-^{\eps},\mathcal{V}\big\rangle\big|> \frac{\delta}{2}\Big\}.
\end{split}
\end{equation*}
By definition of $\mathbf{S}_{\delta,\eps}^N$ in \eqref{eq:S_sets}, it suffices to estimate the probabilities of $\mathbf{A}^N$, $\mathbf{B}^N$, and $\mathbf{C}^N$.

We first consider $\mathbf{A}^N$. Fix $x\in\mathds{B}^N$. For $e\in\{e_1,e_2\}$, let $J_{x,e}^+(T)$ and $J_{x,e}^-(T)$ denote the numbers of jumps across the edge $(x,x+e)$ during $[0,T]$ in the $+e$- and $-e$-directions, respectively. Since each oriented jump rate is bounded by $N^2/2$ according to the generator of SSEP in \eqref{eq:SSEP_generator}, the random variable
\begin{equation*}
\sum_{e\in\{e_1,e_2\}}\big(J_{x,e}^+(T)+J_{x,e}^-(T)\big)
\end{equation*}
is dominated by a Poisson random variable with rate of order $2N^2T$. Moreover the definition of $\mathcal{Q}_{T,x}^N$ in \eqref{eq:empirical_current} implies
\begin{equation*}
N|\mathcal{Q}_{T,x}^N|
\leq
\sum_{e\in\{e_1,e_2\}}\big(J_{x,e}^+(T)+J_{x,e}^-(T)\big).
\end{equation*}
A standard Poisson tail estimate therefore gives constants $C_T,c_T>0$ such that
\begin{equation*}
\mathds{P}_{\gamma,\SS}^N\big[|\mathcal{Q}_{T,x}^N|>N^{10}\big]
\leq C_T e^{-c_TN^5}.
\end{equation*}
Since $|\mathds{B}^N|\leq C N^2$, a union bound yields
\begin{equation}\label{eq:extreme_current}
\mathds{P}_{\gamma,\SS}^N\big[\mathbf{A}^N\big]
\leq
\sum_{x\in\mathds{B}^N}
\mathds{P}_{\gamma,\SS}^N\Big[|\mathcal{Q}_{T,x}^N|>N^{10}\Big]
\leq C_T N^2 e^{-c_TN^5}
\leq C_Te^{-c_TN^5/2}.
\end{equation}
We next consider $\mathbf{B}^N$. For $x\in\mathds{B}^N$ with $\eps\leq |x|\leq 1$, the microscopic averaging radius $\mathfrak{l}_x$ is of order $N$ by \eqref{eq:averaging_radius}. Hence when $\eps\leq |x|\leq 1-2\eps$ with $x$ away from the boundary of $B_1$, the function $ [\rho_x^N]_{\eps} $ is a space-time average over a box containing $O(N^2)$ sites.

Notice that $\mathds{P}_{\gamma,\SS,t}^N$ is the equilibrium law with density $\gamma$ at any time $t\geq 0$. Therefore the large deviation estimates of Bernoulli variables, together with a time-slicing argument controlling the jump of particles, can be applied to the space-time average $ [\rho_x^N]_{\eps} $. For example we can cut $[0,T]$ into $N^{10}$ segments, and within each segment the particles rarely jump, and the space-time average $[\rho_{x}^N]_{\eps}$ can be reduced to spatial averages $\overline{\eta}_{t,x}^{\eps\mathfrak{l}}$ at endpoints of the segments. This yields constants $C_{T,\gamma,\eps},c_{T,\gamma,\eps}>0$ such that
\begin{equation*}
\mathds{P}_{\gamma,\SS}^N\Big[\big|[\rho_x^N]_{\eps}-\gamma\big|>\eps\Big]
\leq C_{T,\gamma,\eps}e^{-c_{T,\gamma,\eps}N^2}
\qquad
\textup{for all }x\in\mathds{B}^N\textup{ with }\eps\leq |x|\leq 1-2\eps.
\end{equation*}
The number of such lattice points is at most $CN^2$. A union bound gives
\begin{equation}\label{eq:boundary_deviation}
\mathds{P}_{\gamma,\SS}^N[\mathbf{B}^N]
\leq
\sum_{x\in\mathds{B}^N:\eps\leq |x|\leq 1-2\eps}
\mathds{P}_{\gamma,\SS}^N\Big[\big|[\rho_x^N]_{\eps}-\gamma\big|>\eps\Big]
\leq C_{T,\gamma,\eps}N^2e^{-c_{T,\gamma,\eps}N^2}
\leq e^{-c_{T,\gamma,\eps}N^2/2}.
\end{equation}
The probability of $\mathbf{C}^N$ is super-exponentially small according to Lemma \ref{lem:LDP_boundary_no_contribution} below. Combining this with \eqref{eq:extreme_current} and \eqref{eq:boundary_deviation}, we obtain a faster exponential decay than $N^2/\log N$
\begin{equation*}
\limsup_{\eps\rightarrow 0}\limsup_{N\rightarrow +\infty}\frac{\log N}{N^2}\log \mathds{P}_{\gamma,\SS}^N\big[\mathbf{S}^N\setminus\mathbf{S}_{\delta,\eps}^N\big]
\leq
\limsup_{\eps\rightarrow 0}\limsup_{N\rightarrow +\infty}\frac{\log N}{N^2}\log\Big(\mathds{P}_{\gamma,\SS}^N[\mathbf{A}^N]+\mathds{P}_{\gamma,\SS}^N[\mathbf{B}^N]+\mathds{P}_{\gamma,\SS}^N[\mathbf{C}^N]\Big)=-\infty.
\end{equation*}
This proves the Lemma \ref{lem:SN_approximation}.
\end{proof}

\begin{lemma}\label{lem:LDP_boundary_no_contribution}
Fix arbitrary parameters $(T,\gamma,\delta,\mathfrak{b})$ with
\begin{equation*}
T>0,\quad \gamma\in (0,1),\quad \delta>0,\quad\mathfrak{b}>0.
\end{equation*}
For any families of vector fields $\{\mathbf{F}_{+}^{N,\eps}\}_{N<+\infty,\eps>0}$  and $\{\mathbf{F}_{-}^{N,\eps}\}_{N<+\infty,\eps>0}$ satisfying
\begin{equation*}
\big|\mathbf{F}_+^{N,\eps}\big|\leq \mathfrak{b}|x|^{-1}\mathds{1}_{N^{-\eps}\leq |x|\leq 1},\quad \big|\mathbf{F}_-^{N,\eps}\big|\leq \mathfrak{b}|x|^{-1}\mathds{1}_{|x|\leq N^{-1+\eps}},
\end{equation*}
we have
\begin{subequations}\label{eq:LDP_boundary_no_contribution}
\begin{align}
&\limsup_{\eps\rightarrow 0}\limsup_{N\rightarrow +\infty}\frac{\log N}{N^2}\log \mathds{P}_{\gamma,\SS}^N\Bigg[\Big|\Big\langle\mathcal{Q}_T^N,\mathbf{F}_{+}^{N,\eps} \Big\rangle\Big|\geq \delta\Bigg]=-\infty, \label{eq:LDP_boundary_no_contribution_+}\\
&\limsup_{\eps\rightarrow 0}\limsup_{N\rightarrow +\infty}\frac{\log N}{N^2}\log \mathds{P}_{\gamma,\SS}^N\Bigg[\Big|\Big\langle\mathcal{Q}_T^N,\mathbf{F}_{-}^{N,\eps} \Big\rangle\Big|\geq \delta\Bigg]=-\infty. \label{eq:LDP_boundary_no_contribution_-}
\end{align}
\end{subequations}
\end{lemma}

\begin{proof}

We only detail the proof of \eqref{eq:LDP_boundary_no_contribution_+}, because \eqref{eq:LDP_boundary_no_contribution_-} can be proved by the same method.

\textit{Step 1: Decompose $\mathbf{F}_+^{N,\eps}$ into a divergence-free component and a gradient component.}

For each $N>0$, apply the discrete Helmholtz decomposition in Proposition \ref{prop:discrete_Helmholtz} to $\mathbf{F}_{+}^{N,\eps}$, where it can be decomposed into a divergence-free component $g^{N,\eps}$ and a gradient component $\nabla_{\textup{disc}}\phi^{N,\eps}$
\begin{equation}\label{eq:discrete_helmholtz_applied}
\begin{split}
&\mathbf{F}_{+}^{N,\eps}=g^{N,\eps}+\nabla_{\textup{disc}}\phi^{N,\eps},\quad \nabla_{\textup{disc}}\cdot g_{x}^{N,\eps}\equiv 0,\\
&|\phi_{x}^{N,\eps}|\leq C_{\mathfrak{b}}(1-\log|x|),\quad \textup{for any $x\in\mathds{B}^N$; }\quad N^{-2}\sum_{x\in\mathds{B}^N}|g^{N,\eps}|^2\leq N^{-2}\sum_{x\in\mathds{B}^N}|\mathbf{F}_+^{N,\eps}|^2=C_{\mathfrak{b}}\eps O(\log N).
\end{split}
\end{equation}
The decomposition in \eqref{eq:discrete_helmholtz_applied} allow us to decompose $\langle \mathcal{Q}_T^N,\mathbf{F}_+^{N,\eps}\rangle$ for each $N$
\begin{equation*}
\langle \mathcal{Q}_T^N,\mathbf{F}_+^{N,\eps}\rangle=\langle \mathcal{Q}_T^N,g^{N,\eps}\rangle +\langle \mathcal{Q}_T^N,\nabla_{\textup{disc}}\phi^{N,\eps}\rangle=\langle \mathcal{Q}_T^N,g^{N,\eps}\rangle+\frac{1}{N^2}\sum_{x\in\mathds{B}^N}\big(\eta_{T,x}-\eta_{0,x}\big)\phi_{x}^{N,\eps},
\end{equation*}
where the last equality is derive by discrete integration by parts. Therefore to prove \eqref{eq:LDP_boundary_no_contribution_+}, it is sufficient to prove
\begin{subequations}\label{eq:boundary_nocontribution_+_decompose}
\begin{align}
&\limsup_{\eps\rightarrow 0}\limsup_{N\rightarrow +\infty}\frac{\log N}{N^2}\log \mathds{P}_{\gamma,\SS}^N\Bigg[\Big|\Big\langle\mathcal{Q}_T^N,g^{N,\eps} \Big\rangle\Big|\geq \frac{\delta}{2}\Bigg]=-\infty,\label{eq:boundary_nocontribution_+_decompose_rot}\\
&\limsup_{\eps\rightarrow 0}\limsup_{N\rightarrow +\infty}\frac{\log N}{N^2}\log \mathds{P}_{\gamma,\SS}^N\Bigg[\Big|\frac{1}{N^2}\sum_{x\in\mathds{B}^N}\big(\eta_{T,x}-\eta_{0,x}\big)\phi_{x}^{N,\eps}\Big|\geq \frac{\delta}{2}\Bigg]=-\infty\label{eq:boundary_nocontribution_+_decompose_grad}
\end{align}
\end{subequations}

\textit{Step 2: Proof of \eqref{eq:boundary_nocontribution_+_decompose_grad} for the contribution of the gradient component $\nabla_{\textup{disc}}\phi^{N,\eps}$.}

The upper bound $|\phi_x^{N,\eps}|\leq C_{\mathfrak{b}}(1-\log|x|)$ in \eqref{eq:discrete_helmholtz_applied} proves 
\begin{equation}\label{eq:nocontribution_gradient_1}
\Big|\frac{1}{N^2}\sum_{x\in\mathds{B}^N:|x|\leq \eps\textup{ or $|x|\geq 1-2\eps$}}\big(\eta_{T,x}-\eta_{0,x}\big)\phi_x^{N,\eps}\Big|= C_{\mathfrak{b}}O(\eps).
\end{equation}
Therefore we only need to estimate the contribution of $\phi_x^{N,\eps}$ with $\eps<|x|< 1-2\eps$. The idea is that we want to replace $\eta$ by its spatial average $\overline{\eta}^{\eps^2\mathfrak{l}}$ (see \eqref{eq:spatial_average}), with $\mathfrak{l}_x$ as the averaging radius $\mathfrak{l}_x=N|x|_{\infty}$. Define $\varphi_{x}^{N,\eps}$ as the spatial average of $\phi_{x}^{N,\eps}$, whose average will be transfered to $\eta$ later 
\begin{equation*}
\varphi_{x}^{N,\eps}:=\frac{1}{(1+2\eps^2\mathfrak{l}_x)^2}\sum_{y\in\mathds{B}^N}\phi_{y}^{N,\eps}\mathds{1}_{|y-x|_{\infty}\leq \eps^2\mathfrak{l}_xN^{-1}}.
\end{equation*}
We have
\begin{equation}
\begin{split}
&\Big|\frac{1}{N^2}\sum_{x\in\mathds{B}^N:\eps<|x|< 1-2\eps}\big(\eta_{T,x}-\eta_{0,x}\big)\phi_{x}^{N,\eps}-\frac{1}{N^2}\sum_{x\in\mathds{B}^N:\eps<|x|< 1-2\eps}\big(\eta_{T,x}-\eta_{0,x}\big)\varphi_{x}^{N,\eps}\Big|\\
&\leq C\frac{\eps^2}{N^2}\sum_{x\in\mathds{B}^N:\frac{\eps}{2}< |x|< 1-\eps}|\nabla_{\textup{disc}}\phi_x^{N,\eps}|\leq C\eps^2\Big(\frac{1}{N^2}\sum_{x\in\mathds{B}^N:\frac{\eps}{2}< |x|< 1-\eps}|\nabla_{\textup{disc}}\phi_{x}^{N,\eps}|^2\Big)^{\frac{1}{2}}\leq C\eps,
\end{split}
\end{equation}
where the last step is due to Proposition \ref{prop:discrete_Helmholtz}. Therefore we can transfer the spatial average from $\phi_x^{N,\eps}$ to $\eta$
\begin{equation}\label{eq:nocontribution_gradient_2}
\begin{split}
\frac{1}{N^2}\sum_{x\in\mathds{B}^N:\eps<|x|< 1-2\eps}\big(\eta_{T,x}-\eta_{0,x}\big)\phi_{x}^{N,\eps}&=\frac{1}{N^2}\sum_{x\in\mathds{B}^N:\eps<|x|< 1-2\eps}\big(\eta_{T,x}-\eta_{0,x}\big)\varphi_{x}^{N,\eps}+C_{\mathfrak{b}}O(\eps)\\
&=\frac{1}{N^2}\sum_{x\in\mathds{B}^N:\eps<|x|< 1-2\eps}\big(\overline{\eta}_{T,x}^{\eps^2\mathfrak{l}}-\overline{\eta}_{0,x}^{\eps^2\mathfrak{l}}\big)\phi_{x}^{N,\eps}+C_{\mathfrak{b}}O(\eps).
\end{split}
\end{equation}
At any time $t\geq 0$, the measure $d\mathds{P}_{\gamma,\SS,t}^N$ is a product Bernoulli measure with density $\gamma$. For any $x\in\mathds{B}^N$ with $\eps< |x|< 1-2\eps$, the configuration $\overline{\eta}_{x}^{\eps^2\mathfrak{l}}$ is the average of $\eta_{\cdot}$ over $\eps^4N^2$ sites. Therefore $\overline{\eta}_{T,x}^{\eps^2\mathfrak{l}}$ and $\overline{\eta}_{0,x}^{\eps^2\mathfrak{l}}$ can always be approximated by $\gamma$ with super-exponential accuracy
\begin{equation}\label{eq:nocontribution_gradient_3}
\limsup_{N\rightarrow +\infty}\frac{\log N}{N^2}\log\mathds{P}_{\gamma,\SS}^N\Bigg[\Big|\frac{1}{N^2}\sum_{x\in\mathds{B}^N:\eps<|x|<1-2\eps}\big(\overline{\eta}_{T,x}^{\eps^2\mathfrak{l}}-\overline{\eta}_{0,x}^{\eps^2\mathfrak{l}}\big)\phi_{x}^{N,\eps}\Big|\geq \frac{\delta}{2}\Bigg]=-\infty.
\end{equation}
Combine \eqref{eq:nocontribution_gradient_1}, \eqref{eq:nocontribution_gradient_2}, and \eqref{eq:nocontribution_gradient_3}, and we can prove \eqref{eq:boundary_nocontribution_+_decompose_grad}.

\textit{Step 3: Proof of \eqref{eq:boundary_nocontribution_+_decompose_rot} for the contribution of the divergence-free component $g^{N,\eps}$.}

For \eqref{eq:boundary_nocontribution_+_decompose_rot}, we only detail proof of the super-exponential estimate of the event $\langle \mathcal{Q}_T^N,g^{N,\eps}\rangle\geq \frac{\delta}{2}$, because the estimate of $\langle \mathcal{Q}_T^N,g^{N,\eps}\rangle\leq -\frac{\delta}{2}$ can be derived by replacing $g^{N,\eps}$ with $-g^{N,\eps}$.

To estimate the cost of $\langle \mathcal{Q}_T^N,g^{N,\eps}\rangle\geq \frac{\delta}{2}$, we tilt the measure $\mathds{P}_{\gamma,\SS}^N$ by an asymmetric forcing $\eps^{-\frac{1}{2}}g^{N,\eps}$. Denote $d\mathds{P}_{\gamma,\eps^{-1/2}g^{N,\eps}}^N$ as the probability measure of trajectories, for the modified dynamics with forcing $\eps^{-\frac{1}{2}}g^{N,\eps}$, see \eqref{eq:modified_generator}
\begin{equation}
\mathds{P}_{\gamma,\SS}^N\Bigg[\Big\langle\mathcal{Q}_T^N,g^{N,\eps} \Big\rangle\geq \delta\Bigg]=\mathds{P}_{\gamma,\eps^{-1/2}g^{N,\eps}}^N\Bigg[\frac{d\mathds{P}_{\gamma,\SS}^N}{d\mathds{P}_{\gamma,\eps^{-1/2}g^{N,\eps}}^N}\mathds{1}_{\langle \mathcal{Q}_T^N,g^{N,\eps}\rangle\geq \delta}\Bigg],\quad \textup{$\frac{d\mathds{P}_{\gamma,\SS}^N}{d\mathds{P}_{\gamma,\eps^{-1/2}g^{N,\eps}}^N}$ is Radon-Nikodym derivative}.
\end{equation}
To compute the Radon-Nikodym derivative, we can apply Proposition \ref{prop:radon_nikodym}. Since $g^{N,\eps}$ is divergence-free according to Proposition \ref{prop:discrete_Helmholtz}, the Radon-Nikodym derivative does not have a divergence term
\begin{equation*}
\begin{split}
&\mathds{P}_{\gamma,\SS}^N\Bigg[\Big\langle\mathcal{Q}_T^N,g^{N,\eps} \Big\rangle\geq \delta\Bigg]\\
&=\mathds{P}_{\gamma,\eps^{-1/2}g^{N,\eps}}^N\Bigg[\exp\Bigg(-\frac{N^2}{\log N}\Big(\langle \mathcal{Q}_T^N,\eps^{-\frac{1}{2}}g^{N,\eps}\rangle-\frac{1}{2}\int_0^T\sum_{\substack{x\in\mathds{B}^N\\ e\in\{e_1,e_2\}}}N^{-2}\frac{\eps^{-1}|g_{x,e}^{N,\eps}|^2}{\log N}\sigma_{x,e}(\eta_{t})dt+o_N(1;\mathfrak{b})\Big)\Bigg)\mathds{1}_{\langle \mathcal{Q}_T^N,g^{N,\eps}\rangle\geq \delta}\Bigg].
\end{split}
\end{equation*}
Because we are looking at the event $\langle \mathcal{Q}_T^N,g^{N,\eps} \rangle\geq \delta$, we always have $\langle \mathcal{Q}_T^N,\eps^{-\frac{1}{2}}g^{N,\eps}\rangle\geq \eps^{-\frac{1}{2}}\delta$. Taking $\eps\rightarrow 0$, we have $\langle \mathcal{Q}_T^N,\eps^{-\frac{1}{2}}g^{N,\eps}\rangle\rightarrow +\infty$. For the term involving $|g_{x,e}^{N,\eps}|^2$, the upper bound of the $L^2$-norm of $g^{N,\eps}$ in \eqref{eq:discrete_helmholtz_applied} implies it is of order $O(\delta)$. Therefore $\langle \mathcal{Q}_T^N,\eps^{-\frac{1}{2}}g^{N,\eps}\rangle\rightarrow +\infty$ dominates, and we conclude the proof.
\end{proof}

\subsection{Coarse-grained Functions and Variational Formulas over Discrete Profiles}\label{sec:coarse_graining}

The next task is to construct a cover  of the set $\mathbf{S}_{\delta,\eps}^N$ in \eqref{eq:S_sets} for each $N>1$. The idea is that, we coarse-grain the discrete configurations $(\rho^N,\mathcal{Q}^N)$ in $\mathbf{S}_{\delta,\eps}^N$ into different groups. To this end, a coarse-grained function $\mathbf{h}_{\eps}^N:=\big(\mathbf{h}_{\eps,\textup{den}}^N,\mathbf{h}_{\eps,\textup{cur}}^N\big)$ is defined as
\begin{equation}\label{eq:coarse_graining_function}
\begin{split}
&\mathbf{h}_{\eps,\textup{den}}^N:\lden\rightarrow [0,1]\cap \frac{\eps}{\log^2N}\mathds{Z},\qquad\mathbf{h}_{\eps,\textup{cur}}^N:\lcur\rightarrow \big[-N^{10},N^{10}\big]^2\cap\frac{\eps}{\log^2N}\mathds{Z}^2,
\end{split}
\end{equation}
with $\lden$ and $\lcur$ as two different sets of points $x\in\mathds{B}^N$
\begin{equation*}
\lcur:=\Big\{x\in\mathds{B}^N\Big|\ |x|\leq N^{-\eps}\Big\},\quad \lden:=\lcur\cup \Big\{x\in\mathds{B}^N\Big|\ x\in \lfloor \sqrt{N}\rfloor \frac{1}{N}\mathds{Z}^2\Big\}.
\end{equation*}
The set $\lcur$ does not contain $\{x\in\mathds{B}^N|\ |x|\geq N^{-\eps}\}$ because the contribution of $\mathcal{Q}_T^N$ in this layer is negligible according to Lemma \ref{lem:LDP_boundary_no_contribution}. The set $\lden$ contains the layer $|x|\geq N^{-\eps}$ because we need to capture the boundary value $\gamma$ of the density profiles, see $\mathbf{S}_{\delta,\eps}^N$ in \eqref{eq:S_sets}. Recall the space-time average $[\rho_{\cdot}^N]_{\eps}$ defined in \eqref{eq:space_time_average}. For any coarse-grained function $\mathbf{h}_{\eps}^N$, we define the associated set $\mathbf{O}_{\eps,\mathbf{h}}^N$ as
\begin{equation}\label{eq:open_set_for_discrete}
\begin{split}
\mathbf{O}_{\eps,\mathbf{h}}^N&:=\Bigg\{(\rho^N,\mathcal{Q}^N)\in D\big([0,T];\mathbf{M}_{\{0,1\}}^N\times \mathbf{M}_{\mathds{R}^2}^N\big)\Bigg|\ \Big|[\rho_x^N]_{\eps}-\gamma\Big|\leq \eps,\ \textup{for all $x\in \mathds{B}^N$ with $\eps\leq |x|\leq 1-2\eps$};\\
&\qquad \Big|[\rho_x^N]_{\eps}-\mathbf{h}_{\eps,\textup{den}}^N(x)\Big|<\frac{\eps}{\log^2N}\textup{ for $x\in\lden$ and }\Big|\mathcal{Q}_{T,x}^N-\mathbf{h}_{\eps,\textup{cur}}^N(x)\Big|<\frac{\eps}{\log^2N}\textup{ for $x\in\lcur$}\Bigg\}.
\end{split}
\tag{Set1}
\end{equation}
We count the size of the relevant sets
\begin{equation*}
\begin{split}
&\#\lden,\#\lcur= O\Big(N^{2(1-\eps)}\Big),\quad \# \Big\{[0,1]\cap \frac{\eps}{\log^2N}\mathds{Z}\Big\}=O\Big(\frac{(\log N)^2}{\eps}\Big),\\
&\#\Big\{\big[-N^{10},N^{10}\big]^2\cap\frac{\eps}{\log^2N}\mathds{Z}^2\Big\}=O\Big(\frac{N^{20}}{\eps^2}\log^4N\Big).
\end{split}
\end{equation*}
Denote $M_N$ as the total number of all possible functions $\mathbf{h}_{\eps}^N$. We can compute the logarithm of $M_N$
\begin{equation}\label{eq:total_number_of_h}
\log M_N\leq \underbrace{CN^{2(1-\eps)}\log\Bigg(C\frac{(\log N)^2}{\eps}\Bigg)}_{\textup{\blue{$\log$ of total number of $\mathbf{h}_{\eps,\textup{den}}^N$ }}}+ \underbrace{CN^{2(1-\eps)}\log \Bigg(C\frac{N^{20}}{\eps^2}\log^4N\Bigg)}_{\textup{\blue{$\log$ of total number of $\mathbf{h}_{\eps,\textup{cur}}^N$}}}\ll \frac{N^2}{\log N}.
\end{equation}
The probability of $\mathbf{S}_{\delta,\eps}^N$ in \eqref{eq:S_sets} can be bounded by the probabilities of all $\doh$ (recall the notation $o_N$ in Remark \ref{rmk:error_term})
\begin{equation}\label{eq:open_cover_max_1}
\begin{split}
\frac{\log N}{N^2}\log\mathds{P}_{\gamma,\SS}^N\Big[\mathbf{S}_{\delta,\eps}^N\Big]\leq \frac{\log N}{N^2}\log\mathds{P}_{\gamma,\SS}^N\Bigg[\bigcup_{\mathbf{h}}\doh\Bigg]&\leq\frac{\log N}{N^2}\log \Bigg(M_N\max_{\mathbf{h}}\mathds{P}_{\gamma,\SS}^N\Big[\doh\Big]\Bigg)\\
&\leq \frac{\log N}{N^2}\log M_N+\frac{\log N}{N^2}\max_{\mathbf{h}}\log\mathds{P}_{\gamma,\SS}^N\Big[\doh\Big]\\
&\leq o_N(1;\eps)+\frac{\log N}{N^2}\max_{\mathbf{h}}\log\mathds{P}_{\gamma,\SS}^N\Big[\doh\Big].
\end{split}
\end{equation}
Later we will prove Lemma \ref{lem:discrete_optimization}, showing that $\log\mathds{P}_{\gamma,\SS}^N[\doh]$ can be bounded from above by a variational formula over forcing $\mathbf{F}$ and discrete profiles $(\rho^N,\mathcal{Q}^N)$. To this end, we define the set $\mathbf{C}^N[\mathfrak{b}]$ of forcings.
\begin{definition}\label{def:set_of_forcing}
For any $(N>1,\mathfrak{b}>0)$, define $\mathbf{C}^N[\mathfrak{b}]$ as the set of forcings $\mathbf{F}:\mathds{B}^N\mapsto \mathds{R}^2$ satisfying conditions \eqref{eq:H1}, \eqref{eq:H2}, and \eqref{eq:disc}. Condition \eqref{eq:H1} is on the regularities of $\mathbf{F}$
\begin{equation}\label{eq:H1}
|\mathbf{F}_x|\leq \mathfrak{b}|x|^{-1},\quad  |\nabla_{\textup{disc}} \mathbf{F}_x|\leq \mathfrak{b}|x|^{-2},\quad  |(\nabla_{\textup{disc}})^2\mathbf{F}_x|\leq \mathfrak{b}|x|^{-3}\quad\textup{for $x\in \mathds{B}^N$}, \tag{H1}
\end{equation}
Condition \eqref{eq:H2} is on the forcing $\mathbf{F}$ being identically $0$ at the central core or near the boundary $\p B_1$
\begin{equation}\label{eq:H2}
|\mathbf{F}_x|\equiv 0\textup{ for $x\in \mathds{B}^N$ with $|x|\leq N^{-1+\frac{\eps}{2}}$ or $|x|\geq \frac{1}{10}$}. \tag{H2}
\end{equation}
The condition \eqref{eq:disc} gives a decomposition of $\mathbf{F}$, which is useful when we consider the inner product $\langle \mathcal{Q}_T^N,\mathbf{F}\rangle$ on $\mathds{B}^N$
\begin{equation}\label{eq:disc}
\begin{split}
&\textup{exists a function $V$ on $\mathds{B}^N$ with $|V_x|\leq \mathfrak{b}^2(1-\log|x|),\ |\nabla_{\textup{disc}}V_x|\leq \mathfrak{b}^2|x|^{-1}$} \\
&\hspace{3cm} \textup{such that $\mathbf{F}_x-\nabla_{\textup{disc}}V_x$ is purely rotational (see Remark \ref{rmk:vector_analysis}) on $x\in\mathds{B}^N$}. 
\end{split}
\tag{Disc}
\end{equation}
We use $\mathbf{F}_x^{\textup{rot}}$ to denote $\mathbf{F}_x-\nabla_{\textup{disc}}V_x$.
\end{definition}

\begin{figure}[h]
    \centering
    \includegraphics[width=0.9\linewidth]{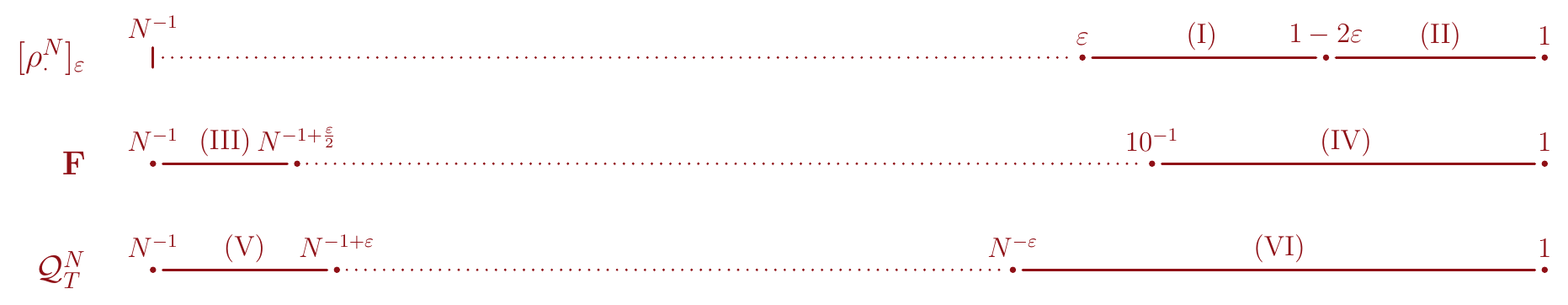}
    \caption{This figure is a summary of different conditions we have on profiles and forcings.}
    \label{fig:truncation_summary}
\end{figure}

\begin{itemize}
    \item In interval (I) the space-time average $[\rho_{\cdot}^N]_{\eps}$ satisfies $|[\rho_x^N]_{\eps}-\gamma|\leq \eps$, and in interval (II) we do not have this boundary condition because $[\rho_{x}^N]_{\eps}$ is not sufficiently averaged due to the presence of the boundary $\p B_1$.
    \item In interval (III) the forcing is taken to be zero $\mathbf{F}_x\equiv 0$, in order to truncate the singularity of forcing $\mathbf{F}\in \mathbf{C}^N[\mathfrak{b}]$. In interval (IV) the forcing is also zero $\mathbf{F}_x\equiv 0$, in order to satisfy the condition \eqref{eq:appendix_A_bound} of Proposition \ref{prop:radon_nikodym}. This also allows us to perform discrete integration by parts to the forcing $\mathbf{F}$. 
    \item In interval (V) and (VI) we have $|\langle \mathcal{Q}_T^N(\mathcal{B}_+^{\eps}+\mathcal{B}_-^{\eps}),\mathcal{V}\rangle|\leq \delta$. This allows us to truncate the empirical current in (V) and (VI). This is necessary because forcing $\mathbf{F}\in \mathbf{C}^N[\mathfrak{b}]$ can not estimate of the cost of the empirical current $\mathcal{Q}_T^N$ in these intervals.
\end{itemize}

\begin{lemma}\label{lem:discrete_optimization}
For arbitrary $(T,N,\eps,\delta,\mathfrak{b})$ with
\begin{equation*}
T>0,\quad N>1,\quad \eps>0,\quad \delta>0,\quad \mathfrak{b}>0,
\end{equation*}
and any set $\doh$ as in \eqref{eq:open_set_for_discrete}, we have
\begin{equation*}
\frac{\log N}{N^2}\log\mathds{P}_{\gamma,\SS}^N\Big[\doh\Big]\leq \Big(-\inf_{\rho^N,\mathcal{Q}^N\in\doh}\sup_{\mathbf{F}\in \mathbf{C}^N[\mathfrak{b}]}\mathbf{R}_{T,\mathbf{F}}^{N,\eps}(\rho^N,\mathcal{Q}^N)\Big)+\err(\mathfrak{b},\delta,\eps,N).
\end{equation*}
Here $\err$ is a generic error described in \eqref{eq:composed_error}. Recall the space-time average $[\rho_{\cdot}^N]_{\eps}$ defined in \eqref{eq:space_time_average} and the empirical current $\mathcal{Q}_T^N$ defined in \eqref{eq:empirical_current}. The functional $\mathbf{R}_{T,\mathbf{F}}^{N,\eps}$ is \begin{equation}\label{eq:discrete_R_functional_eps}
\begin{split}
&\mathbf{R}_{T,\mathbf{F}}^{N,\eps}(\rho^N,\mathcal{Q}^N):=\underbrace{\Big\langle\Psi_x^{\eps}\mathcal{Q}_T^N,\mathbf{F}^{\textup{rot}}\Big\rangle}_{\textup{\blue{truncated current term}}}-\underbrace{T\sum_{x\in\mathds{B}^N,e\in\{e_1,e_2\}}\frac{1}{2N^2}\frac{|\mathbf{F}_{x,e}|^2}{\log N}\sigma\Big([\rho_x^N]_{\eps}\Big)}_{\textup{\blue{nonlinear term of the space-time average}}}-\underbrace{T\sum_{x\in \mathds{B}^N}\frac{1}{2N^2}[\rho_x^N]_{\eps}\big(\ddiv\mathbf{F}_x\big)}_{\textup{\blue{divergence term}}}.
\end{split}
\end{equation}
\end{lemma}

We can compare the functional $\mathbf{R}_{T,\mathbf{F}}^{N,\eps}$ with the functional $\mathbf{R}_{T,\mathbf{F}}^N$ in Definition \ref{def:R_functionals}.
\begin{itemize}
    \item A first difference is the configuration $\eta$ in $\mathbf{R}_{T,\mathbf{F}}^N$ is replaced by the space-time average $[\rho_{\cdot}^N]_{\eps}$ in $\mathbf{R}_{T,\mathbf{F}}^{N,\eps}$.
    \item A second difference is the cost involving the current profile $\mathcal{Q}_T^N$: in $\mathbf{R}_{T,\mathbf{F}}^{N,\eps}$ the current profile $\mathcal{Q}_T^N$ is truncated by $\Psi^{\eps}$ (see \eqref{eq:radial_truncation}), and the current cost is determined by the rotational component $\mathbf{F}^{\textup{rot}}$ of the forcing $\mathbf{F}$. This reduction is important because in Section \ref{sec:true_proof_upper}, it allow us to write the large deviation upper bound as a variational formula over current profiles with only rotational component.
\end{itemize}

\begin{proof}[Proof of Lemma \ref{lem:discrete_optimization}]

\textit{Step 1: Apply the replacement Lemma \ref{lem:replacement_lemma}.}

Take any $\delta>0$ and arbitrary forcing $\mathbf{F}\in \mathbf{C}^N[\mathfrak{b}]$ in Definition \ref{def:set_of_forcing}, define the set $\mathbf{B}_{\delta,\eps,\textup{cur}}^N$ as
\begin{equation*}
\mathbf{B}_{\delta,\eps,\textup{cur}}^N:=\Big\{(\rho^N,\mathcal{Q}^N)\in D\big([0,T];\mathbf{M}_{\{0,1\}}^N\times \mathbf{M}_{\mathds{R}^2}^N\big)\Big|\ \big|\big\langle (\mathcal{B}_{+}^{\eps}+\mathcal{B}_{-}^{\eps})\mathcal{Q}_{T}^N,\mathbf{F}^{\textup{rot}}\big\rangle\big|\leq \delta\Big\},
\end{equation*}
where $\mathbf{F}^{\textup{rot}}$ is given in Definition \ref{def:set_of_forcing}. By Lemma \ref{lem:LDP_boundary_no_contribution}, we have
\begin{equation}\label{eq:super_exp_boundary_current}
\limsup_{\eps\rightarrow 0}\limsup_{N\rightarrow +\infty}\frac{\log N}{N^2}\log\mathds{P}_{\gamma,\SS}^N\big[\Big(\mathbf{B}_{\delta,\eps,\textup{cur}}^{N}\Big)^c\big]=-\infty.
\end{equation}
Define the set $\mathbf{B}_{\delta,\eps}^N\subset \mathbf{B}_{\delta,\eps,\textup{cur}}^{N}$ as
\begin{equation}\label{eq:union_of_replace_error_set}
\begin{split}
\mathbf{B}_{\delta,\eps}^N:=&\mathbf{B}_{\delta,\eps,H_1^N}^{N}\Big(\eta_{x}(1-\eta_{x+e_1})\Big)\cap\mathbf{B}_{\delta,\eps,H_1^N}^{N}\Big(\eta_{x+e_1}(1-\eta_{x})\Big)\\
&\cap \mathbf{B}_{\delta,\eps,H_2^N}^{N}\Big(\eta_{x}(1-\eta_{x+e_2})\Big)\cap\mathbf{B}_{\delta,\eps,H_2^N}^{N}\Big(\eta_{x+e_2}(1-\eta_{x})\Big) \cap \mathbf{B}_{\delta,\eps,H_3^N}^{N}\big(\eta_{x}\big)\cap \mathbf{B}_{\delta,\eps,\textup{cur}}^N,
\end{split}
\end{equation}
where the relevant sets have been defined in \eqref{eq:replacement_error_set}. Here the families of functions $\{H_i^N\}_{N}$ (with $i\in \{1,2,3\}$) are chosen as 
\begin{equation*}
H_{1,x}^N:=\frac{1}{4}\frac{|\mathbf{F}_{x,e_1}|^2}{\log N},\quad H_{2,x}^N:=\frac{1}{4}\frac{|\mathbf{F}_{x,e_2}|^2}{\log N},\quad  H_{3,x}^N:=\frac{1}{2}\ddiv\mathbf{F}_x.
\end{equation*}
Across the proof, some error terms are related to the forcing $\mathbf{F}\in \mathbf{C}^N[\mathfrak{b}]$. However, the errors will only depend on the parameter $\mathfrak{b}$, not on the specific form of $\mathbf{F}$.

To estimate the large deviation upper bound, we first restrict the probability to the set $\mathbf{B}_{\delta,\eps}^N$. Within this set $\mathbf{B}_{\delta,\eps}^N$, we can neglect the current $\mathcal{Q}_{T}^N$ near the outer and inner boundaries, and average the empirical density $\rho^N$.
\begin{equation}\label{eq:restriction_upper}
\begin{split}
&\frac{\log N}{N^2}\log \mathds{E}_{\gamma,\SS}^N\Big[\mathds{1}_{(\rho^N,\mathcal{Q}^N)\in\mathbf{O}_{\eps,\mathbf{h}}^N}\Big]\leq\max\Big\{\frac{\log N}{N^2}\log \mathds{E}_{\gamma,\SS}^N\Big[\mathds{1}_{(\rho^N,\mathcal{Q}^N)\in\mathbf{O}_{\eps,\mathbf{h}}^N}\mathds{1}_{\mathbf{B}_{\delta,\eps}^N}\Big],\frac{\log N}{N^2}\log\mathds{P}_{\gamma,\SS}^N(\mathbf{B}_{\delta,\eps}^{N,c})\Big\}+o_N(1).
\end{split}
\end{equation}
The super-exponential replacement Lemma \ref{lem:replacement_lemma} and the estimate \eqref{eq:super_exp_boundary_current} shows the probability of $\mathbf{B}_{\delta,\eps}^{N,c}$ is super-exponentially small. Therefore it only produces a negligible error
\begin{equation*}
\begin{split}
&\frac{\log N}{N^2}\log \mathds{E}_{\gamma,\SS}^N\Big[\mathds{1}_{(\rho^N,\mathcal{Q}^N)\in\mathbf{O}_{\eps,\mathbf{h}}^N}\Big]\leq \frac{\log N}{N^2}\log \mathds{E}_{\gamma,\SS}^N\Big[\mathds{1}_{(\rho^N,\mathcal{Q}^N)\in\mathbf{O}_{\eps,\mathbf{h}}^N}\mathds{1}_{\mathbf{B}_{\delta,\eps}^N}\Big]+\errtotal,
\end{split}
\end{equation*}
and we are going to estimate the expectation in the RHS using Proposition \ref{prop:radon_nikodym}. For each $\eps$ and $N$, the expectation further equals to the following, with $\mathbf{F}\in\mathbf{C}^N[\mathfrak{b}]$ chosen previously and $\mathbf{R}_{T,\mathbf{F}}^N$ defined in 
\eqref{eq:discrete_R_functional}
\begin{equation*}
\begin{split}
&\frac{\log N}{N^2}\log\mathds{E}_{\gamma,\mathbf{F}}^N\Bigg[\frac{d\mathds{P}_{\gamma,\SS}^N}{d\mathds{P}_{\gamma,\mathbf{F}}^N}\mathds{1}_{(\rho^N,\mathcal{Q}^N)\in\mathbf{O}_{\eps,\mathbf{h}}^N}\mathds{1}_{\mathbf{B}_{\delta,\eps}^N}\Bigg]\\
&\qquad\leq \frac{\log N}{N^2}\log\mathds{E}_{\gamma,\mathbf{F}}^N\Bigg[\exp\Bigg(-\frac{N^2}{\log N}\Big(\underbrace{\mathbf{R}_{T,\mathbf{F}}^{N}(\rho^N,\mathcal{Q}^N)}_{\textup{\blue{dynamical cost}}}+o_N(1;\mathfrak{b})\Big)\Bigg)\mathds{1}_{(\rho^N,\mathcal{Q}^N)\in\mathbf{O}_{\eps,\mathbf{h}}^N}\mathds{1}_{\mathbf{B}_{\delta,\eps}^N}\Bigg].
\end{split}
\end{equation*}
Here the error $o_N(1;\mathfrak{b})$ comes from the Taylor expansion of the Radon-Nikodym derivative (see Proposition \ref{prop:radon_nikodym}). 

Recall that the replacement lemma \ref{lem:replacement_lemma} estimates the error of replacing $\eta$ with the spatial average $\overline{\eta}_x^{\eps\mathfrak{l}}$ given as \eqref{eq:spatial_average}. Also recall $\Psi_{x}^{\eps}$ is the truncation function defined in \eqref{eq:radial_truncation}. Within $\mathbf{B}_{\delta,\eps}^N$, the definition \eqref{eq:union_of_replace_error_set} of $\mathbf{B}_{\delta,\eps}^N$ implies
\begin{equation*}
\begin{split}
&\mathbf{R}_{T,\mathbf{F}}^N(\rho^N,\mathcal{Q}^N)=\frac{1}{N^2}\Bigg(\sum_{x\in\mathds{B}^N}\Psi_x^{\eps}\mathcal{Q}_{T,x}^N\cdot\mathbf{F}_x^{\textup{rot}}-\frac{1}{2}\int_0^T\sum_{x\in \mathds{B}^N}\frac{|\mathbf{F}_x|^2}{\log N}\sigma\big(\bar{\eta}_{t,x}^{\eps \mathfrak{l}}\big)dt-\frac{1}{2}\int_0^T\sum_{x\in \mathds{B}^N}\bar{\eta}_{t,x}^{\eps \mathfrak{l}}\nabla_{\textup{disc}}\cdot \mathbf{F}_xdt\Bigg)+O(\delta).
\end{split}
\end{equation*}
Here to get the $\Psi_{x}^{\eps}\mathcal{Q}_{T,x}^N\cdot\mathbf{F}_x^{\textup{rot}}$, we first use discrete integration by parts to change $\mathcal{Q}_{T,x}^N\cdot\mathbf{F}_x$ into $\mathcal{Q}_{T,x}^N\cdot\mathbf{F}_x^{\textup{rot}}$, similar to the proof of Lemma \ref{lem:LDP_boundary_no_contribution}. Then we use the definition of $\mathbf{B}_{\delta,\eps}^N$ to change $\mathcal{Q}_{T,x}^N\cdot\mathbf{F}_x^{\textup{rot}}$ into $\Psi_{x}^{\eps}\mathcal{Q}_{T,x}^N\cdot\mathbf{F}_x^{\textup{rot}}$.

\textit{Step 2: Space-time average by convexity.}

The convexity of $-\sigma(\cdot)$ further allows us to replace the space average $\overline{\eta}_x^{\eps\mathfrak{l}}$ by the space-time average $[\rho_x^N]_{\eps}$ in \eqref{eq:space_time_average}, leading to the functional $\mathbf{R}_{T,\mathbf{F}}^{N,\eps}$ in \eqref{eq:discrete_R_functional_eps}
\begin{equation}\label{eq:space_time_ave_functional}
\begin{split}
\mathbf{R}_{T,\mathbf{F}}^N(\rho^N,\mathcal{Q}^N)&\geq \frac{1}{N^2}\Bigg(\sum_{x\in\mathds{B}^N}\Psi_{x}^{\eps}\mathcal{Q}_{T,x}^N\cdot \mathbf{F}_x^{\textup{rot}}-\frac{T}{2}\sum_{x\in\mathds{B}^N}\frac{|\mathbf{F}_x|^2}{\log N}\sigma\big([\rho_x^N]_{\eps}\big)-\frac{T}{2}\sum_{x\in\mathds{B}^N}[\rho_x^N]_{\eps}\nabla_{\textup{disc}}\cdot\mathbf{F}_x\Bigg)+O(\delta)\\
&=\mathbf{R}_{T,\mathbf{F}}^{N,\eps}(\rho^N,\mathcal{Q}^N)+O(\delta).
\end{split}
\end{equation}
Consequently for any $\mathbf{F}\in\mathbf{C}^N[\mathfrak{b}]$
\begin{equation*}
\begin{split}
&\mathds{E}_{\gamma,\mathbf{F}}^N\Bigg[\frac{d\mathds{P}_{\gamma,\SS}^N}{d\mathds{P}_{\gamma,\mathbf{F}}^N}\mathds{1}_{(\rho^N,\mathcal{Q}^N)\in\mathbf{O}_{\eps,\mathbf{h}}^N}\mathds{1}_{\mathbf{B}_{\delta,\eps}^N}\Bigg]\leq\mathds{E}_{\gamma,\mathbf{F}}^N\Bigg[\exp\Bigg(-\frac{N^2}{\log N}\Big(\underbrace{\mathbf{R}_{T,\mathbf{F}}^{N,\eps}(\rho^N,\mathcal{Q}^N)}_{\textup{\blue{dynamical cost}}}+\underbrace{O(\delta)+o_N(1;\mathfrak{b})}_{\blue{\substack{\textup{error from replacement}\\ \textup{and expansion}}}}\Big)\Bigg)\mathds{1}_{(\rho^N,\mathcal{Q}^N)\in\mathbf{O}_{\eps,\mathbf{h}}^N}\mathds{1}_{\mathbf{B}_{\delta,\eps}^N}\Bigg].
\end{split}
\end{equation*}
Finally
\begin{equation}\label{eq:sup_inf}
\begin{split}
&\frac{\log N}{N^2}\log \mathds{E}_{\gamma,\SS}^N\Big[\mathds{1}_{(\rho^N,\mathcal{Q}^N)\in\mathbf{O}_{\eps,\mathbf{h}}^N}\mathds{1}_{\mathbf{B}_{\delta,\eps}^N}\Big]\leq -\sup_{\mathbf{F}\in\mathbf{C}^N[\mathfrak{b}]}\inf_{\rho^N,\mathcal{Q}^N\in\mathbf{O}_{\eps,\mathbf{h}}^N}\Big(\underbrace{\mathbf{R}_{T,\mathbf{F}}^{N,\eps}(\rho^N,\mathcal{Q}^N)}_{\textup{\blue{dynamical cost}}}+O(\delta)+o_N(1;\mathfrak{b})\Big).
\end{split}
\end{equation}

\textit{Step 3: Exchange the order of $\sup$ and $\inf$.}

In the definition of $\mathbf{O}_{\eps,\mathbf{h}}^N$ in \eqref{eq:open_set_for_discrete}, when $|x|\leq N^{-\eps}$ we pin the values of $[\rho_x^N]_{\eps}$ and $\mathcal{Q}_{T,x}^N$ around the coarse-grained function $\mathbf{h}$ with error $O(\eps(\log N)^{-2})$. When $|x|\geq N^{-\eps}$, we truncate the current by $\Psi_{x}^{\eps}$ in \eqref{eq:radial_truncation}, and use the spatial regularity of $[\rho_x^N]_{\eps}$ in Lemma \ref{lem:regularity_of_average}. Thus for any two pairs of profiles $(\rho^{N,*},\mathcal{Q}^{N,*})$ and $(\rho^{N,\dagger},\mathcal{Q}^{N,\dagger})$ in $\doh$, we have
\begin{equation*}
\Big|[\rho_{x}^{N,*}]_{\eps}-[\rho_{x}^{N,\dagger}]_{\eps}\Big|\leq O\Big(\frac{\eps}{\log^2N}\Big),\quad \Big|\Psi_{x}^{\eps}\mathcal{Q}_{T,x}^{N,*}-\Psi_{x}^{\eps}\mathcal{Q}_{T,x}^{N,\dagger}\Big|\leq O\Big(\frac{\eps}{\log^2 N}\Big),\quad \textup{for any $x\in\mathds{B}^N$}.
\end{equation*}
Therefore we can choose a fixed profile $(\hat{\rho}^N,\hat{\mathcal{Q}}^N)\in \doh$ such that for every $\mathbf{F}\in \mathbf{C}^N[\mathfrak{b}]$ uniformly (recall the functions of $\mathbf{F}$ in $\mathbf{R}_{T,\mathbf{F}}^{N,\eps}$ are all bounded by $C_{\mathfrak{b}}|x|^{-2}$)
\begin{equation*}
\sup_{(\rho^N,\mathcal{Q}^N)\in \doh}\Big|\mathbf{R}_{T,\mathbf{F}}^{N,\eps}\big(\rho^N,\mathcal{Q}^N\big)-\mathbf{R}_{T,\mathbf{F}}^{N,\eps}\big(\hat{\rho}^N,\hat{\mathcal{Q}}^N\big)\Big|=O\Big(N^{-2}\sum_{x\in\mathds{B}^N}\frac{\eps}{\log^2N}C_{\mathfrak{b}}|x|^{-2}\Big)=\err(\mathfrak{b},\delta,\eps,N).
\end{equation*}
Consequently we can switch the order of $\sup$ and $\inf$
\begin{equation}\label{eq:exchange_order}
\begin{split}
\sup_{\mathbf{F}\in \mathbf{C}^N[\mathfrak{b}]}\inf_{(\rho^N,\mathcal{Q}^N)\in \doh}\mathbf{R}_{T,\mathbf{F}}^{N,\eps}(\rho^N,\mathcal{Q}^N)&=\sup_{\mathbf{F}\in \mathbf{C}^N[\mathfrak{b}]}\mathbf{R}_{T,\mathbf{F}}^{N,\eps}(\hat{\rho}^N,\hat{\mathcal{Q}}^N)+\err(\mathfrak{b},\delta,\eps,N)\\
&\geq \inf_{(\rho^N,\mathcal{Q}^N)\in \doh} \sup_{\mathbf{F}\in \mathbf{C}^N[\mathfrak{b}]}\mathbf{R}_{T,\mathbf{F}}^{N,\eps}(\rho^N,\mathcal{Q}^N)+\err(\mathfrak{b},\delta,\eps,N).
\end{split}
\end{equation}
Applying \eqref{eq:exchange_order} to \eqref{eq:sup_inf}, we conclude the proof of Lemma \ref{lem:discrete_optimization}.
\end{proof}

\subsection{Coarse-grained Functions and Variational Formulas over Continuum Profiles}\label{sec:coarse_graining_cont}

The variational formula given by Lemma \ref{lem:discrete_optimization} is over the discrete profiles $\big(\rho^N,\mathcal{Q}^N\big)\in \mathbf{O}_{\eps,\mathbf{h}}^N\subset D\big([0,T],\mathbf{M}_{\{0,1\}}^N\times \mathbf{M}_{\mathds{R}^2}^N\big)$, where $\rho_t^N$ is a $\{0,1\}$-valued configuration at each time $t$. Since $\{\doh\}$ forms a cover of the set $\mathbf{S}_{\delta,\eps}^N$ in \eqref{eq:S_sets}, by Lemma \ref{lem:discrete_optimization} we have
\begin{equation}\label{eq:upper_bound_tilte_S}
\frac{\log N}{N^2}\log\mathds{P}_{\gamma,\SS}^N\Big[\mathbf{S}_{\delta,\eps}^N\Big]\leq \Big(-\inf_{\rho^N,\mathcal{Q}^N\in\widetilde{\mathbf{S}}_{\delta,\eps}^N}\sup_{\mathbf{F}\in \mathbf{C}^N[\mathfrak{b}]}\mathbf{R}_{T,\mathbf{F}}^{N,\eps}(\rho^N,\mathcal{Q}^N)\Big)+\err(\mathfrak{b},\delta,\eps,N),\quad \widetilde{\mathbf{S}}_{\delta,\eps}^N:=\bigcup_{\mathbf{h}:\doh\cap \mathbf{S}_{\delta,\eps}^N\not =\varnothing}\doh.
\end{equation}
The set $\widetilde{\mathbf{S}}_{\delta,\eps}^N$ can be understood as an approximation of the set $\mathbf{S}_{\delta,\eps}^N$, and as $\eps,\delta\rightarrow 0$ this approximation is more and more accurate. In this section, by defining a continuum counterpart $\mathcal{S}_{\delta,\eps}^N$ of $\widetilde{\mathbf{S}}_{\delta,\eps}^N$, we transform the upper bound in \eqref{eq:upper_bound_tilte_S} into a variational formula over continuum profiles $(\rho,q)$. With $C^{\sharp}$ as a fixed large enough constant and $\mathcal{M}_{\mathds{R}^2}(B_1)$ as the set of $\mathds{R}^2$-valued finite Borel measure on $B_1$, the set $\mathcal{S}_{\delta,\eps}^N$ is defined as
\begin{equation}\label{eq:continuum_S}
\begin{split}
\mathcal{S}_{\delta,\eps}^N:=\Bigg\{(\rho,q)\in \mathcal{C}^1&\big(B_1;[0,1]\big)\times \mathcal{M}_{\mathds{R}^2}(B_1)\Bigg|\ \int_{N^{-1+\eps}\leq |x|\leq N^{-\eps}}q_x\cdot\mathcal{V}_xdx\geq J-C^{\sharp}(\eps+\delta),\ \textup{supp}(q)\subset \big\{N^{-1+\eps}\leq |x|\leq N^{-\eps}\big\}\\
&|\nabla_x\rho|\leq C^{\sharp}\eps^{-1}\min\Big\{|x|^{-1},N\Big\}\textup{ for $x\in B_1$}; \quad |\rho_x-\gamma|\leq C^{\sharp}\eps\textup{ for $x\in B_1$ with $\eps\leq |x|\leq 1$}\Bigg\}.
\end{split}
\end{equation}
Define $\mathcal{C}^N[\mathfrak{b}]$ as the continuum counterpart of $\mathbf{C}^N[\mathfrak{b}]$ in Definition \ref{def:set_of_forcing}. Notice that $\mathbf{C}^N[\mathfrak{b}]$ is a set of forcings on the discrete set $\mathds{B}^N$, and $\mathcal{C}^N[\mathfrak{b}]$ is a set of forcings on the continuum $B_1$
\begin{equation}\label{eq:continuum_forcing_set}
\begin{split}
&\mathcal{C}^N[\mathfrak{b}]:=\Bigg\{\mathcal{F}\in \mathcal{C}^2(B_1;\mathds{R}^2)\textup{ and radially symmetric}\Big|\quad  |\mathcal{F}_x|\equiv 0\textup{ for $x\in B_1$ with $|x|\leq N^{-1+\frac{\eps}{2}}$ or $|x|\geq \frac{1}{10}$}\\
&\hspace{7cm}|\mathcal{F}_x|\leq \mathfrak{b}|x|^{-1},\ |\nabla \mathcal{F}_x|\leq \mathfrak{b}|x|^{-2},\ |\nabla^2\mathcal{F}_x|\leq \mathfrak{b}|x|^{-3}\textup{ for $x\in B_1$}\Bigg\}.
\end{split}
\end{equation}

\begin{lemma}\label{lem:upper_bound_O_continuum}
For any parameters $(\mathfrak{b}>0,\delta>0,\eps>0,N>1)$, there is 
\begin{equation*}
\begin{split}
&\frac{\log N}{N^2}\log \mathds{P}_{\gamma,\SS}^N\Big[\big(\rho^N,\mathcal{Q}^N\big)\in\mathbf{S}_{\delta,\eps}^N\Big]\leq -\Big[\inf_{(\rho,q)\in\mathcal{S}_{\delta,\eps}^N}\sup_{\mathcal{F}\in \mathcal{C}^N[\mathfrak{b}]}\mathcal{R}_{T,\mathcal{F}}^N(\rho,q)\Big]+\textbf{\textup{Err}}(\mathfrak{b},\delta,\eps,N).
\end{split}
\end{equation*}
Here $\err$ is a generic error described in \eqref{eq:composed_error}, and the functional $\mathcal{R}_{T,\mathcal{F}}^N$ is defined in \eqref{eq:continuum_R_cost}.
\end{lemma}

We can compare the functional $\mathcal{R}_{T,\mathcal{F}}^N(\rho,q)$ in Definition \ref{def:R_functionals} with the functional $\mathbf{R}_{T,\mathbf{F}}^{N,\eps}(\rho^N,\mathcal{Q}^N)$ in \eqref{eq:discrete_R_functional_eps}, which also motivates the construction of $\mathcal{S}_{\delta,\eps}^N$ in \eqref{eq:continuum_S}
\begin{itemize}
    \item The cost $\mathbf{R}_{T,\mathbf{F}}^{N,\eps}$ is written in terms of the space-time average $[\rho_{\cdot}^N]_{\eps}$, and in $\mathcal{R}_{T,\mathcal{F}}^N$ the space-time average is replaced by a continuum profile $\rho$. The regularity and boundary condition of $\rho$ in $\mathcal{S}_{\delta,\eps}^N$ is due to the spatial regularity (Lemma \ref{lem:regularity_of_average}) and the boundary condition (see \eqref{eq:S_sets}) of $[\rho_{\cdot}^N]_{\eps}$.
    \item The cost of current profile $\mathcal{Q}_T^N$ in $\mathbf{R}_{T,\mathbf{F}}^{N,\eps}(\rho^N,\mathcal{Q}^N)$ is written in terms of the truncated current profile $\Psi^{\eps}\mathcal{Q}_T^N$. Here $\Psi^{\eps}$ is the truncation function in \eqref{eq:radial_truncation}. In $\mathcal{R}_{T,\mathcal{F}}^N(\rho,q)$, the truncated current profile $\Psi^{\eps}\mathcal{Q}$ is replaced by a finite $\mathds{R}^2$-valued Borel measure $q\in \mathcal{M}_{\mathds{R}^2}(B_1)$. According to the definition of $\mathbf{S}_{\delta,\eps}^N$ in \eqref{eq:S_sets}, the truncation $\Psi^{\eps}$ only produces an error of order $O(\delta)$ to the winding number, therefore in $\mathcal{S}_{\delta,\eps}^N$ we have $\int_{N^{-1+\eps}\leq |x|\leq N^{-\eps}}q_x\cdot \mathcal{V}_xdx\geq J-C^{\sharp}(\eps+\delta)$.
\end{itemize}

\begin{proof}[Proof of Lemma \ref{lem:upper_bound_O_continuum}]

\textit{Step 1: Embedding of $\widetilde{\mathbf{S}}_{\delta,\eps}^N$ into $\mathcal{S}_{\delta,\eps}^N$.}

To prove Lemma \ref{lem:upper_bound_O_continuum}, it suffices to construct a mapping $\mathcal{P}$ from $\widetilde{\mathbf{S}}_{\delta,\eps}^N$ to $\mathcal{S}_{\delta,\eps}^N$ such that for any $(\rho^N,\mathcal{Q}^N)\in \tildeS$
\begin{equation}\label{eq:d_to_c_forcing_matching}
\sup_{\mathbf{F}\in \mathbf{C}^N[\mathfrak{b}]}\mathbf{R}_{T,\mathbf{F}}^{N,\eps}(\rho^N,\mathcal{Q}^N)\geq \sup_{\mathcal{F}\in \mathcal{C}^N[\mathfrak{b}]}\mathcal{R}_{T,\mathcal{F}}^N\big(\mathcal{P}(\rho^N),\mathcal{P}(\mathcal{Q}^N)\big)+\errtotal.
\end{equation}
Intuitionally, Equation \eqref{eq:d_to_c_forcing_matching} implies Lemma \ref{lem:upper_bound_O_continuum} because $\mathcal{P}$ as an embedding of $\widetilde{\mathbf{S}}_{\delta,\eps}^N$ into $\mathcal{S}_{\delta,\eps}^N$, means $\mathcal{S}_{\delta,\eps}^N$ is a larger set than $\tildeS$. Since we take $\inf$ over the set of profiles, a larger set of profiles means a lower cost. Rigorously
\begin{equation*}
\begin{split}
\inf_{\rho^N,\mathcal{Q}^N\in\widetilde{\mathbf{S}}_{\delta,\eps}^N}\sup_{\mathbf{F}\in \mathbf{C}^N[\mathfrak{b}]}\mathbf{R}_{T,\mathbf{F}}^{N,\eps}(\rho^N,\mathcal{Q}^N)&\geq \inf_{\rho^N,\mathcal{Q}^N\in\widetilde{\mathbf{S}}_{\delta,\eps}^N}\underbrace{\sup_{\mathcal{F}\in \mathcal{C}^N[\mathfrak{b}]}\mathcal{R}_{T,\mathcal{F}}^{N}\big(\mathcal{P}(\rho^N),\mathcal{P}(\mathcal{Q}^N)\big)}_{\textup{\blue{derived by \eqref{eq:d_to_c_forcing_matching}}}}+\errtotal\\
&\geq \inf_{(\rho,q)\in \mathcal{P}(\widetilde{\mathbf{S}}_{\delta,\eps}^N)\subset \mathcal{S}_{\delta,\eps}^N}\sup_{\mathcal{F}\in \mathcal{C}^N[\mathfrak{b}]}\mathcal{R}_{T,\mathcal{F}}^{N}\big(\rho,q\big)+\errtotal.
\end{split}
\end{equation*}
For any discrete profile $(\rho^N,\mathcal{Q}^N)\in \widetilde{\mathbf{S}}_{\delta,\eps}^N$, we construct $\mathcal{P}(\mathcal{Q}^N)$ as a $\mathds{R}^2$-valued finite Borel measure, and $\mathcal{P}(\rho^N)$ as a differentiable function on $B_1$
\begin{equation}\label{eq:construction_of_P}
\begin{split}
&\big(\mathcal{P}(\mathcal{Q}^N)\big)_x:=T^{-1}N^{-2}\big(1-\mathcal{B}_{+,x}^{\eps}-\mathcal{B}_{-,x}^{\eps}\big)\mathcal{Q}_{T,x}^N\delta_{x},\\
&\mathcal{P}(\rho^N):=\textup{ one function in $\mathcal{S}_{\delta,\eps}^N$ satisfying }\mathcal{P}(\rho^N)=[\rho_{x}^N]_{\eps}\textup{ for any $x\in \mathds{B}^N$ with $|x|\leq 1-2\eps$}.
\end{split}
\end{equation}
We abbreviate $(\mathcal{P}(\rho^N),\mathcal{P}(\mathcal{Q}^N))$ as $(\rho,q)$. Now we will verify that the profile $(\rho,q)$ constructed in \eqref{eq:construction_of_P} satisfies \eqref{eq:d_to_c_forcing_matching}. 

\textit{Step 2: Embedding of $\mathcal{C}^N[\mathfrak{b}]$ into $\mathbf{C}^N[\mathfrak{b}]$.}

Recall that we abbreviate $(\mathcal{P}(\rho^N),\mathcal{P}(\mathcal{Q}^N))$ as $(\rho,q)$. To prove \eqref{eq:d_to_c_forcing_matching}, for any forcing $\mathcal{F}\in \mathcal{C}^N[\mathfrak{b}]$, we construct a forcing $\mathbf{F}\in \mathbf{C}^N[\mathfrak{b}]$ satisfying
\begin{equation*}
\mathbf{R}_{T,\mathbf{F}}^{N,\eps}(\rho^N,\mathcal{Q}^N)\geq \mathcal{R}_{T,\mathcal{F}}^N(\rho,q)+\errtotal.
\end{equation*}
Decompose the forcing $\mathcal{F}$ into radial and rotation components $\mathcal{F}=\mathcal{F}_r\bm{r}+\mathcal{F}_{\theta}\bm{\theta}$, analogous to the decomposition in \eqref{eq:radial_rotation_decomposition}. The forcing $\mathcal{F}$ at any $x\in B_1$ also equals to
\begin{equation*}
\mathcal{F}=\mathcal{F}_{\theta}\bm{\theta}+\nabla V,\quad V(x):=\int_{N^{-1}}^{|x|}\mathcal{F}_r\Big(L\frac{(x_1,x_2)}{|x|}\Big)dL.
\end{equation*}
Construct the discrete forcing $\mathbf{F}$ at any $x\in\mathds{B}^N$ as
\begin{equation}\label{eq:construction_of_bfF}
    \mathbf{F}_x=\mathcal{F}_{\theta}(x)\bm{\theta}_x+\nabla_{\textup{disc}}V_x.
\end{equation}
Here $\mathcal{F}$ and $\mathbf{F}$ have a common rotation field $\mathcal{F}_{\theta}\bm{\theta}$, while the difference is the $\nabla V$ in $\mathcal{F}$ and the $\nabla_{\textup{disc}}V$ in $\mathbf{F}$. The regularity of $\mathcal{F}$ in \eqref{eq:continuum_forcing_set} implies that the function $V$ satisfies the following, where heuristically the $n$-th gradient of $V$ is the $(n-1)$-th gradient of $\mathcal{F}$
\begin{equation*}
|\nabla V_x|\leq \mathfrak{b}\min\{|x|^{-1},N\},\quad |\nabla^2 V_x|\leq \mathfrak{b}\min\{|x|^{-2},N^2\}.
\end{equation*}
We also have
\begin{equation}\label{eq:error_of_forcing}
\big|\mathcal{F}_x-\mathbf{F}_x\big|=\big|\nabla_{\textup{disc}}V_x-\nabla V_x\big|=O\Big(\frac{1}{N}|x|^{-2}\Big),\quad \big|\nabla_{\textup{disc}}\cdot\mathbf{F}_x-\nabla\cdot \mathcal{F}_x\big|\leq \big|\nabla_{\textup{disc}}\cdot \mathbf{F}^{\textup{rot}}\big|+\big|\Delta_{\textup{disc}}V-\Delta V\big|=O\Big(\frac{1}{N}|x|^{-3}\Big).
\end{equation}
We are going to prove the error estimate from $(\rho^N,\mathcal{Q}^N,\mathbf{F})$ to $(\rho,q,\mathcal{F})$:
\begin{subequations}\label{eq:dc_matching}
\begin{align}
&\big\langle (1-\mathcal{B}_{+}^{\eps}-\mathcal{B}_{-}^{\eps})\mathcal{Q}_T^N,\mathbf{F}^{\textup{rot}}\big\rangle-T\langle q,\mathcal{F}^{\textup{rot}}\rangle=0, \label{eq:dc_matching_a}\\
&\frac{T}{2N^2}\sum_{x\in\mathds{B}^N}\frac{|\mathbf{F}_x|^2}{\log N}\sigma([\rho_{x}^N]_{\eps})-\frac{T}{2N^2}\sum_{x\in\mathds{B}^N}\frac{|\mathcal{F}_x|^2}{\log N}\sigma(\rho_x)=\err(\mathfrak{b},\delta,\eps,N), \label{eq:dc_matching_b}\\
&\frac{T}{2N^2}\sum_{x\in\mathds{B}^N}[\rho_x^N]_{\eps}\nabla_{\textup{disc}}\cdot\mathbf{F}_x-\frac{T}{2N^2}\sum_{x\in\mathds{B}^N}\rho_x\nabla\cdot\mathcal{F}_x=\err(\mathfrak{b},\delta,\eps,N). \label{eq:dc_matching_c}
\end{align}
\end{subequations}
The proof of \eqref{eq:dc_matching_a} is direct, because the construction of $q=T^{-1}(1-\mathcal{B}_+^{\eps}-\mathcal{B}_{-}^{\eps})\mathcal{Q}_T^N$ in \eqref{eq:construction_of_P}. To prove  \eqref{eq:dc_matching_b} and \eqref{eq:dc_matching_c}, we first notice $\rho_x$ strictly equals to $[\rho_x^N]_{\eps}$ at any $x\in\mathds{B}^N$ with $|x|\leq 1-2\eps$ according to \eqref{eq:construction_of_P}. Therefore we only need to measure the error of replacing $\mathbf{F}$ by $\mathcal{F}$, see \eqref{eq:error_of_forcing}
\begin{equation*}
\begin{split}
&\frac{T}{2N^2}\sum_{x\in\mathds{B}^N}\frac{|\mathbf{F}_x|^2}{\log N}\sigma([\rho_{x}^N]_{\eps})-\frac{T}{2N^2}\sum_{x\in\mathds{B}^N}\frac{|\mathcal{F}_x|^2}{\log N}\sigma(\rho_x)\\
&=O\Big(N^{-3}\sum_{x\in\mathds{B}^N:|x|\geq N^{-1+\frac{\eps}{2}}}\frac{|x|^{-3}}{\log N}\Big)=O\Big(\frac{1}{N}\int_{N^{-1+\frac{\eps}{2}}\leq |x|\leq 1}\frac{|x|^{-3}}{\log N}dx\Big)=\errtotal,
\end{split}
\end{equation*}
and the same holds for \eqref{eq:dc_matching_c}.

Next we estimate the error from the discrete summation to the continuum integral
\begin{subequations}\label{eq:dc_matching_cont}
\begin{align}
&\frac{T}{2N^2}\sum_{x\in\mathds{B}^N}\frac{|\mathcal{F}_x|^2}{\log N}\sigma(\rho_x)-\frac{T}{2}\int_{N^{-1}\leq |x|\leq 1}\frac{|\mathcal{F}_x|^2}{\log N}\sigma(\rho_x)dx=\err(\mathfrak{b},\delta,\eps,N), \label{eq:dc_matching_cont_b}\\
&\frac{T}{2N^2}\sum_{x\in\mathds{B}^N}\rho_x\nabla\cdot\mathcal{F}_x-\frac{T}{2}\int_{N^{-1}\leq |x|\leq 1}\rho_x\nabla\cdot\mathcal{F}_xdx=\err(\mathfrak{b},\delta,\eps,N). \label{eq:dc_matching_cont_c}
\end{align}
\end{subequations}
The proof of \eqref{eq:dc_matching_cont_b} and \eqref{eq:dc_matching_cont_c} can be unified as follows. Suppose $\psi:[0,1]\rightarrow \mathds{R}$ is a smooth function, which may either represent $\psi(x)=x$ or $\psi(x)=\sigma(x)$. Let $\mathscr{F}$ be a function satisfying
\begin{equation}\label{eq:scrF_regularity}
\Big|\mathscr{F}_x\Big|\leq C_{\mathfrak{b}}|x|^{-2}\mathds{1}_{|x|\geq N^{-1+\frac{\eps}{2}}},\quad |\nabla \mathscr{F}_x|\leq C_{\mathfrak{b}}|x|^{-3}\mathds{1}_{|x|\geq N^{-1+\frac{\eps}{2}}}.
\end{equation}
The function $\mathscr{F}$ can either represent $|\mathcal{F}_x|^2(\log N)^{-1}$ or $\nabla \cdot\mathcal{F}$. The spatial regularity of $\rho$ in $\mathcal{S}_{\delta,\eps}^N$ is given in \eqref{eq:continuum_S}, and the spatial regularity of $\mathscr{F}$ is given in \eqref{eq:scrF_regularity}. These together derive
\begin{equation*}
\frac{1}{N^2}\sum_{x\in\mathds{B}^N}\psi\big(\rho_x\big)\mathscr{F}_x=\int_{B_1}\psi(\rho_x)\mathscr{F}_xdx+O\Big(N^{-2}\sum_{x\in\mathds{B}^N:N^{-1+\frac{\eps}{2}}\leq |x|\leq 1}C_{\mathfrak{b}}\frac{|x|^{-3}}{N}\Big)=\int_{B_1}\psi(\rho_x)\mathscr{F}_xdx+\err(\mathfrak{b},\delta,\eps,N).
\end{equation*}
This prove that we can go from $\mathbf{R}_{T,\mathbf{F}}^{N,\eps}$ to $\mathcal{R}_{T,\mathcal{F}}^N$ with a negligible error
\begin{equation*}
\frac{\log N}{N^2}\log\mathds{P}_{\gamma,\SS}^N\Big[(\rho^N,\mathcal{Q}^N)\in \mathbf{S}_{\delta,\eps}^N\Big]\leq -\inf_{(\rho,q)\in \mathcal{S}_{\delta,\eps}^N}\sup_{\mathcal{F}\in \mathcal{C}^N[\mathfrak{b}]}\mathcal{R}_{T,\mathcal{F}}^N(\rho,q)+\err(\mathfrak{b},\delta,\eps,N).
\end{equation*}
This concludes the proof of Lemma \ref{lem:upper_bound_O_continuum}.
\end{proof}

\subsection{Proof of Singular Large Deviation Upper Bound in Theorem \ref{th:one-scale-deviation}}\label{sec:true_proof_upper}

\begin{theorem}\label{th:theorem_upper_bound}
For any parameters $\big(J>0,T>0,\gamma\in(0,1)\big)$, we have the large deviation upper bound for the winding number
\begin{equation}\label{eq:upper_bound_theorem}
\begin{split}
&\limsup_{N\rightarrow +\infty}\frac{\log N}{N^2}\log \mathds{P}_{\gamma,\SS}^N\Big[\bl\mathcal{Q}_T^N,\mathcal{V}\br\geq JT\Big]\leq -T\inf_{(\rho,q)}\Bigg(\frac{1}{2}\int_{e^{-1}\leq |x|\leq 1}\frac{|q_x|^2}{\sigma(\rho_x)}dx+\frac{1}{8}\int_{e^{-1}\leq |x|\leq 1}\frac{|\nabla \rho_x|^2}{\sigma(\rho_x)}dx\Bigg),
\end{split}
\end{equation}
where the profile $(\rho,q)$ is taken over
\begin{equation}\label{eq:th4.9_range}
\rho\in \mathcal{C}^1\big(B_1;[0,1]\big),\ \rho\equiv \gamma \textup{ on $\p B_1$};\quad  q\textup{ as function $B_1\mapsto \mathds{R}^2$} \quad \int_{e^{-1}\leq |x|\leq 1}q_x\cdot \mathcal{V}_x dx\geq J.
\end{equation}
This is exactly the large deviation upper bound stated in Theorem \ref{th:one-scale-deviation}.
\end{theorem}

\begin{proof}
\textit{Step 1: Radially symmetric profiles.} 

Recall the decomposition \eqref{eq:radial_rotation_decomposition} of $q$ into a gradient component $q_{r}(x)\bm{r}_x$ and a rotation component $q_{\theta}(x)\boldsymbol{\theta}_x$. Given any density profile $\rho$ and current profile $q$, define $(\overline{\rho},\overline{q})$ as their radial symmetrization
\begin{equation*}
\begin{split}
&\overline{\rho}(x):=\frac{1}{2\pi|x|}\int_{|x'|=|x|}\rho(x')dx',\\
&\overline{q}_{r}(x):=\frac{1}{2\pi |x|}\int_{|x'|=|x|}q_{r}(x')dx',\quad \overline{q}_{\theta}(x):=\frac{1}{2\pi|x|}\int_{|x'|=|x|}q_{\theta}(x')dx',\quad \overline{q}(x)=\overline{q}_{r}(x)\bm{r}_x+\overline{q}_{\theta}(x)\boldsymbol{\theta}_x.
\end{split}
\end{equation*}
Since the functional $\mathcal{R}_{T,\mathcal{F}}^N(\rho,q)$ in \eqref{eq:continuum_R_cost} is convex in $(\rho,q)$, we apply the Jensen inequality $\mathcal{R}_{T,\mathcal{F}}^N(\overline{\rho},\overline{q})\leq \mathcal{R}_{T,\mathcal{F}}^N(\rho,q)$ to Lemma \ref{lem:upper_bound_O_continuum}, which allows us to restrict the variational formula over radially symmetric profiles
\begin{equation}\label{eq:q-theta}
\begin{split}
&\frac{\log N}{N^2}\log \mathds{P}_{\gamma,\SS}^N\Big[\big(\rho^N,\mathcal{Q}^N\big)\in\mathbf{S}_{\delta,\eps}^N\Big]\leq -\Big[\inf_{(\rho,q)\in\overline{\mathcal{S}}_{\delta,\eps}^N}\sup_{\mathcal{F}\in \mathcal{C}^N[\mathfrak{b}]}\mathcal{R}_{T,\mathcal{F}}^N(\rho,q)\Big]+\textbf{\textup{Err}}(\mathfrak{b},\delta,\eps,N),\\
&\overline{\mathcal{S}}_{\delta,\eps}^N:=\Big\{(\rho,q)\Big|\ \textup{exists $(\rho_*,q_*)\in \mathcal{S}_{\delta,\eps}^{N}$ such that $(\overline{\rho}_*,\overline{q}_*)=(\rho,q)$}\Big\}.
\end{split}
\end{equation}

\textit{Step 2: Transformation into self-similar variables.} 

However the set $\overline{\mathcal{S}}_{\delta,\eps}^N$ in \eqref{eq:q-theta} depends on $N$. To prove the $N$-independent upper bound in Theorem \ref{th:theorem_upper_bound}, we need to re-write the variational formula \eqref{eq:q-theta} in a $N$-independent form. Therefore similar to Lemma \ref{lem:C_self_similarity}, we want to write the cost functional and all profiles and forcings in the self-similar variable $|x|=N^{-s}$. In the self-similar variable, functions $(\rho,q_r,q_{\theta})$ is written as $(\uprho,\mathtt{q}_r,\mathtt{q}_{\theta})$, see \eqref{eq:self_similar_rho_q}
\begin{equation}\label{eq:profile_self_similar_sec4}
\mathtt{q}_{r}(s):=\Big(|x|\log N q_r(x)\Big)\Big|_{|x|=N^{-s}},\quad \mathtt{q}_{\theta}(s):=\Big(|x|\log N q_{\theta}(x)\Big)\Big|_{|x|=N^{-s}},\quad \uprho(s):=\rho(x)\big|_{|x|=N^{-s}}.
\end{equation}
Forcing $\mathcal{F}:=(\mathcal{F}_r,\mathcal{F}_{\theta})$ is written as $\mathtt{F}:=(\mathtt{F}_r,\mathtt{F}_{\theta})$
\begin{equation}\label{eq:forcing_self_similar}
\mathtt{F}_{r}(s):=\Big(|x| \mathcal{F}_r(x)\Big)\Big|_{|x|=N^{-s}},\quad \mathtt{F}_{\theta}(s):=\Big(|x| \mathcal{F}_{\theta}(x)\Big)\Big|_{|x|=N^{-s}}.
\end{equation}
Following the proof of Lemma \ref{lem:C_self_similarity} on the self-similarity of $\mathcal{C}^N$-functional in Definition \ref{def:cost_functional_lower}, we get for any radially symmetric profile $(\rho,q)$
\begin{equation}\label{eq:to_self_similar_R_functional}
\mathcal{R}_{T,\mathcal{F}}^N(\rho,q)=T\Bigg(2\pi\int_0^1\mathtt{q}_{\theta}(s)\mathtt{F}_{\theta}(s)ds -\pi\int_{0}^1\Big(\mathtt{F}_r^2(s)+\mathtt{F}_{\theta}^2(s)\Big)\sigma(\uprho(s))ds-\pi\int_0^1\frac{d}{ds}\uprho(s)\mathtt{F}_r(s)ds\Bigg)=:\mathtt{R}_{T,\mathtt{F}}(\uprho,\mathtt{q}),
\end{equation}
In $\mathtt{R}_{T,\mathtt{F}}(\uprho,\mathtt{q})$, the cost related with the current profile $\mathtt{q}$ is only $2\pi T\int_0^1 \mathtt{q}_{\theta}(s)\mathtt{F}_{\theta}(s)ds$, independent of the radial current $\mathtt{q}_r$. This is because in $\mathcal{R}_{T,\mathcal{F}}^N(\rho,q)$ has current-related cost as $\langle q,\mathcal{F}^{\textup{rot}}\rangle=\langle q_{\theta}\boldsymbol{\theta},\mathcal{F}^{\textup{rot}}\rangle$, where we only see the cost of the rotational current.

Equation \eqref{eq:to_self_similar_R_functional} transforms the cost functional $\mathcal{R}_{T,\mathcal{F}}^N$ in physical variables $x\in B_1$, into the cost functional $\mathtt{R}_{T,\mathtt{F}}$ in self-similar variable $s\in [0,1]$. Therefore we can transform the upper bound \eqref{eq:q-theta} into self-similar variables, with profiles and forcings written in self-similar variable $s$
\begin{equation}\label{eq:upper_bound_self_similar_witherror}
\begin{split}
&\frac{\log N}{N^2}\log\mathds{P}_{\gamma,\SS}^N\Big[\langle\mathcal{Q}_T^N,\mathcal{V}\rangle\geq JT\Big]\leq -\inf_{(\uprho,\mathtt{q})\in \mathtt{S}_{\delta,\eps}}\sup_{\mathtt{F}\in \mathtt{C}_{\eps}[\mathfrak{b}]}\mathtt{R}_{T,\mathtt{F}}(\uprho,\mathtt{q})+\err(\mathfrak{b},\delta,\eps,N),
\end{split}
\end{equation}
with set $\mathtt{S}_{\delta,\eps}$ defined as
\begin{equation}\label{eq:upper_bound_self_similar_witherror_S}
\begin{split}
&\mathtt{S}_{\delta,\eps}:=\Big\{(\uprho,\mathtt{q})\Big|\ \underbrace{\uprho\in \mathcal{C}^1\big([0,1];[0,1]\big),\ \Big|\uprho(s)-\gamma\Big|\leq C^{\sharp}\eps\textup{ for $0\leq s\leq \frac{-\log \eps}{\log N}$};}_{\textup{\blue{regularity and boundary condition of the density profile}}}\\
&\hspace{3cm}\underbrace{2\pi\int_{\eps}^{1-\eps}\mathtt{q}_{\theta}(s)ds\geq J-C^{\sharp}(\eps+\delta),\ \textup{supp}(\mathtt{q}_{\theta})\subset [\eps,1-\eps],\ \mathtt{q}_r\equiv 0}_{\textup{\blue{winding number condition and no radial current}}}\Big\}
\end{split}
\end{equation}
and set $\mathtt{C}_{\eps}[\mathfrak{b}]$ as
\begin{equation}\label{eq:upper_bound_self_similar_witherror_F}
\begin{split}
&\mathtt{C}_{\eps}[\mathfrak{b}]:=\Big\{\mathtt{F}=(\mathtt{F}_{r},\mathtt{F}_{\theta})\in \mathcal{C}^2([0,1];\mathds{R}^2)\Big|\ \lVert \mathtt{F}\rVert_{L^{\infty}}\leq \frac{\mathfrak{b}}{10},\quad \Big\lVert \frac{d}{ds}\mathtt{F}(s)\Big\rVert_{L^{\infty}}\leq \frac{\mathfrak{b}}{10}\log N,\quad \Big\lVert \frac{d^2}{ds^2}\mathtt{F}(s)\Big\rVert_{L^{\infty}}\leq \frac{\mathfrak{b}}{10}\log^2 N;\\
&\hspace{8cm}\mathtt{F}_r(s)=\mathtt{F}_{\theta}(s)\equiv 0\textup{ when $1-\frac{\eps}{2}\leq s\leq 1$ or $0\leq s\leq \frac{-\log \frac{1}{10}}{\log N}$} \Big\}.
\end{split}
\end{equation}
Similar to the proof of Lemma \ref{lem:upper_bound_O_continuum}, Equation \eqref{eq:upper_bound_self_similar_witherror} is derived from the following facts: 
\begin{itemize}
    \item \textbf{\textup{[Embedding of profiles]}} For any $(\rho,q)\in \overline{\mathcal{S}}_{\delta,\eps}^N$ in the physical variable $x$, the profile $(\uprho,\mathtt{q})$ (with exceptionally $\mathtt{q}_r\equiv 0$) in the self-similar variable satisfies $(\uprho,\mathtt{q})\in \mathtt{S}_{\delta,\eps}$.
    \item \textbf{\textup{[Embedding of forcings]}} For any forcing $\mathtt{F\in \mathtt{C}_{\eps}[\mathfrak{b}]}$ in the self-similar variable $s$, the forcing $\mathcal{F}$ in the physical variable $x$ satisfies $\mathcal{F}\in \mathcal{C}^N[\mathfrak{b}]$
\end{itemize}
\begin{remark}
The variational formulas in \eqref{eq:q-theta} and \eqref{eq:upper_bound_self_similar_witherror} take infimum over profiles and supremum over forcings. Heuristically, because $\overline{\mathcal{S}}_{\delta,\eps}^N$ is embedded into $\mathtt{S}_{\delta,\eps}$, in \eqref{eq:upper_bound_self_similar_witherror} we take infimum over profiles in a larger set, which makes the cost lower. And because $\mathtt{C}_{\eps}[\mathfrak{b}]$ is embedded into $\mathcal{C}^N[\mathfrak{b}]$, in \eqref{eq:upper_bound_self_similar_witherror} we take supremum over forcings in a smaller set, which makes the cost lower. 
\end{remark}

\textit{Step 2.1: Embedding of profiles.}

We first verify the embedding of profiles. The density profiles satisfy
\begin{equation*}
\begin{split}
 \rho(x)\in [0,1]\quad &\Longrightarrow \quad \uprho(s)=\rho(N^{-s})\in [0,1],\\
|\nabla_x\rho|\leq C^{\sharp}\eps^{-1}\min\{|x|^{-1},N\}\quad &\Longrightarrow \quad \Big|\frac{d}{ds}\uprho(s)\Big|=\Big|\frac{d}{ds}\rho(x)|_{|x|=N^{-s}}\Big|\Longrightarrow \uprho \in \mathcal{C}^1\big([0,1];[0,1]\big)\\
\Big|\rho(x)-\gamma\Big|\leq C^{\sharp}\eps\textup{ for $\eps\leq |x|\leq 1$}\quad &\Longrightarrow\quad  \Big|\uprho(s)-\gamma\Big|\leq C^{\sharp}\eps\textup{ for $0\leq s\leq \frac{-\log\eps}{\log N}$}.
\end{split}
\end{equation*}
The current profiles satisfy
\begin{equation*}
\int_{N^{-1+\eps}\leq |x|\leq N^{-\eps}}q_x\cdot\mathcal{V}_xdx=\int_{N^{-1+\eps}\leq |x|\leq N^{-\eps}}q_{\theta,x}\frac{(-x_2,x_1)}{|x|}\cdot \frac{1}{|x|}\frac{(-x_2,x_1)}{|x|}dx=\int_{N^{-1+\eps}\leq |x|\leq N^{-\eps}}q_{\theta,x}\frac{1}{|x|}dx
\end{equation*}
and
\begin{equation*}
\begin{split}
\int_{N^{-1+\eps}\leq |x|\leq N^{-\eps}}q_{\theta,x}\frac{1}{|x|}dx=2\pi\int_{N^{-1+\eps}}^{N^{-\eps}}q_{\theta,x}|_{|x|=r}dr&=2\pi\int_{\eps}^{1-\eps}q_{\theta,x}|_{|x|=N^{-s}}\Big|\frac{d|x|}{ds}\Big|ds\\
&=2\pi\int_{\eps}^{1-\eps}q_{\theta,x}|_{|x|=N^{-s}}N^{-s}\log Nds=2\pi\int_{\eps}^{1-\eps}\mathtt{q}_{\theta}(s)ds.
\end{split}
\end{equation*}
Therefore 
\begin{equation*}
\int_{N^{-1+\eps}\leq |x|\leq N^{-\eps}}q_x\cdot\mathcal{V}_xdx\geq J-C^{\sharp}(\eps+\delta)\quad \Longrightarrow \quad 2\pi\int_{\eps}^{1-\eps}\mathtt{q}_{\theta}(s)ds\geq J-C^{\sharp}(\eps+\delta).
\end{equation*}

\textit{Step 2.2: Embedding of forcings.}

Now we verify the embedding of forcings. 
\begin{equation}
\frac{d}{d|x|}\mathcal{F}(x)=\frac{ds}{d|x|}\frac{d}{ds}\big(N^{s}\mathtt{F}(s)\big)= -\frac{|x|^{-1}}{\log N} \Big(N^{s}\mathtt{F}(s)\log N+N^{s}\frac{d}{ds}\mathtt{F}(s)\Big)=-|x|^{-2}\Big(\mathtt{F}(s)+\frac{1}{\log N}\frac{d}{ds}\mathtt{F}(s)\Big)
\end{equation}
and
\begin{equation}
\begin{split}
\Big|\frac{d}{d|x|}\Big(\frac{d}{d|x|}\mathcal{F}(x)\Big)\Big|&\leq \Bigg|\frac{ds}{d|x|}\frac{d}{ds}\Bigg(\frac{N^s}{\log N} \Big(N^{s}\mathtt{F}(s)\log N+N^{s}\frac{d}{ds}\mathtt{F}(s)\Big)\Bigg)\Bigg|\\
&\leq \frac{|x|^{-1}}{\log N}\Bigg|2N^{2s}\mathtt{F}(s)\log N+3N^{2s}\frac{d}{ds}\mathtt{F}(s)+\frac{N^{2s}}{\log N}\frac{d^2}{ds^2}\mathtt{F}(s)\Bigg|\\
&\leq 2|x|^{-3}|\mathtt{F}(s)|+\frac{3|x|^{-3}}{\log N}\Big|\frac{d}{ds}\mathtt{F}(s)\Big|+\frac{|x|^{-3}}{\log^2N}\Big|\frac{d^2}{ds^2}\mathtt{F}(s)\Big|.
\end{split}
\end{equation}

Therefore the condition in \eqref{eq:upper_bound_self_similar_witherror_F} for forcings $\mathtt{F}(s)$ in the self-similar variable $s$ implies that, the corresponding forcings in the physical variable $|x|$ satisfy the conditions in \eqref{eq:continuum_forcing_set}.

\textit{Step 3: Absorbing errors and eliminating profiles not absolutely continuous to the Lebesgue measure.}

Next we take $N\rightarrow +\infty$, $\eps\rightarrow 0$, $\delta\rightarrow 0$, and $\mathfrak{b}\rightarrow +\infty$ in order. The error term $\errtotal$ vanishes, and the errors in the definition of $\mathtt{S}_{\delta,\eps}$ and $\mathtt{C}_{\eps}[\mathfrak{b}]$ also disappears 
\begin{equation}\label{eq:self_similar_variational_sec4}
\begin{split}
&\limsup_{N\rightarrow +\infty}\frac{\log N}{N^2}\log\mathds{P}_{\gamma,\SS}^N\Big[\langle\mathcal{Q}_T^N,\mathcal{V}\rangle\geq JT\Big]\leq -\limsup_{\mathfrak{b}\rightarrow +\infty}\inf_{(\uprho,\mathtt{q})\in \mathtt{S}}\sup_{\mathtt{F}\in \mathtt{C}[\mathfrak{b}]}\mathtt{R}_{T,\mathtt{F}}(\uprho,\mathtt{q}),\\
& \mathtt{S}:=\Big\{(\uprho,\mathtt{q})\Big|\ \uprho\in \mathcal{C}^1\big([0,1];[0,1]\big),\quad \uprho(0)=\gamma,\quad 2\pi\int_{0}^{1}\mathtt{q}_{\theta}(s)ds\geq J,\quad\mathtt{q}_r\equiv 0\Big\},\\
& \mathtt{C}[\mathfrak{b}]:=\Big\{\mathtt{F}=(\mathtt{F}_{r},\mathtt{F}_{\theta})\in \mathcal{C}^2([0,1];\mathds{R}^2)\Big|\ \lVert \mathtt{F}\rVert_{L^{\infty}}\leq \mathfrak{b} \Big\}.
\end{split}
\end{equation}
Now we prove
\begin{equation}\label{eq:exchange_order_b_parameter}
\limsup_{\mathfrak{b}\rightarrow +\infty}\inf_{(\uprho,\mathtt{q})\in \mathtt{S}}\sup_{\mathtt{F}\in \mathtt{C}[\mathfrak{b}]}\mathtt{R}_{T,\mathtt{F}}(\uprho,\mathtt{q})=\inf_{(\uprho,\mathtt{q})\in \mathtt{S}}\sup_{\mathtt{F}\in \mathtt{C}[\infty]}\mathtt{R}_{T,\mathtt{F}}(\uprho,\mathtt{q}),
\end{equation}
For the RHS of \eqref{eq:exchange_order_b_parameter}, there exists profiles $(\uprho,\mathtt{q})\in \mathtt{S}$ whose rotational current is not absolutely continuous to the Lebesgue measure on $[0,1]$. Now we prove for such a profile $(\uprho,\mathtt{q})$, its cost is infinite $\sup_{\mathtt{F}\in \mathtt{C}[\infty]}\mathtt{R}_{T,\mathtt{F}}(\uprho,\mathtt{q})=+\infty$. Lebesgue decomposition theorem can decompose the measure $\mathtt{q}_{\theta}$ into an absolutely continuous (to Lebesgue measure) part and an singular part $\mathtt{q}_{\theta}=\mathtt{q}_{\theta}^{\textup{cont}}+\mathtt{q}_{\theta}^{\textup{sing}}$. Jordan decomposition theorem can further decompose $\mathtt{q}_{\theta}^{\textup{sing}}$ into a positive measure and a negative measure $\mathtt{q}_{\theta}^{\textup{sing}}=\mathtt{q}_{\theta}^{\textup{sing},+}+\mathtt{q}_{\theta}^{\textup{sing},-}$. 

Denote $S^+$ as the support of the positive Borel measure $\mathtt{q}_{\theta}^{\textup{sing},+}$. \cite[Theorem 2.18]{rudin1974real} implies there exists a closed set $C^+\subset S^+$ and a constant $c$ such that $\int \mathtt{q}_{\theta}^{\textup{sing},+}(s)\mathds{1}_{C^+}(s)ds>c$. Choose $\psi:[0,+\infty)\rightarrow [0,1]$ as a smooth function with $\psi(0)=1$ and $\textup{supp}(\psi)=[0,1]$. Define forcing $\mathtt{F}^{\eps}$ as the forcing with zero radial component $\mathtt{F}_{r}^{\eps}\equiv 0$ and rotational component as
\begin{equation*}
\mathtt{F}_{\theta}^{\eps}(s)=\psi\Big(\eps^{-1}\textup{dist}(s,C^+)\Big)\textup{ which satisfies $\lim_{\eps\rightarrow 0}\mathtt{F}_{\theta}^{\eps}(s)=\mathds{1}_{C^+}(s)$}.
\end{equation*}
For any $K>0$ as $\eps\rightarrow 0$, the cost in $\mathtt{R}_{T,K\mathtt{F}^{\eps}}(\uprho,\mathtt{q})$ related to the density profile $\uprho$ vanishes, because $|\textup{supp}(K\mathtt{F}^{\eps})|\rightarrow 0$ under the Lebesgue measure as $\eps\rightarrow 0$.  The dominated convergence theorem implies
\begin{equation*}
\lim_{\eps\rightarrow 0}\mathtt{R}_{T,K\mathtt{F}^{\eps}}(\uprho,\mathtt{q})=2\pi TK\lim_{\eps\rightarrow 0}\int \mathtt{q}_{\theta}(s)\mathtt{F}_{\theta}^{\eps}(s)ds=2\pi TK\int \mathtt{q}_{\theta}(s)\mathds{1}_{C^+}(s)ds>2\pi TKc.
\end{equation*}
Taking $K\rightarrow +\infty$, we have $\sup_{\mathtt{F}\in \mathtt{C}[\infty]}\mathtt{R}_{T,\mathtt{F}}(\uprho,\mathtt{q})=+\infty$ and
\begin{equation}\label{eq:eliminate_singular_current}
\inf_{(\uprho,\mathtt{q})\in \mathtt{S}}\sup_{\mathtt{F}\in \mathtt{C}[\infty]}\mathtt{R}_{T,\mathtt{F}}(\uprho,\mathtt{q})= \inf_{(\uprho,\mathtt{q})\in \mathtt{S},\textup{ $\mathtt{q}$ is a function}}\sup_{\mathtt{F}\in \mathtt{C}[\infty]}\mathtt{R}_{T,\mathtt{F}}(\uprho,\mathtt{q}).
\end{equation}

\textit{Step 4: Proof of \eqref{eq:exchange_order_b_parameter}.}

Now to prove \eqref{eq:exchange_order_b_parameter}, it is sufficient to prove
\begin{equation}\label{eq:exchange_order_b_parameter_not_singular}
\limsup_{\mathfrak{b}\rightarrow +\infty}\inf_{(\uprho,\mathtt{q})\in \mathtt{S}}\sup_{\mathtt{F}\in \mathtt{C}[\mathfrak{b}]}\mathtt{R}_{T,\mathtt{F}}(\uprho,\mathtt{q})=\inf_{(\uprho,\mathtt{q})\in \mathtt{S},\textup{ $\mathtt{q}$ is a function}}\sup_{\mathtt{F}\in \mathtt{C}[\infty]}\mathtt{R}_{T,\mathtt{F}}(\uprho,\mathtt{q}),
\end{equation}
Because the functional $\mathtt{R}_{T,\mathtt{F}}(\uprho,\mathtt{q})$ is concave in $\mathtt{F}$ and convex in $(\uprho,\mathtt{q})$, given $(\uprho,\mathtt{q})\in \mathtt{S}$ with $\mathtt{q}$ as a function, the optimal forcing $\mathtt{F}^{\op}$ is characterized by $\frac{\p \mathtt{R}_{T,\mathtt{F}}(\uprho,\mathtt{q})}{\p \mathtt{F}}\equiv 0$
\begin{equation}\label{eq:optimal_forcing}
\mathtt{F}_{\theta}^{\op}(s)=\frac{\mathtt{q}_{\theta}(s)}{\sigma(\uprho(s))},\quad \mathtt{F}_{r}^{\op}(s)=-\frac{1}{2}\frac{\frac{d}{ds}\uprho(s)}{\sigma(\uprho(s))}.
\end{equation}
The forcing $\mathtt{F}^{\op}=(\mathtt{F}_r^{\op},\mathtt{F}_{\theta}^{\op})$ does not necessarily belong to $\mathtt{C}[+\infty]$ due to the fact that $(\uprho,\mathtt{q})$ can have very low regularities. But it can be approximated by forcings in $\mathtt{C}[+\infty]$ with arbitrary accuracy. Inserting the expression \eqref{eq:optimal_forcing} of the optimal forcing into the definition \eqref{eq:to_self_similar_R_functional} of $\mathtt{R}_{T,\mathtt{F}}(\uprho,\mathtt{q})$, it shows the optimal value $\sup_{\mathtt{F}\in \mathtt{C}[\infty]}\mathtt{R}_{T,\mathtt{F}}(\uprho,\mathtt{q})$ equals to the functional $T\mathtt{C}(\uprho,\mathtt{q})$ in Lemma \ref{lem:C_self_similarity}
\begin{equation}\label{eq:exchange_b_1}
\inf_{(\uprho,\mathtt{q})\in \mathtt{S},\textup{ $\mathtt{q}$ is a function}}\sup_{\mathtt{F}\in \mathtt{C}[\infty]}\mathtt{R}_{T,\mathtt{F}}(\uprho,\mathtt{q})=\inf_{(\uprho,\mathtt{q})\in \mathtt{S},\textup{ $\mathtt{q}$ is a function}}\mathtt{R}_{T,\mathtt{F}^{\op}}(\uprho,\mathtt{q})=T\inf_{(\uprho,\mathtt{q})\in \mathtt{S},\textup{ $\mathtt{q}$ is a function}}\mathtt{C}(\uprho,\mathtt{q}).
\end{equation}
Proposition \ref{prop:optimizer} shows that there exists a unique optimizer $(\uprho_*,\mathtt{q}_*)$ of the variational formula in \eqref{eq:exchange_b_1}. We denote $\mathtt{F}_{*}^{\op}$ as the optimal forcing determined by $(\uprho_*,\mathtt{q}_*)$. Since $(\uprho_*,\mathtt{q}_*)$ is the optimizer, we can choose any other profile $(\uprho,\mathtt{q})\in \mathtt{S}$  and define $\eta_1:=\uprho-\uprho_*,\eta_2:=\mathtt{q}-\mathtt{q}_*$
\begin{equation}\label{eq:exchange_b_2}
T\lim_{\eps\rightarrow 0}\eps^{-1}\Big(\mathtt{C}(\uprho_*+\eps\eta_1,\mathtt{q}_*+\eps\eta_2)-\mathtt{C}(\uprho_*,\mathtt{q}_*)\Big)=\lim_{\eps\rightarrow 0}\eps^{-1}\Bigg(\mathtt{R}_{T,\mathtt{F}_*^{\op,\eps}}(\uprho_*+\eps\eta_1,\mathtt{q}_*+\eps\eta_2)-\mathtt{R}_{T,\mathtt{F}_*^{\op}}(\uprho_*,\mathtt{q}_*)\Bigg)\geq 0,
\end{equation}
where $\mathtt{F}_*^{\op,\eps}$ is the optimal forcing of $(\uprho_*+\eps\eta_1,\mathtt{q}_*+\eps\eta_2)$ given by \eqref{eq:optimal_forcing}. Equation \eqref{eq:exchange_b_2} implies
\begin{equation}\label{eq:exchange_b_3}
\lim_{\eps\rightarrow 0}\eps^{-1}\Bigg(\mathtt{R}_{T,\mathtt{F}_*^{\op,\eps}}(\uprho_*,\mathtt{q}_*)-\mathtt{R}_{T,\mathtt{F}_*^{\op}}(\uprho_*,\mathtt{q}_*)\Bigg)+\lim_{\eps\rightarrow 0}\eps^{-1}\Bigg(\mathtt{R}_{T,\mathtt{F}_*^{\op}}(\uprho_*+\eps\eta_1,\mathtt{q}_*+\eps\eta_2)-\mathtt{R}_{T,\mathtt{F}_*^{\op}}(\uprho_*,\mathtt{q}_*)\Bigg)\geq 0
\end{equation}
The fact that $\mathtt{F}_*^{\op}$ is the optimal forcing of $(\uprho_*,\mathtt{q}_*)$ implies
\begin{equation}\label{eq:exchange_b_4}
\lim_{\eps\rightarrow 0}\eps^{-1}\Bigg(\mathtt{R}_{T,\mathtt{F}_*^{\op,\eps}}(\uprho_*,\mathtt{q}_*)-\mathtt{R}_{T,\mathtt{F}_*^{\op}}(\uprho_*,\mathtt{q}_*)\Bigg)=0.
\end{equation}
Insert \eqref{eq:exchange_b_4} into \eqref{eq:exchange_b_3} and we have
\begin{equation*}
\lim_{\eps\rightarrow 0}\eps^{-1}\Bigg(\mathtt{R}_{T,\mathtt{F}_*^{\op}}(\uprho_*+\eps\eta_1,\mathtt{q}_*+\eps\eta_2)-\mathtt{R}_{T,\mathtt{F}_*^{\op}}(\uprho_*,\mathtt{q}_*)\Bigg)\geq 0.
\end{equation*}
This shows $(\uprho_*,\mathtt{q}_*)$ is the critical point of the convex functional $\mathtt{R}_{T,\mathtt{F}_*^{\op}}(\cdot,\cdot)$ with the fixed forcing $\mathtt{F}_*^{\op}$, and thus the minimizer of $\mathtt{R}_{T,\mathtt{F}_*^{\op}}(\cdot,\cdot)$
\begin{equation*}
\limsup_{\mathfrak{b}\rightarrow +\infty}\inf_{(\uprho,\mathtt{q})\in \mathtt{S}}\sup_{\mathtt{F}\in \mathtt{C}[\mathfrak{b}]}\mathtt{R}_{T,\mathtt{F}}(\uprho,\mathtt{q})\geq \inf_{(\uprho,\mathtt{q})\in \mathtt{S}}\mathtt{R}_{T,\mathtt{F}_*^{\op}}(\uprho,\mathtt{q})=\mathtt{R}_{T,\mathtt{F}_*^{\op}}(\uprho_*,\mathtt{q}_*)=\inf_{(\uprho,\mathtt{q})\in\mathtt{S},\textup{$\mathtt{q}$ is a function}}\sup_{\mathtt{F}\in \mathtt{C}[\infty]}\mathtt{R}_{T,\mathtt{F}}(\uprho,\mathtt{q})
\end{equation*}
This proves \eqref{eq:exchange_order_b_parameter} and
\begin{equation}\label{eq:upper_bound_proof_final}
\limsup_{N\rightarrow +\infty}\frac{\log N}{N^2}\log\mathds{P}_{\gamma,\SS}^N\Big[\langle\mathcal{Q}_T^N,\mathcal{V}\rangle\geq JT\Big]\leq -\inf_{(\uprho,\mathtt{q})\in \mathtt{S},\textup{$\mathtt{q}$ is a function}}\sup_{\mathtt{F}\in \mathtt{C}[\infty]}\mathtt{R}_{T,\mathtt{F}}(\uprho,\mathtt{q})=-T\inf_{(\uprho,\mathtt{q})\in \mathtt{S},\textup{$\mathtt{q}$ is a function}}\mathtt{C}(\uprho,\mathtt{q}).
\end{equation}
Apply Lemma \ref{lem:C_self_similarity} (with $N$ taken as $e$) to \eqref{eq:upper_bound_proof_final}, we can write the cost functional and profiles in self-similar variable back to the physical variable. Therefore we conclude the proof of the large deviation upper bound in Theorem \ref{th:theorem_upper_bound}.

\end{proof}

\section{Singular Large Deviation Lower Bound for Winding Number}\label{sec:lower_bound}

This section is devoted to the proof of the large deviation lower bounds (Theorem \ref{th:one-scale-deviation} and Corollary \ref{cor:two_scale}). The proof is based on the relative entropy estimate established in Sections \ref{sec:relative_entropy}, and the concentration estimate established in Section \ref{sec:concentration_RN}.

Throughout Sections \ref{sec:lower_bound}--\ref{sec:concentration_RN}, the symbol $C_*$ denotes a finite positive constant, possibly changing from line to line, which is independent of $N$ and of the lattice sites. It may depend on the parameters fixed such as $\gamma,\zeta,\alpha,J,T,\delta$. When a deterministic function $H^N$ is present, the $C_*$ may also depend on the constant $\mathfrak{b}$, which measures the upper bound of $|H_x^N|$
\begin{equation*}
\big|H_x^N\big|\leq \mathfrak{b}|x|^{-2}\frac{1}{\log N},\quad \textup{ for all $x\in \mathds{B}^N$}.
\end{equation*}
We may also use $o_N^*(1)$ to denote a function converging to $0$ as $N\rightarrow +\infty$, whose convergence rate possibly depends on the parameters $(\gamma,\zeta,\alpha,J,T,\delta,\mathfrak{b})$.

\begin{proof}[Proof of Lower Bound in Theorem \ref{th:one-scale-deviation}]

Fix arbitrary parameters $(\gamma,\zeta,\alpha,J,N)$ with 
\begin{equation*}
\gamma\in (0,1),\quad \zeta\in (0,1),\quad \alpha\in(0,\frac{1}{10}),\quad J>0,\quad N>1.
\end{equation*}
\textit{Step 1: Tilt the dynamics.}

In the proof, we use $(\mathfrak{u},\mathfrak{q})$ to represent the following density and current profiles with explicit expressions in \eqref{eq:u_profile}
\begin{equation*}
\mathfrak{u}(x):=\mathfrak{u}^N[\gamma,\zeta;\alpha](x),\quad \mathfrak{q}(x):=\frac{J}{2\pi\log N}\nabla^{\perp}\log|x|\mathds{1}_{N^{-1}\leq |x|\leq 1},\quad \textup{for $x\in B_1$.}
\end{equation*}
To keep the density profile $\mathfrak{u}$ stationary, and to produce the current profile $\mathfrak{q}$, the forcing should be (up to small errors)
\begin{equation}\label{eq:continuum_forcing_parametrized}
\mathcal{F}^N[\gamma,\zeta,\alpha,J](x)=J\frac{\nabla^{\perp}\log|x|}{2\pi\sigma(\mathfrak{u}_x)}+\frac{1}{2}\log N\frac{\nabla\mathfrak{u}_x}{\sigma(\mathfrak{u}_x)},\quad \textup{for $x\in B_1$.}
\end{equation}
However due to the singularity of $\mathcal{F}^N$, we cannot derive directly a good relative entropy estimate of the measure $d\mathds{P}_{\mathfrak{u},\mathcal{F}^N,t}^N$ of the modified dynamics \eqref{eq:modified_generator}. To this end, in Proposition \ref{prop:general_choice} we will construct a modified forcing $\mathds{F}^N[\gamma,\zeta,\alpha,J]$ such that 
\begin{equation}\label{eq:relative_entropy_inproof}
\mathds{H}^N\Big(d\mathds{P}_{\mathfrak{u},\mathds{F},t}^N\Big|\ d\local{\mathfrak{u}}\Big)\leq C_*N^{4\alpha},\quad \mathds{F}:=\mathds{F}^N[\gamma,\zeta,\alpha,J]=\mathcal{F}^N[\gamma,\zeta,\alpha,J]+C_*O(N^{-1}|x|^{-2})+C_*O\Big(\mathds{1}_{|x|\geq 1-\frac{20}{N}}\Big).
\end{equation}
From now on, we use $\mathcal{F}$ to abbreviate $\mathcal{F}^N[\gamma,\zeta,\alpha,J]$, and use $\mathds{F}$ to abbreviate $\mathds{F}^N[\gamma,\zeta,\alpha,J]$. Notice that $(\mathfrak{u},\mathfrak{q},\mathds{F})$ all depend on $N$. Now we prove the lower bound in Theorem \ref{th:one-scale-deviation}. Based on the relative entropy estimate and concentration estimates derived in Sections \ref{sec:relative_entropy}-\ref{sec:concentration_RN}, the proof is classical.

Given a forcing $\mathcal{F}$ as in \eqref{eq:continuum_forcing_parametrized} and a forcing $\mathds{F}$ as constructed in Proposition \ref{prop:general_choice} satisfying \eqref{eq:relative_entropy_inproof}, for any $\delta>0$ define the set $\mathcal{A}_{\delta}^N$ as
\begin{equation}\label{eq:approximation_A}
\mathcal{A}_{\delta}^N:=\Big\{\{\eta_{s}\}_{0\leq s\leq T}\Big|\  |\mathbf{R}_{T,\mathds{F}}^N(\rho^N,\mathcal{Q}^N)-\mathcal{R}_{T,\mathcal{F}}^N(\mathfrak{u},\mathfrak{q})|<\delta\Big\}.
\end{equation}
Here the functionals $\mathbf{R}_{T,\mathds{F}}^N$ and $\mathcal{R}_{T,\mathcal{F}}^N$ are defined in Definition \ref{def:R_functionals}. Tilt the measure $d\mathds{P}_{\gamma,\SS}^N$ of SSEP into the measure $d\mathds{P}_{\mathfrak{u},\mathds{F}}^N$ of the modified dynamics \eqref{eq:modified_generator} for any $N$
\begin{equation}\label{eq:sec5_tilt}
\begin{split}
&\log \mathds{P}_{\gamma,\SS}^N\Big[\langle \mathcal{Q}_T^N,\mathcal{V}\rangle\geq JT-w\Big]=\log \mathds{E}_{\mathfrak{u},\mathds{F}}^N\Bigg[\frac{d\mathds{P}_{\gamma,\SS}^N}{d\mathds{P}_{\mathfrak{u},\mathds{F}}^N}(\rho^N,\mathcal{Q}^N)\mathds{1}_{\langle \mathcal{Q}_T^N,\mathcal{V}\rangle\geq JT-w}\Bigg].
\end{split}
\end{equation}

\textit{Step 2: Smallness of the initial entropy.}

The initial entropy $\textup{Ent}_{\gamma}(\mathfrak{u})$ corresponds to the cost of the initial deviation from $\local{\gamma}$ to $\local{\mathfrak{u}}$
\begin{equation}\label{eq:initial_entropy_sec5}
\textup{Ent}_{\gamma}(\mathfrak{u})=\int \Big(\mathfrak{u}_x\log\frac{\mathfrak{u}_x}{\gamma_x}+(1-\mathfrak{u}_x)\log\frac{1-\mathfrak{u}_x}{1-\gamma_x}\Big)dx.
\end{equation}
As the density profile $\mathfrak{u}$ in \eqref{eq:u_profile} only differs from $\gamma$ close to origin, we are going to check that the initial entropy in \eqref{eq:initial_entropy_sec5} is small. The fact that $\log y\leq y-1$ for any $y>0$ implies
\begin{equation}\label{eq:inequality_ent}
\mathfrak{u}_x\log\frac{\mathfrak{u}_x}{\gamma_x}+(1-\mathfrak{u}_x)\log\frac{1-\mathfrak{u}_x}{1-\gamma_x}\leq \frac{\mathfrak{u}_x(\mathfrak{u}_x-\gamma_x)}{\gamma_x}+\frac{(1-\mathfrak{u}_x)(\gamma_x-\mathfrak{u}_x)}{1-\gamma_x}\leq (\mathfrak{u}_x-\gamma_x)\Big(\frac{\mathfrak{u}_x}{\gamma_x}-\frac{1-\mathfrak{u}_x}{1-\gamma_x}\Big)=\frac{(\mathfrak{u}_x-\gamma_x)^2}{\gamma_x(1-\gamma_x)}.
\end{equation}
Insert \eqref{eq:inequality_ent} to \eqref{eq:initial_entropy_sec5}, and decompose the integral domain
\begin{equation*}
\textup{Ent}_{\gamma}(\mathfrak{u})\leq \int_{|x|\leq N^{-\alpha}}\frac{(\mathfrak{u}-\gamma_x)^2}{\gamma_x(1-\gamma_x)}dx+\int_{N^{-\alpha}\leq |x|\leq 1}\frac{(\mathfrak{u}-\gamma_x)^2}{\gamma_x(1-\gamma_x)}dx=O(N^{-2\alpha})+\int_{N^{-\alpha}\leq |x|\leq 1}\frac{(\mathfrak{u}-\gamma_x)^2}{\gamma_x(1-\gamma_x)}dx.
\end{equation*}
The expression of $\mathfrak{u}$ in \eqref{eq:u_profile} shows that the difference $|\mathfrak{u}-\gamma_x|$ is bounded above by $|\mathfrak{u}_x-\gamma_x|\leq C_*\frac{-\log |x|}{\log N}$. Consequently $\int_{N^{-\alpha}\leq |x|\leq 1}\frac{(\mathfrak{u}-\gamma_x)^2}{\gamma_x(1-\gamma_x)}dx\leq C_*\int_{N^{-\alpha}\leq |x|\leq 1}\Big(\frac{\log |x|}{\log N}\Big)^2dx\leq C_*\int_{N^{-\alpha}}^1r(\log r)^2(\log N)^{-2}dr=O\Big(\frac{1}{\log^2N}\Big)$. This proves that $\textup{Ent}_{\gamma}(\mathfrak{u})$ is of order $O((\log N)^{-2})$
.

\textit{Step 3: Dynamical cost and concentration estimates.}

According to Proposition \ref{prop:radon_nikodym} and the smallness of the initial entropy $\textup{Ent}_{\gamma}(\mathfrak{u})$, by neglecting a lower order term, we have
\begin{equation*}
\begin{split}
&\liminf_{N\rightarrow +\infty}\frac{\log N}{N^2}\log \mathds{P}_{\gamma,\SS}^N\Big[\langle \mathcal{Q}_T^N,\mathcal{V}\rangle\geq JT-w\Big]\\
&\qquad\qquad\geq \liminf_{N\rightarrow +\infty}\frac{\log N}{N^2}\log \mathds{E}_{\mathfrak{u},\mathds{F}}^N\Bigg[\exp\Bigg(-\frac{N^2}{\log N}\Big(\mathbf{R}_{T,\mathds{F}}^N(\rho^N,\mathcal{Q}^N)+o_N^*(1)\Big)\Bigg)\mathds{1}_{\langle \mathcal{Q}_T^N,\mathcal{V}\rangle\geq JT-w}\Bigg],
\end{split}
\end{equation*}
where $\mathbf{R}_{T,\mathds{F}}^N$ is given as in \eqref{eq:discrete_R_functional}. Adding a restriction by $\mathcal{A}_{\delta
}^N$ as \eqref{eq:approximation_A} and replacing $\mathbf{R}_{T,\mathds{F}}^N(\rho^N,\mathcal{Q}^N)$ with $\mathcal{R}_{T,\mathcal{F}}^N(\mathfrak{u},\mathfrak{q})$, the term above is further greater than
\begin{equation*}
\begin{split}
&\geq \liminf_{N\rightarrow +\infty}\frac{\log N}{N^2}\log \mathds{E}_{\mathfrak{u},\mathds{F}}^N\Bigg[\exp\Bigg(-\frac{N^2}{\log N}\Big(\mathbf{R}_{T,\mathds{F}}^N(\rho^N,\mathcal{Q}^N)+o_N^*(1)\Big)\Bigg)\mathds{1}_{\langle \mathcal{Q}_T^N,\mathcal{V}\rangle\geq JT-w}\mathds{1}_{\mathcal{A}_{\delta}^N}\Bigg]\\
&\geq \liminf_{N\rightarrow \infty}\frac{\log N}{N^2}\log\mathds{E}_{\mathfrak{u},\mathds{F}}^N\Bigg[\exp\Bigg(-\frac{N^2}{\log N}\Big(\mathcal{R}_{T,\mathcal{F}}^N(\mathfrak{u},\mathfrak{q})+o_N^*(1)+O(\delta)\Big)\Bigg)\mathds{1}_{\langle\mathcal{Q}_T^N,\mathcal{V} \rangle\geq JT-w}\mathds{1}_{\mathcal{A}_{\delta}^N}\Bigg],
\end{split}
\end{equation*}
where in the last step we have replaced $\mathbf{R}_{T,\mathds{F}}^N(\rho^N,\mathcal{Q}^N)$ with $\mathcal{R}_{T,\mathcal{F}}^N(\mathfrak{u},\mathfrak{q})$ according to the definition of $\mathcal{A}_{\delta}^N$, contributing an error $O(\delta)$. Consequently
\begin{equation*}
\begin{split}
&\liminf_{N\rightarrow +\infty}\frac{\log N}{N^2}\log \mathds{P}_{\gamma,\SS}^N\Big[\langle \mathcal{Q}_T^N,\mathcal{V}\rangle\geq JT-w\Big]\\
&\qquad\qquad\geq -\limsup_{N\rightarrow \infty}\Big(\mathcal{R}_{T,\mathcal{F}}^N(\mathfrak{u},\mathfrak{q})+O(\delta)+o_N^*(1)\Big)+\liminf_{N\rightarrow \infty}\frac{\log N}{N^2}\log \mathds{E}_{\mathfrak{u},\mathds{F}}^N\Big[\mathds{1}_{\langle \mathcal{Q}_T^N,\mathcal{V} \rangle \geq JT-w}\mathds{1}_{\mathcal{A}_{\delta}^N}\Big].
\end{split}
\end{equation*}
Since the Proposition \ref{prop:relative_entropy} gives a good control of the relative entropy as stated in \eqref{eq:relative_entropy_inproof}, we can prove concentration estimates under the measure $d\mathds{P}_{\mathfrak{u},\mathds{F}}^N$ for arbitrary $\delta>0$
\begin{equation}\label{eq:concentration_statement_sec5}
\lim_{N\rightarrow +\infty}\mathds{P}_{\mathfrak{u},\mathds{F}}^N\Big[\mathcal{A}_{\delta}^N\Big]=1,\quad \lim_{N\rightarrow +\infty}\mathds{P}_{\mathfrak{u},\mathds{F}}^N\Big[\langle \mathcal{Q}_T^N,\mathcal{V} \rangle\geq JT-w\Big]=1,
\end{equation}
with the details of the proof given in Section \ref{sec:concentration_RN}. The first equality in \eqref{eq:concentration_statement_sec5} is proved by Proposition \ref{prop:concentration}, which shows $\mathbf{R}_{T,\mathds{F}}^N(\rho^N,\mathcal{Q}^N)$ converges to $\mathcal{R}_{T,\mathcal{F}}^N(\mathfrak{u},\mathfrak{q})$ as $N\rightarrow +\infty$. The second inequality is proved by Lemma \ref{lem:winding_concentration}.

The derivation of \eqref{eq:concentration_statement_sec5} is a core difficulty of the proof of lower bound, as the forcing $\mathds{F}$ used in the tilting \eqref{eq:sec5_tilt} is very singular. Taking $\delta\rightarrow 0$, there is
\begin{equation*}
\begin{split}
&\liminf_{N\rightarrow +\infty}\frac{\log N}{N^2}\log \mathds{P}_{\gamma,\SS}^N\Big[\langle \mathcal{Q}_T^N,\mathcal{V}\rangle\geq JT-w\Big]\geq -\liminf_{N\rightarrow +\infty}\mathcal{R}_{T,\mathcal{F}}^{N}(\mathfrak{u},\mathfrak{q})=-T\liminf_{N\rightarrow +\infty}\mathcal{C}^N\big(\mathfrak{u},\mathfrak{q}\big).
\end{split}
\end{equation*}
The definition of $\mathcal{C}(\gamma,\zeta,\alpha,J)$ in Definition \ref{def:cost_functional_lower} equals to the cost $T^{-1}\mathcal{R}_{T,\mathcal{F}}^N(\mathfrak{u},\mathfrak{q})$ as $N\rightarrow +\infty$
\begin{equation}\label{eq:lower_bound_near_final}
\begin{split}
&\liminf_{N\rightarrow +\infty}\frac{\log N}{N^2}\log \mathds{P}_{\gamma,\SS}^N\Big[\langle \mathcal{Q}_T^N,\mathcal{V}\rangle\geq JT-w\Big]\geq -T\mathcal{C}(\gamma,\zeta,\alpha,J),\quad\textup{with } \zeta=\frac{1}{2},\quad \alpha=\frac{1}{10}\frac{1}{1+|J|}.
\end{split}
\end{equation}
Equation \eqref{eq:lower_bound_near_final} is true for any $J>0$ and $w>0$. Thus for any $J>0$ and any $w>0$, we have
\begin{equation}
\begin{split}
&\liminf_{N\rightarrow +\infty}\frac{\log N}{N^2}\log \mathds{P}_{\gamma,\SS}^N\Big[\langle \mathcal{Q}_T^N,\mathcal{V}\rangle\geq JT\Big]\geq -T\mathcal{C}(\gamma,\zeta,\alpha,J+w),\quad\textup{with } \zeta=\frac{1}{2},\quad \alpha=\frac{1}{10}\frac{1}{1+|J+w|}.
\end{split}
\end{equation}
Taking $w\rightarrow 0$, the continuity of $\mathcal{C}$-functional in Proposition \ref{prop:C_cost_cubic} implies
\begin{equation}\label{eq:lower_bound_final}
\begin{split}
&\liminf_{N\rightarrow +\infty}\frac{\log N}{N^2}\log \mathds{P}_{\gamma,\SS}^N\Big[\langle \mathcal{Q}_T^N,\mathcal{V}\rangle\geq JT\Big]\geq -T\mathcal{C}(\gamma,\zeta,\alpha,J),\quad\textup{with } \zeta=\frac{1}{2},\quad \alpha=\frac{1}{10}\frac{1}{1+|J|}.
\end{split}
\end{equation}
The proof of the large deviation lower bound in Theorem \ref{th:one-scale-deviation} is concluded.
\end{proof}

In Proposition \ref{prop:C_cost_cubic} we derive an explicit expression of the cost function $\mathcal{C}(\gamma,\zeta,\alpha,J)$ with a negligible error. This enables us to prove Corollary \ref{cor:two_scale} using the asymptotic behavior of $\mathcal{C}(\gamma,\zeta,\alpha,J)$ as $J\rightarrow +\infty$.

\section{Characterization of the Relative Entropy}\label{sec:relative_entropy}

\begin{proposition}\label{prop:relative_entropy}
Take arbitrary constants $(\gamma,\zeta,\alpha,J)$ with
\begin{equation*}
\gamma\in (0,1),\quad \zeta\in (0,1),\quad \alpha\in (0,\frac{1}{3}),\quad J>0.
\end{equation*}
Also take $\mathfrak{u}$ to represent $\mathfrak{u}_x:=\mathfrak{u}^N[\gamma,\zeta;\alpha](x)$ in \eqref{eq:u_profile}, and $\mathds{F}$ as the forcing field $\mathds{F}^N[\gamma,\zeta,\alpha,J]$ defined in  Proposition \ref{prop:general_choice}. For the law $d\mathds{P}_{\mathfrak{u},\mathds{F}}^N$ of the modified dynamics with generator $L_{\mathds{F}}^{N}$, its relative entropy satisfies, for every fixed $T<\infty$ and every $t\in[0,T]$
\begin{equation*}
\mathds{H}^N\Big(d\mathds{P}_{\mathfrak{u},\mathds{F},t}^N\Big|d\local{\mathfrak{u}}\Big)\leq C_*N^{4\alpha}.
\end{equation*}
\end{proposition}
Yau's inequality \cite{yau1991relative} is a classical way to estimate the relative entropy
\begin{equation}\label{eq:Yau}
\frac{d}{dt}\mathds{H}^N\Big(d\mathds{P}_{\mathfrak{u},\mathds{F},t}^N\Big|d\local{\mathfrak{u}}\Big)\leq - \local{\mathfrak{u}}\Big[\Gamma_{\mathds{F}}^N\sqrt{f_t}\Big]+ \local{\mathfrak{u}}\Big[f_tL_{\mathds{F}}^{N,*}\mathds{1}\Big] ,\quad f_t:=\frac{d\mathds{P}_{\mathfrak{u},\mathds{F},t}^N}{d\local{\mathfrak{u}}}.
\end{equation}
Here the dissipation term $\Gamma_{\mathds{F}}^N$ is
\begin{equation}\label{eq:modified_generator_Dirichlet}
\begin{split}
\Gamma_{\mathds{F}}^Nf(\eta)=\frac{N^2}{2}\sum_{x\in\mathds{B}^N}\sum_{e\in\{e_1,e_2\}}\Bigg[&\eta_x(1-\eta_{x+e})\big(f(\eta^{x,x+e})-f(\eta)\big)^2\exp\Big(\frac{\mathds{F}_{x,e}}{N\log N}\Big)\\
&\qquad+\eta_x
(1-\eta_{x-e})\big(f(\eta^{x,x-e})-f(\eta)\big)^2\exp\Big(-\frac{\mathds{F}_{x-e,e}}{N\log N}\Big)\Bigg],
\end{split}
\end{equation}
and $L_{\mathds{F}}^{N,*}$ is the adjoint operator of $L_{\mathds{F}}^N$ with respect to $d\local{\mathfrak{u}}$. The work \cite{Jara_Menezes_2018} establishes a powerful framework of controlling relative entropy, by providing tools to estimate quantitatively the contribution $\local{\mathfrak{u}}\Big[f_tL_{\mathds{F}}^{N,*}\mathds{1}\Big]$ in \eqref{eq:Yau}. Specifically, define $\omega_x$ as
\begin{equation}\label{eq:centered_configuration}
\omega_x:=\frac{\eta_x-\mathfrak{u}_x}{\mathfrak{u}_x(1-\mathfrak{u}_x)}.
\end{equation}
In the current context, \cite{Jara_Menezes_2018} allows us to express the time derivative of relative entropy as
\begin{equation*}
\p_t\mathds{H}^N\Big(d\mathds{P}_{\mathfrak{u},\mathds{F},t}^N\Big|\ d\local{\mathfrak{u}}\Big)\leq -\local{\mathfrak{u}}\Big[\Gamma_{\mathds{F}}^N\sqrt{f_t}\Big]+\textup{term linear in $\omega$}+\textup{term quadratic in $\omega$},
\end{equation*}
whose precise statement will be given in Lemma \ref{lem:linear_quadratic}. The differences between the current paper and \cite{Jara_Menezes_2018} are the following: In \cite{Jara_Menezes_2018}, the coefficients in the linear term are $o_N(1)$, and those in the quadratic term are bounded. Here in this paper, the coefficients in both linear and quadratic terms are singular close to the origin.

To overcome these new difficulties, we build below new methods together with the method of \cite{Jara_Menezes_2018}, to estalish the control of the relative entropy. The estimate of the relative entropy in Proposition \ref{prop:relative_entropy} will be proved in the following sections. In Section \ref{sec:adjoint_operator} we give the explicit expression of the entropy growth by Lemma \ref{lem:linear_quadratic}. In Section \ref{sec:proper_discretization} we construct the forcing $\mathds{F}^N[\gamma,\zeta,\alpha,J]$ (abbreviated $\mathds{F}$) to control the linear terms. In Section \ref{sec:cancellation} we use a cancellation mechanism, combined with the method of Jara-Menezes \cite{Jara_Menezes_2018}, to control the quadratic terms. Based on these ingredients, we prove Proposition \ref{prop:relative_entropy} in Section \ref{sec:proof_entropy}.

\subsection{Exact Formula of the Adjoint Operator}\label{sec:adjoint_operator}

This section contains some definitions and lemmas that apply to general forcing $\mathds{F}$. We denote $\mathbf{E}^N$ as the set of admissible edges
\begin{equation}\label{eq:admissible_edges}
\mathbf{E}^N:=\Big\{(x,e)\Big|\ x\in\mathds{B}^N,\ e\in\{e_1,e_2\},\ x+e\in \mathds{B}^N\Big\},
\end{equation}
and $\p\mathds{B}^N$ as the boundary of $\mathds{B}^N$
\begin{equation}\label{eq:boundary_BN}
\p\mathds{B}^N:=\Big\{x\in\mathds{B}^N\Big|\ \textup{there exists $e\in\{e_1,e_2\}$ such that }x+e\not\in \mathds{B}^N \textup{ or $x-e\not\in \mathds{B}^N$}\Big\}.
\end{equation}

\begin{definition}
For any $x\in\mathds{B}^N$ and $e\in\{e_1,e_2\}$, define the function $\mathbf{Q}_{x,e}$ as
\begin{equation}\label{eq:auxillary_A_usual}
\mathbf{Q}_{x,e}=N\Bigg[\exp\Big(\frac{\mathds{F}_{x,e}}{N\log N}\Big)\mathfrak{u}_x(1-\mathfrak{u}_{x+e})-\exp\Big(\frac{-\mathds{F}_{x,e}}{N\log N}\Big)\mathfrak{u}_{x+e}(1-\mathfrak{u}_{x})\Bigg].
\end{equation}
\end{definition}

\begin{lemma}\label{lem:linear_quadratic}
\textbf{\textup{[Decomposition of the adjoint operator]}} The adjoint operator $L_{\mathds{F}}^{N,*}\mathds{1}$ of the modified dynamics can be decomposed as
\begin{equation*}
L_{\mathds{F}}^{N,*}\mathds{1}(\eta)=\frac{1}{2}\mathds{A}_{\textup{lin}}(\eta)+\frac{1}{2}\mathds{A}_{\textup{qua}}(\eta),
\end{equation*}
with $\mathds{A}_{\textup{lin}}(\eta)$ being a linear-in-$\omega$ term and $\mathds{A}_{\textup{qua}}(\eta)$ being a quadratic-in-$\omega$ term
\begin{equation}\label{eq:linear_quadratic}
\begin{split}
&\mathds{A}_{\textup{lin}}(\eta)=-\sum_{x\in\mathds{B}^N\setminus \p\mathds{B}^N}\omega_x\nabla_{\textup{disc}}\cdot\mathbf{Q}_{x}+N\sum_{x\in \p\mathds{B}^N}\omega_x\sum_{e\in\{e_1,e_2\}}\Big(\mathbf{Q}_{x-e,e}\mathds{1}_{x-e\in\mathds{B}^N}-\mathbf{Q}_{x,e}\mathds{1}_{x+e\in\mathds{B}^N}\Big),\\
&\mathds{A}_{\textup{qua}}(\eta)=\sum_{(x,e)\in\mathbf{E}^N}\mathbf{Q}_{x,e}N(\mathfrak{u}_{x+e}-\mathfrak{u}_x)\omega_{x+e}\omega_{x}.
\end{split}
\end{equation}
\end{lemma}
\begin{proof}

The first step of the proof is to prove the following equality, with details following \cite[Section A.3]{Jara_Menezes_2018}
\begin{equation}\label{eq:adjoint_JM_reduction}
\begin{aligned}
&L_{\mathds{F}}^{N,*}\mathds{1}(\eta)=
\frac{1}{2}\sum_{(x,e)\in\mathbf{E}^N}
N\mathbf{Q}_{x,e}(\omega_{x+e}-\omega_x)+\frac{1}{2}\sum_{(x,e)\in\mathbf{E}^N}
N\mathbf{Q}_{x,e}(\mathfrak{u}_{x+e}-\mathfrak{u}_x)\omega_{x+e}\omega_x.
\end{aligned}
\end{equation}

Define $c_{x,e}^{\mathds{F}}(\eta)$ as the following sum of jump rates
\begin{equation}\label{eq:unoriented_jump_rate}
c_{x,e}^{\mathds{F}}(\eta):=\frac{N^2}{2}\Bigg(\exp\Big(\frac{\mathds{F}_{x,e}}{N\log N}\Big)\eta_x(1-\eta_{x+e})+\exp\Big(-\frac{\mathds{F}_{x,e}}{N\log N}\Big)\eta_{x+e}(1-\eta_x)\Bigg),
\end{equation}
where the first term is the jump rate from $x$ to $x+e$, and the second term is the jump rate from $x+e$ to $x$. With the notion of $c_{x,e}^{\mathds{F}}$ and the notion of admissible edges $\mathbf{E}^N$ in \eqref{eq:admissible_edges}, the generator $L_{\mathds{F}}^N$ defined in \eqref{eq:modified_generator} and its adjoint $L_{\mathds{F}}^{N,*}$ can be written as
\begin{equation}\label{eq:adjoint_original}
    L_{\mathds{F}}^N f(\eta)
    =
    \sum_{(x,e)\in\mathbf{E}^N}
    c_{x,e}^{\mathds{F}}(\eta)
    \Big(f(\eta^{x,x+e})-f(\eta)\Big),\quad L_{\mathds{F}}^{N,*} g(\eta)=\sum_{(x,e)\in\mathbf{E}^N}\left[c_{x,e}^{\mathds{F}}(\eta^{x,x+e})g(\eta^{x,x+e})\frac{\local{\mathfrak{u}}(\eta^{x,x+e})}{\local{\mathfrak{u}}(\eta)}-c_{x,e}^{\mathds{F}}(\eta)g(\eta)\right].
\end{equation}
Replacing the generic function $g(\eta)$ in \eqref{eq:adjoint_original} by $\mathds{1}$, we derive
\begin{equation}\label{eq:adjoint_intermediate_1}
L_{\mathds{F}}^{N,*}\mathds{1}(\eta)=\sum_{(x,e)\in\mathbf{E}^N}\left[c_{x,e}^{\mathds{F}}(\eta^{x,x+e})\frac{\local{\mathfrak{u}}(\eta^{x,x+e})}{\local{\mathfrak{u}}(\eta)}-c_{x,e}^{\mathds{F}}(\eta)\right].
\end{equation}
Evaluate the unoriented jump rate $c_{x,e}^{\mathds{F}}$ at $\eta^{x,x+e}$
\begin{equation}\label{eq:swapped_unoriented_jumprate}
c_{x,e}^{\mathds{F}}(\eta^{x,x+e})=\frac{N^2}{2}\Bigg(\exp\Big(\frac{\mathds{F}_{x,e}}{N\log N}\Big)\eta_{x+e}(1-\eta_x)+\exp\Big(-\frac{\mathds{F}_{x,e}}{N\log N}\Big)\eta_{x}(1-\eta_{x+e})\Bigg).
\end{equation}
The fraction in \eqref{eq:adjoint_intermediate_1} can be evaluated at configurations $\eta$ with $\eta_{x+e}(1-\eta_x)=1$ or $\eta_x(1-\eta_{x+e})=1$
\begin{equation}\label{eq:derivative_in_adjoint}
\eta_{x+e}(1-\eta_x)\frac{\local{\mathfrak{u}}(\eta^{x,x+e})}{\local{\mathfrak{u}}(\eta)}=\eta_{x+e}(1-\eta_x)\frac{\mathfrak{u}_x(1-\mathfrak{u}_{x+e})}{(1-\mathfrak{u}_x)\mathfrak{u}_{x+e}},\quad \eta_x(1-\eta_{x+e})\frac{\local{\mathfrak{u}}(\eta^{x,x+e})}{\local{\mathfrak{u}}(\eta)}=\eta_x(1-\eta_{x+e})\frac{(1-\mathfrak{u}_{x})\mathfrak{u}_{x+e}}{\mathfrak{u}_x(1-\mathfrak{u}_{x+e})}.
\end{equation}
Applying Equations \eqref{eq:swapped_unoriented_jumprate} and \eqref{eq:derivative_in_adjoint} to the expression \eqref{eq:adjoint_intermediate_1} of $L_{\mathds{F}}^{N,*}$, we derive
\begin{equation}\label{eq:adjoint_intermediate_2}
\begin{split}
&c_{x,e}^{\mathds{F}}(\eta^{x,x+e})\frac{\local{\mathfrak{u}}(\eta^{x,x+e})}{\local{\mathfrak{u}}(\eta)}-c_{x,e}^{\mathds{F}}(\eta)\\
&=\frac{N^2}{2}\exp\Big(\frac{\mathds{F}_{x,e}}{N\log N}\Big)\eta_{x+e}(1-\eta_x)\frac{\mathfrak{u}_x(1-\mathfrak{u}_{x+e})}{(1-\mathfrak{u}_{x})\mathfrak{u}_{x+e}}+\frac{N^2}{2}\exp\Big(-\frac{\mathds{F}_{x,e}}{N\log N}\Big)\eta_{x}(1-\eta_{x+e})\frac{(1-\mathfrak{u}_{x})\mathfrak{u}_{x+e}}{\mathfrak{u}_x(1-\mathfrak{u}_{x+e})}\\
&\quad -\frac{N^2}{2}\exp\Big(\frac{\mathds{F}_{x,e}}{N\log N}\Big)\eta_{x}(1-\eta_{x+e})-\frac{N^2}{2}\exp\Big(-\frac{\mathds{F}_{x,e}}{N\log N}\Big)\eta_{x+e}(1-\eta_x).
\end{split}
\end{equation}
Reorganize the terms in \eqref{eq:adjoint_intermediate_2}
\begin{equation*}
\begin{split}
&\frac{N^2}{2}\exp\Big(\frac{\mathds{F}_{x,e}}{N\log N}\Big)\eta_{x+e}(1-\eta_x)\frac{\mathfrak{u}_x(1-\mathfrak{u}_{x+e})}{(1-\mathfrak{u}_{x})\mathfrak{u}_{x+e}}-\frac{N^2}{2}\exp\Big(\frac{\mathds{F}_{x,e}}{N\log N}\Big)\eta_{x}(1-\eta_{x+e})\\
&\qquad\qquad\qquad\qquad\qquad\qquad=\frac{N^2}{2}\exp\Big(\frac{\mathds{F}_{x,e}}{N\log N}\Big)\mathfrak{u}_{x}(1-\mathfrak{u}_{x+e})\Bigg(\frac{\eta_{x+e}(1-\eta_x)}{\mathfrak{u}_{x+e}(1-\mathfrak{u}_x)}-\frac{\eta_x(1-\eta_{x+e})}{\mathfrak{u}_x(1-\mathfrak{u}_{x+e})}\Bigg),
\end{split}
\end{equation*}
and
\begin{equation*}
\begin{split}
&\frac{N^2}{2}\exp\Big(-\frac{\mathds{F}_{x,e}}{N\log N}\Big)\eta_{x}(1-\eta_{x+e})\frac{(1-\mathfrak{u}_{x})\mathfrak{u}_{x+e}}{\mathfrak{u}_x(1-\mathfrak{u}_{x+e})}-\frac{N^2}{2}\exp\Big(-\frac{\mathds{F}_{x,e}}{N\log N}\Big)\eta_{x+e}(1-\eta_x)\\
&\qquad\qquad\qquad\qquad\qquad\qquad=\frac{N^2}{2}\exp\Big(-\frac{\mathds{F}_{x,e}}{N\log N}\Big)(1-\mathfrak{u}_x)\mathfrak{u}_{x+e}\Bigg(\frac{\eta_x(1-\eta_{x+e})}{\mathfrak{u}_x(1-\mathfrak{u}_{x+e})}-\frac{\eta_{x+e}(1-\eta_x)}{\mathfrak{u}_{x+e}(1-\mathfrak{u}_x)}\Bigg).
\end{split}
\end{equation*}
By reorganizing the terms in \eqref{eq:adjoint_intermediate_2}, by the definition of $\mathbf{Q}_{x,e}$ in \eqref{eq:auxillary_A_usual}, we get
\begin{equation*}
\begin{aligned}
&c_{x,e}^{\mathds{F}}(\eta^{x,x+e})\frac{\local{\mathfrak{u}}(\eta^{x,x+e})}{\local{\mathfrak{u}}(\eta)}-c_{x,e}^{\mathds{F}}(\eta)=
\frac{N}{2}\mathbf{Q}_{x,e}\Bigg(\frac{\eta_{x+e}(1-\eta_x)}{\mathfrak{u}_{x+e}(1-\mathfrak{u}_{x})}-\frac{\eta_x(1-\eta_{x+e})}{\mathfrak{u}_x(1-\mathfrak{u}_{x+e})}\Bigg).
\end{aligned}
\end{equation*}
Summing over all admissible edges yields
\begin{equation}\label{eq:adjoint_some_calculation_1}
\begin{aligned}
&L_{\mathds{F}}^{N,*}\mathds{1}(\eta)=
\sum_{(x,e)\in\mathbf{E}^N}
\frac{N}{2}\mathbf{Q}_{x,e}\Bigg(\frac{\eta_{x+e}(1-\eta_x)}{(1-\mathfrak{u}_x)\mathfrak{u}_{x+e}}-\frac{\eta_x(1-\eta_{x+e})}{\mathfrak{u}_x(1-\mathfrak{u}_{x+e})}\Bigg).
\end{aligned}
\end{equation}
\cite[Equation A.3]{Jara_Menezes_2018} gives the equality
\begin{equation}\label{eq:decomposition_into_omega}
\frac{\eta_{x+e}(1-\eta_x)}{(1-\mathfrak{u}_x)\mathfrak{u}_{x+e}}-\frac{\eta_x(1-\eta_{x+e})}{\mathfrak{u}_x(1-\mathfrak{u}_{x+e})}=(\omega_{x+e}-\omega_x)+(\mathfrak{u}_{x+e}-\mathfrak{u}_x)\omega_{x+e}\omega_x.
\end{equation}
Applying \eqref{eq:decomposition_into_omega} to \eqref{eq:adjoint_some_calculation_1}, we can prove \eqref{eq:adjoint_JM_reduction}. Doing integration by parts for the linear term in \eqref{eq:adjoint_JM_reduction}, there is
\begin{equation*}
N\sum_{(x,e)\in\mathbf{E}^N}\mathbf{Q}_{x,e}(\omega_{x+e}-\omega_x)=-\sum_{x\in\mathds{B}^N\setminus \p\mathds{B}^N}\omega_x\nabla_{\textup{disc}}\cdot\mathbf{Q}_{x}+N\sum_{x\in \p\mathds{B}^N}\omega_x\sum_{e\in\{e_1,e_2\}}\Big(\mathbf{Q}_{x-e,e}\mathds{1}_{x-e\in\mathds{B}^N}-\mathbf{Q}_{x,e}\mathds{1}_{x+e\in\mathds{B}^N}\Big)
\end{equation*} 
This concludes the proof of Lemma \ref{lem:linear_quadratic}.
\end{proof}

\subsection{Discrete Differential Calculus}\label{sec:proper_discretization}

The discrete skew-gradient $\nabla_{\textup{disc}}^\perp$ is defined as the following $N$-dependent difference operator, where $g$ is any function on $N^{-1}\mathds{Z}^2$
\begin{equation}\label{eq:discrete_skew_gradient}
\begin{split}
&\big(\nabla_{\textup{disc}}^{\perp}g\big)_x:=\Big(\big(\nabla_{\textup{disc}}^{\perp}g\big)_{x,e_1},\big(\nabla_{\textup{disc}}^{\perp}g\big)_{x,e_2}\Big),\\
&\big(\nabla_{\textup{disc}}^{\perp}g\big)_{x,e_1}:=-N\big(g_{x+e_1+e_2}-g_{x+e_1}\big),\quad \big(\nabla_{\textup{disc}}^{\perp}g\big)_{x,e_2}:=N\big(g_{x+e_1+e_2}-g_{x+e_2}\big).
\end{split}
\end{equation}
A notable property of $\nabla_{\textup{disc}}^{\perp}$ is that for arbitrary function $g$ on $ N^{-1}\mathds{Z}^2$, we have
\begin{equation}\label{eq:discrete_vanish}
\begin{split}
\nabla_{\textup{disc}}\cdot\nabla_{\textup{disc}}^\perp g&=N\Big(\big(\nabla_{\textup{disc}}^{\perp}g\big)_{x,e_1}+\big(\nabla_{\textup{disc}}^{\perp}g\big)_{x,e_2}-\big(\nabla_{\textup{disc}}^{\perp}g\big)_{x-e_1,e_1}-\big(\nabla_{\textup{disc}}^{\perp}g\big)_{x-e_2,e_2}\Big)\\
&=N^2\Big(\big(g_{x+e_1}-g_{x+e_1+e_2}\big)+\big(g_{x+e_1+e_2}-g_{x+e_2}\big)-\big(g_{x}-g_{x+e_2}\big)-\big(g_{x+e_1}-g_x\big)\Big)=0,
\end{split}
\end{equation}
which is the discrete analogy of $\nabla\cdot\nabla^{\perp}\equiv 0$.

Recall that $\mathfrak{u}$ is taken as the abbreviation of $\mathfrak{u}^N[\gamma,\zeta;\alpha]$ given in \eqref{eq:u_profile}, and that $\mathcal{F}$ is taken as the abbreviation of $\mathcal{F}^N[\gamma,\zeta,\alpha,J]$
\begin{equation*}
\mathcal{F}^N[\gamma,\zeta,\alpha,J](x):=J\frac{\nabla^{\perp}\log|x|}{2\pi\sigma(\mathfrak{u}_x)}+\frac{1}{2}\log N\frac{\nabla \mathfrak{u}_x}{\sigma(\mathfrak{u}_x)},\quad \textup{for all $x\in  B_1$.}
\end{equation*}

\begin{proposition}\label{prop:general_choice}
Take arbitrary parameters $(\gamma,\zeta,\alpha,J,N)$ with
\begin{equation*}
\gamma\in (0,1),\quad \zeta\in (0,1),\quad \alpha\in (0,1),\quad J>0,\quad N>1.
\end{equation*}
We are going to construct the forcing field
$\mathds{F}^N[\gamma,\zeta,\alpha,J]$ as a forcing field on $\mathds{B}^N$ such that for any $x\in\mathds{B}^N$ and $e\in\{e_1,e_2\}$ with $x+e\in\mathds{B}^N$
\begin{equation}\label{eq:construction_of_F}
\begin{split}
&\mathbf{Q}_{x,e}(\mathds{F}_{x,e}^N)
=
J\frac{1}{\pi\log N} \Big(\nabla^{\perp}_{\textup{disc}}\big(\log|x|\mathds{1}_{|x|\leq 1-\frac{10}{N}}\big)\Big)_e,
\end{split}
\end{equation}
where $\mathbf{Q}_{x,e}$ is a function of $\mathds{F}_{x,e}^N$ defined in \eqref{eq:auxillary_A_usual}.

Moreover, the forcing field satisfies the following consistency estimate for any $x\in\mathds{B}^N$, where $\mathfrak{u}$ is the abbreviation of $\mathfrak{u}^N[\gamma,\zeta;\alpha]$
\begin{equation}\label{eq:forcing_stability_prop}
\mathds{F}_{x}^{N}[\gamma,\zeta,\alpha,J]=\mathcal{F}^N[\gamma,\zeta,\alpha,J]+C_*O\Big(\frac{1}{N}\frac{1}{|x|^2}\Big)+C_*O\Big(\mathds{1}_{|x|\geq 1-\frac{20}{N}}\Big),\quad \mathds{F}_{x}^N[\gamma,\zeta,\alpha,J]\equiv 0\textup{ when $|x|\geq 1-\frac{5}{N}$}.
\end{equation}
\end{proposition}

\begin{proof}
Since $\mathfrak{u}^N[\gamma,\zeta;\alpha]$ (abbreviated as $\mathfrak{u}$) takes values between $\gamma$ and $\zeta$, there is a constant $c_{\gamma,\zeta}>0$ such that $\sigma(\mathfrak{u}_x)\geq c_{\gamma,\zeta}$ for $x\in\mathds{B}^N$.

\textit{Step 1: Existence of the forcing field $\mathds{F}_{x}^N[\gamma,\zeta,\alpha,J]$.}

At each $x\in\mathds{B}^N$ and $e\in\{e_1,e_2\}$, we regard the term $\mathbf{Q}_{x,e}$ as a function of $\mathds{F}\in\mathds{R}$ (we omit the superscript $N$)
\begin{equation*}
\mathbf{Q}_{x,e}(\mathds{F})=N\exp\Big(\frac{\mathds{F}}{N\log N}\Big)\mathfrak{u}_x(1-\mathfrak{u}_{x+e})-N\exp\Big(-\frac{\mathds{F}}{N\log N}\Big)\mathfrak{u}_{x+e}(1-\mathfrak{u}_x).
\end{equation*}
To prove Proposition \ref{prop:general_choice}, it is sufficient to find $\mathds{F}_{x,e}$ for each $(x,e)$ satisfying the consistency estimate \eqref{eq:forcing_stability_prop} and
\begin{equation}\label{eq:goal_of_forcing}
\mathbf{Q}_{x,e}(\mathds{F}_{x,e})=J\frac{1}{\pi\log N} \Big(\nabla^{\perp}_{\textup{disc}}\big(\log|x|\mathds{1}_{|x|\leq 1-\frac{10}{N}}\big)\Big)_e.
\end{equation}
For any $(x,e)$, the absolute value of RHS in \eqref{eq:goal_of_forcing} is bounded from above by
\begin{equation}\label{eq:bound_A_current}
|\mathbf{Q}_{x,e}(\mathds{F}_{x,e})|=\Big|J\frac{1}{\pi\log N} \Big(\nabla^{\perp}_{\textup{disc}}\big(\log|x|\mathds{1}_{|x|\leq 1-\frac{10}{N}}\big)\Big)_e\Big|\leq C_1^*|x|^{-1}(\log N)^{-1}.
\end{equation}
For the function $\mathbf{Q}_{x,e}(\mathds{F})$, there exists a constant $C_2^*>0$ such that
\begin{equation}\label{eq:intermediate_value}
\mathbf{Q}_{x,e}(C_2^*|x|^{-1})\geq C_1^*|x|^{-1}(\log N)^{-1},\quad \mathbf{Q}_{x,e}(-C_2^*|x|^{-1})\leq -C_1^*|x|^{-1}(\log N)^{-1},
\end{equation}
which shows by intermediate value theorem that, there exists $\mathds{F}_{x,e}\in (-C_2^*|x|^{-1},C_2^*|x|^{-1})$ solving \eqref{eq:goal_of_forcing}. If multiple values solve \eqref{eq:goal_of_forcing}, we choose only one of them, and denote it as $\mathds{F}_{x,e}^N[\gamma,\zeta,\alpha,J]$. This gives us a forcing field for any $(x,e)$, and we show below that this forcing satisfies the consistency estimate \eqref{eq:forcing_stability_prop}.

\textit{Step 2: Consistency estimate.}

When $|x|\geq 1-\frac{5}{N}$ we have $\mathbf{Q}_{x,e}\equiv 0$ according to the construction of $\mathbf{Q}_{x,e}$ in \eqref{eq:construction_of_F}. Furthermore the density profile $\mathfrak{u}_x$ is also constant equal to $\gamma$ for $|x|\geq 1-\frac{5}{N}$ by its definition in \eqref{eq:u_profile}. Expand $\mathbf{Q}_{x,e}$ in \eqref{eq:auxillary_A_usual} and we have $\mathds{F}_{x,e}\equiv 0$ for $|x|\geq 1-\frac{5}{N}$.

Using the bound $C_2^*|x|^{-1}$ of $\mathds{F}_x^N[\gamma,\zeta,\alpha,J]$ and the expression \eqref{eq:u_profile} of $\mathfrak{u}$
\begin{equation}\label{eq:taylor_expansion_u}
\exp\Big(\frac{\mathds{F}_{x,e}^N[\gamma,\zeta,\alpha,J]}{N\log N}\Big)-1=\frac{\mathds{F}_{x,e}^N[\gamma,\zeta,\alpha,J]}{N\log N}+C_*O\Big(\frac{|x|^{-2}}{N^2\log^2N}\Big),\quad \mathfrak{u}_{x+e_i}-\mathfrak{u}_{x}=\frac{1}{N}\nabla_{x_i}\mathfrak{u}_x+C_*O\Big(\frac{|x|^{-2}}{N^2\log N}\mathds{1}_{|x|\geq \frac{1}{2}N^{-\alpha}}\Big),
\end{equation}
we expand the function $\mathbf{Q}_{x,e}$ using Taylor expansion 
\begin{equation}\label{eq:expansion_of_A}
\mathbf{Q}_{x,e_i}=\frac{\mathds{F}_{x,e_i}^N[\gamma,\zeta,\alpha,J]}{\log N}2\sigma(\mathfrak{u}_x)-\nabla_{x_i}\mathfrak{u}_x+C_*O\Big(\frac{|x|^{-2}}{N\log N}\Big),\quad \textup{ for any $x\in\mathds{B}^N$ and $e\in\{e_1,e_2\}$.}
\end{equation}
Expand the RHS of \eqref{eq:goal_of_forcing}
\begin{equation}\label{eq:expansion_of_current}
\mathbf{Q}_{x,e}=J\frac{1}{\pi\log N} \Big(\nabla^{\perp}_{\textup{disc}}\big(\log|x|\mathds{1}_{|x|\leq 1-\frac{10}{N}}\big)\Big)_e=\frac{J}{\pi\log N}\big(\nabla^{\perp}\log|x|\big)_{e_i}+C_*O\Big(\frac{|x|^{-2}}{N\log N}\Big)+C_*O\Big(\frac{1}{\log N}\mathds{1}_{|x|\geq 1-\frac{20}{N}}\Big).
\end{equation}
The Taylor expansion of the LHS of \eqref{eq:goal_of_forcing} is given by \eqref{eq:expansion_of_A}, and the Taylor expansion of the RHS of \eqref{eq:goal_of_forcing} is given by \eqref{eq:expansion_of_current}. Applying the two expansions to \eqref{eq:goal_of_forcing}, we conclude the proof of Proposition \ref{prop:general_choice}.
\end{proof}

\subsection{Microscopic Transport and Cancellation Mechanism}\label{sec:cancellation}

This section proves the following estimate for any $t\in[0,T]$ and $\eps>0$
\begin{equation}\label{eq:cancellation_goal}
\mathds{E}_{\mathfrak{u},\mathds{F},t}^N\Big[\mathds{A}_{\textup{qua}}(\eta)\Big]\leq \eps \local{\mathfrak{u}}\Big[\Gamma_{\SS}^N(\sqrt{f_t})\Big]+\eps^{-1}C_*\mathds{H}^N\Big(d\mathds{P}_{\mathfrak{u},\mathds{F},t}^N\Big|d\local{\mathfrak{u}}\Big)+\eps^{-1}C_*N^{4\alpha}.
\end{equation}
on the quadratic part $\mathds{A}_{\textup{qua}}(\eta)$ of the
adjoint operator $L_{\mathds{F}}^{N,*}\mathds{1}$, which is the last input needed in the proof of Proposition
\ref{prop:relative_entropy}. 
Recall the structure of the quadratic term $\mathds{A}_{\textup{qua}}(\eta)$ in
\eqref{eq:linear_quadratic}
\begin{equation}\label{eq:reorganized_qua_1}
\mathds{A}_{\textup{qua}}(\eta)
=\sum_{x\in\mathds{B}^N}\omega_x\Big[
\mathbf{Q}_{x,e_1}N(\mathfrak{u}_{x+e_1}-\mathfrak{u}_x)\omega_{x+e_1}
+\mathbf{Q}_{x,e_2}N(\mathfrak{u}_{x+e_2}-\mathfrak{u}_x)\omega_{x+e_2}
\Big].
\end{equation}
The difficulty is that the coefficient $N\mathbf{Q}_{x,e}(\mathfrak{u}_{x+e}-\mathfrak{u})$ of
$\omega_x\omega_{x+e}$ is singular near the origin. The crucial observation is that the two directions must be treated together. Define the transported quadratic term by
\begin{equation}\label{eq:after_transport}
\mathds{A}_{\textup{qua}}^{\dagger}(\eta)
:=\sum_{x\in\mathds{B}^N}\omega_x\omega_{x+e_1+e_2}G_{\textup{qua},x}^{\tr},
\qquad
G_{\textup{qua},x}^{\tr}:=\mathbf{Q}_{x,e_1}N(\mathfrak{u}_{x+e_1}-\mathfrak{u}_x)+\mathbf{Q}_{x,e_2}N(\mathfrak{u}_{x+e_2}-\mathfrak{u}_x).
\end{equation}
Thus the variables $\omega_{x+e_1}$ and $\omega_{x+e_2}$ in
\eqref{eq:reorganized_qua_1} have both been moved to the common site
$x+e_1+e_2$. Function $G_{\textup{qua}}^{\dagger}$ has a good upper bound, which lead to a good upper bound of the term $\mathds{A}_{\textup{qua}}^{\dagger}$.

\begin{lemma}\label{lem:A_qua_dagger}
The function $G_{\textup{qua}}^{\tr}$ defined in \eqref{eq:after_transport}
satisfies
\begin{equation}\label{eq:G_qua_bound}
|G_{\textup{qua},x}^{\tr}|\leq C_*\frac{\mathds{1}_{N^{-\alpha}\leq |x|\leq 1}}{N\log^2N}\frac{1}{|x|^3}\textup{ for any $x\in\mathds{B}^N$};\quad \sum_{x\in\mathds{B}^N} |G_{\textup{qua},x}^{\tr}|^2\leq C_*\frac{N^{4\alpha}}{\log^4N}.
\end{equation}
When $\alpha\in (0,1/3)$, we further have $|G_{\textup{qua},x}^{\dagger}|\leq C_*$ for any $x\in\mathds{B}^N$. In addition, we have for any $t\in[0,T]$
\begin{equation}\label{eq:A_dagger_upper_bound}
\mathds{E}_{\mathfrak{u},\mathds{F},t}^N\big[|\mathds{A}_{\textup{qua}}^{\dagger}|\big]\leq C_*\mathds{H}^N\Big(d\mathds{P}_{\mathfrak{u},\mathds{F},t}^N\Big|\ d\local{\mathfrak{u}}\Big)+C_*N^{4\alpha}(\log N)^{-4}.
\end{equation}
\end{lemma}
\begin{proof}

Indeed, by the
Taylor expansion \eqref{eq:expansion_of_current} of $\mathbf{Q}_{x,e_i}$ and the Taylor expansion \eqref{eq:taylor_expansion_u} of $\mathfrak{u}$
\begin{equation*}
\mathbf{Q}_{x,e_i}
=\frac{J}{\pi\log N}(\nabla^\perp\log|x|)_{e_i}
+C_*O\Big(\frac{1}{N\log N}|x|^{-2}\Big)+C_*O\Big(\mathds{1}_{|x|\geq 1-\frac{20}{N}}\Big),\qquad N(\mathfrak{u}_{x+e_i}-\mathfrak{u}_x)
=(\nabla \mathfrak{u}_x)_{e_i}
+C_*O\Big(\frac{1}{N\log N}|x|^{-2}\Big).
\end{equation*}
Consequently the function $G_{\textup{qua}}^{\dagger}$ equals to
\begin{equation*}
\begin{split}
&\mathbf{Q}_{x,e_1}N(\mathfrak{u}_{x+e_1}-\mathfrak{u}_x)
+\mathbf{Q}_{x,e_2}N(\mathfrak{u}_{x+e_2}-\mathfrak{u}_x)=
\frac{J}{\pi\log N}\nabla^\perp\log|x|\cdot\nabla\mathfrak{u}_x
+C_*O\Big(\frac{1}{N\log^2N}|x|^{-3}\Big)+C_*O\Big(\frac{1}{\log N}\mathds{1}_{|x|\geq 1-\frac{20}{N}}\Big).
\end{split}
\end{equation*}
The term $\nabla^\perp\log|x|\cdot\nabla\mathfrak{u}_x$ vanishes because $\mathfrak{u}=\mathfrak{u}^N[\gamma,\zeta;\alpha]$ is
radially symmetric, while $\nabla^\perp\log|x|$ is tangent to the circle $\{|x|=\textup{const.}\}$. This proves the first inequality in \eqref{eq:G_qua_bound}.  Apply this point-wise estimate, and use Lemma \ref{lem:lattice_sum} to estimate the sum of $|x|^{-6}$
\begin{equation*}
\sum_{x\in\mathds{B}^N:|x|\geq N^{-\alpha}}
|G_{\textup{qua},x}^{\tr}|^2\leq \frac{C_*}{N^2\log^4N}\sum_{|x|\geq N^{-\alpha}}\frac{1}{|x|^6}\leq C_*\frac{N^{4\alpha}}{\log^4N}.
\end{equation*}

Now we have proved \eqref{eq:G_qua_bound}.  

Based on estimates \eqref{eq:G_qua_bound}, the next step is to prove \eqref{eq:A_dagger_upper_bound} by applying the entropic inequality in Lemma \ref{lem:entropic_ineq_correlation} with $\xi_x=\omega_x\omega_{x+e_1+e_2}$. For arbitrary $y\in N^{-1}\mathds{Z}^2$ with $|y|_{\infty}\geq 3N^{-1}$, the variable $\xi_{x}$ is independent from $\xi_{x+y}$. Therefore the integer $\ell$ in Lemma \ref{lem:entropic_ineq_correlation} can be chosen as $3$. Now the entropy inequality can be applied to $\mathds{A}_{\textup{qua}}^{\dagger}$ with properly chosen constants $C_*$ to absorb other constants in the statement of Lemma \ref{lem:entropic_ineq_correlation}
\begin{equation*}
\mathds{E}_{\mathfrak{u},\mathds{F},t}^N\big[|\mathds{A}_{\textup{qua}}^{\dagger}|\big]\leq
C_*\mathds{H}^N\Big(d\mathds{P}_{\mathfrak{u},\mathds{F},t}^N\Big|d\local{\mathfrak{u}}\Big)
+\log 2+C_*\max_{k\in\{+1,-1\}}\sum_{x\in\mathds{B}^N}\log\local{\mathfrak{u}}\Big[\exp\Big(kG_{\textup{qua},x}^{\dagger}\omega_x\omega_{x+e_1+e_2}\Big)\Big].
\end{equation*}
Further use the exponential estimate in Lemma \ref{lem:quadratic_subgaussian}
\begin{equation}\label{eq:transported_part_bound}
\mathds{E}_{\mathfrak{u},\mathds{F},t}^N\big[|\mathds{A}_{\textup{qua}}^{\dagger}|\big]\leq
C_*\mathds{H}^N\Big(d\mathds{P}_{\mathfrak{u},\mathds{F},t}^N\Big|d\local{\mathfrak{u}}\Big)+\log 2
+C_*\sum_{x\in\mathds{B}^N}|G_{\textup{qua},x}^{\dagger}|^2.
\end{equation}
The expectation is therefore bounded by
$C_*N^{4\alpha}/\log^4N+C_*\mathds{H}^N$ according to \eqref{eq:G_qua_bound}. This proves \eqref{eq:A_dagger_upper_bound}.
\end{proof}

Subtract the expression \eqref{eq:after_transport} of $\mathds{A}_{\textup{qua}}^{\dagger}$
from the expression \eqref{eq:reorganized_qua_1} of $\mathds{A}_{\textup{qua}}$
\begin{equation}\label{eq:difference_one_step}
\begin{split}
\mathds{A}_{\textup{qua}}-\mathds{A}_{\textup{qua}}^{\dagger}&=\sum_{x\in\mathds{B}^N} \omega_x\mathbf{Q}_{x,e_2}N(\mathfrak{u}_{x+e_2}-\mathfrak{u}_x)
\big(\omega_{x+e_2}-\omega_{x+e_1+e_2}\big)+\sum_{x\in\mathds{B}^N} \omega_x\mathbf{Q}_{x,e_1}N(\mathfrak{u}_{x+e_1}-\mathfrak{u}_x)
\big(\omega_{x+e_1}-\omega_{x+e_1+e_2}\big)\\
&=\sum_{x\in\mathds{B}^N} \omega_{x-e_2}\mathbf{Q}_{x-e_2,e_2}N(\mathfrak{u}_{x}-\mathfrak{u}_{x-e_2})
\big(\omega_{x}-\omega_{x+e_1}\big)+\sum_{x\in\mathds{B}^N} \omega_{x-e_1}\mathbf{Q}_{x-e_1,e_1}N(\mathfrak{u}_{x}-\mathfrak{u}_{x-e_1})
\big(\omega_{x}-\omega_{x+e_2}\big)\\
&=\sum_{x\in\mathds{B}^N,e\in\{e_1,e_2\}}(\omega_{x}-\omega_{x+e})E_{x,e},
\end{split}
\end{equation}
with
\begin{equation}\label{eq:difference_one_step_h}
E_{x,e_1}:=\omega_{x-e_2}\mathbf{Q}_{x-e_2,e_2}N(\mathfrak{u}_x-\mathfrak{u}_{x-e_2}),
\qquad
E_{x,e_2}:=\omega_{x-e_1}\mathbf{Q}_{x-e_1,e_1}N(\mathfrak{u}_x-\mathfrak{u}_{x-e_1}).
\end{equation}
Applying the one-bond integration-by-parts estimate
\eqref{eq:one_bond_ibp_app} to each summand gives, for every $\eps>0$
\begin{equation}\label{eq:transport_error}
\begin{split}
&
\local{\mathfrak{u}}\Big[
f_t\big(\mathds{A}_{\textup{qua}}-\mathds{A}_{\textup{qua}}^{\dagger}\big)
\Big]\leq
\eps \local{\mathfrak{u}}\Big[\Gamma_{\SS}^N(\sqrt{f_t})\Big]
+\frac{C_*}{\eps}
\local{\mathfrak{u}}\big[f_t\mathds{A}_{\textup{qua}}^W\big]
+\local{\mathfrak{u}}\big[f_t\mathds{A}_{\textup{qua}}^Z\big],
\end{split}
\end{equation}
where the constant $\eps_0^2$ in \eqref{eq:one_bond_ibp_app} has been absorbed into $C_*$ and
\begin{equation}\label{eq:AZAW}
    \mathds{A}_{\textup{qua}}^W:=N^{-2}\sum_{x\in\mathds{B}^N,e\in\{e_1,e_2\}}(E_{x,e})^2,
    \qquad
    \mathds{A}_{\textup{qua}}^Z
    :=\sum_{x\in\mathds{B}^N,e\in\{e_1,e_2\}}(\mathfrak{u}_{x+e}-\mathfrak{u}_x)E_{x,e}\omega_x\omega_{x+e}.
\end{equation}
The Dirichlet form in \eqref{eq:transport_error} is controlled by
$\local{\mathfrak{u}}[\Gamma_{\mathds{F}}^N\sqrt {f_t}]$, up to a constant depending
only on $(\gamma,\zeta,\alpha,J)$; this follows from the boundedness from below of
$\mathds{F}_{x,e}/(N\log N)$ on the lattice. It remains to estimate the transport error used in
\eqref{eq:AZAW}.

\begin{lemma}\label{lem:AWAZ}
For the terms $\mathds{A}_{\textup{qua}}^W$ and
$\mathds{A}_{\textup{qua}}^Z$ defined above, and for
$f_t=d\mathds{P}_{\mathfrak{u},\mathds{F},t}^N/d\local{\mathfrak{u}}$,
\begin{equation*}
\local{\mathfrak{u}}\big[f_t |\mathds{A}_{\textup{qua}}^{W}|\big]
\leq
C_*\frac{N^{2\alpha}}{\log^2 N},\quad \local{\mathfrak{u}}\big[f_t|\mathds{A}_{\textup{qua}}^Z|\big]
\leq
C_*
\mathds{H}^N\Big(d\mathds{P}_{\mathfrak{u},\mathds{F},t}^N\Big|d\local{\mathfrak{u}}\Big)
+C_*\frac{N^{4\alpha}}{\log^6N}.
\end{equation*}
\end{lemma}
\begin{proof}
To estimate the coefficient $E_{x,e}$ in \eqref{eq:difference_one_step_h}, we recall the estimate \eqref{eq:bound_A_current} of $\mathbf{Q}_{x,e}$, which gives $|\mathbf{Q}_{x,e}|\leq C_*|x|^{-1}(\log N)^{-1}$. The definition of $\mathfrak{u}$ in \eqref{eq:u_profile} gives $|N(\mathfrak{u}_{x+e}-\mathfrak{u})|\leq C_*|x|^{-1}(\log N)^{-1}\mathds{1}_{|x|\geq N^{-\alpha}}$. Since $\omega$ is uniformly bounded by a constant $C_*$, this gives the pointwise estimate $|E_{x,e}|
\leq
C_*|x|^{-2}(\log N)^{-2}\mathds{1}_{\{|x|\geq N^{-\alpha}\}}$. Then we use Lemma \ref{lem:lattice_sum} to estimate the sum of $|x|^{-4}$
\begin{equation*}
\begin{split}
\local{\mathfrak{u}}\big[f_t|\mathds{A}_{\textup{qua}}^W|\big]
&\leq
C_*
\frac{1}{N^2}\sum_{x\in\mathds{B}^N}
\frac{1}{|x|^4\log^4N}\mathds{1}_{\{|x|\geq N^{-\alpha}\}}\leq C_*\frac{N^{2\alpha}}{\log^4N}.
\end{split}
\end{equation*}
For the term $\mathds{A}_{\textup{qua}}^Z$, we estimate the coefficient $(\mathfrak{u}_{x+e}-\mathfrak{u})E_{x,e}$
\begin{equation*}
\Big|(\mathfrak{u}_{x+e}-\mathfrak{u}_x)E_{x,e}\Big|\leq \Big|\mathfrak{u}_{x+e}-\mathfrak{u}_x\Big|\ \Big| E_{x,e}\Big|\leq C_*\frac{|x|^{-3}}{N\log^3 N}\mathds{1}_{|x|\geq N^{-\alpha}}.
\end{equation*}
Apply the entropic inequality \eqref{eq:entropic_inequality_app} to $\mathds{A}_{\textup{qua}}^Z$ and use Lemma \ref{lem:lattice_sum}
\begin{equation*}
\begin{split}
\local{\mathfrak{u}}\big[f_t|\mathds{A}_{\textup{qua}}^Z|\big]
&\leq
C_*
\mathds{H}^N\Big(d\mathds{P}_{\mathfrak{u},\mathds{F},t}^N\Big|d\local{\mathfrak{u}}\Big)+\log 2
+C_*\sum_{x\in\mathds{B}^N}
\frac{\mathds{1}_{|x|\geq N^{-\alpha}}}{N^2\log^6N}\frac{1}{|x|^6}
\leq C_*\mathds{H}^N\Big(d\mathds{P}_{\mathfrak{u},\mathds{F},t}^N\Big|\ d\local{\mathfrak{u}}\Big)+C_*N^{4\alpha}(\log N)^{-6}.
\end{split}
\end{equation*}
This concludes the proof of Lemma \ref{lem:AWAZ}.
\end{proof}

\subsection{Proof of the Entropy Bound in Proposition \ref{prop:relative_entropy}}\label{sec:proof_entropy}
\begin{proof}[Proof of Proposition \ref{prop:relative_entropy}]
Recall Yau's inequality \eqref{eq:Yau}
\begin{equation}\label{eq:entropy_proof_yau}
\frac{d}{dt}\mathds{H}^N\Big(d\mathds{P}_{\mathfrak{u},\mathds{F},t}^N\Big|d\local{\mathfrak{u}}\Big)
\leq
-\local{\mathfrak{u}}\Big[\Gamma_{\mathds{F}}^N\sqrt{f_t}\Big]
+\local{\mathfrak{u}}\Big[f_t L_{\mathds{F}}^{N,*}\mathds{1}\Big],\quad f_t:=\frac{d\mathds{P}_{\mathfrak{u},\mathds{F},t}^N}{d\local{\mathfrak{u}}}.
\end{equation}
The adjoint formula of Lemma \ref{lem:linear_quadratic} decomposes the last term into a linear contribution and a quadratic contribution
\begin{equation}\label{eq:entropy_proof_decomposition}
L_{\mathds{F}}^{N,*}\mathds{1}
=\frac{1}{2}\mathds{A}_{\textup{lin}}+\frac{1}{2}\mathds{A}_{\textup{qua}},
\end{equation}
where the two terms are given in \eqref{eq:linear_quadratic}. We estimate them separately. We will first prove
\begin{equation}\label{eq:A_lin_recalled}
\mathds{A}_{\textup{lin}}=-\sum_{x\in\mathds{B}^N\setminus \p\mathds{B}^N}\omega_x\nabla_{\textup{disc}}\cdot\mathbf{Q}_{x}+N\sum_{x\in \p\mathds{B}^N}\omega_x\sum_{e\in\{e_1,e_2\}}\Big(\mathbf{Q}_{x-e,e}\mathds{1}_{x-e\in\mathds{B}^N}-\mathbf{Q}_{x,e}\mathds{1}_{x+e\in\mathds{B}^N}\Big)\equiv 0.
\end{equation}
The forcing $\mathds{F}$ is constructed in Proposition \ref{prop:general_choice} so that the field $\mathbf{Q}_{x,e}$ itself is the discrete skew-gradient $\nabla_{\textup{disc}}^{\perp}$ in \eqref{eq:construction_of_F}. Therefore its discrete divergence vanishes. More precisely
\begin{equation}\label{eq:linear_bulk_cancellation}
-\nabla_{\textup{disc}}\cdot\mathbf{Q}_{x}
=-\frac{J}{\pi\log N}\nabla_{\textup{disc}}\cdot\nabla_{\textup{disc}}^{\perp}\big(\log|x|\mathds{1}_{|x|\leq 1-\frac{10}{N}}\big)\equiv 0,
\end{equation}
according to \eqref{eq:discrete_vanish}. The forcing $\mathds{F}$ constructed in Proposition \ref{prop:general_choice} also implies $\mathbf{Q}_{x,e}\equiv 0$ when $|x|\geq 1-\frac{5}{N}$ due to the truncation $\mathds{1}_{|x|\leq 1-\frac{10}{N}}$ in \eqref{eq:construction_of_F}. Therefore we have $\mathds{A}_{\textup{lin}}\equiv 0$ in \eqref{eq:A_lin_recalled}. 

Now we estimate $\mathds{A}_{\textup{qua}}$. Section \ref{sec:cancellation} proves \eqref{eq:cancellation_goal} for any $\eps>0$
\begin{equation}\label{eq:quadratic_entropy_input}
\begin{split}
\mathds{E}_{\mathfrak{u},\mathds{F},t}^N\Big[\mathds{A}_{\textup{qua}}(\eta)\Big]\leq \eps \local{\mathfrak{u}}\Big[\Gamma_{\SS}^N(\sqrt{f_t})\Big]+\eps^{-1}C_*\mathds{H}^N\Big(d\mathds{P}_{\mathfrak{u},\mathds{F},t}^N\Big|d\local{\mathfrak{u}}\Big)+\eps^{-1}C_*N^{4\alpha}.
\end{split}
\end{equation}
Inserting $\eps=\frac{1}{2}$, the positive Dirichlet form in
\eqref{eq:quadratic_entropy_input} is absorbed by the negative Dirichlet form in
\eqref{eq:entropy_proof_yau}.  Insert the estimate \eqref{eq:A_lin_recalled} of $\mathds{A}_{\textup{lin}}$ and the estimate \eqref{eq:quadratic_entropy_input} of $\mathds{A}_{\textup{qua}}$ into Yau's inequality \eqref{eq:entropy_proof_yau}
\begin{equation}\label{eq:entropy_gronwall_differential}
\frac{d}{dt}\mathds{H}^N\Big(d\mathds{P}_{\mathfrak{u},\mathds{F},t}^N\Big|d\local{\mathfrak{u}}\Big)
\leq
C_*\mathds{H}^N\Big(d\mathds{P}_{\mathfrak{u},\mathds{F},t}^N\Big|d\local{\mathfrak{u}}\Big)
+C_*N^{4\alpha}.
\end{equation}
Since the modified dynamics is started from the reference local equilibrium $d\local{\mathfrak{u}}$, we have $\mathds{H}^N\Big(d\mathds{P}_{\mathfrak{u},\mathds{F},0}^N\Big|d\local{\mathfrak{u}}\Big)=0$. Applying Gronwall's inequality to \eqref{eq:entropy_gronwall_differential} on $[0,T]$ gives for $0\leq t\leq T$
\begin{equation*}
\mathds{H}^N\Big(d\mathds{P}_{\mathfrak{u},\mathds{F},t}^N\Big|d\local{\mathfrak{u}}\Big)\leq C_*N^{4\alpha}.
\end{equation*}
This proves Proposition \ref{prop:relative_entropy}.
\end{proof}

\begin{remark}
In Section \ref{sec:concentration_RN}, we also use the estimate of the adjoint operator $L_{\mathds{F}}^{N,*}$, thus we gather a useful inequality here. In the decomposition \eqref{eq:entropy_proof_decomposition}, we have shown $\mathds{A}_{\textup{lin}}\equiv 0$ by \eqref{eq:A_lin_recalled}, and $\mathds{A}_{\textup{qua}}$ is estimated by \eqref{eq:transport_error}. Consequently for any probability density $f$ with respect to $d\local{\mathfrak{u}}$
\begin{equation}\label{eq:sec6.4remark1}
\local{\mathfrak{u}}\Big[f\mathds{L}_{\mathds{F}}^{N,*}\mathds{1}\Big]\leq \eps \local{\mathfrak{u}}\Big[\Gamma_{\SS}^N(\sqrt{f})\Big]+\frac{1}{2}\local{\mathfrak{u}}\Big[f\mathds{A}_{\textup{qua}}^{\dagger}\Big]+\frac{C_*}{\eps}\local{\mathfrak{u}}\Big[f\mathds{A}_{\textup{qua}}^W\Big]+\frac{1}{2}\local{\mathfrak{u}}\Big[f\mathds{A}_{\textup{qua}}^Z\Big].
\end{equation}
When $f_t:=\frac{d\mathds{P}_{\mathfrak{u},\mathds{F},t}^N}{d\local{\mathfrak{u}}}$, Lemma \ref{lem:A_qua_dagger}, \ref{lem:AWAZ}, and Proposition \ref{prop:relative_entropy} show
\begin{equation}\label{eq:sec6.4remark2}
\frac{1}{2}\local{\mathfrak{u}}\Big[f_t|\mathds{A}_{\textup{qua}}^{\dagger}|\Big]+\frac{C_*}{\eps}\local{\mathfrak{u}}\Big[f_t|\mathds{A}_{\textup{qua}}^W|\Big]+\frac{1}{2}\local{\mathfrak{u}}\Big[f_t|\mathds{A}_{\textup{qua}}^Z|\Big]\leq C_*N^{4\alpha}.
\end{equation}
\end{remark}

\section{Concentration Estimate for Radon--Nikodym Derivative}\label{sec:concentration_RN}
Based on the relative entropy estimate in Proposition \ref{prop:relative_entropy}, we now prove the concentration estimate \eqref{eq:concentration_statement_sec5} used in proving the lower bound. We need to show that, after inserting the empirical density and current into the Radon--Nikodym functional $\mathbf{R}_{T,\mathds{F}}^N$ in \eqref{eq:discrete_R_functional}, the whole expression is asymptotically equal to a deterministic value given by the functional $\mathcal{R}_{T,\mathcal{F}}^N$ as in \eqref{eq:continuum_R_cost}. Here we recall that $C_*$ are constants that may change from line to line, depending on parameters $(\gamma,\zeta,\alpha,T,J,\mathfrak{b})$. We also use $\logles$ as $\leq$ neglecting logarithmic factors.

\begin{proposition}\label{prop:concentration}
\textbf{\textup{[Concentration of the Radon--Nikodym functional]}}
Take arbitrary parameters $(\gamma,\zeta,\alpha,J,N)$ with
\begin{equation*}
\gamma\in (0,1),\quad \zeta\in (0,1),\quad \alpha\in \Big(0,\frac{1}{10}\Big),\quad J>0,\quad N>1.
\end{equation*}
Let $\mathfrak{u}$ denote the density profile $\mathfrak{u}^N[\gamma,\zeta;\alpha]$ in \eqref{eq:u_profile}, and let $\mathfrak{q}$ denote the current profile $\mathfrak{q}^N[J]$ in \eqref{eq:u_profile}. Take $\mathds{F}=\mathds{F}^N[\gamma,\zeta,\alpha,J]$ be the microscopic forcing constructed in Proposition \ref{prop:general_choice}.  Let $\mathcal{F}$ denote the forcing $\mathcal{F}^N[\gamma,\zeta,\alpha,J]$ in \eqref{eq:continuum_forcing_parametrized}. Then, for every $\delta>0$ and $T>0$
\begin{equation}\label{eq:RN_concentration_statement}
\mathds{P}_{\mathfrak{u},\mathds{F}}^{N}\Big[
\Big|\mathbf{R}_{T,\mathds{F}}^N(\rho^N,\mathcal{Q}^N)
-\mathcal{R}_{T,\mathcal{F}}^N(\mathfrak{u},\mathfrak{q})\Big|>\delta\Big]=C_*o_N(1).
\end{equation}
\end{proposition}

\begin{proof}
Define a discrete deterministic approximation of $\mathbf{R}_{T,\mathds{F}}^N(\rho^N,\mathcal{Q}^N)$ in Definition \ref{def:R_functionals} by replacing $\mathcal{Q}_T^N$ as $T\mathfrak{q}$ and every occupation variable by $\mathfrak{u}$ in \eqref{eq:discrete_R_functional}
\begin{equation}\label{eq:deterministic_discrete_R}
\begin{split}
\mathbf{R}_{T,\mathds{F}}^{N,\det}(\mathfrak{u},\mathfrak{q})
:=& T\langle \mathfrak{q},\mathds{F}\rangle-T\sum_{x\in\mathds{B}^N,e\in\{e_1,e_2\}}\frac{1}{4N^2}\frac{(\mathds{F}_{x,e})^2}{\log N}
\Big(\mathfrak{u}_x(1-\mathfrak{u}_{x+e})+\mathfrak{u}_{x+e}(1-\mathfrak{u}_x)\Big)-T\sum_{x\in\mathds{B}^N}\frac{1}{2N^2}\mathfrak{u}_x\ddiv\mathds{F}_x.
\end{split}
\end{equation}
We want to compare the difference between the two deterministic functionals $\mathbf{R}_{T,\mathds{F}}^{N,\det}(\mathfrak{u},\mathfrak{q})$ and $\mathcal{R}_{T,\mathcal{F}}^N(\mathfrak{u},\mathfrak{q})$. The forcing $\mathcal{F}$ in the second functional is determined by $(\mathfrak{u},\mathfrak{q})$ as in \eqref{eq:continuum_forcing_parametrized}, and the forcing $\mathds{F}$ in the first functional is a slightly modified version of $\mathcal{F}$ as in Proposition \ref{prop:general_choice}, whose error is quantitatively negligible. Lemma \ref{lem:deterministic_consistency} below proves that
\begin{equation}\label{eq:deterministic_consistency_used}
\mathbf{R}_{T,\mathds{F}}^{N,\det}(\mathfrak{u},\mathfrak{q})
=\mathcal{R}_{T,\mathcal{F}}^N(\mathfrak{u},\mathfrak{q})+C_*o_N(1).
\end{equation}
It remains to show that the random difference $\mathbf{R}_{T,\mathds{F}}^N(\rho^N,\mathcal{Q}^N)
-\mathbf{R}_{T,\mathds{F}}^{N,\det}(\mathfrak{u},\mathfrak{q})$ vanishes in $\mathds{P}_{\mathfrak{u},\mathds{F}}^N$-probability.  According to the definition of $\mathbf{R}_{T,\mathds{F}}^N$ in \eqref{eq:discrete_R_functional}, this difference is
\begin{equation}\label{eq:RN_random_decomposition}
\begin{split}
&\underbrace{\langle\mathcal{Q}_T^N,\mathds{F}\rangle
-T\langle \mathfrak{q},\mathds{F}\rangle}_{\textup{\blue{current term}}}-\underbrace{\int_0^T\sum_{x\in\mathds{B}^N,e\in\{e_1,e_2\}}\frac{1}{2N^2}
\frac{(\mathds{F}_{x,e})^2}{\log N}
\Big[\sigma_{x,e}(\eta_t)-\sigma_{x,e}(\mathfrak{u})\Big]dt}_{\textup{\blue{mobility fluctuation}}}-\underbrace{\int_0^T\sum_{x\in\mathds{B}^N}\frac{1}{2N^2}
(\eta_{t,x}-\mathfrak{u}_x)\ddiv\mathds{F}_xdt}_{\textup{\blue{divergence fluctuation}}},
\end{split}
\end{equation}
where $\sigma_{x,e}$ is defined in \eqref{eq:sigmae}. The three terms in \eqref{eq:RN_random_decomposition} are estimated below. This will complete the proof of Proposition \ref{prop:concentration}.

\textit{The concentration of the mobility fluctuation.}

For the mobility fluctuation in \eqref{eq:RN_random_decomposition}, we use the centered configuration $\omega_x:=\frac{\eta_x-\mathfrak{u}_x}{\mathfrak{u}_x(1-\mathfrak{u}_x)}$. Lemma \ref{lem:centralized_decomposition} decomposes $\sigma_{x,e}(\eta)-\sigma_{x,e}(\mathfrak{u})$ into one quadratic term in $\omega_x\omega_{x+e}$ and two linear terms in $\omega_x,\omega_{x+e}$:
\begin{equation*}
\begin{split}
&\sigma_{x,e}(\eta_t)-\sigma_{x,e}(\mathfrak{u})\\
&=\underbrace{-\mathfrak{u}_{x}(1-\mathfrak{u}_{x})\mathfrak{u}_{x+e}(1-\mathfrak{u}_{x+e})
\omega_{x}\omega_{x+e}}_{\textup{\blue{quadratic term}}}+\underbrace{\frac{1}{2}(1-2\mathfrak{u}_{x+e})\mathfrak{u}_{x}(1-\mathfrak{u}_{x})\omega_{x}
+\frac{1}{2}(1-2\mathfrak{u}_{x})\mathfrak{u}_{x+e}(1-\mathfrak{u}_{x+e})\omega_{x+e}}_{\textup{\blue{linear term}}}.
\end{split}
\end{equation*}
Consequently there exist functions $\mathbf{H}_{1,e}(x)$ and $\mathbf{H}_2(x)$ such that, the mobility fluctuation in \eqref{eq:RN_random_decomposition} can be decomposed into
\begin{equation*}
\begin{split}
&\sum_{x\in\mathds{B}^N,e\in\{e_1,e_2\}}N^{-2}\mathbf{H}_{1,e}(x)\omega_x\omega_{x+e}+\sum_{x\in\mathds{B}^N}N^{-2}\mathbf{H}_2(x)\omega_x,\\
&\hspace{5cm} \textup{ with $|\mathbf{H}_{1,e}(x)|\leq C_*|x|^{-2}(\log N)^{-1}$ and $|\mathbf{H}_2(x)|\leq C_*|x|^{-2}(\log N)^{-1}$}.
\end{split}
\end{equation*}
Lemma \ref{lem:concentration_quadratic} gives the following quadratic concentration estimate for any $e\in\{e_1,e_2\}$
\begin{equation}\label{eq:quadratic_concentration_within}
\mathds{P}_{\mathfrak{u},\mathds{F}}^N\Bigg[\ \Big|\int_0^T\sum_{x\in\mathds{B}^N}N^{-2}\mathbf{H}_{1,e}(x)\omega_{t,x}\omega_{t,x+e}dt\Big|>\delta\ \Bigg]=C_*o_N(1).
\end{equation}
Lemma \ref{lem:concentration_linear} gives the following linear concentration estimate for any $e\in\{e_1,e_2\}$
\begin{equation}\label{eq:linear_concentration_within}
\mathds{P}_{\mathfrak{u},\mathds{F}}^N\Bigg[\ \Big|\int_0^T\sum_{x\in\mathds{B}^N}N^{-2}\mathbf{H}_2(x)\omega_{t,x}dt\Big|>\delta\ \Bigg]=C_*o_N(1).
\end{equation}
The two concentration estimates together derive the concentration of the mobility fluctuation
\begin{equation*}
\mathds{P}_{\mathfrak{u},\mathds{F}}^N\Bigg[\ \Big|\int_0^T\sum_{x\in\mathds{B}^N,e\in\{e_1,e_2\}}\frac{1}{2N^2}
\frac{(\mathds{F}_{x,e})^2}{\log N}
\Big(\sigma_{x,e}(\eta_t)-\sigma_{x,e}(\mathfrak{u})\Big)dt\Big|>\delta \Bigg]=C_*o_N(1).
\end{equation*}
Note that in Section \ref{sec:relative_entropy}, we have estimated similar linear-in-$\omega$ and quadratic-in-$\omega$ terms, which come from the entropy production in Lemma \ref{lem:linear_quadratic}. However the coefficients are different, since here the linear and quadratic terms come from the Radon-Nikodym derivative $\mathbf{R}_{T,\mathds{F}}^N$ instead of entropy production. Here the coefficients of the quadratic terms are always positive, not allowing any cancellation similar to Section \ref{sec:cancellation}. To overcome this difficulty, we adapt the averaging method of \cite{Jara_Menezes_2018} into a multi-scale manner with averaging radius $\ell$ in \eqref{eq:sec8_averaging_radius}, where the averaging radius will be tuned according to the parameter $\alpha$ in \eqref{eq:u_profile} as a measurement of the singularity. Specifically given a fixed parameter $\beta\in(0,1)$, we define $\ell$ as
\begin{equation}\label{eq:sec8_averaging_radius}
\ell:=0.5N^{1-\beta}.
\end{equation}
This $\ell$ will be taken as the averaging radius of $\omega_x$ for $|x|\leq N^{-\beta}$ with properly chosen $\beta$, in the following proof of concentration estimates. Outside $|x|\leq N^{-\beta}$ we will not average, because we are away from the origin and the relevant functions are not singular. Notice that this function $\ell$ is different from the averaging radius $\mathfrak{l}_x$ as \eqref{eq:averaging_radius} used in the large deviation upper bound. When $x$ is close to the origin in the sense of $|x|\leq N^{-\beta}$, the averaging radius $\ell$ is no longer proportional to $N|x|$, which is the case of $\mathfrak{l}_x$. Instead, $\ell$ is fixed as $N\times N^{-\beta}$, which allows us to derive \textit{concentration estimates} even in the region $|x|\leq N^{-\beta}$ where the relevant functions are very singular.

\textit{The concentration of the divergence fluctuation:}

For the divergence fluctuation term in \eqref{eq:RN_random_decomposition}, there exists a function $\mathbf{H}_3(x)$ with
\begin{equation}\label{eq:divergence_omega_decomposition}
\sum_{x\in\mathds{B}^N}\frac{1}{2N^2}
(\eta_{t,x}-\mathfrak{u}_x)\ddiv\mathds{F}_x=\sum_{x\in\mathds{B}^N}N^{-2}\mathbf{H}_3(x)\omega_x,\quad\textup{ with $|\mathbf{H}_3(x)|\leq C_*|x|^{-2}$.}
\end{equation}
Lemma \ref{lem:concentration_linear} gives the concentration estimate of the RHS in \eqref{eq:divergence_omega_decomposition}, therefore derives
\begin{equation*}
\mathds{P}_{\mathfrak{u},\mathds{F}}^N\Bigg[\ \Big|\int_0^T\sum_{x\in\mathds{B}^N}\frac{1}{2N^2}
(\eta_{t,x}-\mathfrak{u}_x)\ddiv\mathds{F}_xdt\Big|>\delta\ \Bigg]=C_*o_N(1).
\end{equation*}

\textit{The concentration of the current term.}

Finally, the concentration estimate of the current term in \eqref{eq:RN_random_decomposition} is proved by Lemma \ref{lem:current_concentration}.

\end{proof}

\begin{lemma}\label{lem:deterministic_consistency}
Under the notations of Proposition \ref{prop:concentration}
\begin{equation}\label{eq:deterministic_consistency_statement}
\left|\mathbf{R}_{T,\mathds{F}}^{N,\det}(\mathfrak{u},\mathfrak{q})
-\mathcal{R}_{T,\mathcal{F}}^N(\mathfrak{u},\mathfrak{q})\right|=C_*o_N(1).
\end{equation}
\end{lemma}

\begin{proof}
It is sufficient to compare the three terms in the expression \eqref{eq:deterministic_discrete_R} of $\mathbf{R}_{T,\mathds{F}}^{N,\textup{det}}$ with the three terms in the continuum functional $\mathcal{R}_{T,\mathcal{F}}^N$.  Equivalently, we prove
\begin{subequations}\label{eq:reduced_concentration}
\begin{align}
& \left|\frac{1}{N^2}\sum_{x\in\mathds{B}^N} \mathfrak{q}_x\cdot\mathds{F}_x-\int_{N^{-1}\leq |x|\leq 1} \mathfrak{q}_x\cdot\mathcal{F}_xdx\right|=C_*o_N(1),\label{eq:reduced_concentration_a}\\
&\left|\sum_{x\in\mathds{B}^N,e\in\{e_1,e_2\}}\frac{1}{4N^2}\frac{|\mathds{F}_{x,e}|^2}{\log N}
\Big(\mathfrak{u}_x(1-\mathfrak{u}_{x+e})+\mathfrak{u}_{x+e}(1-\mathfrak{u}_x)\Big)
-\frac{1}{2}\int_{N^{-1}\leq |x|\leq 1}\frac{|\mathcal{F}_x|^2}{\log N}\sigma(\mathfrak{u}_x)dx\right|=C_*o_N(1),\label{eq:reduced_concentration_b}\\
&\left|\sum_{x\in\mathds{B}^N}\frac{1}{2N^2}\mathfrak{u}_x\ddiv\mathds{F}(x)
-\frac{1}{2}\int_{N^{-1}\leq |x|\leq 1} \mathfrak{u}_x\nabla\cdot\mathcal{F}_xdx\right|=C_*o_N(1).\label{eq:reduced_concentration_c}
\end{align}
\end{subequations}
The consistency estimate \eqref{eq:forcing_stability_prop} says
\begin{equation}\label{eq:F_consistency_for_concentration}
\mathds{F}_x=\mathcal{F}_x+C_*O\Big(\frac{1}{N|x|^2}\Big)+C_*O\Big(\mathds{1}_{|x|\geq 1-\frac{20}{N}}\Big).
\end{equation}
Moreover, from the explicit profile \eqref{eq:u_profile} of $(\mathfrak{u},\mathfrak{q})$ and the formula \eqref{eq:continuum_forcing_parametrized} for $\mathcal{F}$
\begin{equation}\label{eq:regularity_reduced_concen}
|\nabla\mathfrak{u}(x)|
\leq C_*\frac{1}{|x|\log N}\mathds{1}_{|x|\geq N^{-\alpha}},\quad |\nabla\mathcal{F}(x)|\leq C_*\frac{1}{|x|^2} ,\quad |\nabla(\mathfrak{q}\cdot\mathcal{F})(x)|
\leq C_*\frac{1}{\log N}\frac{1}{|x|^3}.
\end{equation}

\textit{Part 1: Proof of \eqref{eq:reduced_concentration_a} on the current term.}

Using $\mathfrak{q}=(J/2\pi\log N)\nabla^\perp\log|x|$ and the consistency estimate \eqref{eq:F_consistency_for_concentration}, we measure the error of replacing $\mathds{F}$ by $\mathcal{F}$
\begin{equation*}
\frac{1}{N^2}\sum_{x\in\mathds{B}^N} \mathfrak{q}_x\cdot\mathds{F}_x
=\frac{1}{N^2}\sum_{x\in\mathds{B}^N} \mathfrak{q}_x\cdot\mathcal{F}(x)
+\frac{1}{N^2}\sum_{x\in\mathds{B}^N}
C_*O\Big(\frac{1}{N\log N}|x|^{-3}\Big)+\frac{1}{N^2}\sum_{x\in\mathds{B}^N}C_*O\Big(\mathds{1}_{|x|\geq 1-\frac{20}{N}}\Big).
\end{equation*}
The error is $C_*o_N(1)$ due to Lemma \ref{lem:lattice_sum} estimating the sum of $|x|^{-3}$
\begin{equation*}
\frac{1}{N^3\log N}\sum_{x\in\mathds{B}^N}|x|^{-3}
\leq
\frac{C}{\log N},\quad \frac{1}{N^2}\sum_{x\in\mathds{B}^N}\mathds{1}_{|x|\geq 1-\frac{20}{N}}\leq CN^{-1}.
\end{equation*}
The Riemann-sum error for $\frac{1}{N^2}\sum_{x\in\mathds{B}^N} \mathfrak{q}_x\cdot\mathcal{F}(x)$ is also $C_*o_N(1)$, by the last estimate in \eqref{eq:regularity_reduced_concen}.  This proves the first line of \eqref{eq:reduced_concentration}.

\textit{Part 2: Proof of \eqref{eq:reduced_concentration_b} on the mobility term.}

From \eqref{eq:F_consistency_for_concentration} together with the bound $|\mathds{F}_{x,e}|\leq C_*|x|^{-1}$ and $|\mathcal{F}_{x,e}|\leq C_*|x|^{-1}$, we derive
\begin{equation*}
|\mathds{F}_{x,e}|^2=|\mathcal{F}_{x,e}(x)|^2+
C_*O\Big(\frac{1}{N}|x|^{-3}\Big)+C_*O\Big(\mathds{1}_{|x|\geq 1-\frac{20}{N}}\Big).
\end{equation*}
Also the regularity \eqref{eq:regularity_reduced_concen} of $\mathfrak{u}$ shows $\mathfrak{u}_x(1-\mathfrak{u}_{x+e})+\mathfrak{u}_{x+e}(1-\mathfrak{u}_x)
=2\sigma(\mathfrak{u}_x)+C_*O\big((N|x|\log N)^{-1}\big)$. Therefore summing the errors, we derive
\begin{equation*}
\begin{split}
\sum_{x\in\mathds{B}^N,e\in\{e_1,e_2\}}\frac{1}{4N^2}\frac{|\mathds{F}_{x,e}|^2}{\log N}
\Big(\mathfrak{u}_x(1-\mathfrak{u}_{x+e})+\mathfrak{u}_{x+e}(1-\mathfrak{u}_x)\Big)&=\frac{1}{2}\sum_{x\in\mathds{B}^N,e\in\{e_1,e_2\}}\frac{1}{N^2}\frac{|\mathds{F}_{x,e}|^2}{\log N}\sigma(\mathfrak{u}_x)+C_*o_N(1)\\
&=\frac{1}{2}\sum_{x\in\mathds{B}^N,e\in\{e_1,e_2\}}\frac{1}{N^2}\frac{|\mathcal{F}_{x,e}|^2}{\log N}\sigma(\mathfrak{u}_x)+C_*o_N(1).
\end{split}
\end{equation*}
Replacing the discrete summation by integrals
\begin{equation*}
\frac{1}{2}\sum_{x\in\mathds{B}^N,e\in\{e_1,e_2\}}\frac{1}{N^2}\frac{|\mathcal{F}_{x,e}|^2}{\log N}\sigma(\mathfrak{u}_x)=\frac{1}{2}\int_{N^{-1}\leq |x|\leq 1}\frac{|\mathcal{F}|^2}{\log N}\sigma(\mathfrak{u})dx+C_*o_N(1).
\end{equation*}
we prove the second line of \eqref{eq:reduced_concentration}.

\textit{Part 3: Proof of \eqref{eq:reduced_concentration_c} on the divergence term.}

The discrete integration-by-parts gives
\begin{equation*}
\sum_{x\in\mathds{B}^N}\frac{1}{2N^2}\mathfrak{u}_x\ddiv\mathds{F}(x)
=-\frac{1}{2}\sum_{x\in\mathds{B}^N,e\in\{e_1,e_2\}}\frac{1}{N^2}N(\mathfrak{u}_{x+e}-\mathfrak{u}_x)\mathds{F}_{x,e}.
\end{equation*}
Using \eqref{eq:F_consistency_for_concentration} and the Taylor expansion of $\mathfrak{u}$, the right-hand side is the Riemann approximation of
\begin{equation*}
-\frac{1}{2}\int_{N^{-1}\leq |x|\leq 1}\nabla\mathfrak{u}\cdot\mathcal{F}dx
=\frac{1}{2}\int_{N^{-1}\leq |x|\leq 1}\mathfrak{u}\nabla\cdot\mathcal{F}dx,
\end{equation*}
This proves the third line of \eqref{eq:reduced_concentration}, and hence \eqref{eq:deterministic_consistency_statement}. Here we perform the integration by parts in order to gain the factor $\frac{1}{\log N}$ from $\nabla\mathfrak{u}_x$, which derives the smallness $C_*o_N(1)$.
\end{proof}

\begin{lemma}\label{lem:centralized_decomposition}
Let
$\omega_x=(\eta_x-\mathfrak{u}_x)/[\mathfrak{u}_x(1-\mathfrak{u}_x)]$.  For any two sites $x_1,x_2$
\begin{equation}\label{eq:centralization}
\begin{split}
&\eta_{x_1}(1-\eta_{x_2})+\eta_{x_2}(1-\eta_{x_1})
-\Big(\mathfrak{u}_{x_1}(1-\mathfrak{u}_{x_2})+\mathfrak{u}_{x_2}(1-\mathfrak{u}_{x_1})\Big)\\
&=\underbrace{-2\mathfrak{u}_{x_1}(1-\mathfrak{u}_{x_1})\mathfrak{u}_{x_2}(1-\mathfrak{u}_{x_2})
\omega_{x_1}\omega_{x_2}}_{\textup{\blue{quadratic term}}}+\underbrace{(1-2\mathfrak{u}_{x_2})\mathfrak{u}_{x_1}(1-\mathfrak{u}_{x_1})\omega_{x_1}
+(1-2\mathfrak{u}_{x_1})\mathfrak{u}_{x_2}(1-\mathfrak{u}_{x_2})\omega_{x_2}}_{\textup{\blue{linear term}}}.
\end{split}
\end{equation}
\end{lemma}
\begin{proof}
Write $\eta_x=\mathfrak{u}_x+\mathfrak{u}_x(1-\mathfrak{u}_x)\omega_x$.  Then
\begin{equation*}
\begin{split}
\eta_{x_1}(1-\eta_{x_2})
&=\mathfrak{u}_{x_1}(1-\mathfrak{u}_{x_2})
+(1-\mathfrak{u}_{x_2})\mathfrak{u}_{x_1}(1-\mathfrak{u}_{x_1})\omega_{x_1}\\
&\quad-\mathfrak{u}_{x_1}\mathfrak{u}_{x_2}(1-\mathfrak{u}_{x_2})\omega_{x_2}
-\mathfrak{u}_{x_1}(1-\mathfrak{u}_{x_1})\mathfrak{u}_{x_2}(1-\mathfrak{u}_{x_2})
\omega_{x_1}\omega_{x_2}.
\end{split}
\end{equation*}
The expansion of $\eta_{x_2}(1-\eta_{x_1})$ is obtained by exchanging $x_1$ and $x_2$.  Adding the two identities and subtracting the deterministic density part gives \eqref{eq:centralization}.
\end{proof}

\subsection{Concentration for the Quadratic Term}\label{sec:concentration_quadratic}
In this subsection, we prove the concentration estimate \eqref{eq:quadratic_concentration_within} for the quadratic term.

\begin{lemma}\label{lem:concentration_quadratic}
Under the notation of Proposition \ref{prop:concentration}, let $\{H_x^N\}_{N<+\infty}$ be any sequence of deterministic functions satisfying the following inequality with a positive constant $\mathfrak{b}>0$ uniformly in $N$
\begin{equation}\label{eq:sec8_H_upper_bound}
\big|H_x^N\big|\leq \mathfrak{b}\frac{1}{|x|^2\log N}\textup{ for any } x\in \mathds{B}^N.
\end{equation}
Then for every $\delta>0$ and $e\in\{e_1,e_2\}$
\begin{equation}\label{eq:quadratic_term_concentration}
\mathds{P}_{\mathfrak{u},\mathds{F}}^N\left[
\left|\int_0^T\sum_{x\in\mathds{B}^N}\frac{1}{N^2}H_x^N\omega_{t,x}\omega_{t,x+e}dt\right|>\delta\right]=C_*o_N(1).
\end{equation}
\end{lemma}

\begin{proof}
It is enough to prove the upper tail, and the lower tail is the same argument with $-H^N$. The Jara-Menezes averaging method in \cite{Jara_Menezes_2018} replaces the local product $\omega_x\omega_{x+e}$ by a product of two mesoscopic averages of $\omega$, with a quantitative estimate of the replacement error. Here we modify the method into a multi-scale manner due to the singularity of $H_x^N$ in \eqref{eq:sec8_H_upper_bound}, using the averaging radius $\ell$ introduced in \eqref{eq:sec8_averaging_radius}. Recall that the parameter $\alpha>0$ determines the singularity of the density profile $\mathfrak{u}$ in \eqref{eq:u_profile}. We need to choose the parameter $\beta\in (0,1)$ in the definition of $\ell$ such that
\begin{equation}\label{eq:beta_alpha_relation}
4\alpha<\beta<1-4\alpha.
\end{equation}
The intuition is that the smaller $\beta$ is, the larger we average, see \eqref{eq:sec8_averaging_radius}. Since $\alpha$ measures the singularity of the density profile $\mathfrak{u}$ in \eqref{eq:u_profile}, it also gives a condition on how large we can average using $\ell$. Too small $\beta$ brings non-neligible error, and too large $\beta$ does not average enough the configuration $\omega_x$ to establish concentration estimates.

\textit{Step 1: layer decomposition.}

We split the domain $\mathds{B}^N$ into an inner layer and an outer layer
\begin{equation*}
\mathds{B}^N=\underbrace{\Big\{x\in \mathds{B}^N\Big|\ |x|\leq N^{-\beta}\Big\}}_{\textup{\blue{inner layer}}}\cup\underbrace{\Big\{x\in\mathds{B}^N\Big|\ N^{-\beta}<|x|\leq 1\Big\}}_{\textup{\blue{outer layer}}}.
\end{equation*}
It is sufficient to prove the inner layer estimate
\begin{equation}\label{eq:inner_layer}
\mathds{P}_{\mathfrak{u},\mathds{F}}^N\left[\int_0^T\sum_{x\in\mathds{B}^N:|x|\leq N^{-\beta}}
\frac{1}{N^2}H_x^N\omega_{t,x}\omega_{t,x+e}dt>\delta\right]=C_*o_N(1),
\end{equation}
and the outer layer estimate
\begin{equation}\label{eq:outer-layer}
\mathds{P}_{\mathfrak{u},\mathds{F}}^N\left[
\int_0^T\sum_{x\in\mathds{B}^N:N^{-\beta}<|x|\leq 1}
\frac{1}{N^2}H_x^N\omega_{t,x}\omega_{t,x+e}dt>\delta\right]=C_*o_N(1).
\end{equation}
The inner layer contains the microscopic singularity and is the main part of the proof.  The outer layer is estimated at the end with the subgaussian bounds from Lemmas \ref{lem:subgaussian_omega} and \ref{lem:subgaussian_H}.

For the proof of the inner layer estimate, define the scaled coefficient $\mathfrak{h}_x^N$ and the inner layer quadratic observable $\mathcal{V}_{\textup{qua}}$ by
\begin{equation}\label{eq:scaled_H}
\mathfrak{h}_x^N:=\frac{N^{4\alpha}\log N}{N^2}H_x^N\mathds{1}_{|x|\leq N^{-\beta}},
\qquad
\mathcal{V}_{\textup{qua}}:=\sum_{x\in\mathds{B}^N:|x|\leq N^{-\beta}}
\mathfrak{h}_x^N\omega_x\omega_{x+e}.
\end{equation}
Then \eqref{eq:inner_layer} is equivalent to
\begin{equation}\label{eq:inner_goal}
\mathds{P}_{\mathfrak{u},\mathds{F}}^N\left[
\int_0^T\mathcal{V}_{\textup{qua}}(\eta_t)dt>\delta N^{4\alpha}\log N\right]=C_*o_N(1),
\end{equation}
and the upper bound \eqref{eq:sec8_H_upper_bound} of $|H_x^N|$ becomes
\begin{equation}\label{eq:asymp_H}
\big|\mathfrak{h}_x^N\big|\leq \mathfrak{b} N^{4\alpha-2}|x|^{-2}\mathds{1}_{|x|\leq N^{-\beta}}.
\end{equation}

\textit{Step 2: replacement by averaged configurations for $\mathcal{V}_{\textup{qua}}$.}

The averaging estimate \cite[Eq. (3.8)]{Jara_Menezes_2018} shows there exists a constant $C_*^{(1)}>0$ such that, for every density $f$ with respect to $d\local{\mathfrak{u}}$ and any $\eps>0$, we can turn $\mathcal{V}_{\textup{qua}}$ into averaged quantities with the cost of introducing a positive Dirichlet form. Define $\avedisk$ as the following cube, where we will average several quantities
\begin{equation}\label{eq:average_disk}
\avedisk:=\Big\{y\in N^{-1}\mathds{Z}_{\geq 0}^2\Big|\ |y|_{\infty}\leq N^{-1}(\ell-1)\Big\}.
\end{equation}
The term $\local{\mathfrak{u}}\big[\mathcal{V}_{\textup{qua}}f\big]$ is bounded from above by the following
\begin{equation*}
\local{\mathfrak{u}}\big[\mathcal{V}_{\textup{qua}}f\big]
\leq \local{\mathfrak{u}}\big[\mathcal{V}^{\ell}f\big]
+\frac{\eps }{2}\local{\mathfrak{u}}\Big[\Gamma_{\SS}^N(\sqrt{f})\Big]
+\eps^{-1}C_*^{(1)}\local{\mathfrak{u}}\big[\mathcal{W}^{\ell}f\big]
+\local{\mathfrak{u}}\big[\mathcal{Z}^{\ell}f\big].
\end{equation*}
The three terms produced by the replacement are
\begin{subequations}\label{eq:decomposition_components}
\begin{align}
\mathcal{V}^{\ell}
&:=\sum_{x\in\mathds{B}^N:|x|\leq CN^{-\beta}}\lomega_x^{\ell}\romega_{x+e}^{\ell}\qquad\qquad\qquad\quad \textup{ with } \lomega_x^{\ell}:=\sum_{y\in \avedisk}
\omega_{x-y}\mathfrak{h}_{x-y}^N\ell^{-2}\textup{ and } \romega_{x+e}^{\ell}:=\sum_{z\in \avedisk}
\omega_{x+e+z}\ell^{-2},\label{eq:decomposition_components_a}
\\
\mathcal{W}^{\ell}
&:=\frac{1}{N^2}\sum_{e'\in\{e_1,e_2\}}\sum_{x\in\mathds{B}^N:|x|\leq CN^{-\beta}}
\Big(\mathcal{H}_{x,e'}^N\Big)^2\qquad \textup{with }\mathcal{H}_{x,e'}^N:=
\sum_{y\in\avediskplus}
\phi_{\ell}(y;e')\omega_{x-y}\mathfrak{h}_{x-y}^N,\label{eq:decomposition_components_b}\\
\mathcal{Z}^{\ell}
&:=\sum_{e'\in\{e_1,e_2\}}\sum_{x\in\mathds{B}^N:|x|\leq CN^{-\beta}}
(\mathfrak{u}_{x+e'}-\mathfrak{u}_x)\mathcal{H}_{x-e,e'}^N\omega_x\omega_{x+e'}. \label{eq:decomposition_components_c}
\end{align}
\end{subequations}
Notice that the edge $e'$ in the definitions of $\mathcal{W}^{\ell}$ and $\mathcal{Z}^{\ell}$ is different from the edge $e$ in the statement of Lemma \ref{lem:concentration_quadratic}. The edge $e$ refers to the quadratic term $\omega_x\omega_{x+e}$ we want to average, and the edge $e'$ refers to the part of transport error we want to estimate. In this paper, since the density profile $\mathfrak{u}$ is taken as $\mathfrak{u}^N[\gamma,\zeta;\alpha]$ and $\beta>\alpha$ according to \eqref{eq:beta_alpha_relation}, we have $\mathfrak{u}_{x}\equiv \zeta$ for $|x|\leq N^{-\beta}$
and therefore $\mathcal{Z}^{\ell}\equiv 0$
\begin{equation}\label{eq:JM_replacement_quadratic}
\local{\mathfrak{u}}\big[\mathcal{V}_{\textup{qua}}f\big]
\leq \local{\mathfrak{u}}\big[\mathcal{V}^{\ell}f\big]
+\frac{\eps}{2}\local{\mathfrak{u}}\big[\Gamma_{\SS}^N(\sqrt{f})\big]
+\eps^{-1}C_*^{(1)}\local{\mathfrak{u}}\big[\mathcal{W}^{\ell}f\big].
\end{equation}

\textit{Step 3: absorb the Dirichlet form via Feynman-Kac formula.}

By Lemma \ref{lem:feynman_kac_variational_app} there is a constant $C_*^{(2)}>0$ such that for every 
observable $G(\eta_t)$
\begin{equation}\label{eq:FK_variational_quadratic}
\log\mathds{E}_{\mathfrak{u},\mathds{F}}^N
\Bigg[\exp\Bigg(\int_0^T G(\eta_t)dt\Bigg)\Bigg]
\leq
T\sup_f
\left\{
\local{\mathfrak{u}}[Gf]
-C_*^{(2)}\local{\mathfrak{u}}[\Gamma_{\SS}^N(\sqrt f)]
+\frac{1}{2}\local{\mathfrak{u}}[fL_{\mathds{F}}^{N,*}\mathds{1}]
\right\},
\end{equation}
where the supremum is over densities $f$ with respect to
$\local{\mathfrak{u}}$. Take $\mathcal{U}$ as a function that will be defined later in \eqref{eq:U_function_sec8}. To prove \eqref{eq:inner_goal}, the Markov inequality shows that it suffices to prove
\begin{subequations}\label{eq:inner_decomposed_goal}
\begin{align}
&\log\mathds{E}_{\mathfrak{u},\mathds{F}}^N\Bigg[\exp\Big(\int_0^T\Big(\mathcal{V}_{\textup{qua}}(\eta_t)-\mathcal{U}(\eta_t)\Big)dt\Big)\Bigg]-\delta N^{4\alpha}\log N\rightarrow -\infty,\quad \textup{as $N\rightarrow +\infty$}, \label{eq:inner_decomposed_goal_a}\\
&\mathds{E}_{\mathfrak{u},\mathds{F}}^N\Bigg[\Big|\int_0^T\mathcal{U}(\eta_t)dt\Big|\Bigg]\Big(\delta N^{4\alpha}\log N\Big)^{-1}=C_*o_N(1). \label{eq:inner_decomposed_goal_b}
\end{align}
\end{subequations}
For \eqref{eq:inner_decomposed_goal_a}, we apply Feynman-Kac formula \eqref{eq:FK_variational_quadratic} with $G=\mathcal{V}_{\textup{qua}}-\mathcal{U}$, and we want to show the LHS of \eqref{eq:applied_FK_sec8} is negative
\begin{equation}\label{eq:applied_FK_sec8}
\log \mathds{E}_{\mathfrak{u},\mathds{F}}^N\Bigg[\exp\Bigg(\int_0^T\Big(\mathcal{V}_{\textup{qua}}(\eta_t)-\mathcal{U}(\eta_t)\Big)dt\Bigg)\Bigg]\leq T\sup_{f}\Bigg\{\local{\mathfrak{u}}\big[\mathcal{V}_{\textup{qua}}f\big]-C_*^{(2)}\local{\mathfrak{u}}\big[\Gamma_{\SS}^N(\sqrt{f})\big]+\frac{1}{2}\local{\mathfrak{u}}\big[fL_{\mathds{F}}^{N,*}\mathds{1}\big]-\local{\mathfrak{u}}\big[\mathcal{U}f\big]\Bigg\}.
\end{equation}
We apply the averaging formula \eqref{eq:JM_replacement_quadratic} to the term $\local{\mathfrak{u}}[\mathcal{V}_{\textup{qua}}f]$, where we get averaged terms having better estimates
\begin{equation}\label{eq:JM_applied_to_Vqua_sec8}
\local{\mathfrak{u}}\big[\mathcal{V}_{\textup{qua}}f\big]\leq \local{\mathfrak{u}}\big[\mathcal{V}^{\ell}f\big]
+\frac{\eps}{2}\local{\mathfrak{u}}\big[\Gamma_{\SS}^N(\sqrt{f})\big]
+\eps^{-1}C_*^{(1)}\local{\mathfrak{u}}\big[\mathcal{W}^{\ell}f\big].
\end{equation}
For the adjoint operator, in Section \ref{sec:relative_entropy} we have proved \eqref{eq:sec6.4remark1}
\begin{equation}\label{eq:JM_applied_to_adjoint_sec8}
\frac{1}{2}\local{\mathfrak{u}}\Big[f\mathds{L}_{\mathds{F}}^{N,*}\mathds{1}\Big]\leq \frac{\eps}{2} \local{\mathfrak{u}}\Big[\Gamma_{\SS}^N(\sqrt{f})\Big]+\frac{1}{4}\local{\mathfrak{u}}\Big[f\mathds{A}_{\textup{qua}}^{\dagger}\Big]+\frac{1}{2}\frac{C_*}{\eps}\local{\mathfrak{u}}\Big[f\mathds{A}_{\textup{qua}}^W\Big]+\frac{1}{4}\local{\mathfrak{u}}\Big[f\mathds{A}_{\textup{qua}}^Z\Big],
\end{equation}
where $\mathds{A}_{\textup{qua}}^{\dagger}$ is defined in \eqref{eq:after_transport}, and $\mathds{A}_{\textup{qua}}^{W}$ as well as $\mathds{A}_{\textup{qua}}^Z$ are defined in \eqref{eq:AZAW}. The Dirichlet form in \eqref{eq:JM_applied_to_Vqua_sec8} and \eqref{eq:JM_applied_to_adjoint_sec8} can be absorbed by the negative Dirichlet form in \eqref{eq:applied_FK_sec8} by choosing $\eps$ small enough. We construct $\mathcal{U}$ as
\begin{equation}\label{eq:U_function_sec8}
\mathcal{U}=\mathcal{V}^{\ell}+\eps^{-1}C_*^{(1)}\mathcal{W}^{\ell}+\frac{1}{4}\mathds{A}_{\textup{qua}}^{\dagger}+\frac{1}{2}\eps^{-1}C_*\mathds{A}_{\textup{qua}}^W+\frac{1}{4}\mathds{A}_{\textup{qua}}^Z.
\end{equation}
This choice of $\mathcal{U}$ makes the RHS of \eqref{eq:applied_FK_sec8} negative, and proves \eqref{eq:inner_decomposed_goal_a}. To prove the estimate of $\mathcal{U}$ in \eqref{eq:inner_decomposed_goal_b}, we need to show
\begin{equation*}
\mathds{E}_{\mathfrak{u},\mathds{F}}^N\Big[\int_0^T|\mathcal{V}^{\ell}(\eta_t)|dt\Big]=C_*O(N^{4\alpha}),\quad \mathds{E}_{\mathfrak{u},\mathds{F}}^N\Big[\int_0^T|\mathcal{W}^{\ell}(\eta_t)|dt\Big]=C_*O(N^{4\alpha}),
\end{equation*}
which will be proved by Lemma \ref{lem:V} and Lemma \ref{lem:W}, by estimating the upper bound of $\mathcal{V}^{\ell}(\eta_t)$ and $\mathcal{W}^{\ell}(\eta_t)$ for each $t\in[0,T]$. For terms $(\mathds{A}_{\textup{qua}}^{\dagger},\mathds{A}_{\textup{qua}}^W,\mathds{A}_{\textup{qua}}^Z)$, we have the same type of estimate as in \eqref{eq:sec6.4remark2}. This completes \eqref{eq:inner_decomposed_goal_b} and the proof of Lemma \ref{lem:concentration_quadratic}.
\end{proof}

\begin{lemma}\label{lem:V}
Under the notation of Proposition \ref{prop:concentration}, the term $\mathcal{V}^{\ell}$ in \eqref{eq:decomposition_components} satisfies
\begin{equation*}
\mathds{E}_{\mathfrak{u},\mathds{F},t}^N\Big[|\mathcal{V}^{\ell}|\Big]=C_*O(N^{4\alpha}),\qquad 0\leq t\leq T.
\end{equation*}
\end{lemma}

\begin{proof}
Recall $\mathcal{V}^{\ell}=\sum_{x\in\mathds{B}^N:|x|\leq CN^{-\beta}}\lomega_x^{\ell}\romega_{x+e}^{\ell}$, where $\ell$ is the averaging radius defined in \eqref{eq:sec8_averaging_radius}. The left average $\lomega_x^{\ell}$ contains the singular coefficient $\mathfrak{h}^N$; the right average $\romega_{x+e}^{\ell}$ does not.  Lemma \ref{lem:subgaussian_omega} shows that the subgaussian order of $\lomega_x^{\ell}$ is $N^{8\alpha-4(1-\beta)}$, and the subgaussian order of $\romega_{x+e}^{\ell}$ is $N^{-2(1-\beta)}$. Applying the entropic inequality with correlations in Lemma \ref{lem:entropic_ineq_correlation} to $\mathcal{V}^{\ell}$
\begin{equation}\label{eq:vli_entropic_inequal}
\begin{split}
\mathds{E}_{\mathfrak{u},\mathds{F},t}^N \big[|\mathcal{V}^{\ell}|\big]
&\leq
C\mathds{H}^N\Big(d\mathds{P}_{\mathfrak{u},\mathds{F},t}^N\Big|d\local{\mathfrak{u}}\Big)+\log 2+\max_{k\in\{+1,-1\}}
\frac{C}{\ell^2}\sum_{x\in\mathds{B}^N:|x|\leq C N^{-\beta}}
\log \local{\mathfrak{u}}
\left[\exp\Big(kC \ell^2\lomega_x^{\ell}\romega_{x+e}^{\ell}\Big)\right].
\end{split}
\end{equation}
Next we use Lemma \ref{lem:quadratic_subgaussian} to estimate the $\log\local{\mathfrak{u}}$ term, and we need to verify the condition $\theta\sigma_1\sigma_2\leq 1/4$ of Lemma \ref{lem:quadratic_subgaussian}. Recall that $\mathcal{SG}$ is the subgaussian order defined in \eqref{eq:def_subgaussian}. Take $\theta$ as $\ell^2$, $\sigma_1$ as $\sqrt{\mathcal{SG}(\lomega_x^{\ell})}$, and $\sigma_2$ as $\sqrt{\mathcal{SG}(\romega_{x+e}^{\ell})}$, and use $\beta<1-4\alpha$ in \eqref{eq:beta_alpha_relation}
\begin{equation}\label{eq:pointwise_SG_Vquadra}
\ell^2\sqrt{\mathcal{SG}(\lomega_x^{\ell})\mathcal{SG}(\romega_{x+e}^{\ell})}
\leq C_* N^{2(1-\beta)}N^{4\alpha-3(1-\beta)}\log N
=C_*N^{4\alpha-(1-\beta)}\log N\leq \frac{1}{4}\quad\textup{when $N$ is large enough}.
\end{equation}
Therefore Lemma \ref{lem:quadratic_subgaussian} and Equation \eqref{eq:pointwise_SG_Vquadra} imply as $\ell=0.5N^{1-\beta}$
\begin{equation*}
\mathds{E}_{\mathfrak{u},\mathds{F},t}^N\big[|\mathcal{V}^{\ell}|\big]\leq
\mathds{H}^N\Big(d\mathds{P}_{\mathfrak{u},\mathds{F},t}^N\Big|d\local{\mathfrak{u}}\Big)
+\log 2+C\ell^{-2}\sum_{x\in\mathds{B}^N:|x|\leq CN^{-\beta}}
\underbrace{\ell^{2}\sqrt{\mathcal{SG}(\lomega_x^{\ell})\mathcal{SG}(\romega_{x+e}^{\ell})}}_{\textup{\blue{smaller than $1$}}}\leq C_*O(N^{4\alpha}),
\end{equation*}
where we used Proposition \ref{prop:relative_entropy} in the last inequality to control the entropy.
\end{proof}

\begin{lemma}\label{lem:W}
Under the notation of Proposition \ref{prop:concentration}, the term $\mathcal{W}^{\ell}$ in \eqref{eq:decomposition_components} satisfies
\begin{equation*}
\mathds{E}_{\mathfrak{u},\mathds{F},t}^N\Big[|\mathcal{W}^{\ell}|\Big]=C_*O(N^{4\alpha}),\qquad 0\leq t\leq T.
\end{equation*}
\end{lemma}
\begin{proof}
Recall
\begin{equation*}
\mathcal{W}^{\ell}=\frac{1}{N^2}\sum_{e'\in\{e_1,e_2\}}
\sum_{x\in\mathds{B}^N:|x|\leq CN^{-\beta}}\Big(\mathcal{H}_{x,e'}^N\Big)^2.
\end{equation*}
By Lemma \ref{lem:subgaussian_H}, in the inner layer the subgaussian order of $\mathcal{H}_{x,e'}^N$ is $N^{8\alpha}\log N$. Apply the entropic inequality with correlations \ref{lem:entropic_ineq_correlation} to $\mathcal{W}^{\ell}$
\begin{equation*}
\mathds{E}_{\mathfrak{u},\mathds{F},t}^N[ |\mathcal{W}^{\ell}|]\leq C\mathds{H}^N\Big[d\mathds{P}_{\mathfrak{u},\mathds{F},t}^N\Big|\ d\local{\mathfrak{u}}\Big]+\log 2+\max_{k\in\{+1,-1\}}\frac{C}{\ell^2}\sum_{e'\in\{e_1,e_2\}}\sum_{x\in\mathds{B}^N:|x|\leq CN^{-\beta}}\log\local{\mathfrak{u}}\Big[\exp\Big(k\ell^2N^{-2}\big(\mathcal{H}_{x,e'}^N\big)^2\Big)\Big].
\end{equation*}
Verify the condition $\theta\sigma_1\sigma_2\leq 1/4$ of the subgaussian estimate Lemma \ref{lem:quadratic_subgaussian}, with $\theta=\ell^2N^{-2}$ and $\sigma_1=\sigma_2=\sqrt{\mathcal{SG}(\mathcal{H}_{x,e'}^N)}$
\begin{equation*}
\ell^2N^{-2}\mathcal{SG}\big(\mathcal{H}_{x,e'}^N\big)\leq C^* N^{2(1-\beta)}N^{-2}N^{8\alpha}\log N=C_*N^{8\alpha-2\beta}\log N\leq \frac{1}{4},
\end{equation*}
which follows from $\beta>4\alpha$ in \eqref{eq:beta_alpha_relation}. Now we can apply Lemma \ref{lem:quadratic_subgaussian} 
\begin{equation*}
\begin{split}
\mathds{E}_{\mathfrak{u},\mathds{F},t}^N[\mathcal{W}^{\ell}]&\leq C\mathds{H}^N\Big[d\mathds{P}_{\mathfrak{u},\mathds{F},t}^N\Big|\ d\local{\mathfrak{u}}\Big]+\log 2+\frac{C}{\ell^2}\sum_{e'\in\{e_1,e_2\}}\sum_{x\in\mathds{B}^N:|x|\leq CN^{-\beta}}\underbrace{\ell^2N^{-2}\mathcal{SG}\big(\mathcal{H}_{x,e'}^N\big)}_{\textup{\blue{smaller than $1$}}}\leq C_*O(N^{4\alpha}),
\end{split}
\end{equation*}
where we have used Proposition \ref{prop:relative_entropy} to control the entropy.
\end{proof}

\textit{Step 4: the outer layer estimate \eqref{eq:outer-layer}.}

For the proof of the outer later estimate, denote $\mathfrak{h}_x^N$ as
\begin{equation}\label{eq:asymp_H_linear}
\mathfrak{h}_{x}^N:=H_x^N\frac{N^{4\alpha}\log N}{N^2}\mathds{1}_{|x|>N^{-\beta}},\quad \big|\mathfrak{h}_x^N\big|\leq \mathfrak{b}N^{4\alpha-2}|x|^{-2}\mathds{1}_{|x|>N^{-\beta}}.
\end{equation}
It is sufficient to prove
\begin{equation*}
\mathds{P}_{\mathfrak{u},\mathds{F}}^N\Bigg[\ \int_0^T\sum_{x\in\mathds{B}^N:|x|> N^{-\beta}}\mathfrak{h}_x^N\omega_{t,x}\omega_{t,x+e}dt>\delta N^{4\alpha}\log N\Bigg]=C_*o_N(1).
\end{equation*}
Using Markov inequality, it is further sufficient to prove the outer layer estimate
\begin{equation}\label{eq:lin_concentration_markov}
\frac{1}{\delta N^{4\alpha}\log N}\mathds{E}_{\mathfrak{u},\mathds{F}}^N\Bigg[\int_0^T\ \Big|\sum_{x\in\mathds{B}^N:|x|> N^{-\beta}}\mathfrak{h}_x^N\omega_{t,x}\omega_{t,x+e}\Big|dt\Bigg]=C_*o_N(1).
\end{equation}
The entropic inequality in Lemma \ref{lem:entropic_ineq_correlation} implies for any $t\in [0,T]$
\begin{equation*}
\begin{split}
&\mathds{E}_{\mathfrak{u},\mathds{F}}^N\Bigg[\ \Big|\sum_{x\in\mathds{B}^N:|x|> N^{-\beta}}\mathfrak{h}_x^N\omega_{t,x}\omega_{t,x+e}\Big|\Bigg]\\
&\leq C_* \mathds{H}^N\Big(d\mathds{P}_{\mathfrak{u},\mathds{F},t}^N\Big|\ d\local{\mathfrak{u}}\Big)+\log 2+\max_{k\in\{+1,-1\}}C_*\sum_{x\in\mathds{B}^N:|x|>N^{-\beta}}\log\local{\mathfrak{u}}\Big[\exp\Big(kC_*\mathfrak{h}_x^N\omega_x\omega_{x+e}\Big)\Big]\\
&\leq C_* \mathds{H}^N\Big(d\mathds{P}_{\mathfrak{u},\mathds{F},t}^N\Big|\ d\local{\mathfrak{u}}\Big)+\log 2+C_*\sum_{x\in\mathds{B}^N:|x|> N^{-\beta}}|\mathfrak{h}_x^N|^2.
\end{split}
\end{equation*}
Recall the upper bound for $|\mathfrak{h}_x^N|$ in \eqref{eq:asymp_H_linear}, and use Lemma \ref{lem:lattice_sum} to calculate the sum of $|x|^{-4}$
\begin{equation*}
C_*\sum_{x\in\mathds{B}^N:|x|> N^{-\beta}}|\mathfrak{h}_x^N|^2\leq C_*\sum_{x\in\mathds{B}^N:|x|> N^{-\beta}}\frac{N^{8\alpha}}{N^4}|x|^{-4}\leq C_*\frac{N^{8\alpha}}{N^4}N^{2+2\beta}\leq C_*\frac{N^{8\alpha+2\beta}}{N^2}.
\end{equation*}
The fact that $\beta<1-4\alpha$ in \eqref{eq:beta_alpha_relation} implies $N^{8\alpha+2\beta-2}=o(N^{4\alpha})$. Proposition \ref{prop:relative_entropy} bounds the relative entropy $\mathds{H}^N$ from above by $C_*N^{4\alpha}$. This time-wise upper bound concludes the proof of \eqref{eq:lin_concentration_markov} by an integral in time

\subsection{Concentration for the Linear Term}\label{sec:concentration_linear}

\begin{lemma}\label{lem:concentration_linear}
Under the notation of Proposition \ref{prop:concentration}, let $\{H_x^N\}_{N<+\infty}$ be a sequence of  deterministic functions satisfying the bound with parameter $\mathfrak{b}>0$
\begin{equation}\label{eq:linear_H_bound}
\big|H_x^N\big|\leq \mathfrak{b}|x|^{-2}\quad \textup{ for any } x\in \mathds{B}^N \textup{ and $N<+\infty$}.
\end{equation}
Then for every $\delta>0$
\begin{equation}\label{eq:concentration_linear}
\mathds{P}_{\mathfrak{u},\mathds{F}}^N\left[
\left|\int_0^T \sum_{x\in\mathds{B}^N}\frac{1}{N^2}H_x^N\omega_{t,x}dt\right|>\delta\right]=C_*o_N(1).
\end{equation}
\end{lemma}
\begin{proof}
Again it suffices to prove the upper tail, since the lower tail follows by replacing $H^N$ with $-H^N$.  Use the scaled coefficient $\mathfrak{h}^N:=H_x^NN^{4\alpha-2}\log N$ from \eqref{eq:scaled_H}, with the following upper bound due to the condition \eqref{eq:linear_H_bound}
\begin{equation}\label{eq:h_bound_linear_sec8}
|\mathfrak{h}_{x}^N|\leq \mathfrak{b}N^{4\alpha-2}|x|^{-2}\log N.
\end{equation}
The upper-tail estimate of \eqref{eq:concentration_linear} is equivalent to
\begin{equation}\label{eq:linear_scaled_goal}
\mathds{P}_{\mathfrak{u},\mathds{F}}^N\left[
\int_0^T\sum_{x\in\mathds{B}^N}\mathfrak{h}_x^N\omega_{t,x}dt>
\delta N^{4\alpha}\log N\right]=C_*o_N(1).
\end{equation}
The proof of \eqref{eq:linear_scaled_goal} can be split into the inner layer and the outer layer, where we need to prove
\begin{subequations}\label{eq:linear_inner_outer}
\begin{align}
&\mathds{P}_{\mathfrak{u},\mathds{F}}^N\left[
\int_0^T\sum_{x\in\mathds{B}^N:|x|\leq N^{-\beta}}\mathfrak{h}_x^N\omega_{t,x}dt>
\delta N^{4\alpha}\log N\right]=C_*o_N(1), \label{eq:linear_inner_outer_a}\\
&\mathds{P}_{\mathfrak{u},\mathds{F}}^N\left[
\int_0^T\sum_{x\in\mathds{B}^N:N^{-\beta}<|x|\leq 1}\mathfrak{h}_x^N\omega_{t,x}dt>
\delta N^{4\alpha}\log N\right]=C_*o_N(1). \label{eq:linear_inner_outer_b}
\end{align}
\end{subequations}
To prove the outer layer part \eqref{eq:linear_inner_outer_b}, we directly apply the entropic inequality in Lemma \ref{lem:entropic_inequality_app} for any time $t\in[0,T]$, and use the fact that $\omega_x$ is a subgaussian variable in the sense of \eqref{eq:def_subgaussian} and that $\mathfrak{h}_x^N$ satisfies \eqref{eq:h_bound_linear_sec8}
\begin{equation}\label{eq:linear_entropic_estimate_naive}
\begin{split}
&\mathds{E}_{\mathfrak{u},\mathds{F}}^N\Big[\Big|\sum_{x\in\mathds{B}^N:N^{-\beta}<|x|\leq 1}\mathfrak{h}_x^N\omega_{t,x}\Big|\Big]\\
&\leq \mathds{H}^N\Big(d\mathds{P}_{\mathfrak{u},\mathds{F},t}^N\Big|\ d\local{\mathfrak{u}}\Big)+\log 2+\max_{k\in\{+1,-1\}}\sum_{x\in\mathds{B}^N:N^{-\beta}<|x|\leq 1}\log\local{\mathfrak{u}}\Big[\exp\Big(k\mathfrak{h}_x^N\omega_x\Big)\Big]\\
&\leq \mathds{H}^N\Big(d\mathds{P}_{\mathfrak{u},\mathds{F},t}^N\Big|\ d\local{\mathfrak{u}}\Big)+\log 2+C_*\sum_{x\in\mathds{B}^N:N^{-\beta}<|x|\leq 1}N^{8\alpha-4}|x|^{-4}\log^2N\leq C_*O(N^{4\alpha}).
\end{split}
\end{equation}
Here the last inequality is due to the control of the relative entropy in Proposition \ref{prop:relative_entropy} and the condition of $\beta<1-4\alpha$ in \eqref{eq:beta_alpha_relation}
\begin{equation*}
\mathds{H}^N\Big(d\mathds{P}_{\mathfrak{u},\mathds{F},t}^N\Big|\ d\local{\mathfrak{u}}\Big)\leq C_*O(N^{4\alpha}),\quad \sum_{x\in\mathds{B}^N:N^{-\beta}<|x|\leq 1}N^{8\alpha-4}|x|^{-4}\log^2N\leq C_*N^{8\alpha-4}N^{2+2\beta}\log^2N\leq C_*O(N^{4\alpha}).
\end{equation*}

The inner layer estimate \eqref{eq:linear_inner_outer_a} can not be derived using the naive entropic estimate as in \eqref{eq:linear_entropic_estimate_naive}. For example, by applying the entropic inequality in Lemma \ref{lem:entropic_inequality_app} and using Lemma \ref{lem:lattice_sum} to estimate the sum of $|x|^{-4}$, we have 
\begin{equation*}
\mathds{E}_{\mathfrak{u},\mathds{F}}^N\Bigg[\Big|\sum_{x\in\mathds{B}^N:|x|\leq N^{-\beta}}\mathfrak{h}_x^N\omega_{t,x}\Big| \Bigg]\leq \mathds{H}^N\Big(d\mathds{P}_{\mathfrak{u},\mathds{F},t}^N\Big|\ d\local{\mathfrak{u}}\Big)+C_*\sum_{x\in\mathds{B}^N:|x|\leq N^{-\beta}}N^{8\alpha-4}|x|^{-4}\log^2N=C_*O(N^{8\alpha}\log^2N),
\end{equation*}
which is much larger than the desired estimate $C_*O(N^{4\alpha})$. This failure of the previous estimate is due to the singularity of $\mathfrak{h}_x^N$ in the inner layer $|x|\leq N^{-\beta}$. 

Instead, to prove \eqref{eq:linear_inner_outer_a}, we will perform the averaging method as in proof of Lemma \ref{lem:concentration_quadratic}. First we add an indicator function and change the position index to the definition of $\mathfrak{h}_x^N$ in the proof for the inner layer estimate \eqref{eq:linear_inner_outer_a}, which will be convenient later
\begin{equation}
\mathfrak{h}_x^N:=H_{x+e_1}^NN^{4\alpha-2}\log N\mathds{1}_{|x|\leq N^{-\beta}}.
\end{equation}
It is sufficient to estimate $\sum_{x\in\mathds{B}^N:|x|\leq N^{-\beta}}\mathfrak{h}_{x}^N\omega_{x+e_1}$. This expression is very similar to $\sum_{x\in\mathds{B}^N:|x|\leq N^{-\beta}}\mathfrak{h}_{x}^N\omega_x\omega_{x+e_1}$ we estimated in the proof of Lemma \ref{lem:concentration_quadratic}, and will be the same if we formally replace $\omega_x$ by $1$. By applying Feynman-Kac formula in Lemma \ref{lem:feynman_kac_variational_app} as we have done in the \textit{step 3} of the proof of Lemma \ref{lem:concentration_quadratic}, it is sufficient to estimate the following terms corresponding to \eqref{eq:decomposition_components}
\begin{equation}\label{eq:linear_decomposition_components}
\begin{split}
&\mathcal{V}_{\lin}^\ell
:=\sum_{x\in\mathds{B}^N:|x|\leq CN^{-\beta}}\overleftarrow{\mathfrak{h}}_x^{\ell}\romega_{x+e_1}^{\ell},\qquad  \mathcal{W}_{\lin}^{\ell}:=\frac{1}{N^2}\sum_{e'\in\{e_1,e_2\}}\sum_{x\in\mathds{B}^N:|x|\leq CN^{-\beta}}
\big(\mathscr{H}_{x,e'}^N\big)^2\\
&\textup{ with } \overleftarrow{\mathfrak{h}}_x^{\ell}
:=\ell^{-2}\sum_{y\in\avedisk}
\mathfrak{h}_{x-y}^N,\qquad \mathscr{H}_{x,e'}^N
:=\sum_{y\in \avediskplus}
\phi_{\ell}(y;e')\mathfrak{h}_{x-y}^N
\end{split}
\end{equation}
The only difference with \eqref{eq:decomposition_components} is that, we do not have random term $\omega$ in the definitions of $\overleftarrow{\mathfrak{h}}_x^{\ell}$ and $\mathscr{H}_{x,e'}^N$, because we are dealing with linear terms. The $\mathcal{Z}^{\ell}$-term also does not exist here, because $\mathfrak{u}$ is constant in the inner layer $|x|\leq N^{-\beta}$. We can compare the different terms involved in \eqref{eq:decomposition_components} and \eqref{eq:linear_decomposition_components}
\begin{equation*}
\left\{
\begin{array}{lr}
\mathcal{V}^{\ell}:=\sum_{x\in\mathds{B}^N:|x|\leq N^{-\beta}}\Big(\sum_{y\in\avedisk}\ell^{-2}\omega_{x-y}\mathfrak{h}_{x-y}^N\Big)\Big(\sum_{y\in\avedisk}\ell^{-2}\omega_{x+e_1+y}\Big),    &  \\
\mathcal{V}_{\textup{lin}}^{\ell}:=\sum_{x\in\mathds{B}^N:|x|\leq N^{-\beta}}\Big(\sum_{y\in\avedisk}\ell^{-2}\mathfrak{h}_{x-y}^N\Big) \Big(\sum_{y\in\avedisk}\ell^{-2}\omega_{x+e_1+y}\Big),    & 
\end{array}
\right.
\end{equation*}
and
\begin{equation*}
\left\{
\begin{array}{lr}
\mathcal{H}_{x,e'}^N:=\sum_{y\in\avediskplus}\phi_{\ell}(y;e')\omega_{x-y}\mathfrak{h}_{x-y}^N,     & \textup{in the quadratic case}, \\
\mathscr{H}_{x,e'}^N:=\sum_{y\in\avediskplus}\phi_{\ell}(y;e')\mathfrak{h}_{x-y}^N,      &  \textup{in the linear case} .
\end{array}
\right.
\end{equation*}

The same Feynman--Kac argument used in \eqref{eq:JM_replacement_quadratic}--\eqref{eq:applied_FK_sec8} reduces \eqref{eq:linear_inner_outer_a} to proving, uniformly in $t\in[0,T]$
\begin{equation*}
\mathds{E}_{\mathfrak{u},\mathds{F},t}^N[|\mathcal{V}_{\lin}^\ell|]
+\mathds{E}_{\mathfrak{u},\mathds{F},t}^N[|\mathcal{W}_{\lin}^{\ell}|]\leq C_*N^{4\alpha}.
\end{equation*}
This will be proved by Lemmas \ref{lem:linear_V} and \ref{lem:linear_W}.
\end{proof}

\begin{lemma}\label{lem:linear_V}
Under the notation above
\begin{equation*}
\mathds{E}_{\mathfrak{u},\mathds{F},t}^N[|\mathcal{V}_{\lin}^\ell|]=C_*O(N^{4\alpha}),\qquad 0\leq t\leq T.
\end{equation*}
\end{lemma}
\begin{proof}
We first derive an upper bound of $|\overleftarrow{\mathfrak{h}}_{x}^{\ell}|$ when $|x|\leq N^{-\beta}$. In this case, the set $x-\mathbf{D}_{\ell}^N$ is always within the domain $|x|\leq CN^{-\beta}$ with a universal constant $C$. Recall that within the inner layer, $\ell=0.5N^{1-\beta}$
\begin{equation*}
|\overleftarrow{\mathfrak{h}}_x^{\ell}|\leq \ell^{-2}\sum_{y\in\mathbf{D}_{\ell}^N}|\mathfrak{h}_{x-y}^N|\leq C_*N^{-2(1-\beta)}\Big(\sum_{y\in\mathds{B}^N:|y|\leq CN^{-\beta}}N^{4\alpha-2}|y|^{-2}\log^2N\Big)=C_*O\Big(N^{4\alpha+2\beta-2}\log^3 N\Big).
\end{equation*}
By \eqref{eq:decomposition_components_a}, $\romega_x^{\ell}$ is an average over $\omega_x$. Therefore in $\mathcal{V}_{\textup{lin}}^{\ell}:=\sum_{x\in\mathds{B}^N:|x|\leq CN^{-\beta}}\overleftarrow{\mathfrak{h}}_x^{\ell}\romega_{x+e_1}^{\ell}$, this average over $\omega$ can be transferred into an average over $\overleftarrow{\mathfrak{h}}_x^{\ell}$ by change of variables
\begin{equation*}
\begin{split}
\sum_{x\in\mathds{B}^N:|x|\leq CN^{-\beta}}\overleftarrow{\mathfrak{h}}_x^{\ell}\romega_{x+e_1}^{\ell}&=\sum_{x\in\mathds{B}^N:|x|\leq CN^{-\beta}}\sum_{y\in\avedisk}\ell^{-2}\overleftarrow{\mathfrak{h}}_x^{\ell}\omega_{x+e_1+y}\\
&=\sum_{x\in\mathds{B}^N:|x|\leq CN^{-\beta}}\Bigg(\ell^{-2}\sum_{y\in\avedisk}\overleftarrow{\mathfrak{h}}_{x-e_1-y}^{\ell}\Bigg)\omega_x=\sum_{x\in\mathds{B}^N:|x|\leq C N^{-\beta}}\omega_x\Tilde{\mathfrak{h}}_x^{\ell}.
\end{split}
\end{equation*}
with function $\Tilde{\mathfrak{h}}_x^{\ell}$ as
\begin{equation}\label{eq:tilde_h_upper_bound}
\Tilde{\mathfrak{h}}_x^{\ell}:=\sum_{y\in\mathbf{D}_{\ell}^N}\ell^{-2}\overleftarrow{\mathfrak{h}}_{x-e_1-y}^{\ell},\quad |\Tilde{\mathfrak{h}}_x^{\ell}|=C_*O\Big(N^{4\alpha+2\beta-2}\log^3 N\Big).
\end{equation}
Now we apply the entropic inequality in Lemma \ref{lem:entropic_inequality_app}. Each $\log\local{\mathfrak{u}}\big[\exp\big(\Tilde{\mathfrak{h}}_x^{\ell}\omega_x\big)\big]$ is bounded from above by $C_*|\Tilde{\mathfrak{h}}_x^{\ell}|^2$ according to \eqref{eq:def_subgaussian}, which is $C_*N^{8\alpha-4(1-\beta)}\log^6N$ by \eqref{eq:tilde_h_upper_bound}. The summation over $|x|\leq CN^{-\beta}$ contributes a factor $N^{2(1-\beta)}$ 
\begin{equation*}
\begin{split}
\mathds{E}_{\mathfrak{u},\mathds{F},t}^N\Big[|\mathcal{V}_{\textup{lin}}^{\ell}|\Big]&\leq \mathds{H}^N\Big(d\mathds{P}_{\mathfrak{u},\mathds{F},t}^N\Big|\local{\mathfrak{u}}\Big)+\log 2+\max_{k\in\{+1,-1\}}\sum_{x\in\mathds{B}^N:|x|\leq CN^{-\beta}}\log\local{\mathfrak{u}}\Big[\exp\Big(k\tilde{\mathfrak{h}}_x^{\ell}\omega_x\Big)\Big]\\
&\leq C_*O(N^{4\alpha})+C_*N^{2(1-\beta)}N^{8\alpha-4(1-\beta)}\log^6N=C_*O(N^{4\alpha})+C_*N^{8\alpha-2(1-\beta)}\log^6N\leq C_*O(N^{4\alpha}),
\end{split}
\end{equation*}
where last inequality uses the condition of $\beta<1-4\alpha$ in \eqref{eq:beta_alpha_relation}. This concludes the proof.
\end{proof}

\begin{lemma}\label{lem:linear_W}
Under the notation above
\begin{equation*}
\mathds{E}_{\mathfrak{u},\mathds{F},t}^N[|\mathcal{W}_{\lin}^{\ell}|]=C_*O(N^{4\alpha}),\qquad 0\leq t\leq T.
\end{equation*}
\end{lemma}
\begin{proof}
Recall the definitions of $\mathcal{W}_{\textup{lin}}^{\ell}$ and $\mathscr{H}_{x,e'}^N$ in \eqref{eq:linear_decomposition_components}
\begin{equation*}
\mathcal{W}_{\lin}^{\ell}:=\frac{1}{N^2}\sum_{e'\in\{e_1,e_2\}}\sum_{x\in\mathds{B}^N:|x|\leq N^{-\beta}}
\big(\mathscr{H}_{x,e'}^N\big)^2,\quad \mathscr{H}_{x,e'}^N
:=\sum_{y\in \avediskplus}
\phi_{\ell}(y;e')\mathfrak{h}_{x-y}^N.
\end{equation*}
Since the term $\mathcal{W}_{\lin}^{\ell}$ is deterministic, the estimate will be performed in a deterministic way. For the inner layer $|x|\leq N^{-\beta}$, Cauchy--Schwarz and the
square-summability of the flow in Lemma
\ref{lem:two_dimensional_flow_lemma_self_contained} give
\begin{equation*}
\big(\mathscr{H}_{x,e'}^N\big)^2
\leq
\left(\sum_{y\in\mathbf{D}_{2\ell-1}^N}\phi_{\ell}(y;e')^2\right)
\left(\sum_{y\in\mathbf{D}_{2\ell-1}^N}(\mathfrak{h}_{x-y}^N)^2\right)
\logles N^{8\alpha}.
\end{equation*}
Here the last square sum is obtained directly from \eqref{eq:h_bound_linear_sec8}:
\begin{equation*}
\begin{split}
\sum_{y\in\mathbf{D}_{2\ell-1}^N}(\mathfrak{h}_{x-y}^N)^2
&\logles N^{8\alpha-4}\sum_{z\in\mathds{B}^N:|z|\leq CN^{-\beta}}|z|^{-4}\logles N^{8\alpha}.
\end{split}
\end{equation*}
There are $O(N^{2(1-\beta)})$ inner-layer sites, hence
\begin{equation*}
\frac{1}{N^2}\sum_{e'\in\{e_1,e_2\}}\sum_{|x|\leq N^{-\beta}}\big(\mathscr{H}_{x,e'}^N\big)^2
\logles N^{8\alpha-2\beta}.
\end{equation*}
Use the condition of $\beta>4\alpha$ in \eqref{eq:beta_alpha_relation}, we conclude the proof.
\end{proof}

\subsection{Concentration for the Current Term}\label{sec:concentration_current}

The current term $\langle \mathcal{Q}_T^N,\mathds{F}\rangle$ is the only term in \eqref{eq:RN_random_decomposition} not expressed directly through the $\eta$ variables. Now we use martingale methods to prove the concentration estimate of $\langle \mathcal{Q}_T^N,\mathds{F}\rangle-T\langle \mathfrak{q},\mathds{F}\rangle$ in Proposition \ref{prop:concentration}. Let
\begin{equation}\label{eq:RN_bond_set}
\mathbf{E}^N
:=
\big\{(x,e):x\in\mathds{B}^N,\ e\in\{e_1,e_2\},\ x+e\in\mathds{B}^N\big\}
\end{equation}
be the set of nearest-neighbor bonds.  For a bond
$(x,e)\in\mathbf{E}^N$, write
\begin{equation*}
\mathbf{N}_t^{x,e,+}:=\#\{\textup{jumps }x\to x+e\textup{ during }[0,t]\},
\qquad
\mathbf{N}_t^{x,e,-}:=\#\{\textup{jumps }x+e\to x\textup{ during }[0,t]\}.
\end{equation*}
Thus the empirical current \eqref{eq:empirical_current} can be written equivalently as
\begin{equation}\label{eq:RN_current_counting}
\mathcal{Q}_{t,x,e}^N
=N^{-1}\big(\mathbf{N}_t^{x,e,+}-\mathbf{N}_t^{x,e,-}\big),
\qquad
\langle\mathcal{Q}_t^N,\mathds{F}\rangle
=\sum_{(x,e)\in\mathbf{E}^N}\frac{1}{N^2}\mathcal{Q}_{t,x,e}^N\mathds{F}_{x,e}=\sum_{(x,e)\in\mathbf{E}^N}\frac{1}{N^3}\big(\mathbf{N}_{t}^{x,e,+}-\mathbf{N}_{t}^{x,e,-}\big)\mathds{F}_{x,e}.
\end{equation}
We also recall the discrete divergence operator $\nabla_{\textup{disc}}\cdot$ 
\begin{equation}\label{eq:RN_divergence_convention}
\nabla_{\textup{disc}}\cdot\mathds{F}(x)
:=N\sum_{e\in\{e_1,e_2\}}
\Big(\mathds{F}_{x,e}
-\mathds{F}_{x-e,e}\Big).
\end{equation}
According to the generators $L_{\SS}^{N}$ and $L_{\mathds{F}}^{N}$ defined in \eqref{eq:SSEP_generator} and \eqref{eq:modified_generator}, the SSEP rates and the modified dynamics' rates on the bond $(x,x+e)$ are
\begin{equation}\label{eq:RN_rates}
\begin{array}{lll}
r_{x,e}^{0,+}(\eta)=\frac{N^2}{2}\eta_x(1-\eta_{x+e}),
&\qquad&
r_{x,e}^{\mathds{F},+}(\eta)=r_{x,e}^{0,+}(\eta)\exp\!\Big(\dfrac{\mathds{F}_{x,e}}{N\log N}\Big),\\[0.9em]
r_{x,e}^{0,-}(\eta)=\frac{N^2}{2}\eta_{x+e}(1-\eta_x),
&&
r_{x,e}^{\mathds{F},-}(\eta)=r_{x,e}^{0,-}(\eta)\exp\!\Big(-\dfrac{\mathds{F}_{x,e}}{N\log N}\Big).
\end{array}
\end{equation}
The sign $+$ records the jump rate in the chosen orientation $x\to x+e$, and the sign $-$ records the jump rate in the opposite orientation $x+e\to x$.

\begin{lemma}\label{lem:current_concentration}
Under the notation of Proposition \ref{prop:concentration}, for every $\delta>0$,
\begin{equation}\label{eq:current_concentration_statement}
\mathds{P}_{\mathfrak{u},\mathds{F}}^N\left[
\left|\langle\mathcal{Q}_T^N,\mathds{F}\rangle-T\langle \mathfrak{q},\mathcal{F}\rangle\right|>\delta\right]=C_*o_N(1).
\end{equation}
\end{lemma}

\begin{proof}
The following processes are martingales for any $(x,e)\in\mathbf{E}^N$
\begin{equation*}
\mathbf{m}_{t}^{x,e,+}:=\frac{\mathds{F}_{x,e}}{N^3}\Big(\mathbf{N}_{t}^{x,e,+}-\int_0^tr_{x,e}^{\mathds{F},+}(\eta_s)ds\Big),\quad \mathbf{m}_{t}^{x,e,-}:=\frac{\mathds{F}_{x,e}}{N^3}\Big(\mathbf{N}_{t}^{x,e,-}-\int_0^tr_{x,e}^{\mathds{F},-}(\eta_s)ds\Big).
\end{equation*}
We can calculate the quadratic variations $\langle \mathbf{m}_{\cdot}^{x,e,+}\rangle$ and $\langle \mathbf{m}_{\cdot}^{x,e,-} \rangle$
\begin{equation}\label{eq:m_quadratic_variation}
\mathds{E}_{\mathfrak{u},\mathds{F}}^N\Big[\langle \mathbf{m}_{\cdot}^{x,e,+}\rangle_t\Big]=C_*O\Big(\frac{\mathds{F}_{x,e}^2}{N^4}\Big),\quad \mathds{E}_{\mathfrak{u},\mathds{F}}^N\Big[\langle \mathbf{m}_{\cdot}^{x,e,-}\rangle_t\Big]=C_*O\Big(\frac{\mathds{F}_{x,e}^2}{N^4}\Big).
\end{equation}
Equation \eqref{eq:RN_current_counting} implies we can define $\mathcal{M}_t^N$ as a martingale
\begin{equation*}
\mathcal{M}_{t}^N:=\sum_{(x,e)\in\mathbf{E}^N}\Big(\mathbf{m}_t^{x,e,+}-\mathbf{m}_t^{x,e,-}\Big)=\langle \mathcal{Q}_t^N,\mathds{F}\rangle-\frac{1}{N^3}\sum_{(x,e)\in\mathbf{E}^N}\int_0^t\mathds{F}_{x,e}\Big(r_{x,e}^{\mathds{F},+}(\eta_s)-r_{x,e}^{\mathds{F},-}(\eta_s)\Big)ds.
\end{equation*}
The quadratic variation of $\mathcal{M}_t^N$ can be calculated by \eqref{eq:m_quadratic_variation} as
\begin{equation*}
\mathds{E}_{\mathfrak{u},\mathds{F}}^N\Big[\langle \mathcal{M}_{\cdot}^N\rangle_t\Big]=C_*O\Big(\sum_{x\in\mathds{B}^N}\frac{1}{N^4}|x|^{-2}\Big)=C_*O\Big(\frac{\log N}{N^2}\Big).
\end{equation*}
Consequently $\mathcal{M}_{T}^N$ converges to $0$ in the $L^2$-sense, and $\langle \mathcal{Q}_T^N,\mathds{F}\rangle$ converges to $N^{-3}\sum_{(x,e)\in\mathbf{E}^N}\int_0^T\mathds{F}_{x,e}\big(r_{x,e}^{\mathds{F},+}(\eta_s)-r_{x,e}^{\mathds{F},-}(\eta_s)\big)ds$ with
\begin{equation*}
\begin{split}
r_{x,e}^{\mathds{F},+}(\eta)-r_{x,e}^{\mathds{F},-}(\eta)&=\frac{N^2}{2}\Bigg(\eta_{x}(1-\eta_{x+e})\exp\Big(\frac{\mathds{F}_{x,e}}{N\log N}\Big)-\eta_{x+e}(1-\eta_{x})\exp\Big(-\frac{\mathds{F}_{x,e}}{N\log N}\Big)\Bigg)\\
&=\frac{N^2}{2}\Bigg(\eta_x-\eta_{x+e}+\frac{\mathds{F}_{x,e}}{N\log N}\Big(\eta_{x}(1-\eta_{x+e})+\eta_{x+e}(1-\eta_{x})\Big)+O\Big(\frac{|x|^{-2}}{N^2\log^2N};\mathfrak{b}\Big)\Bigg)
\end{split}
\end{equation*}
Thus $\langle \mathcal{Q}_T^N,\mathds{F}\rangle$ converges to
\begin{equation}\label{eq:current_concentration_final}
\begin{split}
&N^{-3}\sum_{(x,e)\in\mathbf{E}^N}\int_0^T\mathds{F}_{x,e}\big(r_{x,e}^{\mathds{F},+}(\eta_s)-r_{x,e}^{\mathds{F},-}(\eta_s)\big)ds=N^{-2}\sum_{(x,e)\in\mathbf{E}^N}\int_0^T\mathds{F}_{x,e}\Bigg(-\frac{N}{2}(\eta_{x+e,s}-\eta_{x,s})+\frac{\mathds{F}_{x,e}}{\log N}\sigma_{x,e}(\eta_s)\Bigg)ds+o_N(1;\mathfrak{b}).
\end{split}
\end{equation}
The concentration estimates in Lemma \ref{lem:concentration_quadratic} and \ref{lem:concentration_linear} implies the RHS of \eqref{eq:current_concentration_final} converges in probability to the following as $N\rightarrow +\infty$
\begin{equation*}
T\int_{N^{-1}\leq |x|\leq 1}\mathcal{F}_{x}\cdot\Big(-\frac{1}{2}\nabla\mathfrak{u}_x+\frac{\mathcal{F}_x}{\log N}\sigma(\mathfrak{u}_x)\Big)dx=T\langle \mathfrak{q},\mathcal{F}\rangle.
\end{equation*}
This concludes the proof of Lemma \ref{lem:current_concentration}.
\end{proof}

For the concentration estimate in Proposition \ref{prop:concentration}, by Equation \eqref{eq:RN_random_decomposition} we have decomposed it into three sub-components, and they are proved separately in Lemma \ref{lem:concentration_quadratic}, Lemma \ref{lem:concentration_linear}. and \ref{lem:current_concentration}. Therefore we have fully established the proof of Proposition \ref{prop:concentration}.

In the same manner of proving Lemma \ref{lem:current_concentration}, we can prove the following lemma. It is useful because in the proof of large deviation lower bound, we need to show $\langle \mathcal{Q}_T^N,\mathcal{V}\rangle\geq JT-w$ with probability $1$ as $N\rightarrow +\infty$, see \eqref{eq:concentration_statement_sec5}. The following Lemma proves the concentration of the winding number $\langle \mathcal{Q}_T^N,\mathcal{V}\rangle$.

\begin{lemma}\label{lem:winding_concentration}
Under the notation of Proposition \ref{prop:concentration} and for any $w>0$, we have
\begin{equation}\label{eq:winding_concentration_statement}
\lim_{N\rightarrow +\infty}\mathds{P}_{\mathfrak{u},\mathds{F}}^N\Big[|\langle \mathcal{Q}_T^N,\mathcal{V} \rangle-T\langle \mathfrak{q},\mathcal{V}\rangle|>\frac{w}{2}\Big]=0.
\end{equation}
\end{lemma}

\appendix

\section{General Forcing and the Radon-Nikodym Derivative}\label{app:radon_nikodym}

In this section, we use the formula in \cite[Appendix A]{bertini2002macroscopic} to establish the Radon--Nikodym derivative between the modified dynamics \eqref{eq:modified_generator} and the SSEP process \eqref{eq:SSEP_generator}. 

\begin{proposition}\label{prop:radon_nikodym}
For the forcing $\mathds{F}$, assume there exists a constant $\mathfrak{b}>0$ such that
\begin{equation}\label{eq:appendix_A_bound}
\big|\mathds{F}_{x}\big|\leq \mathfrak{b}|x|^{-1}\textup{ for any $x\in\mathds{B}^N$},\quad |\mathds{F}_x|\equiv 0 \textup{ for any $|x|\geq 1-\frac{3}{N}$}.
\end{equation}
Then the Radon--Nikodym derivative of $\mathds{P}_{\rho,\mathds{F}}^N$ with respect to the SSEP law $\mathds{P}_{\rho,\SS}^N$, both started from $d\local{\rho}$, is
\begin{equation}\label{eq:RN_derivative_general_forcing}
\begin{split}
\frac{d\mathds{P}_{\rho,\mathds{F}}^N}{d\mathds{P}_{\rho,\SS}^N}
=\exp\Bigg\{\frac{N^2}{\log N}\Bigg[
\langle \mathcal{Q}_T^N,\mathds{F}\rangle-\int_0^T\sum_{(x,e)\in\mathbf{E}^N}\frac{1}{4N^2}\frac{(\mathds{F}_{x,e})^2}{\log N}
\Big(&\eta_{t,x}(1-\eta_{t,x+e})+\eta_{t,x+e}(1-\eta_{t,x})\Big)dt\\
&-\int_0^T\sum_x\frac{1}{2N^2}\eta_{t,x}\ddiv\mathds{F}_xdt
+o_N(1;\mathfrak{b},T)\Bigg]\Bigg\}.
\end{split}
\end{equation}
\end{proposition}

\begin{proof}
For each bond $(x,e)\in\mathbf{E}^N$, set $\theta_{x,e}$ as
\begin{equation*}
\theta_{x,e}:=\frac{\mathds{F}_{x,e}}{N\log N}.
\end{equation*}
The Radon-Nikodym derivative for two Markov jump processes with rates
$r^{\mathds{F},\pm}_{x,e}$ and $r^{0,\pm}_{x,e}$, see \cite[Equation A.6]{bertini2002macroscopic}, gives
\begin{equation}\label{eq:RN_girsanov_raw}
\begin{split}
\log\frac{d\mathds{P}_{\rho,\mathds{F}}^N}{d\mathds{P}_{\rho,\SS}^N}
=&\underbrace{\sum_{(x,e)\in\mathbf{E}^N}
\int_0^T\log\frac{r_{x,e}^{\mathds{F},+}(\eta_{t-})}{r_{x,e}^{0,+}(\eta_{t-})}d\mathbf{N}_t^{x,e,+}
 +
\sum_{(x,e)\in\mathbf{E}^N}
\int_0^T\log\frac{r_{x,e}^{\mathds{F},-}(\eta_{t-})}{r_{x,e}^{0,-}(\eta_{t-})}d\mathbf{N}_t^{x,e,-}}_{\textup{\blue{term (A)}}}\\
&\underbrace{-\int_0^T\sum_{(x,e)\in\mathbf{E}^N}
\Big[
r_{x,e}^{\mathds{F},+}(\eta_t)-r_{x,e}^{0,+}(\eta_t)
+r_{x,e}^{\mathds{F},-}(\eta_t)-r_{x,e}^{0,-}(\eta_t)
\Big]dt}_{\textup{\blue{term (B)}}}.
\end{split}
\end{equation}
For term (A), the definition \eqref{eq:RN_rates} of jump rates shows
\begin{equation*}
\log\frac{r_{x,e}^{\mathds{F},+}}{r_{x,e}^{0,+}}=\frac{\mathds{F}_{x,e}}{N\log N},\quad \log\frac{r_{x,e}^{\mathds{F},-}}{r_{x,e}^{0,-}}=-\frac{\mathds{F}_{x,e}}{N\log N}.
\end{equation*}
Consequently, the term (A) in \eqref{eq:RN_girsanov_raw} is exactly the inner product between $\mathcal{Q}_T^N$ and $\mathds{F}$ with factor $\frac{N^2}{\log N}$
\begin{equation*}
\sum_{(x,e)\in\mathbf{E}^N}\frac{\mathds{F}_{x,e}}{N\log N}
\big(\mathbf{N}_T^{x,e,+}-\mathbf{N}_T^{x,e,-}\big)
=\frac{N^2}{\log N}\langle\mathcal{Q}_T^N,\mathds{F}\rangle,
\end{equation*}
according to the relation \eqref{eq:RN_current_counting} between $\mathcal{Q}_{T,x,e}^N$ and $N_{T}^{x,e,\pm}$.

It remains to expand the term (B) in \eqref{eq:RN_girsanov_raw}.  Apply Taylor expansion and use the bound \eqref{eq:appendix_A_bound} of $\mathds{F}$
\begin{equation*}
\begin{split}
&\exp\Big(\frac{\pm \mathds{F}_{x,e}}{N\log N}\Big)=1\pm \frac{\mathds{F}_{x,e}}{N\log N}+\frac{1}{2}\frac{|\mathds{F}_{x,e}|^2}{(N\log N)^2}+o_N(1;\mathfrak{b})|x|^{-3}N^{-3}\big(\log N\big)^{-2},
\end{split}
\end{equation*}
For a fixed edge $(x,e)\in\mathbf{E}^N$
\begin{equation}\label{eq:RN_one_bond_compensator}
\begin{split}
 &r_{x,e}^{\mathds{F},+}(\eta)-r_{x,e}^{0,+}(\eta)
+r_{x,e}^{\mathds{F},-}(\eta)-r_{x,e}^{0,-}(\eta)\\
&=\frac{N^2}{2}\Bigg(\exp\Big(\frac{\mathds{F}_{x,e}}{N\log N}\Big)\eta_x(1-\eta_{x+e})+\exp\Big(-\frac{\mathds{F}_{x,e}}{N\log N}\Big)\eta_{x+e}(1-\eta_{x})-\eta_{x}(1-\eta_{x+e})-\eta_{x+e}(1-\eta_x)\Bigg)\\
&=\underbrace{\frac{N}{2\log N}\mathds{F}_{x,e}(\eta_x-\eta_{x+e})}_{\textup{\blue{linear term}}}
+\underbrace{\frac{(\mathds{F}_{x,e})^2}{4\log^2N}
\Big(\eta_x(1-\eta_{x+e})+\eta_{x+e}(1-\eta_x)\Big)}_{\textup{\blue{quadratic term}}}
+\underbrace{o_N(1;\mathfrak{b})|x|^{-3}N^{-1}\big(\log N\big)^{-2}}_{\textup{\blue{error term}}}.
\end{split}
\end{equation}
Summing \eqref{eq:RN_one_bond_compensator} over $\mathbf{E}^N$, the linear term becomes a discrete divergence
\begin{equation}\label{eq:RN_discrete_ibp}
\sum_{(x,e)\in\mathbf{E}^N}N\mathds{F}_{x,e}(\eta_x-\eta_{x+e})
=
\sum_{x\in\mathds{B}^N}\eta_x\nabla_{\textup{disc}}\cdot\mathds{F}(x).
\end{equation}
The linear term further becomes the following by separating the space normalization
\begin{equation*}
-\int_0^T\frac{1}{2\log N}\sum_x\eta_t(x)\nabla_{\textup{disc}}\cdot\mathds{F}(x)dt
=-\frac{N^2}{\log N}\int_0^T
\sum_x\frac{1}{2N^2}\eta_t(x)\ddiv\mathds{F}(x)dt.
\end{equation*}
The quadratic term is similarly
\begin{equation*}
\begin{split}
&-\int_0^T\sum_{(x,e)\in\mathbf{E}^N}
\frac{(\mathds{F}_{x,e})^2}{4\log^2N}
\Big(\eta_{t,x}(1-\eta_{t,x+e})+\eta_{t,x+e}(1-\eta_{t,x})\Big)dt\\
&=-\frac{N^2}{\log N}
\int_0^T\sum_{(x,e)\in\mathbf{E}^N}
\frac{1}{4N^2}\frac{(\mathds{F}_{x,e})^2}{\log N}\Big(\eta_{t,x}(1-\eta_{t,x+e})+\eta_{t,x+e}(1-\eta_{t,x})\Big)dt,
\end{split}
\end{equation*}
These two identities give the two terms in \eqref{eq:RN_derivative_general_forcing} involving $\eta$.

Finally the error term in \eqref{eq:RN_one_bond_compensator} gives 
\begin{equation*}
\sum_{(x,e)\in\mathbf{E}^N}o_N(1;\mathfrak{b})|x|^{-3}N^{-1}\big(\log N\big)^{-2}=\frac{N^2}{\log^2N}o_N(1;\mathfrak{b}).
\end{equation*}
which concludes the proof of Proposition \ref{prop:radon_nikodym}.
\end{proof}

\section{Analysis of Cost Functionals}
\subsection{The density profile $\mathfrak{u}^N[\gamma,\zeta;\alpha]$ }\label{app:optimization_singular_profiles}

\begin{proposition}\label{prop:C_cost_cubic}
Take arbitrary parameters $(\gamma,\zeta,\alpha,J)$ with
\begin{equation*}
\gamma\in (0,1),\quad \zeta\in (0,1),\quad \alpha\in (0,1),\quad J\in\mathds{R}.
\end{equation*}
Recall the density profile $\mathfrak{u}^N[\gamma,\zeta;\alpha]$ and the current profile $\mathfrak{q}^N[J]$ defined in \eqref{eq:u_profile}. The following equality has a $N$-independent leading term and a small error
\begin{equation}\label{eq:explicit_C_cost}
\begin{split}
&\mathcal{C}^N\Big(\mathfrak{u}^N[\gamma,\zeta;\alpha],\mathfrak{q}^N[J]\Big)=\frac{J^2}{4\pi}\Big(\frac{1-\alpha}{\zeta(1-\zeta)}+\alpha A_0(\gamma,\zeta)\Big)+\frac{9\pi}{4}\frac{(\gamma-\zeta)^2}{\alpha}A_4(\gamma,\zeta)+C_*o_N(1).
\end{split}
\end{equation}
Here the $N$-independent functions $A_0(\gamma,\zeta)$ and $A_4(\gamma,\zeta)$ are defined as
\begin{equation}\label{eq:function_D}
A_0(\gamma,\zeta):=\int_0^1 \frac{ds}{D_{\gamma,\zeta}(s)},\quad A_4(\gamma,\zeta):=\int_0^1\frac{s^4ds}{D_{\gamma,\zeta}(s)},\quad \textup{ with }D_{\gamma,\zeta}(s)
:=
\bigl(\zeta+(\gamma-\zeta)s^3\bigr)
\bigl(1-\zeta-(\gamma-\zeta)s^3\bigr).
\end{equation}
Consequently the function $\mathcal{C}(\gamma,\zeta,\alpha,J)$ in \eqref{eq:C_functional_parametrized} has the explicit expression
\begin{equation*}
\mathcal{C}(\gamma,\zeta,\alpha,J)=\frac{J^2}{4\pi}\Big(\frac{1-\alpha}{\zeta(1-\zeta)}+\alpha A_0(\gamma,\zeta)\Big)+\frac{9\pi}{4}\frac{(\gamma-\zeta)^2}{\alpha}A_4(\gamma,\zeta).
\end{equation*}

\end{proposition}

\begin{proof}
The presence of $\varphi$ in the definition of $\mathfrak{u}$ in \eqref{eq:u_profile} is smooth enough, only producing an error $C_*o_N(1)$. Therefore in the following calculation, we will omit the presence of $\varphi$ for simplicity. Recall the logarithmic variable $s$
\begin{equation}\label{eq:r_s_relation}
s:=\frac{-\log |x|}{\log N}\in[0,1]\quad \textup{with}\quad r=|x|=N^{-s}\quad \textup{and}\quad \frac{dr}{r}=-\log Nds.
\end{equation}
Define $\mathtt{u}(s)$ as a function on $[0,1]$, being the equivalent of $\mathfrak{u}^N[\gamma,\zeta;\alpha]$ in the logarithmic variable $s$
\begin{equation}\label{eq:u_r_relation}
\mathtt{u}(s)=\mathfrak{u}^N[\gamma,\zeta;\alpha](N^{-s}).
\end{equation}
The explicit expression \eqref{eq:u_profile} of $\mathfrak{u}^N$ implies
\begin{equation*}
\mathtt{u}(s)=
\begin{cases}
\displaystyle
\zeta+(\gamma-\zeta)\left(1-\frac{s}{\alpha}\right)^3, &0\leq s<\alpha,\\[1em]
\zeta, &\alpha\leq s\leq 1.
\end{cases}
\end{equation*}
Lemma \ref{lem:C_self_similarity} implies
\begin{equation*}
\mathcal{C}^N\Big(\mathfrak{u}^N[\gamma,\zeta;\alpha],\mathfrak{q}^N[J]\Big)=\frac{J^2}{4\pi}\int_0^1\frac{1}{\mathtt{u}(s)(1-\mathtt{u}(s))}ds+\frac{\pi}{4}\int_0^1\frac{|\mathtt{u}'(s)|^2}{\mathtt{u}(s)(1-\mathtt{u}(s))}ds
\end{equation*}

We compute the two terms separately. On the interval $0\leq s\leq \alpha$, set $s_*:=1-\frac{s}{\alpha}$. Then $ds=-\alpha ds_*$, and
\begin{equation*}
\mathtt{u}(s)=\zeta+(\gamma-\zeta)s_*^3.
\end{equation*}
For the current term
\begin{equation*}
\begin{aligned}
\int_0^1\frac{ds}{\mathtt{u}(s)(1-\mathtt{u}(s))}&=\int_0^\alpha\frac{ds}{\mathtt{u}(s)(1-\mathtt{u}(s))}+\int_\alpha^1\frac{ds}{\zeta(1-\zeta)}=\alpha\int_0^1\frac{ds_*}{D_{\gamma,\zeta}(s_*)}+\frac{1-\alpha}{\zeta(1-\zeta)}=\alpha A_0(\gamma,\zeta)+\frac{1-\alpha}{\zeta(1-\zeta)}.
\end{aligned}
\end{equation*}

For the gradient term, $\mathtt{u}$ is constant on
$[\alpha,1]$, while for $0<s<\alpha$
\begin{equation*}
\mathtt{u}'(s)=-\frac{3(\gamma-\zeta)}{\alpha}
\left(1-\frac{s}{\alpha}\right)^2,\quad |\mathtt{u}'(s)|^2=\frac{9(\gamma-\zeta)^2}{\alpha^2}s_*^4.
\end{equation*}
Consequently
\begin{equation*}
\begin{aligned}
\int_0^1
\frac{|\mathtt{u}'(s)|^2}{\mathtt{u}(s)(1-\mathtt{u}(s))}ds&=\int_0^\alpha\frac{|\mathtt{u}'(s)|^2}{\mathtt{u}(s)(1-\mathtt{u}(s))}ds=\frac{9(\gamma-\zeta)^2}{\alpha}\int_0^1\frac{s_*^4ds_*}{D_{\gamma,\zeta}(s_*)}=\frac{9(\gamma-\zeta)^2}{\alpha}A_4(\gamma,\zeta).
\end{aligned}
\end{equation*}

Combining the current and gradient contributions gives \eqref{eq:explicit_C_cost} and completes the proof of Proposition \ref{prop:C_cost_cubic}.
\end{proof}

\begin{proof}[Proof of Lower Bound in Corollary \ref{cor:two_scale}]

It is sufficient to prove $\lim_{a\rightarrow +\infty}a^{-2}F^{\downarrow}(aJ)=\pi^{-1} J^2$. In Theorem \ref{th:one-scale-deviation}, together with the explicit expression of $\mathcal{C}$ in Proposition \ref{prop:C_cost_cubic}, the cost function $F^{\downarrow}(aJ)$ is explicitly given as
\begin{equation*}
\begin{split}
&F^{\downarrow}(aJ)\\
&=\mathcal{C}\Big(\gamma,\frac{1}{2},\frac{1}{10}\frac{1}{1+|Ja|},Ja\Big)\\
&=\frac{J^2a^2}{4\pi}\Bigg[4\Big(1-10^{-1}(1+|Ja|)^{-1}\Big)+10^{-1}(1+|Ja|)^{-1}A_0(\gamma,1/2)\Bigg]+\frac{9\pi}{4}(\gamma-1/2)^2A_4(\gamma,\frac{1}{2})10(1+|Ja|).
\end{split}
\end{equation*}
Consequently
\begin{equation*}
\frac{1}{a^2}F^{\downarrow}(aJ)=\frac{J^2}{4\pi}\underbrace{\Bigg[4\Big(1-\frac{1}{10}\frac{1}{1+|Ja|}\Big)+10^{-1}(1+|Ja|)^{-1}A_0(\gamma,1/2)\Bigg]}_{\textup{\blue{goes to $4$ as $a\rightarrow +\infty$}}}+\underbrace{\frac{9\pi}{4}(\gamma-1/2)^2A_4(\gamma,\frac{1}{2})\frac{10(1+|Ja|)}{a^2}}_{\textup{\blue{goes to $0$ as $a\rightarrow +\infty$}}}.
\end{equation*}
This concludes the proof.
\end{proof}

\subsection{Euler-Lagrange Equation for the Optimizer}\label{app:optimizer_euler_lagrange}

This appendix records the variational calculation of the cost functional $\mathcal{C}^N$ defined in \eqref{eq:C_cost_functional}.

\begin{proposition}\label{prop:optimizer}
Consider the variational problem 
\begin{equation*}
\inf_{(\rho,q)}\mathcal{C}^N(\rho,q),\quad \textup{with $\rho\equiv \gamma$ on $\p B_1$ and $\rho\in \mathcal{C}^1\big(B_1;(0,1)\big)$ and $\int_{N^{-1}\leq |x|\leq 1}q_x\cdot\mathcal{V}_xdx\geq J$}.
\end{equation*}
The variational problem
admits a radially symmetric minimizer $(\rho_*,q_*)$. Use \eqref{eq:self_similar_rho_q} to write $(\rho_*,q_*)$ into the self-similar coordinate $s=-\log_N |x|\in[0,1]$ as $(\uprho_*,\mathtt{q}_{*,\theta},\mathtt{q}_{*,r})$, the optimal radial density profile satisfies
\begin{equation}\label{eq:optimizer_first_order_ode}
\frac{d}{ds}\uprho_*(s)=\textup{sgn}(\uprho_*(1)-\gamma)\frac{J}{\pi\int_0^1\sigma(\uprho_*(s_1))ds_1}\sqrt{\sigma(\uprho_*(s))\Big(\sigma(\uprho_*(1))-\sigma(\uprho_*(s))\Big)},\qquad\uprho_*(0)=\gamma,\qquad \uprho_*'(1)=0.
\end{equation}
For this density, the optimal current $\mathtt{q}_*$ is
\begin{equation}\label{eq:optimizer_current}
\mathtt{q}_{*,\theta}(s)=\Lambda\sigma(\uprho_*(s)),\quad \mathtt{q}_{*,r}(s)\equiv 0,\quad \Lambda:=\frac{J}{2\pi\int_0^1\sigma(\uprho_*(s_1))ds_1}.
\end{equation}
\end{proposition}

\begin{proof}

The functional $\mathcal{C}^N(\rho,q)$ in \eqref{eq:C_cost_functional} is jointly convex in $(q,\rho,\nabla\rho)$. Thus for any profile $(\rho,q)$, we have $\mathcal{C}^N(\rho,q)\geq \mathcal{C}^N(\overline{\rho},\overline{q})$ with $(\overline{\rho},\overline{q})$ as the radial average. Thus the variational problem of $\mathcal{C}^N(\rho,q)$ can be restricted to radially symmetric profiles. Recall the self-similar variable $s$ with $|x|=N^{-s}$. Lemma \ref{lem:C_self_similarity} shows
\begin{equation}\label{eq:self_similar_appendix}
\inf_{(\uprho,\mathtt{q})}\Bigg(\pi\int_0^1\frac{|\mathtt{q}_{\theta}(s)|^2+|\mathtt{q}_{r}(s)|^2}{\uprho(s)(1-\uprho(s))}ds+\frac{\pi}{4}\int_0^1\frac{|\uprho'(s)|^2}{\uprho(s)(1-\uprho(s))}ds\Bigg),\quad \uprho\in C^1\big([0,1];(0,1)\big),\quad \uprho(0)=\gamma,\quad 2\pi\int_0^1 \mathtt{q}_{\theta}(s)ds=J.
\end{equation}
Fix any admissible density profile $\uprho$, we minimize the cost related with the current $\mathtt{q}_{\theta}$ satisfying the constrain \eqref{eq:self_similar_appendix}. Cauchy-Schwarz inequality implies
\begin{equation}\label{eq:current_cauchy_schwarz}
\frac{J^2}{4\pi^2}=\Bigg(\int_0^1 \frac{\mathtt{q}_{\theta}(s)}{\sqrt{\sigma(\uprho(s))}}\sqrt{\sigma(\uprho(s))} ds\Bigg)^2\leq \Bigg(\int_0^1 \frac{|\mathtt{q}_{\theta}(s)|^2}{\sigma(\uprho(s))}ds\Bigg)\Bigg(\int_0^1\sigma(\uprho(s))ds\Bigg),\quad \int_0^1\frac{|\mathtt{q}_{\theta}(s)|^2}{\sigma(\uprho(s))}ds\geq \frac{J^2}{4\pi^2}\frac{1}{\int_0^1\sigma(\uprho(s))ds}
\end{equation}
\eqref{eq:current_cauchy_schwarz} reaches equality when
\begin{equation}\label{eq:ansatz_current_givendensity}
\mathtt{q_{\theta}}(s)=\Lambda \sigma(\uprho(s)),\quad \Lambda:=\frac{J}{2\pi\int_0^1\sigma(\uprho(s_1))ds_1}.
\end{equation}
Therefore we only need to optimize the following functional $\mathcal{F}(\uprho)$ of $\uprho$
\begin{equation}\label{eq:radial_functional_for_EL}
\mathcal{F}(\uprho):=\underbrace{\frac{J^2}{4\pi\int_0^1\sigma(\uprho)ds}}_{\textup{\blue{current cost}}}+\underbrace{\frac{\pi}{4}\int_0^1\frac{(\uprho')^2}{\sigma(\uprho)}ds}_{\textup{\blue{gradient cost}}} .
\end{equation}
Let $\eta$ be a perturbation and $\uprho_\varepsilon=\uprho+\varepsilon\eta$, where
$\eta(0)=0$.  The optimizer $\uprho$ must have $\uprho'(1)=0$, otherwise we can properly choose the perturbation $\eta$ to reduce the cost $|\uprho'|^2$ around $s=1$. Therefore we only consider $\uprho'(1)=0$. Denote $\Lambda:=J(2\pi\int_0^1\sigma(\uprho)ds)^{-1}$ as in \eqref{eq:optimizer_current}. Compute the variation of the functional $\mathcal{F}$, and the optimizer must satisfy for any perturbation
\begin{equation}\label{eq:variation_of_F}
\frac{d}{d\eps}\mathcal{F}(\uprho+\eps \eta)|_{\eps=0}=\underbrace{\frac{\pi}{2}\int_0^1\frac{\uprho'}{\sigma(\uprho)}\eta'ds-\frac{\pi}{4}\int_0^1\frac{(\uprho')^2\sigma'(\uprho)}{\sigma(\uprho)^2}\eta ds}_{\textup{\blue{variation of the gradient cost}}}-\underbrace{\frac{J^2}{4\pi(\int_0^1\sigma(\uprho)ds)^2}\int_0^1\sigma'(\uprho)\eta ds}_{\textup{\blue{variation of the current cost}}}\equiv 0
\end{equation}
Consequently the optimizer must satisfy
\begin{equation*}
    2\frac{d}{ds}\left(\frac{\uprho'}{\sigma(\uprho)}\right)+\frac{(\uprho')^2\sigma'(\uprho)}{\sigma(\uprho)^2}+4\Lambda^2\sigma'(\uprho)=0.
\end{equation*}
Multiplying by $\uprho'$, this equation becomes the conservation law
\begin{equation*}
\frac{d}{ds}\left(\frac{(\uprho')^2}{\sigma(\uprho)}+4\Lambda^2\sigma(\uprho)\right)=0.
\end{equation*}
Since $\uprho'(1)=0$ and $m=\uprho(1)$, the conserved value is
$4\Lambda^2\sigma(m)$.  Hence
\begin{equation*}
(\uprho')^2=4\Lambda^2\sigma(\uprho)\bigl(\sigma(\uprho(1))-\sigma(\uprho)\bigr).
\end{equation*}
When $\uprho(0)=\gamma<\frac{1}{2}$, we take the increasing branch with  $\uprho'(s)\geq 0$, and the density profile $\uprho(s)$ goes towards $\frac{1}{2}$ as $s\rightarrow 1$. In this case we have \eqref{eq:optimizer_first_order_ode}. When $\uprho(0)=\gamma>\frac{1}{2}$, we take the decreasing branch with  $\uprho'(s)\leq 0$, and the density profile $\uprho(s)$ also goes towards $\frac{1}{2}$ as $s\rightarrow 1$. This concludes the proof of \eqref{eq:optimizer_first_order_ode}. Given the optimal density profile $\uprho_*$, the optimal current profile $\mathtt{q}_*$ is given by \eqref{eq:ansatz_current_givendensity}, and it proves \eqref{eq:optimizer_current}.
\end{proof}

\section{Auxiliary Entropy, Integration by Parts, and Subgaussian Estimates}
\label{app:integration_subgaussian}

This appendix collects the elementary results used in the entropy and
concentration estimates in Sections \ref{sec:relative_entropy}-\ref{sec:concentration_RN}. First we state an elementary estimate on the summation of singular functions on $\mathds{B}^N$. This estimate is frequently used across the whole paper.

\begin{lemma}\label{lem:lattice_sum}
For any $s_1,s_2\in[0,1]$ with $s_1> s_2$, and any integer $m\geq 3$, when $N$ is large enough, we have
\begin{equation*}
\sum_{x\in\mathds{B}^N:N^{-s_1}\leq |x|\leq N^{-s_2}}|x|^{-2}\leq C|s_1-s_2|N^2\log N,\quad \sum_{x\in\mathds{B}^N:N^{-s_1}\leq |x|\leq N^{-s_2}}|x|^{-m}\leq CN^2N^{(m-2)s_1}.
\end{equation*}
\end{lemma}

\begin{proof}
The summation over lattice points $x\in\mathds{B}^N$ will be bounded by integrals
\begin{equation*}
\begin{split}
\sum_{x\in\mathds{B}^N:N^{-s_1}\leq |x|\leq N^{-s_2}}|x|^{-2}\leq N^{2}\sum_{x\in\mathds{B}^N:N^{-s_1}\leq |x|\leq N^{-s_2}}N^{-2}|x|^{-2}&\leq  CN^2\int_{N^{-s_1}\leq |x|\leq N^{-s_2}}|x|^{-2}dx\\
&\leq CN^2\int_{N^{-s_1}}^{N^{-s_2}}r^{-1}dr\leq C|s_1-s_2|N^2\log N,
\end{split}
\end{equation*}
and in the same way
\begin{equation*}
\begin{split}
\sum_{x\in\mathds{B}^N:N^{-s_1}\leq |x|\leq N^{-s_2}}|x|^{-m}&\leq  CN^2\int_{N^{-s_1}\leq |x|\leq N^{-s_2}}|x|^{-m}dx\\
&\leq CN^2\int_{N^{-s_1}}^{N^{-s_2}}r^{-{m-1}}dr\leq CN^2N^{(m-2)s_1}\leq CN^{2+(m-2)s_1}.
\end{split}
\end{equation*}
\end{proof}

\subsection{Entropy Inequalities}

\begin{lemma}\label{lem:entropic_inequality_app}
Let $\rho$ be a density profile with values in $(0,1)$ and let $f$ be a density with
respect to $d\local{\rho}$ on $\{0,1\}^{\mathds{B}^N}$.  Then
\begin{equation}\label{eq:entropy_variational_app}
\local{\rho}[f\log f]=\sup_V\left\{\local{\rho}[fV]-\log\local{\rho}[e^V]\right\},
\end{equation}
where the supremum is over all real-valued functions $V$ on
$\{0,1\}^{\mathds{B}^N}$.  Consequently, for every $\gamma>0$ and every $V$,
\begin{equation}\label{eq:entropic_inequality_app}
\local{\rho}[Vf]\leq\frac{1}{\gamma}\local{\rho}[f\log f]+\frac{1}{\gamma}\log\local{\rho}[e^{\gamma V}] .
\end{equation}
Moreover, for every event $A$ with
$0<\local{\rho}(A)<1$,
\begin{equation}\label{eq:entropy_event_inequality_app}
\local{\rho}[f\mathds{1}_A]\leq\frac{\local{\rho}[f\log f]+\log2}{\log\local{\rho}(A)^{-1}} .
\end{equation}
\end{lemma}

Given a family of random variables $\{\xi_x:x\in\mathds{B}^N\}$, we say it is $\ell$-independent with respect to $\local{\rho}$ for a positive integer $\ell$, if for any set $B\subset \mathds{B}^N$ with
\begin{equation}
\inf_{x,y\in B,x\not=y}|x-y|\geq N^{-1}\ell,
\end{equation}
the random variables $\{\xi_x:x\in B\}$ are jointly independent.

Using Lemma \ref{lem:entropic_inequality_app} and following \cite[Appendix F.1]{Jara_Menezes_2018}, we can prove the following lemma.

\begin{lemma}\label{lem:entropic_ineq_correlation}
Let $\{\xi_x:x\in \mathds{B}^N\}$ be real-valued random variables. Assume there exists an integer $\ell$ such that $\{\xi_x:x\in\mathds{B}^N\}$ is $\ell$-independent. Then for every 
$\gamma>0$
\begin{equation}\label{eq:entropic_coloring_l2_app}
\begin{split}
&\mathds{E}_{\rho,\mathds{F},t}^N\left[\sum_{x\in \mathds{B}^N}\xi_x\right]\leq\frac{C}{\gamma}\mathds{H}^N\Big[d\mathds{P}_{\rho,\mathds{F},t}^N\Big|d\local{\rho}\Big]+\frac{C}{\gamma \ell^2}\sum_{x\in \mathds{B}^N}\log\local{\rho}[e^{\gamma \ell^2\xi_x}] \\
&\mathds{E}_{\rho,\mathds{F},t}^N\left[\Big|\sum_{x\in \mathds{B}^N}\xi_x\Big|\right]\leq\frac{C}{\gamma}\mathds{H}^N\Big[d\mathds{P}_{\rho,\mathds{F},t}^N\Big|d\local{\rho}\Big]+\frac{C\log 2}{\gamma }+\max_{k\in\{+1,-1\}}\frac{C}{\gamma \ell^2}\sum_{x\in \mathds{B}^N}\log\local{\rho}[e^{k\gamma \ell^2\xi_x}]
\end{split}
\end{equation}
\end{lemma}

\subsection{Feynman-Kac Formula and Dirichlet Forms}
Here we recall the statement of the Feynman-Kac formula, in a form adapted to the simple exclusion process, see \cite[Lemma A.2]{Jara_Menezes_2018}. Recall that $\Gamma_{\SS}^N$ is the Dirichlet form given in \eqref{eq:full_dirichley}, and $\Gamma_{\mathds{F}}^N$ is the Dirichlet form given in \eqref{eq:modified_generator_Dirichlet}.
\begin{lemma}[Feynman--Kac variational bound]\label{lem:feynman_kac_variational_app}
Consider any positive integer $N>1$. For any density profile $\rho$ and any forcing $\mathds{F}$ on $\mathds{B}^N$, and every real-valued function
$G$ on $\{0,1\}^{\mathds{B}^N}$
\begin{equation}\label{eq:FK_variational_app}
\begin{split}
&\log\mathds{E}_{\rho,\mathds{F}}^N\left[\exp\left\{\int_0^T G(\eta_t)dt\right\}\right]\leq T\sup_f\left\{\local{\rho}[Gf]-\frac{1}{2}\local{\rho}\big[\Gamma_{\mathds{F}}^N(\sqrt f)\big]+\frac{1}{2}\local{\rho}[fL_{\mathds{F}}^{N,*}\mathds{1}]
\right\},
\end{split}
\end{equation}
where the supremum is over densities $f$ with respect to
$d\local{\rho}$. For parameter $\mathfrak{b}>0$ with $|\mathds{F}_x|\leq \mathfrak{b}|x|^{-1}$ for $x\in\mathds{B}^N$, we also have
\begin{equation}\label{eq:FK_variational_app_comparison}
\begin{split}
&\log\mathds{E}_{\rho,\mathds{F}}^N\left[\exp\left\{\int_0^T G(\eta_t)dt\right\}\right]\leq T\sup_f\left\{\local{\rho}[Gf]-C_*\local{\rho}\big[\Gamma_{\SS}^N(\sqrt f)\big]+\frac{1}{2}\local{\rho}[fL_{\mathds{F}}^{N,*}\mathds{1}]\right\},
\end{split}
\end{equation}
where $C_*$ is a constant determined by $\mathfrak{b}$.
\end{lemma}

Here we recall the statement of \cite[Lemma E.3]{Jara_Menezes_2018}
\begin{lemma}\label{lem:integration_by_parts_app}
Let $x,y\in\mathds{B}^N$ be two adjacent sites and $\rho$ a density profile with value within $[\eps_0,1-\eps_0]$ for a constant $\eps_0>0$. Let
$f:\{0,1\}^{\mathds{B}^N}\to[0,\infty)$ be a density with respect to
$\local{\rho}$, and let
$h:\{0,1\}^{\mathds{B}^N}\to \mathds{R}$ satisfy
\begin{equation*}
    h(\eta^{x,y})=h(\eta),
    \qquad \eta\in\{0,1\}^{\mathds{B}^N}.
\end{equation*}
Then, for every $\delta>0$,
\begin{equation}\label{eq:one_bond_ibp_app}
\begin{split}
\local{\rho}\Big[h(\omega_y-\omega_x)f\Big]
&\leq
\delta \local{\rho}\big[\Gamma_{\SS}^N(\sqrt{f})\big]+\frac{4}{\delta\eps_0^2 N^2}\local{\rho}\Big[h^2f\Big]
-(\rho_y-\rho_x)\local{\rho}\Big[h\omega_x\omega_y f\Big].
\end{split}
\end{equation}
\end{lemma}

\begin{proof}
This is the one-bond integration-by-parts estimate of Appendix E in
\cite{Jara_Menezes_2018}, written with the normalization used here.
\end{proof}

\subsection{Subgaussian Estimates}
We say a random variable $\xi$ is sub-gaussian of order $\sigma^2$ if
\begin{equation}\label{eq:def_subgaussian}
\log\mathds{E}\Big[e^{\theta \xi}\Big]\leq \frac{1}{2}\sigma^2\theta^2,\quad \textup{for all $\theta\in\mathds{R}$, and we write $\mathcal{SG}(\xi)$ for $\sigma^2$}.
\end{equation}
Given a sequence $\{\xi_{i}\}_{1\leq i\leq n}$ of independent subgaussian variables, and a sequence $\{a_i\}_{1\leq i\leq n}$ of coefficients
\begin{equation*}
\log\mathds{E}\Big[\exp\Big(\theta\sum_{1\leq i\leq n}a_i\xi_i\Big)\Big] =\sum_{1\leq i\leq n}\log\mathds{E}\Big[e^{\theta a_i\xi_i}\Big]\leq \frac{1}{2}\Big(\sum_{1\leq i\leq n}|a_i|^2\mathcal{SG}(\xi_i)\Big)\theta^2.
\end{equation*}
Therefore $\sum_{1\leq i\leq n}a_i\xi_i$ is also sub-gaussian and
\begin{equation}\label{eq:subgaussian_sum}
\mathcal{SG}\big(\sum_{1\leq i\leq n}a_i\xi_i\big)\leq \sum_{1\leq i\leq n}|a_i|^2\mathcal{SG}(\xi_i).
\end{equation}
We recall the statement of \cite[Lemma F.7]{Jara_Menezes_2018} and \cite[Lemma F.8]{Jara_Menezes_2018}.
\begin{lemma}\label{lem:quadratic_subgaussian}
Let $\xi$ be subgaussian of order $\sigma^2$.  If
$0\leq \theta\leq (4\sigma^2)^{-1}$, then
\begin{equation*}
\log \mathds{E}\big[e^{\theta \xi^2}\big]\lesssim\theta\sigma^2 .
\end{equation*}
Let $\xi_1$ and $\xi_2$ be subgaussian variables of orders
$\sigma_1^2$ and $\sigma_2^2$, respectively.  If
$0\leq \theta\leq (4\sigma_1\sigma_2)^{-1}$, then
\begin{equation*}
\log\mathds{E}\big[e^{\theta\xi_1\xi_2}\big]\lesssim\theta\sigma_1\sigma_2 .
\end{equation*}
\end{lemma}

We recall the terms entering Lemmas \ref{lem:subgaussian_omega} and
\ref{lem:subgaussian_H}, so that the square-sum estimates below can be read
without returning to Section \ref{sec:concentration_quadratic}.  The centered
coordinate used in the decomposition \eqref{eq:decomposition_into_omega} is $\omega_x:=\frac{\eta_x-\mathfrak{u}_x}{\mathfrak{u}_x(1-\mathfrak{u}_x)}$. The singular coefficient is first rescaled in \eqref{eq:scaled_H}:
\begin{equation}\label{eq:appD_recalled_scaled_coefficient}
\mathfrak{h}_x^N:=\frac{N^{4\alpha}\log N}{N^2}H_x^N\mathds{1}_{|x|\leq N^{-\beta}},\qquad|\mathfrak{h}_x^N|\lesssim N^{4\alpha-2}|x|^{-2}\mathds{1}_{|x|\leq N^{-\beta}},
\end{equation}
where the bound is exactly \eqref{eq:asymp_H}.  For the mesoscopic length
$\ell$, the left and right averages from \eqref{eq:decomposition_components_a} are
\begin{equation}\label{eq:appD_recalled_averages}
\begin{split}
\lomega_x^{\ell}&:=\sum_{y\in\avedisk}\omega_{x-y}\mathfrak{h}_{x-y}^N\ell^{-2},\quad \romega_{x+e}^{\ell}:=\sum_{z\in\avedisk}\omega_{x+e+z}\ell^{-2}.
\end{split}
\end{equation}
Finally, for each fixed direction $e'$, the flow coefficient from
\eqref{eq:decomposition_components_b} is
\begin{equation}\label{eq:appD_recalled_flow_coefficient}
\mathcal{H}_{x,e'}^N=\sum_{y\in\avediskplus}\phi_{\ell}(y;e')\omega_{x-y}\mathfrak{h}_{x-y}^N .
\end{equation}
The weights $\phi_{\ell}(\cdot;e')$ are the two-dimensional flow weights whose $\ell^2$ estimate is recorded in
Lemma \ref{lem:two_dimensional_flow_lemma_self_contained}.

\begin{lemma}\label{lem:subgaussian_omega}
The sub-gaussian orders in \eqref{eq:def_subgaussian} of the left average $\lomega_x^{\ell}$ and the right average $\romega_x^{\ell}$ in \eqref{eq:appD_recalled_averages} satisfies the following inequality, for any $x\in\mathds{B}^N$ with $|x|\leq N^{-\beta}$ (recall $\ell=0.5N^{1-\beta}$)
\begin{equation*}
\mathcal{SG}(\lomega_x^{\ell})\leq C_*N^{8\alpha-4(1-\beta)},\quad \mathcal{SG}(\romega_x^\ell)\leq C_*\ell^{-2}.
\end{equation*}
\end{lemma}
\begin{proof}
Since the densities stay between $\gamma$ and $\zeta$, the centered variables
$\omega_z$ are uniformly subgaussian.
Therefore the subgaussian order of a linear combination is controlled by the
sum of the squared coefficients, see \eqref{eq:subgaussian_sum}.  For $\lomega_x^{\ell}$, the coefficient of
$\omega_{x-y}$ is $\ell^{-2}\mathfrak{h}_{x-y}^N$ by
\eqref{eq:appD_recalled_averages}.  Hence
\begin{equation}\label{eq:lomega_sg_base}
\begin{split}
\mathcal{SG}(\lomega_x^{\ell})\leq C_*\sum_{y\in\mathbf{D}_{\ell}^N}\Big( \ell^{-2}\mathfrak{h}_{x-y}^N\Big)^2\quad \textup{with }\Big|\ell^{-2}\mathfrak{h}_{x-y}^N\Big|\leq C\ell^{-2}N^{4\alpha-2}|x-y|^{-2}.
\end{split}
\end{equation}
Since $|\ell^{-2}\mathfrak{h}_{x-y}^N|^2$ is bounded from above by $\ell^{-4}N^{8\alpha-4}|x-y|^{-4}$, we use Lemma \ref{lem:lattice_sum} to bound the sum of $|x-y|^{-4}$ over $y\in \mathds{B}^N$ by $N^4$, and then we use the fact that $\ell=0.5N^{1-\beta}$
\begin{equation*}
\mathcal{SG}(\lomega_x^{\ell})\leq C\ell^{-4}N^{8\alpha-4}\sum_{y\in \avedisk}|x-y|^{-4}\leq C\ell^{-4}N^{8\alpha}\leq CN^{8\alpha-4(1-\beta)}.
\end{equation*}
Finally, $\romega_x^{\ell}$ has no singular coefficient: by
\eqref{eq:appD_recalled_averages} it is an average of $\ell^2$ centered
variables with weights $\ell^{-2}$.  Its square sum is therefore $\sum_{|z|\leq N^{-1}\ell}(\ell^{-2})^2\asymp\ell^2\ell^{-4}=\ell^{-2}$ which gives $\mathcal{SG}(\romega_x^{\ell})\asymp \ell^{-2}$.
\end{proof}
\begin{lemma}\label{lem:subgaussian_H}
The random coefficient $\mathcal{H}_{x,e'}^N$ in \eqref{eq:appD_recalled_flow_coefficient} satisfies $\mathcal{SG}(\mathcal{H}_{x,e'}^N)\leq C_*N^{8\alpha}\log N$ for any $x\in\mathds{B}^N$ with $|x|\leq N^{-\beta}$ and any $e'\in\{e_1,e_2\}$.
\end{lemma}
\begin{proof}
By the recalled definition \eqref{eq:appD_recalled_flow_coefficient} of $\mathcal{H}_{x,e'}^N$, the
coefficient of the centered variable $\omega_{x-y}$ in $\mathcal{H}_{x,e'}^N$
is $\phi_{\ell}(y;e')\mathfrak{h}_{x-y}^N$.  By \eqref{eq:subgaussian_sum}, the
subgaussian order is therefore bounded by the square sum of these coefficients. For arbitrary $y\in\avedisk$, the function $(\mathfrak{h}_{x-y}^N)^2$ is bounded from above by $N^{8\alpha}$ according to \eqref{eq:appD_recalled_scaled_coefficient}. Consequently by Lemma \ref{lem:two_dimensional_flow_lemma_self_contained}
\begin{equation*}
\mathcal{SG}\big(\mathcal{H}_{x,e'}^N\big)\leq C_*\sum_{y\in\mathbf{D}_{2\ell-1}^N}\phi_{\ell}(y;e')^2\big(\mathfrak{h}_{x-y}^N\big)^2\leq C_*N^{8\alpha}\sum_{y\in\avediskplus}\phi_{\ell}(y;e')^2\leq C_*N^{8\alpha}\log N.
\end{equation*}
\end{proof}

\section{Two-dimensional Flow Lemma and Discrete Helmholtz Decomposition}
\label{app:two_dimensional_flow_lemma}

This appendix records the two-dimensional form of the flow lemma used in the
Jara--Menezes averaging argument, \cite[Lemma G.2]{Jara_Menezes_2018}.  We write the statement directly on the
scaled lattice $N^{-1}\mathds{Z}^2$, and we use the set of edges $e_1=(N^{-1},0)$ and $e_2=(0,N^{-1})$. For $\ell\geq1$, recall the averaging square
\begin{equation*}
\mathbf{D}_\ell^N:=\{0,N^{-1},\dots,(\ell-1)N^{-1}\}^2\subset N^{-1}\mathds{Z}^2,
\end{equation*}
Set
\begin{equation*}
p_\ell(x):=\ell^{-2}\mathds{1}_{\mathbf{D}_\ell^N}(x),\qquad q_\ell(x):=\big(p_\ell*p_\ell\big)(x)=\sum_{y\in N^{-1}\mathds{Z}^2}p_\ell(y)p_\ell(x-y).
\end{equation*}
Then $q_\ell$ is supported in $\mathbf{D}_{2\ell-1}^N$. A flow is a function $\phi:N^{-1}\mathds{Z}^2\times\{e_1,e_2\}\to\mathds{R}$. We say that $\phi$ connects two probability measures $p$ and $q$ on
$N^{-1}\mathds{Z}^2$ if, for every $z\in N^{-1}\mathds{Z}^2$
\begin{equation}\label{eq:flow_connection_d2}
p(z)-q(z)=\sum_{e\in\{e_1,e_2\}}\big(\phi(z;e)-\phi(z-e;e)\big).
\end{equation}
The support of $\phi$ is contained in a set $A\subset N^{-1}\mathds{Z}^2$ if
$\phi(x;e)\neq0$ implies $x,x+e\in A$.

\begin{lemma}\label{lem:two_dimensional_flow_lemma_self_contained}
For every $\ell\geq1$,
there exists a flow function $\phi_\ell:N^{-1}\mathds{Z}^2\times\{e_1,e_2\}\to\mathds{R}$ which connects $\delta_0$ to $q_\ell$, and is supported in
$\mathbf{D}_{2\ell-1}^N$. It satisfies
\begin{equation*}
\sum_{\substack{x\in N^{-1}\mathds{Z}^2\\ e\in\{e_1,e_2\}}}\phi_\ell(x;e)^2\leq C\log\ell,\qquad \ell\geq 2.
\end{equation*}
\end{lemma}

\begin{proof}
We want to find a function $\Psi$ supported on $\mathbf{D}_{2\ell-1}^N$ such that, the discrete gradient $N^{-1}\nabla_{\textup{disc}}\Psi$ restricted within $\mathbf{D}_{2\ell-1}^N$ as a flow function $\phi_{\ell}(x;e)=\Psi_{x+e}-\Psi_x$ solves \eqref{eq:flow_connection_d2}, or equivalently
\begin{equation}\label{eq:discrete_potential_equation}
\delta_0(x)-q_{\ell}(x)=\sum_{y\in\mathbf{D}_{2\ell-1}^N:y\sim x}\big(\Psi_{y}-\Psi_{x}\big),\quad\textup{ for all $x\in \mathbf{D}_{2\ell-1}^N$},
\end{equation}
which is actually the discrete potential equation on $\mathbf{D}_{2\ell-1}^N$. We denote $h(x):=\delta_0(x)-q_{\ell}(x)$. Taking $G_{\ell}(x,y)$ as the two-dimensional Green function on $\mathbf{D}_{2\ell-1}^N$, we know that
\begin{equation}\label{eq:green_function}
\Psi_{x}:=\sum_{y\in \mathbf{D}_{2\ell-1}^N}G_{\ell}(x,y)h(y)
\end{equation}
solves \eqref{eq:discrete_potential_equation}. We integrate by parts
\begin{equation}\label{eq:appendix_D_integration_by_parts}
\frac{1}{2}\sum_{x,y\in\mathbf{D}_{2\ell-1}^N:x\sim y}|\Psi_{y}-\Psi_{x}|^2=-\sum_{x\in\avediskplus}\Psi_x\sum_{y\in\mathbf{D}_{2\ell-1}^N:y\sim x}\big(\Psi_{y}-\Psi_{x}\big)=-\sum_{x\in\avediskplus}\Psi_xh_x.
\end{equation}
Apply \eqref{eq:green_function} to \eqref{eq:appendix_D_integration_by_parts}, and use the fact that $\phi_{\ell}(x;e)=\Psi_{x+e}-\Psi_{x}$ and $|G_{\ell}|$ is always bounded from above by $C\log\ell$ (for example see \cite[Proposition 2.2]{garban2023statistical})
\begin{equation*}
\begin{split}
&\sum_{x\in N^{-1}\mathds{Z}^2,e\in\{e_1,e_2\}}\phi_\ell(x;e)^2\\
&\leq \frac{1}{2}\sum_{x,y\in\mathbf{D}_{2\ell-1}^N:x\sim y}|\Psi_{x}-\Psi_{y}|^2\leq \sum_{x\in\mathbf{D}_{2\ell-1}^N}\Big|\Psi_{x}h_x\Big|\leq\sum_{x,y\in\mathbf{D}_{2\ell-1}^N}\Big|h_xG_{\ell}(x,y)h_y\Big|\leq C\log\ell\Big(\sum_{x\in\mathbf{D}_{2\ell-1}^N}|h(x)|\Big)^2\leq C\log\ell.
\end{split}
\end{equation*}
This concludes the proof.
\end{proof}

\begin{proposition}\label{prop:discrete_Helmholtz}
\textup{\textbf{[Discrete Helmholtz decomposition]}} Given any vector field $\mathbf{F}$ on $\mathds{B}^N$, it can be decomposed into $\mathbf{F}=g+\nabla_{\textup{disc}}\phi$, with  
\begin{equation}
\tilde{\nabla}_{\textup{disc}}\cdot g_x\equiv 0,\quad \tilde{\nabla}_{\textup{disc}}\cdot g_x:=\sum_{e\in\{e_1,e_2\}}\Big(\mathds{1}_{x+e\in \mathds{B}^N}g_{x,e}-\mathds{1}_{x-e\in \mathds{B}^N}g_{x-e,e}\Big),
\end{equation}
and
\begin{equation}\label{eq:discrete_helmholtz_bound}
\begin{split}
&|\phi_x|\leq C\underbrace{\frac{1}{N^2}\sum_{e\in\{e_1,e_2\}}\sum_{y\in \mathds{B}^N}\frac{|\mathbf{F}_{y,e}|}{N^{-1}+|x-y|}}_{\textup{\blue{convolution with the gradient of Green function}}},\quad \big|\nabla_{\textup{disc}}\phi_x\big|\leq C\frac{1}{N^2}\underbrace{\sum_{e\in\{e_1,e_2\}}\sum_{y\in \mathds{B}^N}\frac{|\mathbf{F}_{y,e}|}{N^{-2}+|x-y|^2}}_{\textup{\blue{convolution with the $\nabla^2$ of Green function}}}\\
&N^{-2}\sum_{x\in \mathds{B}^N}|g_x|^2\leq N^{-2}\sum_{x\in \mathds{B}^N}|\mathbf{F}_x|^2
\end{split}
\end{equation}
Furthermore, for any vector field $\mathbf{F}$ defined on $N^{-1}\mathds{Z}^2$ satisfying for any $x\in \mathds{B}^N$ and any $e\in\{e_1,e_2\}$
\begin{equation}\label{eq:dis_helm_condition}
|\mathbf{F}_{x,e}|\leq \mathfrak{b}|x|^{-1}\mathds{1}_{N^{-\eps}\leq |x|\leq 1},
\end{equation}
we have 
\begin{equation}\label{eq:phi_upper_bound}
\begin{split}
&|\phi_x|\leq C_{\mathfrak{b}}(1-\log|x|)\textup{ for any $x\in\mathds{B}^N$};\Big(N^{-2}\sum_{x\in\mathds{B}^N:\frac{\eps}{2}<|x|<1-\eps}|\nabla_{\textup{disc}}\phi_{x}|^2\Big)^{\frac{1}{2}}\leq C_{\mathfrak{b}}\big(1+|\log\eps|\big)\\
& N^{-2}\sum_{x\in \mathds{B}^N}|g_x|^2\leq N^{-2}\sum_{x\in \mathds{B}^N}|\mathbf{F}_x|^2=C_{\mathfrak{b}}\eps O(\log N)
\end{split}
\end{equation}
\end{proposition}
\begin{proof}
The decomposition $\mathbf{F}=g+\nabla_{\textup{disc}}\phi$ and the upper bound in \eqref{eq:discrete_helmholtz_bound} is classical, for example see \cite[Lemma 12]{ehrlacher2016analysis} and \cite[Equation (3.5)]{lim2020hodge}. Now we prove \eqref{eq:phi_upper_bound}. Given a vector field $\mathbf{F}$ supported on $\mathds{B}^N$ satisfying \eqref{eq:dis_helm_condition}, by \eqref{eq:discrete_helmholtz_bound} we have the upper bound on $|\phi_x|$
\begin{equation}\label{eq:phi_first}
|\phi_x|\leq C_{\mathfrak{b}}N^{-2}\sum_{e\in\{e_1,e_2\}}\sum_{y\in \mathds{B}^N}\frac{1}{N^{-1}+|x-y|}\frac{1}{|y|}\leq C\int_{|y|\leq 1}\frac{1}{(N^{-1}+|x-y|)|y|}dy.
\end{equation}
Decompose the integral into two regions
\begin{equation*}
A_1:=\Big\{y\in B_1\Big|\ |y|\leq 2|x|\Big\},\quad A_2:=\Big\{y\in B_1\Big|\ |y|\geq 2|x|\Big\}.
\end{equation*}
Estimate the integral in \eqref{eq:phi_first} within the region $A_1$
\begin{equation*}
\int_{|y|\leq 2|x|}\frac{1}{(N^{-1}+|x-y|)}\frac{1}{|y|}dy\leq \int_{|y|\leq 2|x|}\frac{1}{|x-y|}\frac{1}{|y|}dy\leq \int_{|y'|\leq 2}\frac{1}{|(0,1)-y'|}\frac{1}{|y'|}dy'\leq C,
\end{equation*}
where in the last step we have changed variables from $(x,y)$ to $(x|x|^{-1},y')$ with $y'=y|x|^{-1}$. Now we estimate the integral in \eqref{eq:phi_first} within the region $A_2$, where $|y|\leq C|x-y|$
\begin{equation*}
\int_{1\geq |y|\geq 2|x|}\frac{1}{(N^{-1}+|x-y|)}\frac{1}{|y|}dy\leq C\int_{1\geq |y|\geq 2|x|}\frac{1}{|y|^2}dy\leq C\log\big(\frac{1}{|x|}\big).
\end{equation*}
Putting the estimates for $A_1$ and $A_2$ together, we prove the upper bound of $|\phi_x|$ in \eqref{eq:phi_upper_bound}. The $L^2$-bound of $\nabla_{\textup{disc}}\phi_x$ follows from the discrete Caccioppoli estimate, for example see \cite[Equation (3.4)]{gloria2026moment}. The $L^2$-bound of $g$ is derived by the orthogonality between $g$ and $\nabla_{\textup{disc}}\phi$.
\end{proof}

\paragraph{AI Statement} The initial formulation of this project and the mathematical contents are entirely produced by the author, starting from 2024, without AI assistance. The writing of this paper is also entirely produced by the author except for: AI tools were used for collecting literature references on classical contents, proofreading, and transforming the figures hand-written by the author into latex-drawn figures.

\bibliographystyle{alpha}
\bibliography{bibIT.bib}

\end{document}